\documentclass[11pt]{article}

\usepackage{amssymb,latexsym,amsmath,amsthm,verbatim}
\usepackage{graphicx,epsfig,epstopdf,amssymb,color}
\usepackage[normalem]{ulem}

\newtheorem{theorem}{Theorem}[section]
\newtheorem{lemma}[theorem]{Lemma}

\newtheorem{corollary}[theorem]{Corollary}
\newtheorem{remark}[theorem]{Remark}

\catcode`\@=11
\@addtoreset{equation}{section}

\def\beq{\begin{equation}}
\def\eeq{\end{equation}}
\def\beqa{\begin{eqnarray}}
\def\eeqa{\end{eqnarray}}

\def\eps{\varepsilon}

\def\RR{\mathbb{R}}
\def\CC{\mathbb{C}}
\def\ZZ{\mathbb{Z}}

\def\vnab{\vec{\nabla}}
\def\vmuta{\vmu_{\rm t}^\ast}

\def\A{{\cal A}}
\def\B{{\cal B}}
\def\C{{\cal C}}
\def\D{{\cal D}}

\def\H{{\cal H}}

\def\L{{\cal L}}
\def\M{{\cal M}}
\def\N{{\cal N}}
\def\O{{\cal O}}

\def\Q{{\cal Q}}
\def\R{{\cal R}}
\def\S{{\cal S}}

\def\V{{\cal V}}
\def\W{{\cal W}}

\def\al{\alpha}
\def\be{\beta}
\def\de{\delta}
\def\ga{\gamma}
\def\la{\lambda}
\def\si{\sigma}

\def\Ga{\Gamma}

\def\tc{\tilde{c}}
\def\tf{\tilde{f}}
\def\tF{\tilde{F}}
\def\tG{\tilde{G}}

\def\tM{\tilde{M}}
\def\tN{\tilde{N}}

\def\tQ{\tilde{Q}}
\def\tu{\tilde{u}}
\def\tv{\tilde{v}}

\def\tX{\tilde{X}}
\def\tY{\tilde{Y}}
\def\tmu{\tilde{\mu}}
\def\tla{\tilde{\lambda}}
\def\tsi{\tilde{\sigma}}

\def\tth{\tilde{\theta}}

\def\bu{\bar{u}}
\def\bU{\bar{U}}
\def\bv{\bar{v}}
\def\bV{\bar{V}}

\def\oU{\overline{U}}
\def\oV{\overline{V}}
\def\oF{\overline{F}}
\def\oG{\overline{G}}
\def\oh{\overline{h}}

\def\vmu{\vec{\mu}}

\def\gas{\gamma_{\rm s}}
\def\ps{p_{\rm s}}
\def\qs{q_{\rm s}}
\def\us{u_{\rm s}}
\def\Us{U_{\rm s}}
\def\vs{v_{\rm s}}
\def\Vs{V_{\rm s}}
\def\Ls{{\cal L}_{\rm s}}
\def\vsx{v_{{\rm s},X}}

\def\La{{\cal L}_{\ast}}
\def\fa{f_{\ast}}
\def\Fa{F^{\ast}}
\def\Ga{G^{\ast}}
\def\Ma{M_{\ast}}

\def\qa{q_{\ast}}
\def\ua{u_{\ast}}
\def\va{v_{\ast}}
\def\vax{v_{{\ast},X}}
\def\vaxx{v_{{\ast},XX}}

\def\vb{v_{\rm b}}
\def\vbx{v_{{\rm b},X}}
\def\vu{v_{\rm u}}

\def\Xh{X_{\rm h}}
\def\Yh{Y_{\rm h}}

\def\che{c_{\rm het}}

\def\qhe{q_{\rm het}}

\def\Uhe{U_{\rm het}}
\def\vhe{v_{\rm het}}
\def\Vhe{V_{\rm het}}
\def\lahe{\lambda_{\rm het}}

\def\cho{c_{\rm hom}}
\def\pho{p_{\rm hom}}
\def\qho{q_{\rm hom}}
\def\uho{u_{\rm hom}}
\def\uhox{u_{{\rm hom},X}}
\def\Uho{U_{\rm hom}}
\def\vho{v_{\rm hom}}
\def\vhox{v_{{\rm hom},X}}
\def\Vho{V_{\rm hom}}
\def\laho{\lambda_{\rm hom}}

\def\hI{\hat{I}}
\def\hJ{\hat{J}}
\def\hN{\hat{N}}

\def\hu{\hat{u}}
\def\hv{\hat{v}}

\def\stackdef{\stackrel{\rm def}{=}}

\begin{document}

\title{Slow localized patterns in singularly perturbed 2-component reaction-diffusion equations}

\author{Arjen Doelman\footnotemark[1]}

\maketitle
\renewcommand{\thefootnote}{\fnsymbol{footnote}}
\footnotetext[1]{Mathematisch Instituut, Universiteit Leiden, the Netherlands, doelman@math.leidenuniv.nl.}
\renewcommand{\thefootnote}{\arabic{footnote}}

\begin{abstract}
Localized patterns in singularly perturbed reaction-diffusion equations typically consist of slow parts -- in which the associated solution follows an orbit on a slow manifold in a reduced spatial dynamical system -- alternated by fast excursions -- in which the solution jumps from one slow manifold to another, or back to the original slow manifold. In this paper we consider the existence and stability of stationary and traveling slow localized patterns that do not exhibit such jumps, i.e. that are completely embedded in a slow manifold of the singularly perturbed spatial dynamical system. These patterns have rarely been considered in the literature, for two reasons: (i) in the classical Gray-Scott/Gierer-Meinhardt type models that dominate the literature, the flow on the slow manifold is typically linear and thus cannot exhibit homoclinic pulse or heteroclinic front solutions; (ii) the slow manifolds occurring in the literature are typically `vertical' -- i.e. given by $u \equiv u_0$, where $u$ is the fast variable -- so that the stability problem is determined by a simple (decoupled) scalar equation. The present research concerns a general system of singularly perturbed reaction-diffusion equations and is motivated by several explicit ecosystem models that do give rise to non-vertical normally hyperbolic slow manifolds on which the flow may exhibit both homoclinic and heteroclinic orbits -- that correspond to either stationary or traveling localized patterns. The associated spectral stability problems are at leading order given by a nonlinear, but scalar, eigenvalue problem with Sturm-Liouville characteristics and we establish that homoclinic pulse patterns are typically unstable, while heteroclinic fronts can either be stable or unstable. However, we also show that homoclinic pulse patterns that are asymptotically close to a heteroclinic cycle may be stable. This result is obtained by explicitly determining the leading order approximations of 4 critical asymptotically small eigenvalues. By this analysis -- that involves several orders of magnitude in the small parameter -- we also obtain full control over the nature of the bifurcations -- saddle-node, Hopf, global, etc. -- that determine the existence and stability of the (stationary and/or traveling) heteroclinic fronts and/or homoclinic pulses. Finally, we show that heteroclinic orbits may correspond to stable (slow) interfaces in 2-dimensional space, while the homoclinic pulses must be unstable as localized stripes -- even when they are stable in 1 space dimension.
\end{abstract}

\section{Introduction}
\label{s:Intro}
\begin{figure}[t]
\centering
	\begin{minipage}{.45\textwidth}
		\centering
		\includegraphics[width =\linewidth]{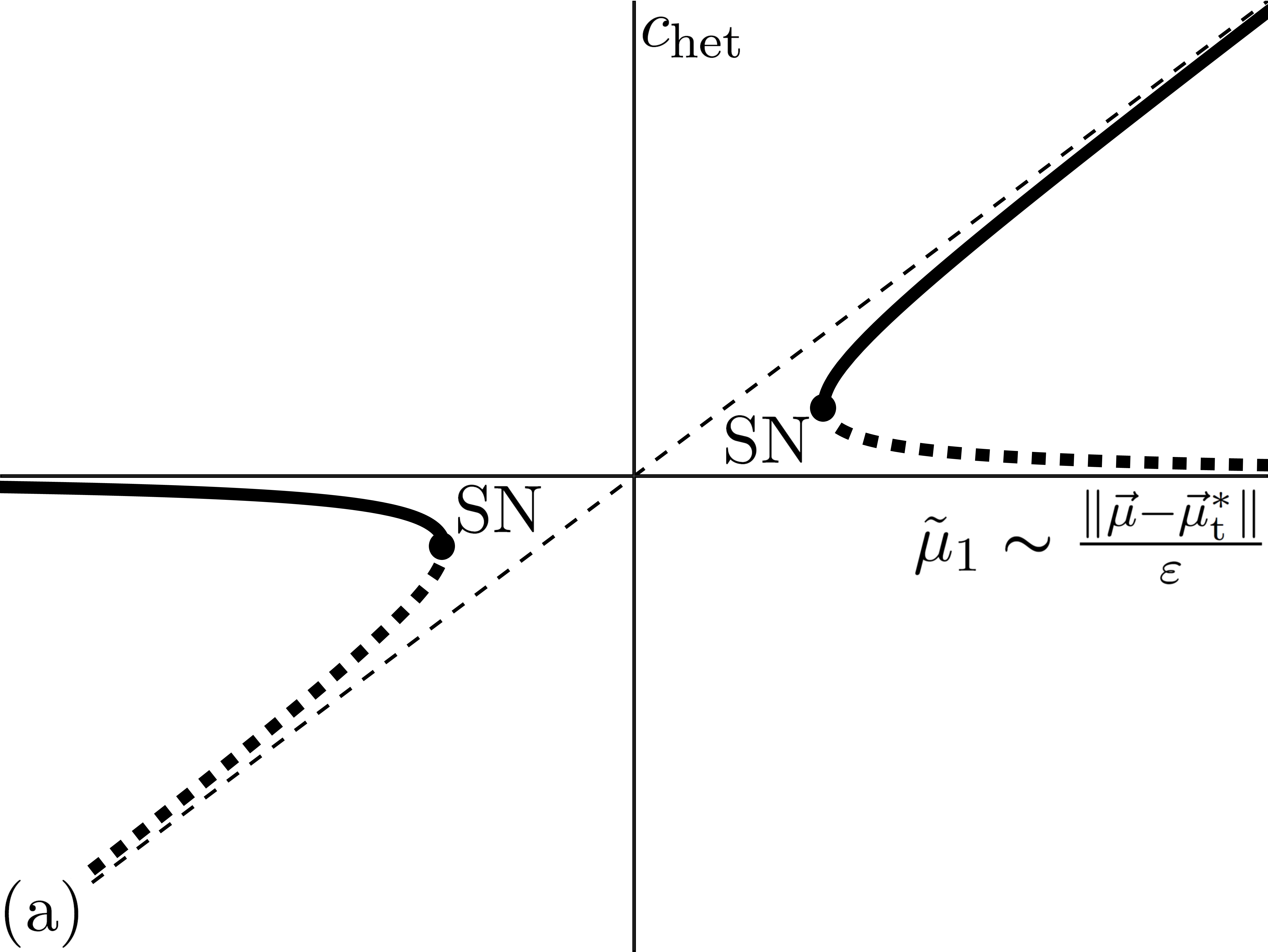}
	\end{minipage}%
	\hspace{.5cm}
	\begin{minipage}{0.45\textwidth}
		\includegraphics[width=\linewidth]{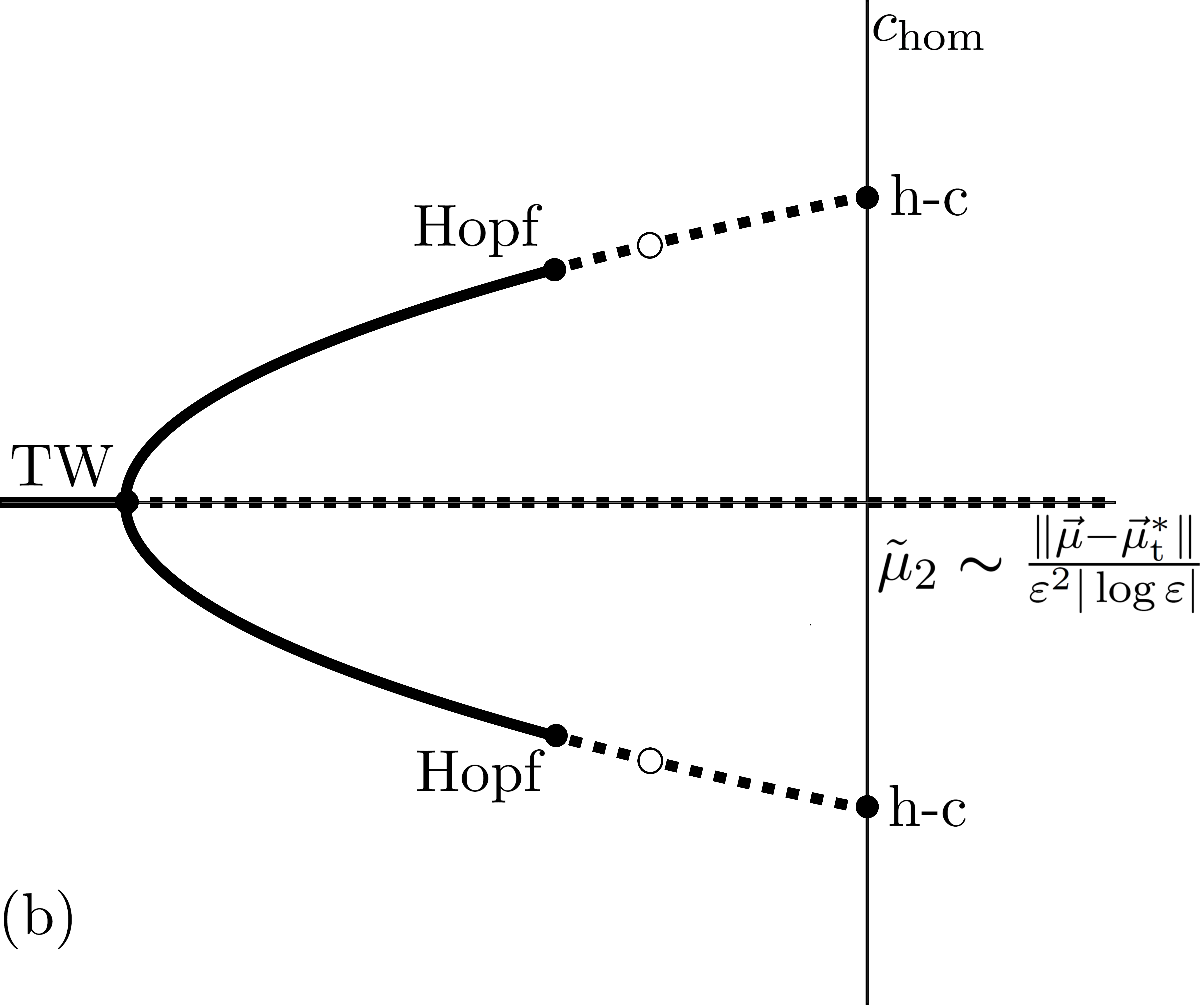}
	\end{minipage}
\caption{\small{Two examples of bifurcation diagrams associated to the existence and stability of slow localized patterns, necessarily for values of $\vmu$ such that $\|\vmu - \vmuta\|$ is asymptotically small (with $\vmuta$ defined by $\Ma(\vmuta) = 0$). (a) Heteroclinic fronts represented by their speeds $c_{\rm het}$ as function of scaled parameter $\tmu_1$, with $\tmu_1 \thicksim \|\vmu - \vmuta\|/\eps$: there are 2 saddle-node bifurcations $O(\eps)$ close to $\vmu = \vmuta$ at which a stable front merges with and unstable front. (b) Nearly heteroclinic pulses given by $c_{\rm hom} = c_{\rm hom}(\tmu_2)$ with $\tmu_2 \thicksim \|\vmu - \vmuta\|/(\eps^2 |\log \eps|)$: the stationary pulse looses its stability as a pair of stable traveling pulses bifurcates (at the point indicated by TW), these pulses destabilize by a Hopf bifurcation as $\tmu_2$ increases further and eventually split up into 2 pairs of traveling fronts (at the point indicated by h-c); the open bullet indicates the moment at which the spectrum associated to the stability of the pulse becomes real. See section \ref{sss:neartravstab} and Fig. \ref{f:hethombifs-stab} for the exact (open) conditions on the parameters for which these bifurcation diagrams occur and (\ref{d:tmu1}), (\ref{d:tmu2}) for the precise definitions of $\tmu_1$ and $\tmu_2$.}}
\label{f:HetHomBifsIntro}
\end{figure}

In this paper, we study both the existence and the stability of localized patterns in a general class of 2-component singularly perturbed reaction-diffusion equations
\begin{equation}
\label{e:RDE}
\left\{	
\begin{array}{rcrcl}
\tau U_t &=& \Delta U & + & F(U,V;\vmu)\\
V_t &=& \frac{1}{\varepsilon^2} \Delta V & + & G(U,V;\vmu)
\end{array}
\right.
\end{equation}
with $U(x,y,t), V(x,y,t): \RR^2 \times \RR^+ \to \RR$, $F(U,V;\vmu)$ and $G(U,V;\vmu)$ sufficiently smooth as function of $U$, $V$ and $\vmu$ and $\tau > 0$ and $\vmu \in \RR^m$, $m \geq 2$ parameters. Localized solutions in systems of this type have been studied extensively in the mathematical literature, however, typically in the setting of explicit models -- usually of Gray-Scott or Gierer-Meinhardt type, see \cite{CW09,DV15,KWW09,SD17,Ward18} and the references therein. Moreover, to our knowledge, the kind of patterns considered here -- slow patterns -- have not been considered, for the very good reason that these patterns do not occur in the standard models in the literature. However, this is an artifact of these models, the present study is partly motivated by explicit ecosystem models within which slow localized patterns may occur -- see section \ref{ss:Eco}.
\\
\begin{figure}[t]
\label{f:M0genGS}
\centering
	\begin{minipage}{.48\textwidth}
		\centering
		\includegraphics[width =\linewidth]{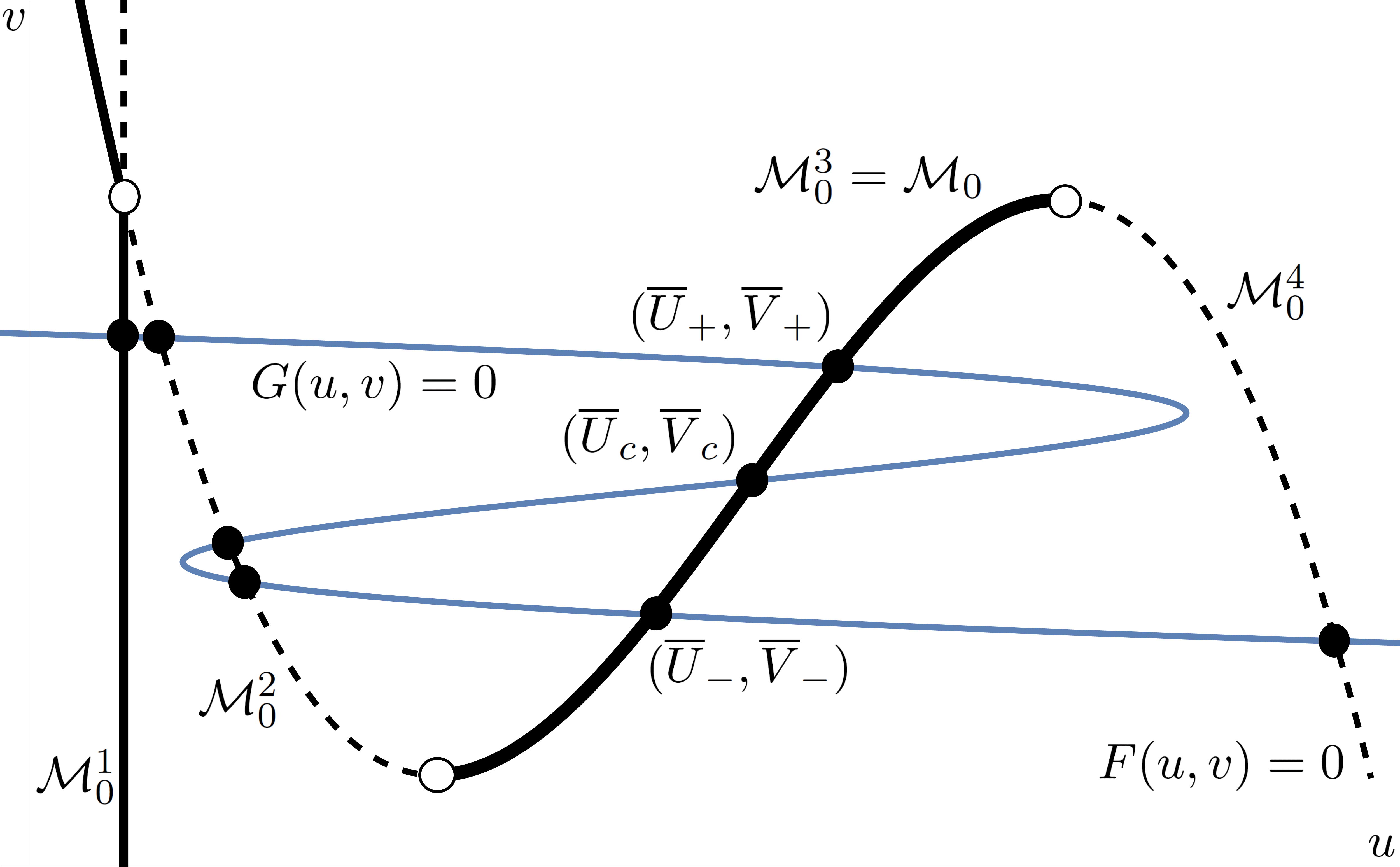}
	\end{minipage}%
	\hspace{.2cm}
	\begin{minipage}{0.48\textwidth}
		\includegraphics[width=\linewidth]{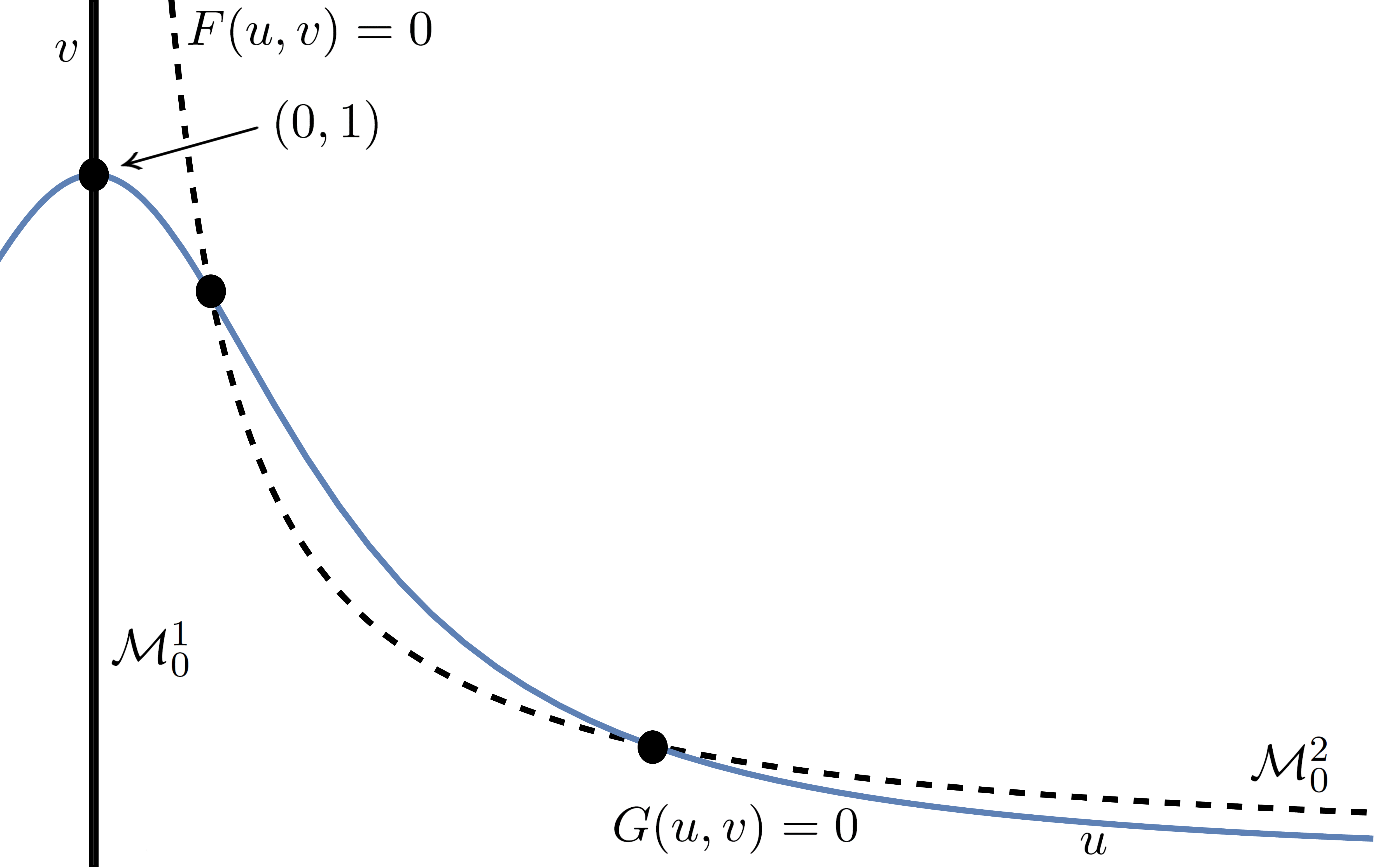}
	\end{minipage}
\caption{\small{Two examples of slow manifolds for system (\ref{e:DS}), determined by $F(u,v) = 0$, combined with critical points of (\ref{e:DS}) on these slow manifolds (i.e. the trivial background states of (\ref{e:RDE})) -- given by $F(u,v) = G(u,v) = 0$. Continuous, respectively dashed, curves indicate the normally hyperbolic, resp. elliptic, (sub)manifolds $\M^j_0$ (\ref{d:Mj0}); the open bullets represent the transition zones where normal hyperbolicity breaks down. Left: a general setting (partly based on \cite{JDCBM20,ZMB15}) with $\M_0 = \M^3_0$ such that $f'(v) \not \equiv 0$ and that the reduced flow (\ref{e:RedSF}) has 3 critical points on $\M_0$. Right: the Gray-Scott case with only a vertical normally hyperbolic slow manifold $\M^1_0$ ($F(U,V; \vmu) = -\mu_1 U + U^2V$, $G(U,V; \vmu) = \mu_2(1-V) - U^2V$ with $\mu_1, \mu_2 = \O(1)$ -- see Remark \ref{r:naturepertslowflow-1}.}}
\end{figure}
\\
We consider the most simple patterns potentially exhibited by (\ref{e:RDE}): uniformly traveling waves that do not depend on the $y$ direction. Naturally, by going into a traveling frame -- and thus by introducing $\xi = x - ct$, $U(x,y,t) = u(\xi)$, $V(x,y,t) = v(\xi)$ -- (\ref{e:RDE}) is reduced to a system of coupled second order ODEs in $\xi$,
\begin{equation}
\label{e:ODE}
\left\{	
\begin{array}{rclcccc}
u_{\xi\xi} & + & c \tau u_\xi & + & F(u,v) & = & 0\\
\frac{1}{\eps^2}v_{\xi\xi} & + & c v_\xi & + & G(u,v) & = & 0
\end{array}
\right.
\end{equation}
where we have chosen the magnitude (in orders of $\eps$) of the speed of the traveling wave such that it appears as leading order term in the (fast) $u$-equation: this is the most natural scaling in which patterns may `jump' from one slow manifold of singularly perturbed system (\ref{e:ODE}) to another, while the level (in the slow variable $v$) at which this jump occurs is directly controlled by $c$ -- as is for instance the case in the classical FitzHugh-Nagumo setting \cite{Jon84}, but also in the more recent paper \cite{JDCBM20} that (partly) inspired the present research. System (\ref{e:ODE}) can be written as a 4-dimensional spatial dynamical system, that reads in its slow form,
\begin{equation}
\label{e:DS}
\begin{array}{rclcrcl}
\eps u_X & = & p & \hspace{0.5cm} & v_X & = & q\\
\eps p_X & = & - F(u,v) - c \tau p & \hspace{0.5cm} & q_X & = & - G(u,v)- \eps c q
\end{array}
\end{equation}
in $X = \eps \xi$. By taking the $\eps \to 0$ limit, we find that the 2-dimensional (reduced) slow manifolds are determined by  $F(u,v_0) = 0$ (and $p =0$, $(v_0,q_0) \in \RR^2$). In general, this determines $J \geq 1$ branches, locally given by graphs,
\beq
\label{d:Mj0}
\M^j_0 = \{(u,p,v,q) \in \mathbb{R}^4: u=f^j(v), p=0\}, \; j = 1,2, ..., J,
\eeq
with $f_j(v)$ such that $F(f^j(v),v) \equiv 0$ -- see Fig. \ref{f:M0genGS}. For those (parts of) $\M^j_0$ that are normally hyperbolic, $\M^j_0$ persists as $\M^j_\eps$, with $\M^j_\eps$ typically depending on $c$: $\M^j_\eps = \M^j_\eps(c)$. For notational convenience, we drop the $j$-dependence of $\M_{0,\eps}^j$ and $f^j(v)$ in the forthcoming analysis and explicitly note that $F(f(v),v) = 0$ implies that,
\beq
\label{e:f'}
F_u(f(v),v) f'(v) + F_v(f(v),v) = 0: f'(v) = - \frac{F_v(f(v),v)}{F_u(f(v),v)}
\eeq
Naturally, the reduced slow flow on $\M_0$ is given by,
\beq
\label{e:RedSF}
v_{XX} + G(f(v),v) = 0
\eeq
In this paper, we study the existence and stability of bounded solutions $(\us(X,\ps(X),\vs(X),\qs(X)) \subset \M_\eps(c)$ of (\ref{e:DS}) that limit (as $\eps \to 0$) on orbits $(f(v_0(X)),0,v_0(X),v_{0,X}(X)) \subset \M_0$ in which $v_0(X)$ is a bounded solution of (\ref{e:RedSF}). Thus, the corresponding patterns $(U(x,y,t),V(x,y,t)) = (\Us(X),\Vs(X))$ in (\ref{e:RDE}) -- with $(\Us(X),\Vs(X))$ at leading order given by $(f(v_0(X)),v_0(X))$ -- are indeed slowly varying in the $x$-direction (and trivial in the $y$-direction): unlike in essence all nontrivial patterns exhibited by singularly perturbed reaction-diffusion systems studied in the literature, the patterns we consider do not consist of alternating slow parts and fast jumps, $(\us(X,\ps(X),\vs(X),\qs(X))$ never takes off from $\M_\eps(c)$, hence the terminology `slow patterns' -- see Remark \ref{r:literaturestableslowpatterns}.
\\ \\
It was already noted in \cite{DV15} that the slow reduced flow for the models studied in the literature almost without exception is linear in $v$ -- except for FitzHugh-Nagumo type models that have 1-dimensional slow manifolds, see \cite{Jon84,CRS16,CS19} and Remark \ref{r:literaturestableslowpatterns}. Naturally, this is a severe limitation and hence the terms `slowly linear' -- for Gray-Scott and Gierer-Meinhardt type models -- and `slowly nonlinear' -- for the more general class of models -- were coined in \cite{DV15}. Moreover, the reduced slow flows associated to the Gray-Scott/Giere-Meinhartd type models considered in the literature -- see \cite{CW09,DV15,KWW09,SD17,Ward18} and the references therein -- do not have bounded solutions $v_0(X)$ (except for critical points) due to their linearity: Gray-Scott/Gierer-Meinhardt type models cannot exhibit the patterns considered here (see however Remark \ref{r:Champneys}). Moreover, even in the slowly nonlinear class of model introduced in \cite{DV15}, slow patterns will typically not be observed, since it is assumed in \cite{DV15} that the normally hyperbolic manifold $\M_0$ is vertical in the $(u,v)$-plane, i.e. that $f(v) \equiv 0$ and thus that $f'(v) \equiv 0$ (like in Gray-Scott/Gierer-Meinhardt cases, see Fig. \ref{f:M0genGS}) -- as we shall argue below.
\\
\begin{figure}[t]
\centering
	\begin{minipage}{.45\textwidth}
		\centering
		\includegraphics[width =\linewidth]{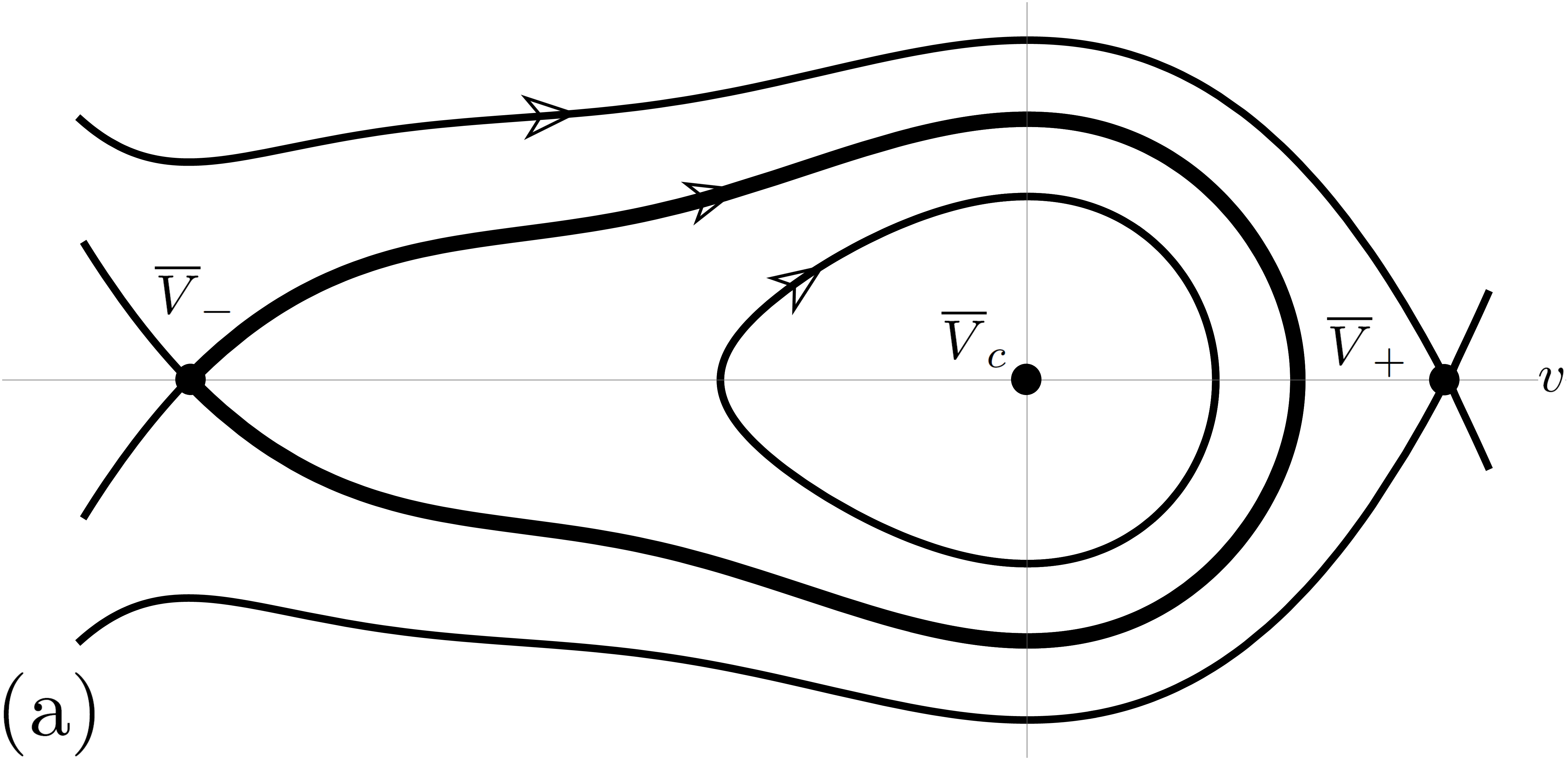}
	\end{minipage}%
	\hspace{.2cm}
	\begin{minipage}{0.45\textwidth}
		\includegraphics[width=\linewidth]{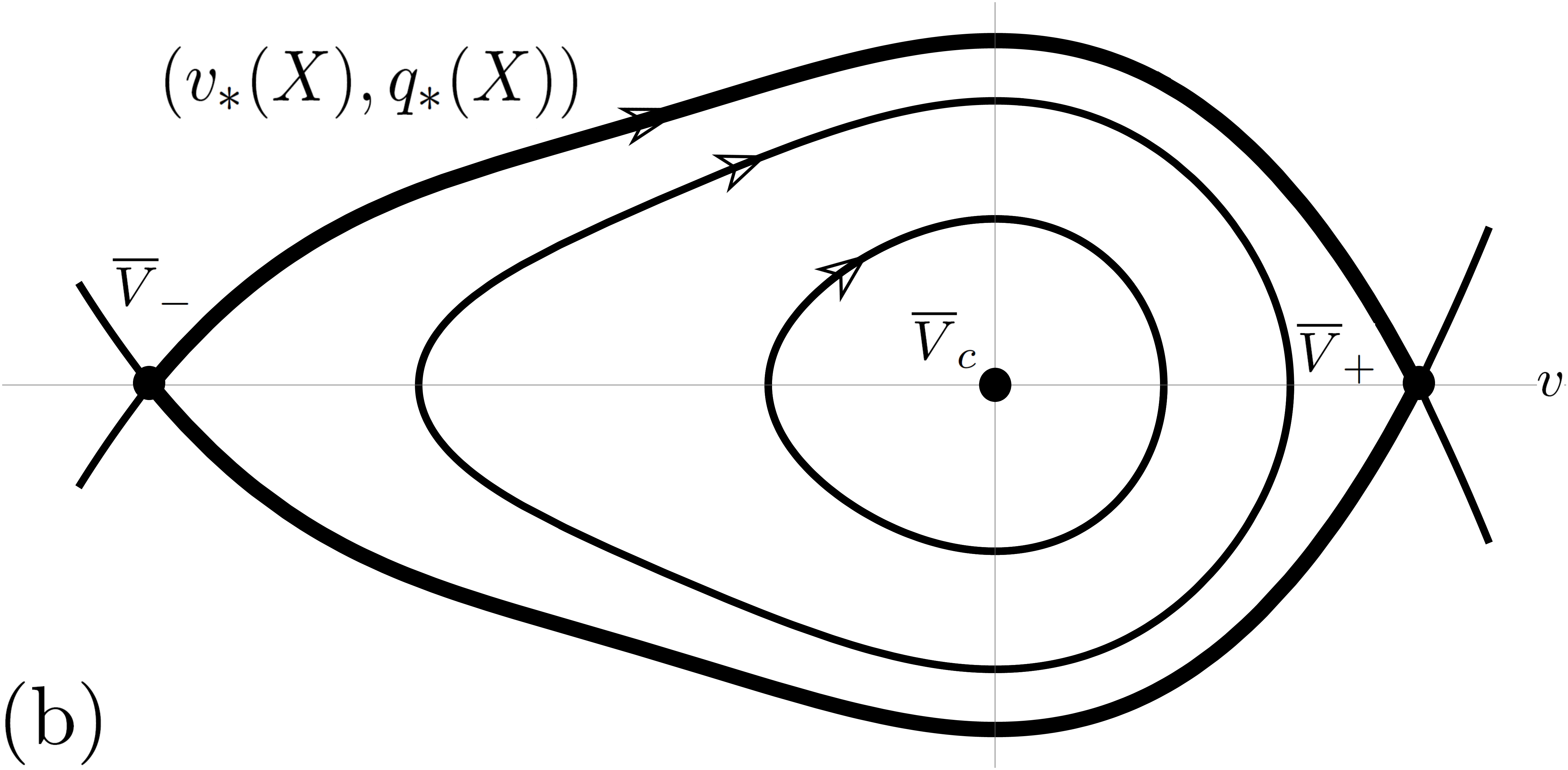}
	\end{minipage}
\caption{\small{The slow flows on $\M_0$ given by (\ref{e:RedSF}) under the assumption that the potential $\W_0(v)$ (\ref{d:HW0}) is of double well type with minima at $v = \oV_\pm$ and maximum at $\oV_c$. (a) The case of unequal wells: $\W_0(\oV_-) > \W_0(\oV_+)$. (b) The equal well case, i.e. $\W_0(\oV_-) = \W_0(\oV_+)$, with heteroclinic orbit $(\va(X),\qa(X))$ connecting the saddle $(\oV_-,0)$ to $(\oV_+,0)$.}}
\label{f:SlowFlows}
\end{figure}
\\
A priori, one may thus question the relevance of the slow patterns considered here. We claim however that slow patterns are both relevant and interesting and that the fact that these patterns have not been considered yet in the quite extended literature on singularly perturbed reaction-diffusion equations mostly is an indication of the restricted `spectrum' of models considered so far. Based on the dryland and savanna ecosystem models (\ref{e:RDE-exZelnik}) and (\ref{e:RDE-exvLangevelde}) introduced in section \ref{ss:Eco} -- and on general considerations -- we assume throughout this paper -- also for simplicity -- that the reduced slow flow on $\M_0$ has 3 critical points: two (non-degenerate) saddle points $(\oV_\pm,0)$ and a center $(\oV_c,0)$ with $\oV_- < \oV_c < \oV_+$. Thus, we assume that the unperturbed integral associated to (\ref{e:RedSF})
\beq
\label{d:HW0}
\H_0(v,q) = \frac12 q^2 + \int_{\oV_-}^v G(f(\tv),\tv) \, d \tv \stackdef \frac12 q^2 - \W_0(v)
\eeq
has a potential $\W_0(v)$ of double well type, with 2 (non-degenerate) minima at $v = \oV_\pm$ and a (local, non-degenerate) maximum at $v = \oV_c$. This is the most simple setting in which (\ref{e:RedSF}) may have localized solutions of both homoclinic pulse/stripe and of heteroclinic front/interface type -- see Fig. \ref{f:SlowFlows}. Moreover, we show in section \ref{ss:Eco} that there are open regions in parameter space for which both ecosystem models (\ref{e:RDE-exZelnik}) and (\ref{e:RDE-exvLangevelde}) have normally hyperbolic slow manifolds with reduced slow flows of double well type.
\\
\begin{figure}[t]
\centering
	\begin{minipage}{.45\textwidth}
		\centering
		\includegraphics[width =\linewidth]{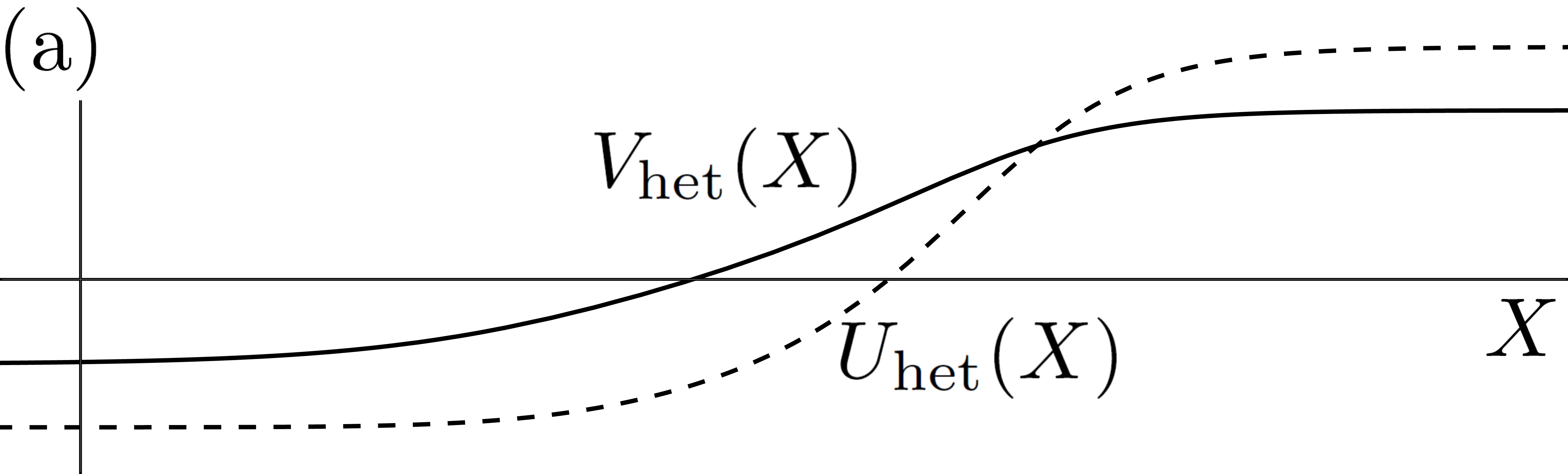}
	\end{minipage}%
	\hspace{.2cm}
	\begin{minipage}{0.45\textwidth}
		\includegraphics[width=\linewidth]{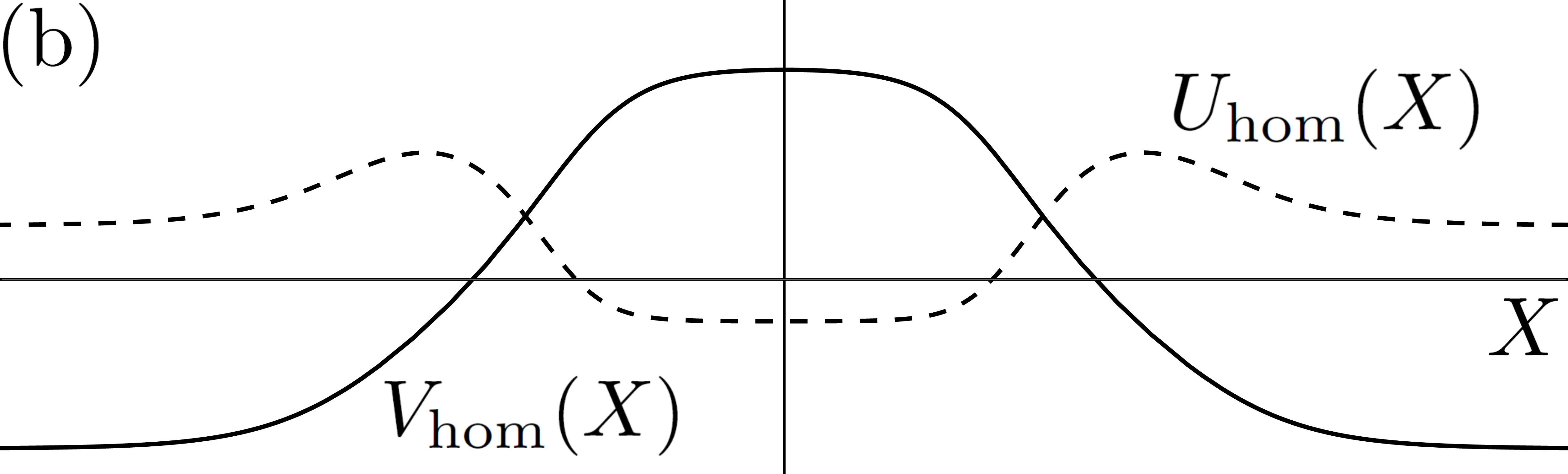}
	\end{minipage}
\caption{\small{Example of a slow heteroclinic front $(\Uhe(X),\Vhe(X))$ (a) and of a slow, stationary, symmetric, nearly double front, homoclinic pulse $(\Uho(X),\Vho(X))$ (b) -- both as 1-dimensional localized solutions of system (\ref{e:RDE})/(\ref{e:RDE-S}) (i.e. for $x/X \in \RR$). For the front, $f(v)$ is chosen to be monotonous ($f'(v) > 0$, as in Fig. \ref{f:M0genGS}), $f'(v)$ changes sign in the pulse example. The front $(\Uhe(X),\Vhe(X))$ connects the trivial background state $(\oU_-,\oV_-)$ to $(\oU_+,\oV_+)$, the pulse $(\Uho(X),\Vho(X))$ is biasymptotic to $(\oU_-,\oV_-)$ and passes closely along $(\oU_+,\oV_+)$ (where $\Vho(X)$ attains its maximal value).}}
\label{f:frontspulses}
\end{figure}
\\
Although periodic patterns are at least as interesting and relevant -- certainly also from the ecological point of view -- we focus here on localized patterns $(\Us(X),\Vs(X))$ as fronts (in $\RR^1$) or interfaces (in $\RR^2$) or pulses/fronts -- Fig. \ref{f:frontspulses} --  i.e. on heteroclinc or homoclinic orbits $(\us(X,\ps(X),\vs(X),\qs(X)) \subset \M_\eps(c)$ in (\ref{e:DS}) -- except for our first persistence result (Theorem \ref{t:E-Pers}) and the discussion in section \ref{ss:Projects}. Naturally, we only consider orbits that limit on saddle points $(\oU,0_\pm,\oV_\pm,0) \in \M_\eps(c)$ of (\ref{e:DS}) that correspond to stable trivial patterns $(U(x,y,t),V(x,y,t) \equiv (\oU_\pm,\oV_\pm)$ of (\ref{e:RDE}) (and note that the patterns $(U(x,y,t),V(x,y,t) \equiv (\oU_c,\oV_c)$ are necessarily unstable \cite{Doe19}).
\\ \\
As a first indication of the impact of the `geometry' of the slow manifold $\M_\eps(c)$, i.e. of the fact that $f'(v) \not \equiv 0$, we note that the full flow on $\M_\eps(c)$ is given by,
\beq
\label{e:SF-I}
v_{XX} + G(f(v),v) + \eps\left[c v_X \tG_{1c}(v) \right] + \eps^2 \left[ c^2 \tau^2 \tG_{2cc}(v,v_X^2) + \tG_{2}(v,v_X^2) \right] = \O(\eps^3)
\eeq
with,
\beq
\label{d:tG1c-I}
\tG_{1c}(v)  = 1- \tau f'(v) \frac{G_u(f(v),v)}{F_u(f(v),v)}
\eeq
(see (\ref{d:tFs}), (\ref{d:tGs}) in section \ref{ss:E-SM}). Thus, in the vertical Gray-Scott/Gierer-Meinhardt-type case with $f'(v) \equiv 0$, the leading order perturbation is a linear friction term: there cannot be traveling localized structures in this case (Fig. \ref{f:pertSlowFlows}(b)). This is very different in the general case, since here the leading order perturbation term has the nature of a nonlinear friction term that may change its sign -- see the example of Fig. \ref{f:pertSlowFlows}(a) in which the unperturbed heteroclinic orbit $(\va(X),\qa(X))$ of Fig. \ref{f:SlowFlows}(b) persists as traveling front (thus for a well-defined value of speed $c = c_{\rm het}(\vmu)$ -- see Theorem \ref{t:exhets}(ii)). On the other hand, it should be noted that this difference only exists for traveling patterns, i.e. for $c \neq 0$. The stationary patterns of (\ref{e:RedSF}) naturally persist under the perturbations of (\ref{e:SF-I}), since (\ref{e:DS}) -- and thus (\ref{e:SF-I}) -- is reversible for $c=0$, see Theorems \ref{t:E-Pers}(i), \ref{t:exhets}(i) and \ref{t:exhoms}(i).
\\
\\
To motivate another distinction between the present study and the literature, we need to briefly sketch the first steps of the (spectral) stability analysis of a slow pattern $(\Us(X),\Vs(X))$. To do so, we consider (\ref{e:RDE}) in its equivalent slow form,
\begin{equation}
\label{e:RDE-S}
\left\{	
\begin{array}{rcrcrcl}
\tau U_t &=& \eps^2 \Delta U & + & \eps c \tau U_X & + & F(U,V)\\
V_t &=& \Delta V & + & \eps c V_X & + & G(U,V)
\end{array}
\right.
\end{equation}
for $(X,Y) \in \RR^2$ with $X = \eps(x-ct)$ and $Y = \eps y$. Since the existence problem could -- at leading order -- be reduced to the simple (integrable) planar system (\ref{e:RedSF}), it is tempting assume that the stability of the slow solutions is determined by the scalar equation naturally associated to (\ref{e:RedSF}) and thus to introduce
\beq
\label{e:RedPDE}
V_t = \Delta V + G(f(V),V)
\eeq
as possible (scalar) `slow reduced PDE'. However, the leading order linearization
\beq
\label{e:specstabexpO1}
(U(x,y,t), V(x,y,t)) = (u_0(X) + \bu_0(X) e^{\la_0 t + i L Y}, v_0(X) + \bv_0(X) e^{\la_0 t + i L Y})
\eeq
about a localized pattern $(\Us(X),\Vs(X)$ at leading order given by $(u_0(X),v_0(X))$ -- with $v_0(X)$ a homoclinic/heteroclinic orbit of (\ref{e:RedSF}) on $\M_0$ and $u_0(X) = f(v_0(X))$ -- is, again at leading order, given by
\begin{equation}
\label{e:ODE-lin-0}
\left\{	
\begin{array}{rcl}
\tau \la_0 \bu_0 &=& F_u(u_0,v_0) \bu_0 + F_v(u_0,v_0) \bv_0\\
\la_0 \bv_0 &=& \bv_{0,XX} - L^2 \bv_0 + G_u(u_0,v_0) \bu_0 + G_v(u_0,v_0) \bv_0
\end{array}
\right.
\end{equation}
Thus, as in the existence problem, $\bu_0(X)$ can be expressed in terms of $\bv_0(X)$,
\beq
\label{e:bu0inbv0}
\bu_0 = \frac{F_v(u_0,v_0) \bv_0}{\tau \la_0 - F_u(u_0,v_0)} = - \frac{f'(v_0) F_u (u_0,v_0)}{\tau \la_0 - F_u(u_0,v_0)} \bv_0 =
f'(v_0) \left[1 - \frac{\tau \la_0}{\tau \la_0 - F_u(u_0,v_0)} \right] \bv_0
\eeq
(\ref{e:f'}). By substitution of (\ref{e:bu0inbv0}) into the $\bv_0$-equation of (\ref{e:ODE-lin-0}), we arrive at the nonlinear eigenvalue problem
\beq
\label{e:O1EigPb}
\left[\L_s - \frac{\tau \la_0 f'(v_0) G_u(f(v_0),v_0)}{\tau \la_0 - F_u(f(v_0),v_0)} - (\la_0 + L^2) \right] \bv_0 = 0
\eeq
where we have introduced the operator $\Ls = \Ls(X)$
\beq
\label{d:Ls}
\L_s(X) = \frac{d^2}{dX^2} + \left[ f'(v_0(X)) G_u(f(v_0(X)),v_0(X)) + G_v(f(v_0(X)),v_0(X)) \right]
\eeq
Now, we notice that the operator $\L_s(X)$ by itself determines the spectral stability of the localized pattern $V(X,t) = v_0(X)$ in the scalar slow reduced PDE (\ref{e:RedPDE}). Thus, the leading order spectral stability of the (localized) slow pattern $(\Us(X),\Vs(X)$ in (\ref{e:RDE-S}) indeed can be seen as being determined by the scalar slow reduction (\ref{e:RedPDE}) if $f'(v_0(X)) G_u(f(v_0(X)),v_0(X)) \equiv 0$ in (\ref{e:O1EigPb}) and thus in particular if $f'(v) \equiv 0$.
\\
\begin{figure}[t]
\centering
	\begin{minipage}{.45\textwidth}
		\centering
		\includegraphics[width =\linewidth]{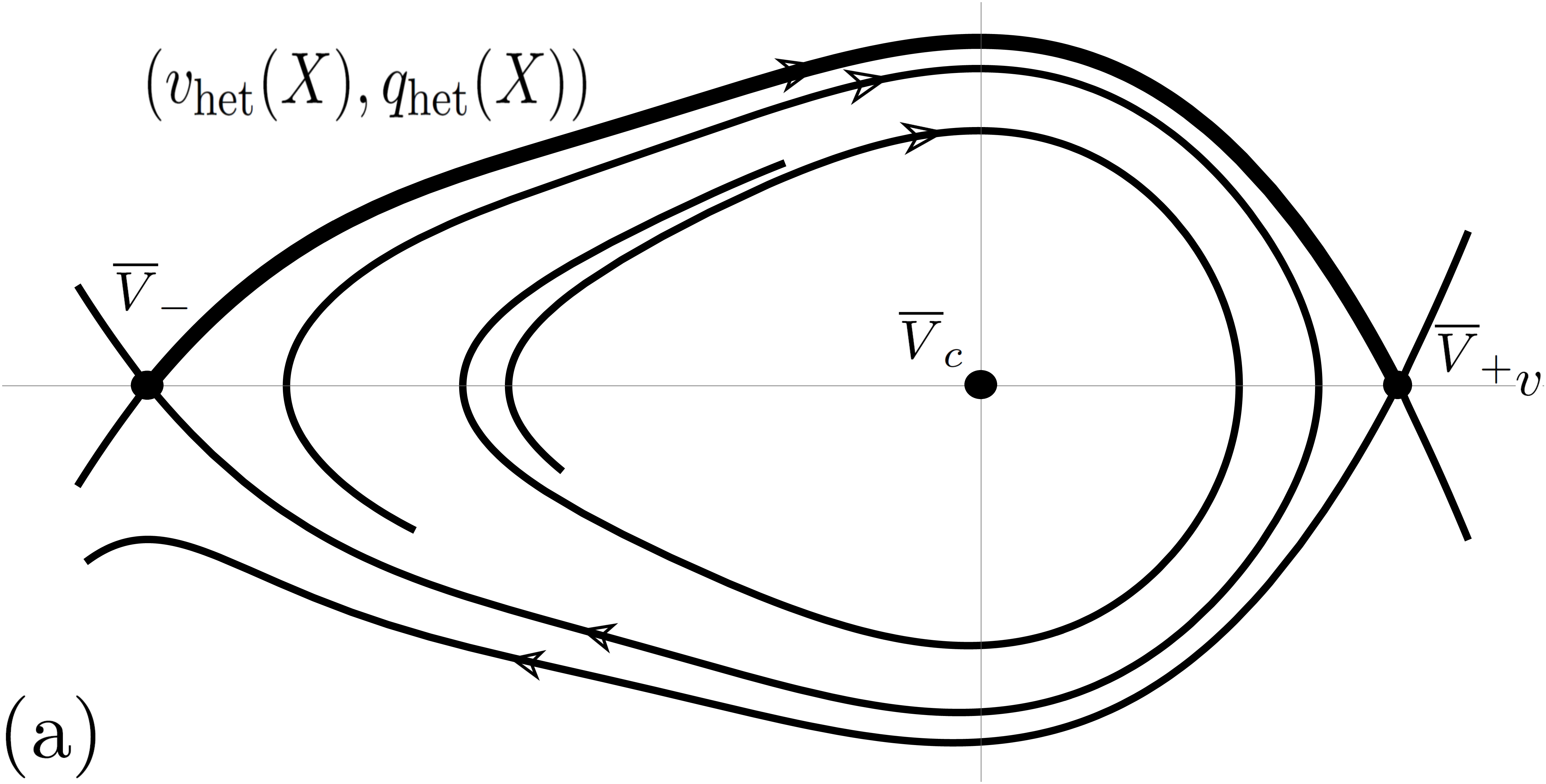}
	\end{minipage}%
	\hspace{.2cm}
	\begin{minipage}{0.45\textwidth}
		\includegraphics[width=\linewidth]{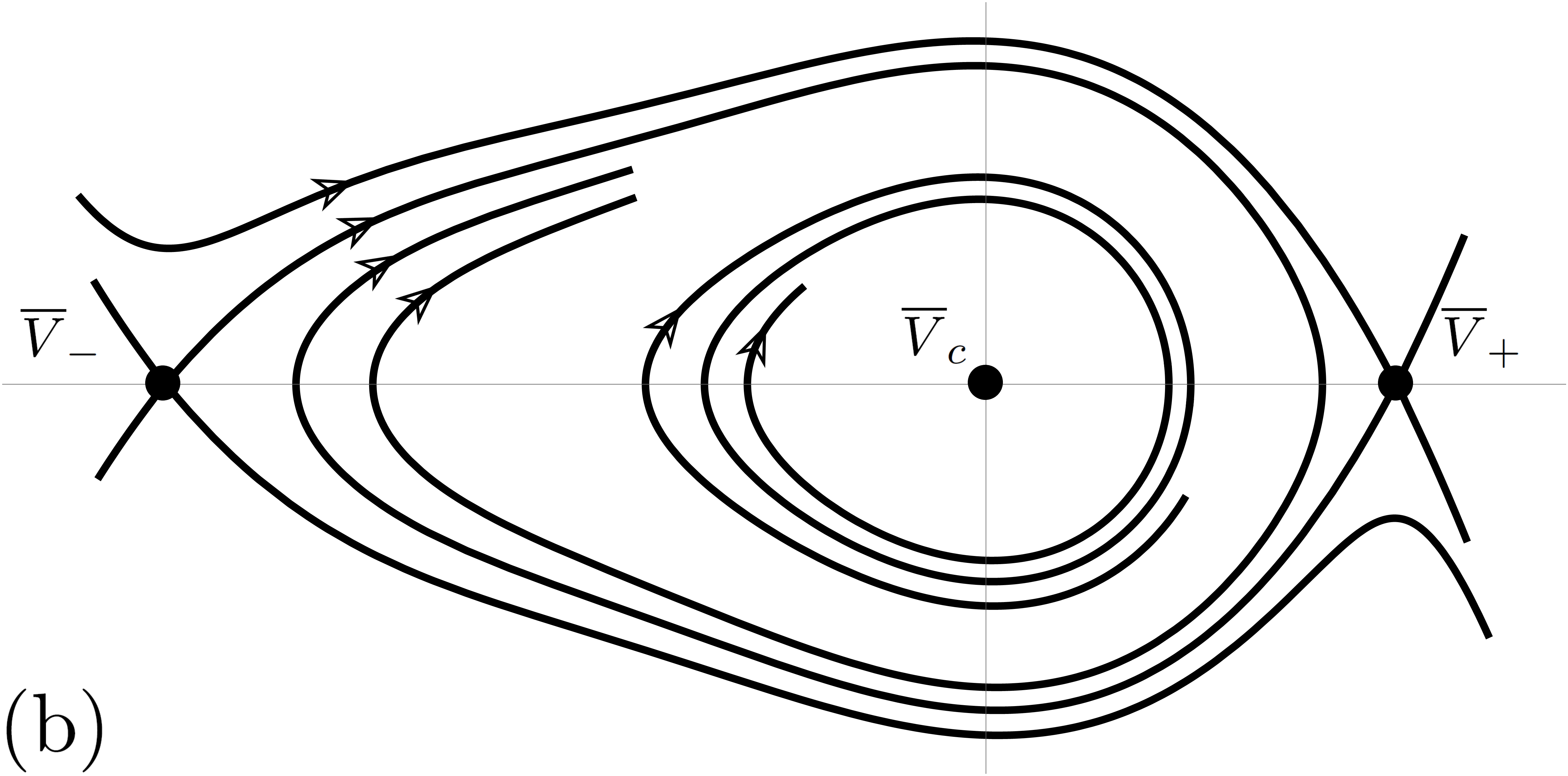}
	\end{minipage}
\caption{\small{Sketches of the flows (\ref{e:SF-I}) on the slow manifold $\M_\eps(c)$ ($c \neq 0$) in the case of (unperturbed) wells of equal depth (cf. Fig. \ref{f:SlowFlows}(b)) . (a) The general case $f'(v) \not \equiv 0$ with $c=c_{\rm het}$ such that there is a heteroclinic connection $(v_{\rm het}(X),q_{\rm het}(X))$ between the saddles $(\oV_\pm,0)$ (Theorem \ref{t:exhets}(ii)). (b) The vertical linear friction case $f'(v) \equiv 0$: there cannot be (slow) homoclinic or heteroclinic localized structures.}}
\label{f:pertSlowFlows}
\end{figure}
\\
In the (Klausmeier-)Gray-Scott/Gierer-Meinhardt type 2-component models considered in the literature, including the slowly nonlinear models of \cite{DV15,Vee15,VD13}, the slow manifolds $\M_0$ from which a slow-fast orbit jumps are given by $\{u = {\rm const.}, p = 0\}$ (typically also without any $\O(\eps)$ corrections) -- see \cite{CW09,DV15,KWW09,SD17,Ward18}  and the reference therein and the Gray-Scott example of Fig. \ref{f:M0genGS} (see however also Remark \ref{r:Champneys}). Thus, these manifolds $\M_0$ have $f'(v) \equiv 0$ and the stability of a slow pattern is at leading order indeed governed by reduced PDE (\ref{e:RedPDE}). As a consequence, the stability is (at leading order) governed by the spectral problem $[\L_s(X) - (\la_0 + L^2)] \bv_0 = 0$, which is a Sturm-Liouville problem for localized $v_0(X)$ (homoclinic or heteroclinic) and a Hill's equation for periodic $v_0(X)$. Hence, it follows in a straightforward fashion that pulse and periodic patterns must be unstable \cite{KP13,MW66}. Clearly, the situation is very different for non-vertical slow manifolds $\M_0$ for which the leading order stability of localized patterns is governed by nonlinear eigenvalue problem (\ref{e:O1EigPb}).
\\ \\
Nevertheless, the first stability result of the paper, Theorem \ref{t:InstabRegHom}, establishes that all slow homoclinic pulse solutions of (\ref{e:DS}) in which the wells of $\W_0(v)$ (\ref{d:HW0}) are of unequal depth -- that can both be stationary and traveling (Theorem \ref{t:E-Pers}) -- are unstable. However, unlike for Sturm-Liouville problems, there may be more than 1 $\O(1)$ unstable eigenvalue -- see the condition in Corollary \ref{c:1laast}(ii) under which there are at least 2 $\O(1)$ unstable eigenvalues and the sketch in Fig. \ref{f:Lambdajrho} that indicates the possibility of 6 or more. On the other hand, one can also formulate an explicit condition on $F(U,V)$ and $G(U,V)$ -- see (\ref{d:lajrhodecreases}) in Corollary \ref{c:1laast}(i) -- under which all $\O(1)$ unstable eigenvalues become asymptotically small as the wells of $\W_0(v)$ become of equal depth. Thus, in that case a localized pattern can potentially be stable through the impact of the many -- so far neglected -- asymptotically small terms. Therefore, the major part of the analysis developed in this paper considers the `nearly heteroclinic' case of an unperturbed potential $\W_0(v)$ with wells of equal depth. And thus, a central role will be played by the unperturbed heteroclinic connections $(\va(X),\qa(X))$ that exist in the reduced slow flows of such double well potentials (Figs. \ref{f:SlowFlows}(b), \ref{f:frontspulses}).
\\ \\
To establish the stability of the nearly heteroclinic homoclinic pulse/stripe and heteroclinic front/interface patterns, it will be necessary to perform an accurate and somewhat subtle asymptotic analysis. In fact, we will find that the stability is only settled at the $\O(\eps^2|\log \eps|)$ level for the pulses, where we note that the $|\log \eps|$ factor -- that comes up throughout the existence and stability analysis of the nearly heteroclinic pulses -- originates from the fact that the nearly heteroclinic homoclinic orbit passes asymptotically close to the saddle $(\oV_+,0) \in \M_\eps(c)$ and that such a passage takes $\O(|\log \eps|)$ `time' (Fig. \ref{f:frontspulses}(b)).
\\ \\
As is usual in these types of asymptotic stability problems, the existence analysis also needs to go into deeper asymptotic `details' than in the more straightforward case of Theorem \ref{t:E-Pers} (that settles the persistence of stationary pulses and an associated bifurcation into traveling pulses for a potential well $\W_0(v)$ with wells of non-equal depth. In this analysis, a central role is played by the Melnikov expression,
\beq
\label{d:Ma}
M_{\ast}(\vmu) = \int_{-\infty}^{\infty} \left[1- \tau f'(\va(\tX)) \frac{G_u(f(\va(\tX)),\va(\tX))}{F_u(f(\va(\tX)),\va(\tX))} \right] \vax^2(\tX) \, d \tX
\eeq
(cf. (\ref{d:tG1c-I})) and the (co-dimension 1 manifold of) critical value(s) $\vmu = \vmuta$ for which $M_{\ast}(\vmu) =0$ (note that $M_{\ast}(\vmu) > 0$ in the vertical case $f'(v) \equiv 0$). Naturally, the analysis should be straightforward, since it only concerns the perturbed integrable planar `nonlinear oscillator' (\ref{e:SF-I}). However, it turns out to be more subtle than expected, especially since the, a priori unknown, speed of a traveling pattern appears as a pre-factor in the leading order perturbation of (\ref{e:SF-I}). As a consequence, traveling fronts and pulses appear in `vertical' bifurcations, in the sense that at leading order, these fronts and pulses exist for any $c$ (of $\O(1)$) exactly at $\vmu = \vmuta$. This is of course only a leading order effect, a higher order analysis -- incorporating the $\O(\eps^2)$ terms of (\ref{e:SF-I}) (from which $c$ cannot be factored out) -- reveals that the (leading order) bifurcation curves have the expected local quadratic shapes. For the traveling heteroclinic fronts, the unfolding takes place $\O(\eps)$ close to $\vmuta$ -- see Fig. \ref{f:HetHomBifsIntro}(a) for an example and Theorem \ref{t:exhets}(ii) and Corollary \ref{c:hetbifs} for the general case. The analysis is more subtle for the homoclinic pulses, for which one has to zoom in $\O(\eps^2 | \log \eps|)$ close to $\vmuta$ to establish the nature of the various bifurcations that may take place -- see for instance Fig. \ref{f:HetHomBifsIntro}(b) in which an example is presented of the bifurcation of a symmetric pair of traveling pulses from the stationary pulse. Both branches end in a global bifurcation in which the homoclinic orbit merges with a heteroclinic cycle that subsequently breaks up into distinct heteroclinic connections (for different values of $c$): the traveling pulse splits up into a pair of traveling fronts that travel with (initially slightly) different speeds -- see Corollary \ref{c:hombifs}.
\\
\\
The main challenge of the spectral stability analysis is to develop an approach by which the discrete spectrum associated to the fronts and pulses against 1-dimensional perturbations can be approximated -- i.e. setting $L=0$ in (\ref{e:O1EigPb}) and in its higher order corrections. By the reversibility symmetry of (\ref{e:RDE}), this is a relatively direct enterprise for the stationary fronts and pulses with $\Ma(\vmu) \neq 0$ and $\O(1)$ with respect to $\eps$. In fact, guided by a Sturm-Liouville-based intuition, one may expect that the pulses with, additional to the translational eigenvalue $\la = 0$, one asymptotically small unstable eigenvalue are unstable and the fronts, with only to the translational eigenvalue $\la = 0$ asymptotically close to the origin, stable. This is indeed the case when $\Ma(\vmu) > 0$ -- although one needs to impose extra conditions to control potential additional $\O(1)$ unstable, perhaps even complex valued, eigenvalues, see Theorem \ref{t:nearstabstat}(i-b). Naturally, this is due to the fact that the stability problem is not Sturm-Liouville, which is immediately apparent when $\Ma(\vmu) < 0$ and $\O(1)$, i.e. the situation that only occurs for non-vertical slow manifolds $\M_\eps$ (\ref{d:Ma}). In that case, it is established that also the standing fronts are unstable and that the stationary nearly heteroclinic pulses have $\O(1)$ unstable eigenvalues (Theorem \ref{t:nearstabstat}).
\\ \\
The situation is more complex for the (bifurcating) traveling fronts and pulses, i.e. for $\vmu$ asymptotically close to $\vmuta$ and thus for $|\Ma(\vmu)| \ll 1$. In this analysis, the higher order Melnikov-type expression $\hN_{2 \la \la}(\vmuta)$ -- see upcoming equation (\ref{d:hN2ll}) -- plays a decisive role. Apart from the translational eigenvalue $\la^1_{\rm het} = 0$, the spectrum associated to the traveling fronts has a second eigenvalue near $0$, $\la^2_{\rm het} = \O(\eps)$, it can be determined by an expression that depends on the speed $c_{\rm het}$ of the front and that has a factor $1/\hN_{2 \la \la}(\vmuta)$ (Lemma \ref{l:asymptevhet}). Naturally, the sign of $\la^2_{\rm het}$ -- and thus that of $\hN_{2 \la \la}(\vmuta)$ -- determines the stability of the traveling front, see Fig. \ref{f:HetHomBifsIntro}(a) for an example and Theorem \ref{t:nearstabtravfrontspulses}(i) for the general case.
\\ \\
The traveling pulses have 4 asymptotically small eigenvalues that vary as function of their speed $c_{\rm hom} = c_0 + \O(\eps)$ -- that is determined by the parameter(s) $\vmu$ (Corollary \ref{c:hombifs}, Fig. \ref{f:HetHomBifsIntro}(b)). In (the proof of) Lemma \ref{l:asymptevhom}, we develop a scheme by which explicit leading order expressions can be determined: $\la^{1,1}_{\rm hom}(c_0) \equiv 0$, $|\la^{1,2}_{\rm hom}(c_0)| \ll \eps$, $\la^{1,\pm}_{\rm hom} (c_0) = \O(\eps)$ and $\la^{1,-}_{\rm hom}(c_0) = - \la^{1,+}_{\rm hom} (c_0)$ at leading order in $\eps$. It follows that if $\hN_{2 \la \la}(\vmuta) < 0$, then $\la^{1,\pm}_{\rm hom}(c_0) \in \RR$ so that the pulses must be unstable. However if $\hN_{2 \la \la}(\vmuta) > 0$,  then the leading order expressions for $\la^{1,\pm}_{\rm hom}(c_0)$ are purely imaginary, i.e. $|{\rm Im}(\la^{1,\pm}_{\rm hom}(c_0))| = \O(\eps)$ and $|{\rm Re}(\la^{1,\pm}_{\rm hom}(c_0))| \ll \eps$, for open regions of $c_0$-values: the stability of the traveling pulses is decided by the next order terms in the approximations of $\la^{1,2}_{\rm hom}(c_0)$ and $\la^{1,\pm}_{\rm hom}(c_0)$, where
\[
\la^{1,2}_{\rm hom}(c_0) = \O(\eps^2 |\log \eps|), \; \; {\rm Re}(\la^{1,+}_{\rm hom}(c_0)) = {\rm Re}(\la^{1,-}_{\rm hom}(c_0)) = \O(\eps^2 |\log \eps|)
\]
are determined in Lemma \ref{l:asymptevhom-hot}. Theorem \ref{t:nearstabtravfrontspulses}(ii) gives a precise description of the (open) conditions under which stable traveling pulses exist. In fact, by considering $c_0 \to 0$, a similar, but much simpler, stability result is obtained for the standing pulses for $\vmu$ $\O(\eps^2 |\log \eps|)$ close to $\vmuta$, naturally also under the condition that $\hN_{2 \la \la}(\vmuta) > 0$ (Corollary \ref{c:stabstandfrontpulseMsmall}); see Fig. \ref{f:HetHomBifsIntro}(b) for an example of the outcome of the combined existence and (spectral) stability analysis.
\\ \\
The stability of the fronts and pulses as interfaces and (localized) stripes for $(X,Y) \in \RR^2$ follows along similar lines. In scalar equations as (\ref{e:RedPDE}) for which the stability of localized patterns is determined by a Sturm-Liouville operator, the $Y$-independent perturbations are `the most harmful', so that a stable front on $\RR$ `automatically' is stable as (trivially extended) interface on $\RR^2$. This is also the case for the stable stationary fronts and for the traveling fronts that have $\hN_{2 \la \la}(\vmuta) < 0$ (Theorem \ref{t:stabstatinterfaces}). However, traveling fronts with $\hN_{2 \la \la}(\vmuta) > 0$ are destabilized by $Y$-dependent perturbations (of side band type) and thus cannot be stable as traveling interfaces on $\RR^2$ (Theorem \ref{t:stabstatinterfaces}(ii)). Since the stability results for pulses against perturbations that only depend on $X \in \RR$  require that $\hN_{2 \la \la}(\vmuta) > 0$, this same mechanism implies that there cannot be stable homoclinic stripes (Theorem \ref{t:unstabtravstripes}).
\\ \\
Finally we remark that although the theme of this paper is the existence -- section \ref{s:Ex} -- and spectral stability -- section \ref{s:Stab} -- of localized patterns in (\ref{e:RDE})/(\ref{e:RDE-S}) that remain on the slow manifold $\M_\eps$ for all $X \in \RR$ (as solutions of spatial dynamical system (\ref{e:DS})), the present research only is a first step towards understanding the relevance of slow patterns on the dynamics of (\ref{e:RDE})/(\ref{e:RDE-S}). Systems that have several normally hyperbolic slow manifolds such as the dryland ecosystem model (\ref{e:RDE-exZelnik}) or that have a `return mechanism' so that orbits may take off and subsequently touch down again on the slow manifold -- such as the Gray-Scott/Gierer-Meinhardt type models -- typically exhibit slow-fast patterns. Such patterns are only (exponentially) close to a slow pattern for $X$ in a half-space -- $X \in (X_{\rm down}, \infty)$ or $X \in (-\infty, X_{\rm off})$ -- in a bounded interval -- $X \in (X_{\rm down}, X_{\rm off})$ -- or during several passages along $\M_\eps$ -- $X \in \cup_{k=1}^{K}(X^k_{\rm down}, X^k_{\rm off})$. In the latter two cases, the stability characteristics of the slow orbits on $\M_\eps$ traced by the full pattern are perhaps not always relevant (in fact, these orbits typically are not even bounded, cf. the classical Gray-Scott/Gierer-Meinhardt type patterns of \cite{CW09,DGK98,DGK01,IWW01,KWW09,SD17,WW03}), but in the former case the stability of the traced slow pattern will matter. Moreover, if the orbit on $\M_\eps$ that is traced is a closed orbit, there may be large families of slow-fast patterns that touch down on the orbit in distinct manners -- see for instance the families of traveling front patterns constructed in \cite{JDCBM20} that all touch down on a homoclinic orbit on $\M_\eps$: if the full orbit on $\M_\eps$ is stable, then (many elements in) this family may be stable (and thus relevant for the dynamics of the modeled (eco)system). Moreover, in the alternative case that there is no jump or return mechanism -- such as in savanna ecosystem model (\ref{e:RDE-exvLangevelde}) or in the class of models constructed in Remark \ref{r:larhosin} -- the present analysis only scratched the surface: the slow patterns most likely play a central role in this case, but our choice to scale their speed as in (\ref{e:ODE}) is very restrictive (for instance, it is natural to expect fronts/interfaces that travel with speeds $|c_0| \gg 1$ -- thereby extending the region in $\vmu$-space in which these patterns occur -- also for reduced slow flows (\ref{e:RedSF}) with double well potentials $\W_0(v)$ (\ref{d:HW0}) with wells of unequal depth). In section \ref{s:Disc} we will discuss several aspects of the research -- spatially periodic patterns, bifurcations, slow-fast orbits -- that may build on the present methods and results.

\begin{remark}
\label{r:naturepertslowflow-1}
\rm 
The reaction terms $F(U,V;\vmu)$ and $G(U,V;\vmu)$ typically do depend explicitly on $\eps$ in the prototypical singularly perturbed 2-component reaction-diffusion models considered in the literature -- such as the (Klausmeier-)Gray-Scott and (generalized) Gierer-Meinhardt models -- and that this typically has a major impact in the set-up of the (asymptotic) analysis (see for instance \cite{CW09,DGK01,KWW09,SD17,Ward18,WW03} and the references therein). Nevertheless, we have for simplicity chosen to not take this into account -- it would yield `additional layers of technicalities' in the upcoming, already quite technical, analysis -- see also Remark \ref{r:naturepertslowflow-2} and discussion section \ref{ss:Projects}. Moreover, it is a natural assumption from the ecological point of view -- see section \ref{ss:Eco}. However, to built a bridge between the present analysis of slow localized patterns and the extended literature of slow/fast localized patterns in Gray-Scott/Gierer-Meinhardt type problems, it may be natural to introduce an explicit $\eps$-dependence in (especially) $G(U,V)$ -- see section \ref{ss:Projects}.
\end{remark}

\begin{remark}
\label{r:literaturestableslowpatterns}
\rm
In the literature on localized patterns in singularly perturbed reaction-diffusion equations, the emphasis has been on patterns with a slow/fast structure -- with the classical FitzHugh-Nagumo pulse(s) as most well-known and well-studied example, see \cite{Jon84} and \cite{CRS16,CS19} and the references therein for recent developments. Like the FitzHugh-Nagumo model, the Klausmeier model \cite{CD18,Kla99,Sher10} and its generalization in \cite{BCD19} is a degenerate singularly perturbed 2-component reaction-diffusion model in the sense that the diffusion coefficient of one of the components vanishes. As a consequence, the associated traveling wave problems are 3-dimensional (instead of 4-dimensional system (\ref{e:DS})), with 2 fast and 1 slow direction. Thus, in these systems, both the (typically non-vertical) slow manifolds and the reduced slow flows are 1-dimensional, which implies that there cannot be slow homoclinic pulses. In principle, such systems could have slow fronts -- if the (normally hyperbolic) slow manifold $\M_\eps$ contains 2 critical points -- whose existence follows immediately from Fenichel theory \cite{Fen79,Jon95}. We are not aware of any studies in which the stability of orbits of this type is considered (but note that it may be done with the methods developed in \cite{BCD19,CRS16}). A similar situation appears in the study of weak shocks in $\la-\omega$ systems \cite{Kap91} for which the existence problem also reduces to a 3-dimensional singularly perturbed system with 2 fast and 1 slow direction. The stability results of \cite{Kap96} on heteroclinic orbits for the Ginzburg-Landau equation, a $\la-\omega$ system, can be seen as a predecessor of the present work, although it should be remarked that the most important challenge for these systems is to establish the nonlinear stability of the slow (or weak) fronts: the spectral stability problem associated to these heteroclinic orbits has essential spectrum touching the imaginary axis (see also \cite{BNSZ14}) -- which is not the case here.
\end{remark}

\subsection{Two motivating ecological models}
\label{ss:Eco}
As direct motivation for the research in this paper, we briefly discuss two explicit singularly perturbed 2-component reaction-diffusion systems from the ecological literature. Both examples model ecosystems that are known to exhibit a rich variety of patterns -- see \cite{dKetal08,GvdVvL17,JDCBM20,Mer15,vLetal03,ZMB15} and the references therein -- and both models are of the type studied here, i.e. both have non-vertical normally hyperbolic slow manifolds while there is an open region in parameter space for which the reduced slow flow on the slow manifold is of the double well type assumed throughout this paper -- as we shall briefly show below.
\\ \\
The first model considers the interaction of vegetation and water, it is the 2-dimensional (in space) version of the 1-dimensional model studied by the methods of geometric singular perturbation theory in \cite{JDCBM20},
\beq
\label{e:RDE-exZelnik}
\left\{	
\begin{array}{rcl}
U_t &=& \hspace{0.375cm} \Delta U  - \left[ \mu_1 - \mu_2 (1-\mu_3 U)(1+\mu_4 U)^2 V \right] U \\[1mm]
V_t &=& \frac{1}{\varepsilon^2} \Delta V + \mu_5 - \left[ \mu_6 (1- \mu_3 \mu_7 U) + \mu_8 U (1+\mu_4 U)^2 \right] V
\end{array}
\right.
\eeq
\cite{Mer15,ZMB15} in which $U(x,y,t)$ models the vegetation biomass and $V(x,y,t)$ the much faster diffusing water content. We refrain from going into the details of all specific terms -- see \cite{JDCBM20,Mer15,ZMB15} -- we just notice that $\vmu \in \RR^8$ with $\mu_j \geq 0$, $j=1,...,8$ (naturally, this can be scaled back to an essential 5-dimensional parameter space -- see \cite{JDCBM20}). The second example is the savanna grass/woodland ecosystem model of \cite{dKetal08,GvdVvL17,vLetal03},
\beq
\label{e:RDE-exvLangevelde}
\left\{	
\begin{array}{rcl}
U_t &=& \hspace{0.375cm} \Delta U  + \mu_1 + \left[\frac{\mu_2}{\mu_3 U + V + \mu_4} - \mu_5 - \mu_6 V \right] U \\[1mm]
V_t &=& \frac{1}{\varepsilon^2} \Delta V + \left[\frac{\mu_7}{\mu_3 U + V + \mu_4} - \mu_8 \right] V
\end{array}
\right.
\eeq
(in the form of \cite{vLetal03} and slightly simplified) in which $U(x,y,t)$ models the woody biomass and $V(x,y,t)$ the grass biomass -- that in this case is the (relatively) fast diffuser -- and which also has $\vmu \in \RR^8$ and $\mu_j \geq 0$, $j=1,...,8$. Again, we refer to the ecological literature \cite{dKetal08,GvdVvL17,vLetal03} for the nature of each of the terms in (\ref{e:RDE-exvLangevelde}) and the interpretation of its parameters.
\\ \\
To show that dryland model (\ref{e:RDE-exZelnik}) indeed has the nature of the systems considered here, we first bring (the existence problem associated to) (\ref{e:RDE-exZelnik}) in the form (\ref{e:DS}) and note that -- except for the vertical slow manifold $\{u=p=0\}$ -- the slow manifolds are determined by the cubic equation (in $u$),
\beq
\label{e:cubicZelnik}
(1-\mu_3 u)(1+\mu_4 u)^2 = \frac{\mu_1}{\mu_2 v},
\eeq
in which $v$ is a parameter. Thus, there a priori are 3 branches $u = f^j(v)$ that can be ordered (for $v$ such that all 3 exist): $f^1(v) \leq f^2(v) \leq f^3(v)$ (and clearly $f^1(v) < -1/\mu_4 < 0$: this branch thus is ecologically irrelevant). The position of the local the maximum of the cubic in (\ref{e:cubicZelnik}) is given by $u_M = (2 \mu_4 - \mu_3)/(3 \mu_3 \mu_4)$. Assuming that $2 \mu_4 - \mu_3 > 0$, we define for $u > u_M > 0$ the normally hyperbolic slow manifold $\M_0$ as the third branch: $\M_0 = \{u=f(v) = f^3(v), p= 0\}$. Next we note that apart from those with $u=0$, the critical points/trivial background states of (\ref{e:RDE-exZelnik}) are determined by another cubic equation in $u$,
\[
\mu_2 \mu_5 (1-\mu_3 u)(1+\mu_4 u)^2 - \mu_1 \mu_6 (1- \mu_3 \mu_7 u) - \mu_1 \mu_8 u (1+\mu_4 u)^2 = 0
\]
(after elimination of $v$). Clearly, $\mu_1, ..., \mu_8$ can be tuned such that all 3 solutions merge at a $u_d > u_M$. For such a (singular) combination of $\mu_j$-values, there thus is 1 degenerate critical point on $\M_0$ that can be unfolded into 3 critical points $(\oV_-,0), (\oV_c,0), (\oV_+,0) \in \M_0$ with $v$-coordinates $\oV_- < \oV_c < \oV_+$: a center in between 2 saddles. Thus, the reduced slow flow on $\M_0$ indeed can have the double well character assumed throughout this paper (for $\vmu$ in an open set in $\RR^8$).
\\ \\
The slow manifolds of the savanna grass/woodland model (\ref{e:RDE-exvLangevelde}) are determined by a quadratic equation (with $v$ as parameter),
\[
\mu_1(\mu_3 u + v + \mu_4) + \mu_2 u - (\mu_5 + \mu_6 v)(\mu_3 u + v + \mu_4) u = 0,
\]
that always has 2 real solutions $u = f^{1,2}(v)$ with $f^1(v) < 0 < f^2(v)$. Thus, there is only one ecologically relevant slow manifold $\M_0 = \{u = f(v) = f^2(v), p=0\}$. For any given $u>0$, the fast flow of the spatial problem associated to (\ref{e:RDE-exvLangevelde}) has a homoclinic solution attached to the critical point on $\M_0$ that necessarily  enters into negative $u$-space, since it orbits around the center point determined by $u = f^1(v) < 0$. (This may seem counter-intuitive, but it is natural from an ecological point of view: trees may regrow from a burned savanna \cite{dKetal08,GvdVvL17,vLetal03}.) Thus, $\M_0$ is normally hyperbolic and more importantly (for the setting of the present paper): there is no fast jumping mechanism, there is only one slow manifold -- $\M_0$ -- and solutions cannot jump from $\M_0$ back to itself without becoming ecologically unrealistic. In other words, the only ecologically meaningful patterns (of traveling wave type that do not vary in the $y$-direction) are the slow patterns -- i.e. solutions on $\M_0$ (in the limit $\eps \to 0$) -- considered here. Moreover, the reduced slow flow on $\M_0$ is given by
\beq
\label{e:RedSF-exvLangevelde}
v_{XX} - \mu_8  v +  \frac{\mu_7 v}{\mu_3 f(v) + v + \mu_4} = 0.
\eeq
Apart from the critical point at the origin, (\ref{e:RedSF-exvLangevelde}) can have 2 additional critical points. It can be checked that there is a large open region in $\vmu$-space for which (\ref{e:RedSF-exvLangevelde}) indeed has the double well character with 3 critical points -- a center $(\oV_c,0)$ in between 2 saddles $(\oV_-,0) = (0,0)$ and $(\oV_+,0)$ -- assumed throughout this paper. Finally, we note that the slow heteroclinic front solution that will be studied in detail in the upcoming sections represents the ecologically crucial interface between the grassless homogeneous woody state $(U(x,y,t), V(x,y,t)) \equiv (f(0),0)$ and the grassy state $(U(x,y,t), V(x,y,t)) \equiv (f(\oV_+),\oV_+)$.
\begin{remark}
\label{r:Champneys}
\rm
The recently in \cite{Cetal20} introduced general class of models,
\beq
\label{e:RDE-Champneys}
\left\{	
\begin{array}{rcl}
U_t &=& \hspace{0.375cm} \Delta U  + \mu_1 - \mu_2 U + \mu_3 V + U^2V \\[1mm]
V_t &=& \frac{1}{\varepsilon^2} \Delta V +  \mu_4 + \mu_5 U - \mu_6 V - U^2V
\end{array}
\right.
\eeq
with $\vec{\mu} \in \RR^6$ and $\mu_j \geq 0$ ($j=1,...,6$), includes various versions of the well-studied Brusselator, Gray-Scott and Schnakenberg models (for instance, the Gray-Scott example of Fig. \ref{f:M0genGS} has $\mu_1 = \mu_3 = \mu_5 = 0$ and $\mu_4 = \mu_6$). The fast reduced flow associated to (\ref{e:RDE-Champneys}) typically has homoclinic orbits, which implies that there is a normally hyperbolic center manifold $\M_0$ determined by $u = f(v)$ with in general -- i.e. for $(\mu_1,\mu_3) \neq (0,0)$ -- $f'(v) \not\equiv 0$. Although it is studied in \cite{Cetal20} for $\vec{\mu}$ such that there is 1 unique background state, system (\ref{e:RDE-Champneys}) can have up to 3 critical points/trivial background states. It can be shown along the above lines (for system (\ref{e:RDE-exZelnik})) that the slow reduced flow on $\M_0$ can have the double well character considered here: there is an open region in $\vec{\mu}$-space for which also (\ref{e:RDE-Champneys}) has the nature of the models considered here.
\end{remark}

\section{Existence}
\label{s:Ex}

We define the critical values  $H_{0,\pm}$ of the unperturbed Hamiltonian (\ref{d:HW0}),
\beq
\label{d:H0pm}
H_{0,-} = \H_0(\oV_-,0) = -\W_0(\oV_-) = 0, \; H_{0,+} = \H_0(\oV_+,0) = -\W_0(\oV_+).
\eeq
By the assumptions on (the double well) potential $\W_0(v)$, we know that the reduced slow system (\ref{e:RedSF}) has  a homoclinic orbit for $H_{0,+} \neq  H_{0,-}$ and a pair of a heteroclinic connections if $H_{0,+} = H_{0,-}$ -- see Fig. \ref{f:SlowFlows}. Without loss of generality, we may only consider cases in which $H_{0,+} > 0$ -- i.e. $\W_0(\oV_+) < \W_0(\oV_-) = 0$ -- together with the limiting case $H_{+,0} = 0$. Thus, we assume that the well around $\oV_+$ is either deeper than the one around $\oV_+$, in which (\ref{e:RedSF}) has an orbit homoclinic to $(\oV_-,0)$ (Fig. \ref{f:SlowFlows}(a)), or that the wells are of equal depth, in which the homoclinic orbit has merged into a heteroclinic cycle between the saddles $(\oV_-,0)$ to $(\oV_+,0)$ (Fig. \ref{f:SlowFlows}(b)).
\\ \\
In this section we study both the persistence of the unperturbed stationary orbits and the appearance of new traveling (localized) patterns for the full slow flow on $\M_\eps = \M_\eps(c)$. Although we focus on localized solutions in this paper, we do consider in (sub)section \ref{ss:E-Pers} the fate under perturbation of both the (unperturbed) periodic and homoclinic orbits of (\ref{e:RedSF}) for $H_{0,+} > 0$, i.e. of all bounded solutions $(v_0(X),q_0(X))$ of (\ref{e:RedSF}) on level sets $\{H_0(v,q) = H_0\}$ with $H_0 \in (H_{0,c}, 0]$ -- where $H_{0,c} = \H_0(\oV_c,0) = \W_0(\oV_c) < 0$ (and we recall that $\W_0(v)$ has a local maximum at $\oV_c$ and that $(v_0(X),q_0(X))$ thus is periodic on $\{H_0(v,q) = H_0\}$ with $H_0 \in (H_{0,c}, 0)$ and homoclinic on $\{H_0(v,q) = 0\}$) -- see Theorem \ref{t:E-Pers}. In (sub)section \ref{ss:E-NearHet}, we consider the limiting, more subtle, nearly heteroclinic case of the (perturbed) double well potential with wells of equal depth (i.e. $H_{0,+} = H_{0,-} = 0$). Here, we focus on only the localized patterns that may bifurcate from the unperturbed heteroclinic cycle: the stationary and traveling heteroclinic fronts of Theorem \ref{t:exhets} and the stationary and traveling homoclinic pulses of Theorem \ref{t:exhoms}. However, we first set the stage by determining explicit approximations of $\M_\eps$ and the flow on $\M_\eps$.

\subsection{The slow manifold and slow reduced flow}
\label{ss:E-SM}
By Fenichel's First Theorem, we know that if $\M_0$ is normally hyperbolic, it persists as slow manifold $\M_\eps$ for $\eps > 0$ and sufficiently small \cite{Fen79,Jon95}; $\M_\eps$ can be explicitly approximated by,
\beq
\label{e:Meps}
\M_\eps(c) = \{(u,p,v,q): u = f(v) + \eps f_1(v,q) + \eps^2 f_2(v,q) + \O(\eps^3), \; p = \eps p_1(v,q) + \eps^2 p_2(v,q) + \O(\eps^3)\}.
\eeq
By also expanding $c$ into $c = c_0 + \eps c_1 + \O(\eps^2)$ and denoting $\frac{\partial F}{\partial u}(u,v)$ by $F_u(u,v)$, etc., we find,
\beq
\label{d:fjpj}
\begin{array}{ll}
f_1(v,q) = - c_0 \tau q \tf_1(v)
&
f_2(v,q) = -c_1 \tau q \tf_1(v) - c_0^2 \tau^2 \tF_{2cc}(v,q^2) - \tF_2(v,q^2)
\\
p_1(v,q) = q f'(v)
&
p_2(v,q) = -c_0 \tau \left[ q^2 \tf_1'(v) - \tf_1(v) G(f(v),v) \right]
\end{array}
\eeq
where,
\beq
\label{d:tFs}
\begin{array}{ccc}
\tf_1(v) & = & \frac{f'(v)}{F_u(f(v),v)}
\\
\tF_{2cc}(v,q^2) & = & \frac{\frac12 q^2 \tf_1^2(v)F_{uu}(f(v),v) - q^2\tf_1'(v) + \tf_1(v) G(f(v),v)}{F_u(f(v),v)}
\\
\tF_{2}(v,q^2) & = & \frac{q^2 f''(v) - f'(v) G(f(v),v)}{F_u(f(v),v)}
\end{array}
\eeq
Note that we have made it explicit in the notation that $\tF_{2cc}$ and $\tF_{2}$ are functions of $q^2$ (and thus symmetric under the transformation $q \to -q$). The assumption that $\M_\eps(c)$ is normally hyperbolic imposes a sign condition on $F_u(f(v),v)$, therefore the appearance of $F_u(f(v),v)$ as denominator in the above expressions is not a problem. Throughout this paper, we assume that $\M_\eps(c)$ is uniformly normally hyperbolic for all values of $v$ covered by the patterns we study. Thus, we assume that for some $\kappa > 0$,
\beq
\label{e:NormHyp}
F_u(f(v),v) < - \kappa < 0 \; \; \forall \; \; \oV_- \leq v  \leq \oV_+.
\eeq
Naturally, the slow flow on $\M_\eps(c)$ is now given by,
\beq
\label{e:SF}
v_{XX} + G(f(v),v) + \eps\left[ c_0 v_X \tG_{1c}(v) \right] + \eps^2 \left[ c_1 v_X \tG_{1c}(v) + c_0^2 \tau^2 \tG_{2cc}(v,v_X^2) + \tG_{2}(v,v_X^2) \right] = \O(\eps^3),
\eeq
where,
\beq
\label{d:tGs}
\begin{array}{ccl}
\tG_{1c}(v) & = & 1- \tau \tf_1(v) G_u(f(v),v)
\\
\tG_{2cc}(v,q^2) & = & - \tF_{2cc}(v,q^2) G_u(f(v),v)  + \frac12 q^2 \tf_1^2(v)G_{uu}(f(v),v)
\\
\tG_{2}(v,q^2) & = & - \tF_{2}(v,q^2) G_u(f(v),v)
\end{array}
\eeq
(\ref{d:tFs}) -- where we explicitly note that (\ref{e:SF}) has inherited the symmetry,
\beq
\label{d:symm}
X \to - X, \; p \to -p, \; q \to -q, \; c \to - c
\eeq
of the full system (\ref{e:DS}) and that we re-obtain its slightly more natural version (\ref{e:SF-I}) by re-writing $c_0 + \eps c_1 + \O(\eps^2)$ as $c$. Moreover, we note that $\M_\eps(c)$ is uniquely determined, since we have assumed that the 3 trivial patterns $(U(x,y,t),V(x,y,t) \equiv (\oU_\pm,\oV_\pm), (\oU_c,\oV_c)$ of (\ref{e:RDE}) correspond to critical points $(\oV_\pm, 0), (\oV_c, 0)$ on $\M_\eps(c)$. In fact, we confirm by (\ref{d:tFs}), (\ref{d:tGs}) that indeed all $(\oV,0) \in \RR^2$ with $G(f(\oV),\oV) = 0$ are critical points of the (full) planar slow system (\ref{e:SF}) and that the assumption that $(\oV_\pm,0) \in \M_\eps(c)$ are (non-degenerate) saddle points implies by linearization that,
\beq
\label{e:Gfvvlin}
\frac{d}{dV} \left[G(f(v),v)\right]|_{v=\oV_\pm} = G_u(\oU_\pm,\oV_\pm) f'(\oV_\pm) + G_v(\oU_\pm,\oV_\pm) =
\frac{\oF^\pm_u \oG^\pm_v - \oF_v^\pm \oG^\pm_u}{\oF^\pm_u} < 0
\eeq
(\ref{e:f'}) -- see also Lemma \ref{l:sigmaess} -- where we have introduced the notation,
\beq
\label{d:oFpmetc}
\oF^\pm_u = \frac{\partial F}{\partial u}(\oU_\pm,\oV_\pm), \; \oG^\pm_v = \frac{\partial G}{\partial v}(\oU_\pm,\oV_\pm), \; \; {\rm etc.}
\eeq
\begin{remark}
\label{r:naturepertslowflow-2}
\rm
Note that the fact that the leading order perturbation term in (\ref{e:SF}) has $c_0$ as pre-factor is a direct consequence of the assumption that $F(U,V)$ and/or $G(U,V)$ do depend directly on $\eps$ (Remark \ref{r:naturepertslowflow-1}). If either $F(U,V)$ or $G(U,V)$ (or both) depend(s) on $\eps$, there are $\O(\eps)$ perturbations without pre-factor $c_0$ in (\ref{e:SF}), which would immediately induce the appearance of traveling waves. Thus, the special role played by the stationary patterns and the associated vertical nature of the bifurcations into traveling waves -- see Fig. \ref{f:HetHomBifsIntro} and the upcoming sections -- also is a consequence of this choice.
\end{remark}

\subsection{$H_{0,+} > 0$: the persistence of the slow reduced patterns}
\label{ss:E-Pers}
To establish the persistence of solutions of (\ref{e:RedSF}) into solutions of (\ref{e:SF-I}), we need to distinguish between the stationary and traveling cases. When $c=0$, the full system (\ref{e:DS}) is reversible (\ref{d:symm}) which implies that the planar slow flow system (\ref{e:SF-I}) is also reversible (note that it thus is crucial/natural that $\tG_{2}$ indeed is  function of $q^2$ (\ref{d:tFs}), (\ref{d:tGs})). Thus, for $c=0$, (\ref{e:SF-I}) remains integrable in the sense that all periodic orbits and their homoclinic limit persist -- see Theorem \ref{t:E-Pers} below for a more precise statement. In fact, near a level set $\{H_0(v,q) = H_0\}$, a local $\O(\eps^2)$ correction of the $\O(1)$ integral $\H_0$ can be constructed explicitly by writing (\ref{e:SF-I}) as,
\beq
\label{e:SFc=0H0}
v_{XX} + G(f(v),v) + \eps^2 \left[\tG_{2}\left(v,2H_0 + 2\int_{\oV_-}^v G(f(\tv),\tv) \, d \tv\right)  \right] = \O(\eps^3)
\eeq
(\ref{d:HW0}). Clearly, (\ref{e:SFc=0H0}) is integrable (a priori up to $\O(\eps^2)$), with (local, approximate) integral $\H(v,q; H_0) = \frac12 q^2 - \W_0(v) - \eps^2 \W_2(v; H_0)$ and
\beq
\label{d:H2}
\W_2(v; H_0) = - \int_{\oV_-}^{v} \tG_{2} \left( \hat{v},2H_0 + 2\int_{\oV_-}^{\hat{v}} G(f(\tv),\tv) \, d \tv \right) \, d \hat{v}.
\eeq
Away from the heteroclinic limit, this $\O(\eps^2)$ correction of (\ref{d:HW0}) does not have a leading order impact on the (non-degenerate) double well character of $\W_0(v)$. However, in the heteroclinic limit -- i.e. if $H_0 = H_{0,+} = H_{0,-} = 0$ -- this correction term does make a difference. By construction, $\H(\oV_-,0; 0) = H_{0,-} - \eps^2 \W_2(\oV_-; 0) = 0$, while
\beq
\label{e:W2=0}
\H(\oV_+,0; 0) = - \eps^2 \W_2(\oV_+; 0) = \eps^2 \int_{\oV_-}^{\oV_+} \tG_{2} \left( \hat{v}, 2\int_{\oV_-}^{\hat{v}} G(f(\tv),\tv) \, d \tv \right) \, d \hat{v} \neq 0,
\eeq
in general. Thus, the saddles $(\oV_-,0)$ and $(\oV_+,0)$ of (\ref{e:SF-I}) generally are no longer on the same level set: the heteroclinic connections typically break and become homoclinic loops -- either to $(\oV_-,0)$ or to $(\oV_+,0)$. Especially since the associated `almost double front' homoclinic patterns may have asymptotically small eigenvalues -- and thus potentially could be stable -- we will study this limit in more detail in section \ref{ss:E-NearHet}.
\\ \\
For $c \neq 0$, the leading order perturbation of (\ref{e:SF-I}) is $\O(\eps)$. More importantly, (\ref{e:SF-I}) is no longer reversible (nor integrable). In fact,
\beq
\label{e:dH0dX}
\frac{d}{dX} \H_0(v,q) = -\eps c \tG_{1c}(v) q^2 - \eps^2 \left[\tG_{2}(v,q^2) + c^2 \tau^2 \tG_{2cc}(v,q^2) \right] q + \O(\eps^3).
\eeq
Thus, if $(v_0(X),q_0(X))$ is a periodic or homoclinic solution of (\ref{e:RedSF}) as defined above -- i.e. $(v_0(X),q_0(X)) \subset \{H_0(v,q) = H_0\}$ with $H_0 \in (H_{0,c},0]$ and initial condition $(v_0(0),0)$ with $v_0(0) \in (\oV_c,\oV_+)$ -- then it follows by a classical application of the Melnikov method (see for instance \cite{GH83}) that it persists as solution of (\ref{e:SF-I}) if,
\beq
\label{e:DeltaH0}
\Delta \H_0(H_0) = \int_{0}^{T} \frac{d}{dX} \H_0(v(X),q(X)) \, dX = -\eps c \int_{0}^{T} \tG_{1c}(v_0) q_0^2 \, dX + \O(\eps^2) = 0,
\eeq
where $T=T(H_0)$ is the period of $(v_0(X),q_0(X))$ -- with $T(H_0) \to \infty$ in the homoclinic limit $H_0 \uparrow 0$. By defining $A(H_0) = v_0(\frac12 T)$ and $B(H_0) = v_0(0)$ -- so that $\oV_- \leq A(H_0) < \oV_c < B(H_0)$ -- we thus find the following (leading order) condition on the persistence of $(v_0,q_0)$,
\beq
\label{e:DHmu=0}
\tilde{\Delta} \H_0(H_0, \vmu) = 2 c \int_{A(H_0)}^{B(H_0)}  \left[1- \tau \frac{f'(v;\vmu) G_u(f(v),v;\vmu)}{F_u(f(v),v;\vmu)} \right] \sqrt{2H_0 - 2\int_{\oV_-}^v G(f(\tv),\tv; \vmu)} \, dv = 0
\eeq
(\ref{d:tFs}), (\ref{d:tGs}), (\ref{e:dH0dX}) and $\tilde{\Delta} \H_0(H_0, \vmu) =$ the leading order approximation of $-\Delta \H_0(H_0)/\eps$ -- where we have (temporarily) reintroduced the parameter dependence of $F(u,v;\vmu)$ and $G(u,v;\vmu)$. Note that the pre-factor $c$ in (\ref{e:DHmu=0}) confirms the persistence of all orbits $(v_0,q_0)$ in the stationary problem (at leading order in $\eps)$.
\begin{theorem}
\label{t:E-Pers}
Assume that $H_{0,+} = \H_0(\oV_+,0) > \H_0(\oV_-,0) = H_{0,-} = 0$ and that $\eps > 0$ is sufficiently small. Let $(v_0(X),q_0(X))$ be a periodic or homoclinic solution of (\ref{e:RedSF}) on the level set $\{H_0(v,q) = H_0\}$ with $H_0 \in (H_{0,c},0]$ and with initial condition $(v_0(0),0)$ such that $v_0(0) \in (\oV_c,\oV_+)$.
\\
{\bf (i) Stationary patterns.} If $c=0$, then $(v_0(X),q_0(X))$ persists as a periodic or homoclinic solution $\gas(X) = (\us(X),\ps(X),\vs(X),\qs(X)) \subset \M_\eps(0)$ of (\ref{e:DS}) with $(\vs(X),\qs(X))$ a solution of (\ref{e:SF-I}) with $(\vs(0),\qs(0)) = (v_0(0),0)$ such that $\|(\vs(X),\qs(X)) - (v_0(X),q_0(X))\| < C \eps^2$ for all $X \in [0,T(H_0))$ and some $C > 0$ -- where $T(H_0)$ is the period of $(v_0(X),q_0(X))$ (and $T(H_0) \to \infty$ in the homoclinic limit $H_0 \uparrow 0$). The $(u,p)$-components of $\gas(X)$ are given by $\us(X) = f(v_0(X)) + \O(\eps^2)$, $\ps(X) = \eps q_0(X) f'(v_0(X)) + \O(\eps^3)$ (\ref{e:Meps}), (\ref{d:fjpj}).
\\
{\bf (ii) Traveling patterns.} If $c\neq 0$, then for any $c \in \RR$ (of $\O(1)$ w.r.t. $\eps$) $(v_0(X),q_0(X))$ persists as a solution $\gas(X) \subset \M_\eps(c)$ of (\ref{e:DS}) if there is a $\vmu_{\rm t}$ such that $\tilde{\Delta} \H_0(H_0, \vmu) = 0$ (\ref{e:DHmu=0}). Again, $(\vs(X),\qs(X))$ is a solution of (\ref{e:SF-I}) with $(\vs(0),\qs(0)) = (v_0(0),0)$ but only of $\O(\eps)$ accuracy: $\|(\vs(X),\qs(X)) - (v_0(X),q_0(X))\| < C \eps$ for all $X \in [0,T(H_0))$; likewise $\us(X) = f(v_0(X)) + \O(\eps)$, $\ps(X) = \eps q_0(X) f'(v_0(X)) + \O(\eps^2)$.
\\
In both cases, the homoclinic orbits of (\ref{e:DS}) correspond to (traveling or stationary) localized homoclinic pulse/stripe patterns $(U(x,y,t),V(x,y,t)) = (\Us(X),\Vs(X))$ of (\ref{e:RDE}) and the periodic solutions to (planar) periodic wave trains in (\ref{e:RDE}).
\end{theorem}
\noindent
Note that this theorem also establishes that there is an open region $\R_{\rm per}$ in parameter space for which system (\ref{e:RDE}) has slow traveling spatially periodic patterns and a co-dimensional 1 manifold $\S_{\rm hom} \subset \partial \R_{\rm per}$ for which it has slow traveling pulses, determined by the condition
\beq
\label{d:M0vmuW0}
M_{{\rm hom}}(\vmu; \W_0) \stackdef \int_{\oV_-}^{v_0^{\rm max}}  \left[1- \tau \frac{f'(v;\vmu) G_u(f(v;\vmu),v;\vmu)}{F_u(f(v;\vmu),v;\vmu)} \right] \sqrt{2 \W_0(v; \vmu)} \, dv = 0,
\eeq
where $v_0^{\rm max} = B(0) \in (\oV_c,\oV_+)$ exists as solution of $\W_0(v_0^{\rm max}) =0$ by the assumed double well nature of $\W_0(v)$ (\ref{d:HW0}) and indeed is the maximal value of the homoclinic orbit $v_0$: $v_0(X) \in (\oV_-, v_0^{\rm max}]$ (this argument implicitly uses the assumption that $H_{0,+} > H_{0,-} = 0$, i.e. that $\W_0(\oV_+) < \W_0(\oV_-) = 0$ (\ref{d:H0pm}): the $\oV_+$-well of $\W_0(v)$ is deeper than the $\oV_-$-well).
\\ \\
{\bf Proof of Theorem \ref{t:E-Pers}.} The main part of the proof of this theorem follows by the above standard, leading order, Melnikov arguments. However, the proof of case (ii) needs a higher order analysis. The fact that $c$ only appears as pre-factor of the $\O(\eps)$ perturbation terms naturally motivates the `{\it for any $c \in \RR$}' claim in (ii), nevertheless, since $c$ does not factor out from the $\O(\eps^2)$ perturbation terms in (\ref{e:SF-I}), one may expect to encounter conditions on $c$ when one considers higher order effects. This is not the case, although one needs to consider $\O(\eps^3)$ effects to establish this. We do not go into the -- technical and quite involved -- details here, but refer to Remark \ref{r:verttopar} for a sketch of this procedure (which is very much along the lines of the proof of Theorem \ref{t:exhoms}(ii), although in that case the higher order analysis does yield a condition on `allowable' $c$-values -- see also Corollary \ref{c:hombifs}). The higher order analysis sketched in Remark \ref{r:verttopar} can also be interpreted as establishing the natural parabolic nature of the bifurcation into traveling waves that occurs near values of $\vmu$ for which (\ref{e:DHmu=0}) holds: the present `{\it for any $c \in \RR$}' formulation gives its vertical leading order approximation. \hfill $\Box$

\begin{remark}
\label{r:ecopatt}
\rm
Theorem \ref{t:E-Pers} is a generalized version of Theorem 2.4 in \cite{JDCBM20} on the existence of slow periodic and homoclinic patterns in the explicit ecological model considered in \cite{JDCBM20} (although the slow reduced flow is not of double well type in \cite{JDCBM20}).
\end{remark}

\begin{remark}
\label{r:f'0ex}
\rm
In the vertical case $f'(v) \equiv 0$, the reduced slow manifold $\M_0$ also is invariant for the full system , i.e. $\M_\eps(c) equiv \M_0$, so that the only perturbation term in the flow on $\M_\eps(c)$ indeed is a linear friction term -- as already noted in the Introduction -- that directly originates from the `$-\eps c q$'-term in (\ref{e:DS}) (more specifically, $\tG_{1c}(v) \equiv 1$ and $\tG_{2cc}(v,v_X^2) \equiv \tG_{2}(v,v_X^2) \equiv 0$ (etc.) in (\ref{e:SF})). Hence, in this case only stationary patterns exist (which is covered by Theorem \ref{t:E-Pers} since $\tilde{\Delta} \H_0(H_0, \vmu)$ cannot be $0$ (\ref{e:DHmu=0})) -- see Fig. \ref{f:pertSlowFlows}(b).
\end{remark}

\subsection{$H_{0,+} =0$: nearly heteroclinic fronts and pulses}
\label{ss:E-NearHet}
In this section we consider the limit case -- from the point of the previous section -- in which the wells of the double well potential $\W_0(v)$ (\ref{d:HW0}) are of equal depth, i.e. the case in which $H_{0,+} = 0$ (\ref{d:H0pm}). By construction/assumption, (\ref{e:RedSF}) has 2 heteroclinic connections between the saddle points $(\oV_-,0)$ and $(\oV_+,0)$: $\va(X)$ and $\va(-X)$ with $\lim_{X \to -\infty} \va(X) = \oV_-$, $\lim_{X \to \infty} \va(X) = \oV_+$, $\vax(X) = \qa(X) > 0$ and $\lim_{X \to \pm \infty} \qa(X) = 0$. We assume for simplicity that varying the parameter $\vmu$ does not change the fact that the wells of $\W_0(v)$ are of equal depth: the connecting orbit $\va(X)$ is assumed to exist for all choices of $\vmu \in \RR^m$. This may be considered to be slightly unnatural, since one typically would prefer to have the freedom to change the relative depth of the wells by varying a parameter. However, allowing this additional (1-dimensional) freedom is not a fundamental issue, it would just add `another layer of technicalities' -- see also Remark \ref{r:naturepertslowflow-1}.
\\ \\
Unlike in the previous section, we will explicitly construct orbits $\gas(X) = (\us(X),\ps(X),\vs(X),\qs(X)) \subset \M_\eps(c)$ of (\ref{e:DS}) that are either homoclinic to the critical point $(\oU_-,0,\oV_-,0)$ -- i.e. $\gas(X) \subset W^u((\oU_-,0,\oV_-,0)) \cap W^s((\oU_-,0,\oV_-,0)) \cap \M_\eps(c)$ -- or heteroclinic between $(\oU_-,0,\oV_-,0)$ and $(\oU_+,0,\oV_+,0)$ -- i.e. $\gas(X) \subset W^u((\oU_-,0,\oV_-,0)) \cap W^s((\oU_+,0,\oV_+,0)) \cap \M_\eps(c)$. The analysis will naturally be done in the setting of the slow flow on $\M_\eps(c)$, thus we will construct solutions $(\vs(X),\qs(X))$ of (\ref{e:SF-I})/(\ref{e:SF}) with $\lim_{X \to -\infty} (\vs(X),\qs(X)) = (\oV_-,0)$ and $\lim_{X \to \infty}(\vs(X),\qs(X)) =$ either $(\oV_-,0)$ or $(\oV_+,0)$. Clearly, in both cases the orbits $(\vs(X),\qs(X))$ will be asymptotically close to the heteroclinic cycle between $(\oV_-,0)$ and $(\oV_+,0)$ spanned by $(\va(X),\qa(X))$ and $(\va(-X),-\qa(-X))$: the heteroclinic orbit $(\va(X),\qa(X))$ of (\ref{e:RedSF}) is the foundation of the upcoming  asymptotic analysis. Clearly, there may also be periodic patterns asymptotically close to the cycle spanned by $(\va(X),\qa(X))$ and $(\va(-X),-\qa(-X))$, however, these patterns will not be considered here -- see section \ref{ss:Projects}.
\\ \\
The perturbation analysis naturally starts out with the classical expansion,
\beq
\label{d:v12q12}
\vs(X) = \va(X) + \eps v_1(X) + \eps^2 v_2(X) + \O(\eps^3), \;
\qs(X) = \qa(X) + \eps q_1(X) + \eps^2 q_2(X) + \O(\eps^3),
\eeq
where we note that the validity of this expansion -- and the control over the (magnitude of the) correction terms -- is settled by Poincar\'e's Expansion Theorem -- see for instance \cite{Ver96} and the upcoming discussions below. Substitution of these expansions into (\ref{e:SF}) yields at leading order in $\eps$,
\beq
\label{e:Lav1}
v_{1,XX} + \left[ G_u(f(\va),\va) f'(\va) + G_v(f(\va),\va) \right] v_1 = - c_0 \qa \tG_{1c}(\va)
\eeq
(\ref{d:tG1c-I}). Therefore, we define the (Sturm-Liouville) operator
\beq
\label{d:La}
\La = \La(X; \vmu) = \frac{d^2}{dX^2} + \left[ G_u(f(\va(X)),\va(X)) f'(\va(X)) + G_v(f(\va(X)),\va(X)) \right].
\eeq
The following (technical) lemma considers the homogeneous and inhomogeneous problems associated to $\La(X)$, i.e.
\beq
\label{e:Lah}
\La(X) v = h(X),
\eeq
where we a priori only assume that $h: \RR \to \RR$ is sufficiently smooth. This simple lemma is crucial to both the upcoming existence and (spectral) stability analysis. However, at several places in the text we will need a somewhat more refined result: Lemma \ref{l:Lah-sharp}, that is formulated and proven in Appendix \ref{a:refined}.
\begin{lemma}
\label{l:Lah}
Let $\va(X)$ be the increasing heteroclinic solution of (\ref{e:RedSF}) written as $v_{XX} - \W'_0(v) = 0$ (\ref{d:HW0}), where $\W_0(v)$ is a double well potential with equal (non-degenerate) wells at $v = \oV_\pm$ -- i.e. $\H_0(\oV_\pm) = 0$ -- and a local maximum at $\oV_c \in (\oV_-,\oV_+)$: $\lim_{X \to \pm \infty} (\va(X), \vax(X)) = (\oV_\pm,0)$, $\vax(X) > 0$ and (by choice) $\va(0) = \oV_c$. Define,
\beq
\label{d:muEnpm}
\al_{\pm} = \W''_0(\oV_{\pm}) = -\frac{\oF^\pm_u \oG^\pm_v - \oF_v^\pm \oG^\pm_u}{\oF^\pm_u} > 0, \; \be_\pm > 0, \; E_{\pm}(X) =  e^{\mp \sqrt{\al_\pm} X} \;
\eeq
(\ref{e:Gfvvlin}), (\ref{d:oFpmetc}), i.e. $E_{\pm}(X) \to 0$ as $X \to \pm \infty$, so that for $\pm X \gg 1$,
\beq
\label{e:vaXgg1}
\va(X) = \oV_{\pm} \mp \be_{\pm} e^{\mp \sqrt{\al_\pm} X}\left(1 + \O(E_{\pm}(X))\right)
\eeq
(where $\be_\pm$ are determined by our choice $\va(0) = \oV_c$). Then, the homogeneous problem (\ref{e:Lah}) (with $h(X) \equiv 0$) has 2 independent solutions, $\vb(X)$ -- bounded and even as function of $X$ -- and $\vu(X)$ -- unbounded and odd --
\beq
\label{d:vbu}
\vb(X) = \vax(X) \, (= \qa(X)), \; \; \vu(X) = \vax(X) \int_0^X \frac{d \tX}{\vax^2(\tX)}
\eeq
with,
\beq
\label{e:propsvbu}
\begin{array}{llll}
\vb(0) = \sqrt{2 \W_0(\oV_c)} & \vbx(0) = 0 & \vb(X) = \be_{\pm} \sqrt{\al_\pm} e^{\mp \sqrt{\al_\pm} X}(1 + \O(E_{\pm})) & {\rm for} \; \pm X \gg 1
\\
\vu(0) = 0 & \vbx(0) = \frac{1}{\sqrt{2 \W_0(\oV_c)}} & \vu(X) = \frac{\pm 1}{2\al_\pm \be_{\pm}} e^{\pm \sqrt{\al_\pm} X}(1 + \O(E_{\pm})) & {\rm for} \; \pm X \gg 1
\end{array}
\eeq
If $h(X)$ decays exponentially as $X\to -\infty$, i.e. if there is a $\al > 0$ such that $\lim_{X \to -\infty} h(X) e^{-\sqrt{\al}X}$ exists, then the (unique) solution $v(X)$ of the inhomogeneous problem (\ref{e:Lah}) that converges to 0 as $X \to -\infty$ with initial condition $v(0) = 0$, is given by
\beq
\label{e:homsolinh}
v(X) = -\left[\int_0^X h(\tX) \vu(\tX) \, d\tX\right] \vb(X) + \left[\int_{-\infty}^X h(\tX) \vb(\tX) \, d\tX\right] \vu(X)
\eeq
where $v(X)$, $v_X(X)$ both decay exponentially to $0$ as $X \to -\infty$. If it is additionally assumed that $h(x)$ remains bounded for $X \gg 1$ -- i.e. that there is a $C > 0$ such that $|h(X)| < C$ for all $X \in \RR$ -- then for $X \gg 1$,
\beq
\label{e:vgg1}
v(X) = \frac{\int_{-\infty}^\infty h(\tX) \vb(\tX) \, d\tX}{2\al_+ \be_{+}} e^{\sqrt{\al_+} X}\left(1 + \O(E_{+})\right) \stackdef \frac{M_h}{2\al_+ \be_{+}} e^{\sqrt{\al_+} X}\left(1 + \O(E_{+})\right),
\eeq
with $|M_h| < \infty$.
\end{lemma}
\noindent
{\bf Proof.} By construction we have $\vb(X) = \vax(X)$ so that $\vu(X)$ follows by standard methods. Likewise, properties (\ref{e:propsvbu}) of $\vb(X)$ and $\vu(X)$ follow directly from (\ref{d:muEnpm}), (\ref{e:vaXgg1}), (\ref{d:vbu}). Moreover, the general solution to the inhomogeneous problem (\ref{e:Lah}) is given by
\beq
\label{e:gensolinh}
v(X) = \left[C_{\rm b} - \int_0^X h(\tX) \vu(\tX) \, d\tX\right] \vb(X) + \left[C_{\rm u} + \int_0^X h(\tX) \vb(\tX) \, d\tX\right] \vu(X),
\eeq
with $C_{{\rm b},{\rm u}} \in \RR$ -- where we have used that the Wronskian associated to (\ref{e:Lah}) $\equiv 1$ by (\ref{e:propsvbu}). These same properties also imply that the condition $v(0) = 0$ is satisfied by setting $C_{\rm b} = 0$. If $\lim_{X \to -\infty} h(X) e^{-\sqrt{\al}X}$ exists for some $\al > 0$, then it follows again from (\ref{e:propsvbu}) that the first term in (\ref{e:gensolinh}) goes to $0$ as $X \to -\infty$ and that we need to set
\[
C_{\rm u} = \int_{-\infty}^0 h(\tX) \vb(\tX) \, d\tX
\]
to counteract the growth of $\vu(X)$ for $X \ll -1$. Hence, we recover (\ref{e:homsolinh}) from (\ref{e:gensolinh}) and we may conclude by a direct check -- again using (\ref{e:propsvbu}) -- that indeed $\lim_{X \to -\infty} v(X) = \lim_{X \to -\infty} v_X(X) = 0$. Finally, we note that for bounded $h(x)$, the first term of (\ref{e:homsolinh}) remains bounded as $X \to \infty$ and that the integral in the second term indeed converges to $M_h$ with $|M_h| < \infty$, so that (\ref{e:vgg1}) follows from the leading order growth behavior of $\vu(X)$. See also (the proof of) Lemma \ref{l:Lah-sharp} for more details.
\hfill $\Box$
\\ \\
The application of Lemma \ref{l:Lah} to (\ref{e:Lav1}) yields explicit (leading order) control over the unstable manifold of the saddle $(\oV_-,0)$ of (\ref{e:SF}), i.e. over the solution $(\vs(X),\qs(X)) \subset W^u((\oV_-,0))$. We therefore introduce
$\tv_1(X)$ as the unique solution of
\beq
\label{d:tv1}
\La \tv_1 = - \qa \tG_{1c}(\va) = - \left[1- \tau \frac{f(\va)G_u(f(\va),\va)}{F_u(f'(\va),\va)}\right] \vax
\eeq
with $v_1(X) = c_0 \tv_1(X)$ (\ref{e:Lav1}), that decays to $0$ as $X \to - \infty$ and has initial condition $\tv_1(0) = 0$ (note that this initial condition just fixes the position of the point $X=0$ through $\vs(0) = \oV_c$ (\ref{d:v12q12})). Since the inhomogeneous term of (\ref{d:tv1}) clearly decays as $X \to \pm \infty$, it follows from Lemma \ref{l:Lah} that $v_1(X)$ -- the first term in the expansion of $\vs(X)$ (\ref{d:v12q12}) -- can be written as
\beq
\label{e:tv1}
v_1(X) =  c_0 \tv_1(X) = c_0\left\{\left[\int_0^X \tG_{1c}(\va) \vax \vu d\tX\right] \vax(X) - \left[\int_{-\infty}^X \tG_{1c}(\va) \vax^2 d\tX\right] \vu(X)\right\},
\eeq
and that for $X \gg 1$,
\beq
\label{e:tv1gg1-0}
\tv_1(X) = - \frac{M_{\ast}(\vmu)}{2 \al_+ \be_{+}} e^{\sqrt{\al_+} X}(1 + \O(E_{+})),
\eeq
with $\Ma(\vmu)$ as defined in (\ref{d:Ma}), where we note that $M_{\ast}(\vmu)$ can also be written as
\[
M_{\ast}(\vmu) = \int_{\oV_-}^{\oV_+} \left[1- \tau \frac{f'(v;\vmu) G_u(f(v;\vmu),v;\vmu)}{F_u(f(v;\vmu),v;\vmu)} \right] \sqrt{2\W_0(v; \vmu))} \, dv,
\]
which naturally coincides with the limit of $M_{{\rm hom}}(\vmu; \W_0)$ (\ref{d:M0vmuW0}) in which the wells of $\W_0(v)$ have become of equal depth -- i.e. the limit $H_{0,+} \downarrow 0$ (\ref{d:H0pm}). To be homoclinic to $(\oV_-,0)$, $\vs(X)$ needs to cross through the $\{q=0\}$-axis, i.e. there has to be a $X_{\rm h}$ such that $\vsx(\Xh) = 0$. Clearly, $\Xh$ cannot be $\O(1)$ since $\vax(X) = \vb(X) = \qa(X) > 0$ for all $X$. However, it follows from (\ref{e:propsvbu}) and (\ref{e:tv1gg1-0}) that for $X \gg 1$
\beq
\label{e:vsgg1}
\vsx(X) = \left[\be_{+} \sqrt{\al_+} e^{-\sqrt{\al_+} X} - \eps c_0 \frac{M_{\ast}(\vmu)}{2\be_{+} \sqrt{\al_+}} e^{\sqrt{\al_+} X}\right](1 + \O(E_{+})) + \O(\eps^2 v_2).
\eeq
A priori, Poincar\'e's Expansion Theorem establishes the validity of regular expansions like (\ref{d:v12q12}) on $X$-intervals of $\O(1)$ length (w.r.t. $\eps$), which would be insufficient for solving $\vsx(\Xh) = 0$ (\ref{e:vsgg1}). However in the present setting, the interval of validity can be extended -- with an associated loss of accuracy (see also \cite{DV15}). First we note that by assuming that $(\vs(X),\qs(X)) \subset W^u((\oV_-,0))$, it necessary follows that $\lim_{X \to -\infty} (v_n(X),q_n(X)) = (0,0)$ (\ref{d:v12q12}) for all $n \geq 1$. As a consequence, the interval of validity of Poincar\'e's Expansion Theorem may be straightforwardly extended to $(-\infty, \O(1))$. Based on the explicit control on the growth of $\vb(x)$ and $\vu(X)$ for $X \gg 1$ (by (\ref{e:propsvbu}), we may extend the $X$-interval of validity of (\ref{d:v12q12}) to values of $X$ such that $E^{-1}_+(X) = \O(1/\eps^\si)$ for $\si > 0$ -- which means that $X$ is of $\O(|\log \eps|)$. Naturally, all equations for $v_n(X)$ are of the form (\ref{e:Lah}) and it follows iteratively that the inhomogeneous terms $h_n(X)$ in general grow exponentially as $E_+^{-n}(X)$ for $X \gg 1$. To see this, we first note that the inhomogeneous term $h_2(X)$ of the equation for $v_2$ must contain terms with $v_1^2(X)$ (cf. (\ref{e:Lav2}) in which $v_1(X)$ is represented by $\tv_1(X)$ (\ref{e:tv1})): the (maximal) growth of $h_2(X)$ indeed is of the form $E_+^{-2}(X)$ for $X \gg 1$ (\ref{e:tv1gg1-0}). Thus, it follows that $v_2(X)$ also grows as $E_+^{-2}(X)$ for $X \gg 1$ -- either by a direct interpretation of (\ref{e:homsolinh}) or more precisely by (\ref{d:approxsvX-sharp}) for $j=2$ in Lemma \ref{l:Lah-sharp}. Using this as input for the $v_3$ equation, we see that $h_3(X)$ grows as $E_+^{-3}(X)$ for $X \gg 1$ so that -- again by Lemma \ref{l:Lah-sharp} -- also $v_3(X)$ grows as $E_+^{-3}(X)$ for $X \gg 1$, etc.. Thus indeed, both the inhomogeneous term $h_n(X)$ and the solution $v_n(X)$ of (\ref{e:Lah}) with $h(X) = h_n(X)$ grows as $E_+^{-n}(X)$ for $X \gg 1$ (a priori for $n \geq 2$, but we already found by Lemma \ref{l:Lah} that this is also the case for $n=1$). Therefore, for $n \geq 1$,
\beq
\label{e:extPoincare}
\eps^n |v_n(X)| = \O\left(\eps^{n(1 - \si)}\right) \; \; {\rm for} \; \; X \in \left(-\infty, \frac{\si}{\sqrt{\al_+}}|\log \eps| + \O(1)\right),
\eeq
i.e. for $X$ up to values such that $E^{-1}_+(X) = \O(1/\eps^\si)$ -- which provides the desired extension of the (asymptotic) size of the interval of validity in Poincar\'e's Extension Theorem with associated loss of accuracy, under the assumption that $0 < \si < 1$. Since $v_{0,X}(X) = \vax(X) = \O(\eps^\si)$ for $E^{-1}_+(X) = \O(1/\eps^\si)$ (\ref{e:propsvbu}), it follows that $\vsx(\Xh) = 0$ can be solved from (\ref{e:vsgg1}) by setting $\si=\frac12$, and thus also that (\ref{e:vsgg1}) is asymptotically correct up to  $\O(\eps^1)$ correction terms ($1 =  2(1-\frac12)$),
\beq
\label{e:Xh-0}
e^{\sqrt{\al_+}\Xh} = \frac{\be_+ \sqrt{2\al_+}}{\sqrt{c_0 M_{\ast}(\vmu)}} \frac{1}{\sqrt{\eps}}
\eeq
under the assumption that $c_0 M_{\ast} >  0$ (\ref{d:Ma}). Moreover, it follows that the $v$-coordinate of $W^u((\oV_-,0)) \cap \{q=0\}$ is given by
\beq
\label{e:vaXh}
\vs(\Xh) =
\oV_{+} - \sqrt{\frac{2 c_0 M_{\ast}(\vmu)}{\al_+}}  \sqrt{\eps}  + \O(\eps),
\eeq
which yields that $c_0 \M_{\ast}(\vmu) = 0$ is the leading order approximation of the condition for which $\vs \subset W^u((\oV_-,0)) \cap W^s((\oV_+,0))$: the heteroclinic case that separates between situations in which $\W^u((\oV_-,0))$ does intersect the $\{q=0\}$-axis and may turn around back towards the saddle $(\oV_-,0)$ and cases in which $\vs$ keeps on increasing and the orbit $(\vs(X),\qs(X))$ passes along the saddle $(\oV_+,0)$.
\\ \\
However, it follows by symmetry (\ref{d:symm}) that the stable manifold $W^s((\oV_-,0))$ of the slow flow (\ref{e:SF}) with $c_0$ replaced by $-c_0$ intersects the $\{q=0\}$-axis as it passes along $(\oV_+,0)$ -- in backwards `time' -- with $v$-coordinate given by (\ref{e:vaXh}) with $c_0$ replaced by $-c_0$. Thus, by changing $c_0$ back to $-c_0$ again, we conclude that $W^s((\oV_-,0))$ only intersects the  $\{q=0\}$-axis in the original system (in backwards time) if $c_0 M_\ast \leq 0$ (and with $v$-coordinate $\oV_{+} - \sqrt{2 |c_0 M_{\ast}|/\al_+}  \sqrt{\eps}  + \O(\eps)$). The situation is thus similar to that of Theorem \ref{t:E-Pers}(ii): traveling heteroclinic or homoclinic orbits can only exist if $\Ma(\vmu) = 0$. More precisely, traveling patterns can a priori only exist for $\vmu$ in an asymptotically small neighborhood of (the critical co-dimension 1 manifold $\vmu=$) $\vmu_{\rm t}^\ast$ for which $\Ma(\vmu_{\rm t}^\ast) = 0$.
\\ \\
Thus, we conclude that expression (\ref{e:Xh-0}) for $\Xh$ cannot be correct, it needs to be re-derived in a setting in which $|c_0 M_{\ast}(\vmu)|$ is asymptotically small -- see (\ref{e:Xh-1}). (Note that this includes the case of stationary fronts or pulses, but that there is no condition on $\Ma(\vmu)$ for stationary patterns -- as in Theorem \ref{t:E-Pers}(i).) The analysis so far can be interpreted as providing a leading order insight in the persistence of the unperturbed heteroclinic orbit $\va(X)$, however, to establish existence results on traveling localized patterns, we need to zoom in near $\vmu = \vmu_{\rm t}^\ast$: we introduce $\vmu_{\tsi}$ by $\vmu = \vmu_{\rm t}^\ast + \eps^{\tsi} \vmu_{\tsi}$ and expand
\beq
\label{d:vmutsi}
\Ma(\vmu) = \Ma(\vmu_{\rm t}^\ast + \eps^{\tsi} \vmu_{\tsi}) = \vec{\nabla} \Ma(\vmu_{\rm t}^\ast) \cdot \vmu_{\tsi} \, \eps^{\tsi} + \O(\eps^{\tsi + 1}) \; \; {\rm with} \; \; \tsi > 0,
\eeq
where we note that $\tsi$ will have different values for the upcoming heteroclinic and homoclinic cases. Before we can formulate our main result on the persistence of the heteroclinic connection $\va(X)$ (Theorem \ref{t:exhets}) we need to introduce some notation. First, we note that for $\vmu = \vmu_{\rm t}^\ast + \eps \vmu_1$, i.e. for $\tsi = 1$, expansion (\ref{d:vmutsi}) yields a sharper version of (\ref{e:tv1gg1-0}): for $X \gg1$,
\beq
\label{e:tv1gg1}
\tv_1(X) = - \eps  \frac{\vec{\nabla} \Ma(\vmu_{\rm t}^\ast) \cdot \vmu_1}{2 \al_+ \be_+} e^{\sqrt{\al_+} X} + \O(X E_+(X)),
\eeq
where the nature of the leading order correction term follows from Lemma \ref{l:Lah-sharp}, the refined version of Lemma \ref{l:Lah}. (More precise: (\ref{d:tv1}) is of the form (\ref{e:Lah}) with $h(X)$ as in (\ref{d:oh01j}) with $j = -1$, hence (\ref{d:approxsvX-sharp}) holds with $j=-1$ and $M_h = \Ma(\vmu_{\rm t}^\ast) = 0$.) We can now define the (Melnikov-type) expressions $\tM_{2cc}(\vmu_{\rm t}^\ast)$ and $\tM_{2}(\vmu)$,
\beq
\label{d:tM2c2}
\begin{array}{ccl}
\tM_{2cc}(\vmu_{\rm t}^\ast) & = & \int_{-\infty}^{\infty} [\tG_{1c,v}^\ast \vax \tv_1 - \tG_{1c}^\ast \tv_{1,X} -  \frac12 ((f'_\ast)^2 G^\ast_{uu}  + G^\ast_{vv} + 2 f'_\ast G^\ast_{uv} + f''_\ast G^\ast_{u})\tv_1^2 + \tau^2 \tG_{2cc}^\ast] \vax dX
\\
\tM_{2}(\vmu) & = & \int_{-\infty}^{\infty} \tG_{2}^\ast \, \vax \, dX
\end{array}
\eeq
(\ref{e:SF}), (\ref{d:tGs}), (\ref{e:tv1}) in which
\beq
\label{d:fFGast}
f'_\ast = f'(\va), \; \tG_{1c}^\ast = \tG_{1c}(\va), \; \tG_{2}^\ast = \tG_{2}(\va,\qa^2), \; G^\ast_{uu} = G_{uu}(f(\va),\va), \; {\rm etc.}
\eeq
Note that for general $\vmu$, the terms containing $\tv_1$ in $\tM_{2cc}(\vmu)$ prevent the integral in (\ref{d:tM2c2}) from converging (by (\ref{e:propsvbu}), (\ref{e:tv1gg1-0})). However, $\Ma(\vmu_{\rm t}^\ast) = 0$, so that the integral does converge at $\vmu= \vmu_{\rm t}^\ast$ and thus $|\tM_{2cc}(\vmu_{\rm t}^\ast)| < \infty$ (cf. (\ref{e:tv1gg1}) with $\vmu_1 = 0$). Note also that the situation for $\tM_{2}(\vmu)$ is less singular: $|\tM_{2}(\vmu)| < \infty$ for general $\vmu$ -- as will be needed in the upcoming theorem.
\begin{theorem}
\label{t:exhets}
Let $\W_0(v)$ (\ref{d:HW0}) be a double well potential with wells of equal depth, i.e. let $H_{0,+} = 0$ (\ref{d:H0pm}) and let $\eps > 0$ be sufficiently small.
\\
{\bf (i) Stationary fronts.} Let $\vmu_s(0) \in \RR^m$ be such that $\tM_{2}(\vmu_s(0)) = 0$ (\ref{d:tM2c2}), then there is a uniquely determined co-dimension 1 manifold in parameter space determined by $\vmu_s(\eps)$, such that for $\vmu = \vmu_s(\eps) = \vmu_s(0) + \O(\eps)$ there exists a heteroclinic solution $(\vhe(X), \qhe(X)) \subset W^u((\oV_-,0))\cap W^s((\oV_+,0))$ of (\ref{e:SF}) with $c=0$.
\\
{\bf (ii) Traveling fronts.} Let $\vmuta$ be such that $\Ma(\vmu_{\rm t}^\ast) = 0$, let $\vmu = \vmu_{\rm t}^\ast + \eps \vmu_1$ and consider (\ref{d:vmutsi}) with $\tsi = 1$. Let $c_0$ be a solution of,
\beq
\label{e:exhets}
c_0 \vec{\nabla} \Ma(\vmu_{\rm t}^\ast) \cdot \vmu_1 + c_0^2 \tM_{2cc}(\vmu_{\rm t}^\ast) + \tM_{2}(\vmu_{\rm t}^\ast) = 0
\eeq
(\ref{d:tM2c2}) and impose the non-degeneracy condition $\vec{\nabla} \Ma(\vmu_{\rm t}^\ast) \neq \vec{0}$. Then, there exists for $c_{\rm het}(\eps) = c_{\rm het}(0) + \O(\eps) = c_0 + \O(\eps)$ a heteroclinic solution $(\vhe(X), \qhe(X)) \subset W^u((\oV_-,0))\cap W^s((\oV_+,0))$ of (\ref{e:SF}).
\\
In both cases, the orbit $(\vhe(X), \qhe(X))$ corresponds to a slow heteroclinic orbit of (\ref{e:DS}) on $\M_\eps(c)$ -- with $c=0$ for (i) and $c=c_{\rm het}(\eps)$ for (ii) -- between the critical points $(\oU_-,0,\oV_-,0)$ and $(\oU_+,0,\oV_+,0)$ and to
a front/interface $(U(x,y,t),V(x,y,t)) = (\Uhe(X),\Vhe(X))$ pattern of (\ref{e:RDE})/(\ref{e:RDE-S}) -- with $\Vhe(X) = \vhe(X)$ and $\Uhe(X) = f(\vhe(X)) + \eps f_1(\vhe(X),\qhe(X)) + \eps^2 f_2(\vhe(X),\qhe(X)) + \O(\eps^3)$ (\ref{e:Meps}) (\ref{d:fjpj}) that is either stationary (i) or travels with speed $c_{\rm het}(\eps)$ (ii) and that connects the homogeneous states $(U(x,y,t),V(x,y,t)) \equiv (\oU_-,\oV_-)$ (as $X \to -\infty$) to $(U(x,y,t),V(x,y,t)) \equiv (\oU_+,\oV_+)$ ($X \to +\infty$) -- see Fig. \ref{f:frontspulses}(a).
\end{theorem}
\noindent
By symmetry (\ref{d:symm}) (with $c = 0$), Theorem \ref{t:exhets}(i) establishes the persistence of both stationary heteroclinic front orbits of (\ref{e:RedSF}) under the condition $\tM_{2}(\vmu_s)=0$, i.e. of the full heteroclinic cycle between the saddles $(\oV_-,0)$ and $(\oV_+,0)$ on $\M_\eps(0)$. Therefore, this condition is equivalent to $\W_2(\oV_+; 0) = 0$ (\ref{e:W2=0}).
\begin{corollary}
\label{c:hetbifs} Assume $\tM_{2cc}(\vmu_{\rm t}^\ast) \neq 0$, define
\beq
\label{d:tmu1}
\tilde{\mu}_1 = \vec{\nabla} \Ma(\vmu_{\rm t}^\ast) \cdot \vmu_1 \in \RR
\eeq
((\ref{d:vmutsi}) with $\si = 1$), and consider the solutions $c = c_{\rm het}^\pm(\tilde{\mu}_1; \eps)$ as determined (at leading order in $\eps$) by (\ref{e:exhets}) in Theorem \ref{t:exhets}(ii).
\\
{\bf (i)} If $\tM_{2}(\vmu_{\rm t}^\ast) = 0$, then the lines $c_{\rm het}^+(\tilde{\mu}_1) \equiv 0$ and $c_{\rm het}^-(\tilde{\mu}_1) = -\tilde{\mu}_1/ \tM_{2cc}$ intersect at $\tilde{\mu}_1 = 0$: (\ref{e:exhets}) describes (at leading order in $\eps$) a transcritical bifurcation at which a branch of heteroclinic orbits with $c_{\rm het}^-(\tilde{\mu}_1) \neq 0$ intersects the branch of stationary heteroclinic orbits.
\\
{\bf (ii)} If $\tM_{2}(\vmu_{\rm t}^\ast) \tM_{2cc}(\vmu_{\rm t}^\ast)< 0$, then both $c_{\rm het}^\pm(\tilde{\mu}_1)$ exist for all $\tilde{\mu}_1 = \O(1)$: there are no bifurcations. If $\tmu_1 \to 0$, i.e. as $\| \vmu - \vmu_{\rm t}^\ast \| \ll \eps$, then $c_{\rm het}^+(\tmu_1) \to  - c_{\rm het}^-(\tmu_1) =  \sqrt{-\tM_{2}/\tM_{2cc}}$.
\\
{\bf (iii)} If $\tM_{2}(\vmu_{\rm t}^\ast) \tM_{2cc}(\vmu_{\rm t}^\ast) > 0$, then the curves $c_{\rm het}^\pm(\tilde{\mu}_1)$ merge at the saddle-node bifurcation points $\tmu_{{\rm het}-SN}^\pm = \pm 2 \sqrt{\tM_{2}\tM_{2cc}}$ (with $c_{\rm het}(\tilde{\mu}_{{\rm het}-SN}^\pm) = -\tmu_{{\rm het}-SN}^\pm/(2 \tM_{2cc}) = \mp \sqrt{\tM_{2}/\tM_{2cc}}$); there are no heteroclinic fronts for $\tmu_{{\rm het}-SN}^- < \tmu_1 < \tmu_{{\rm het}-SN}^+$.
\end{corollary}
\begin{figure}[t]
\centering
	\begin{minipage}{.30\textwidth}
		\centering
		\includegraphics[width =\linewidth]{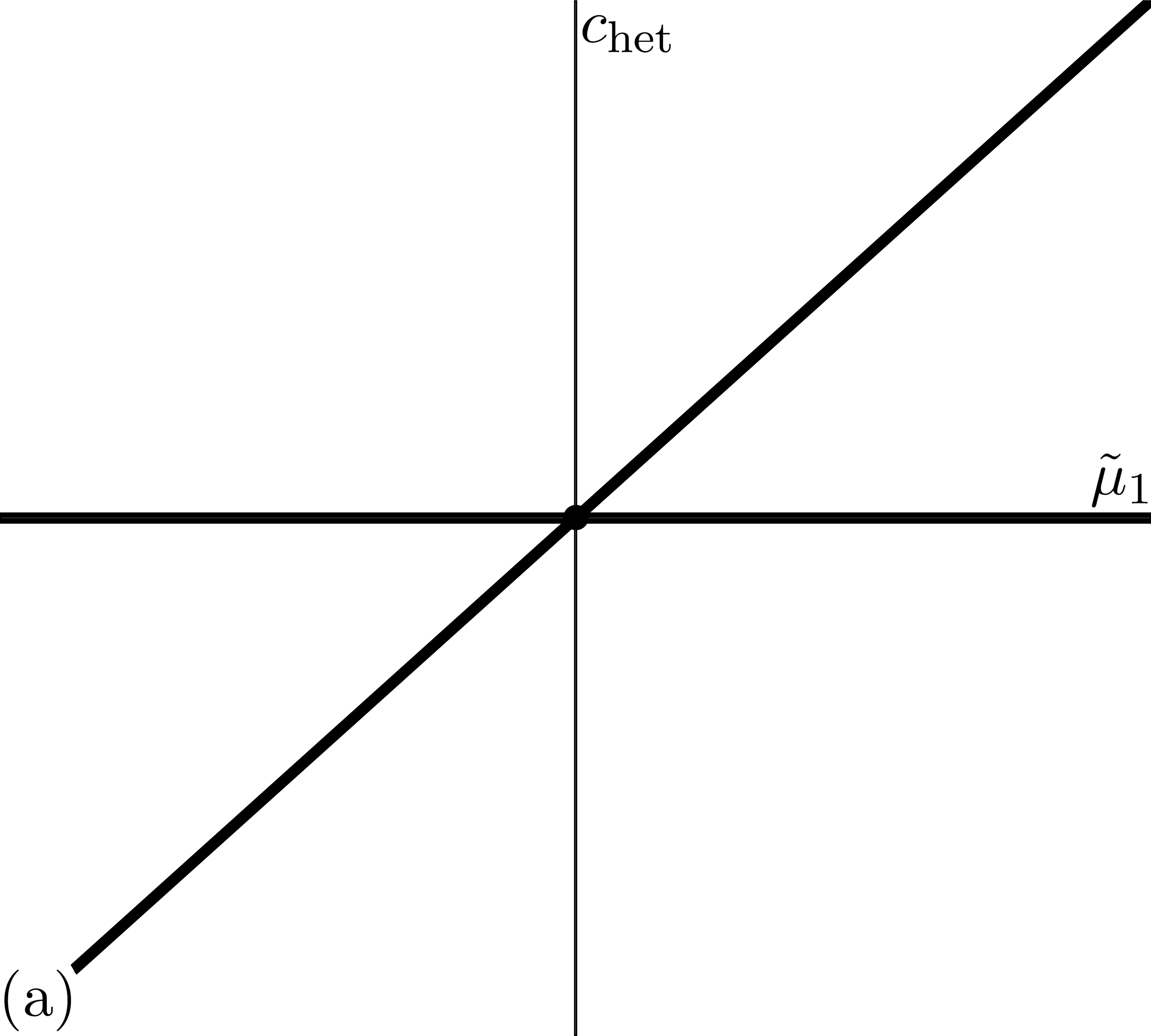}
	\end{minipage}%
	\hspace{.5cm}
	\begin{minipage}{0.30\textwidth}
		\includegraphics[width=\linewidth]{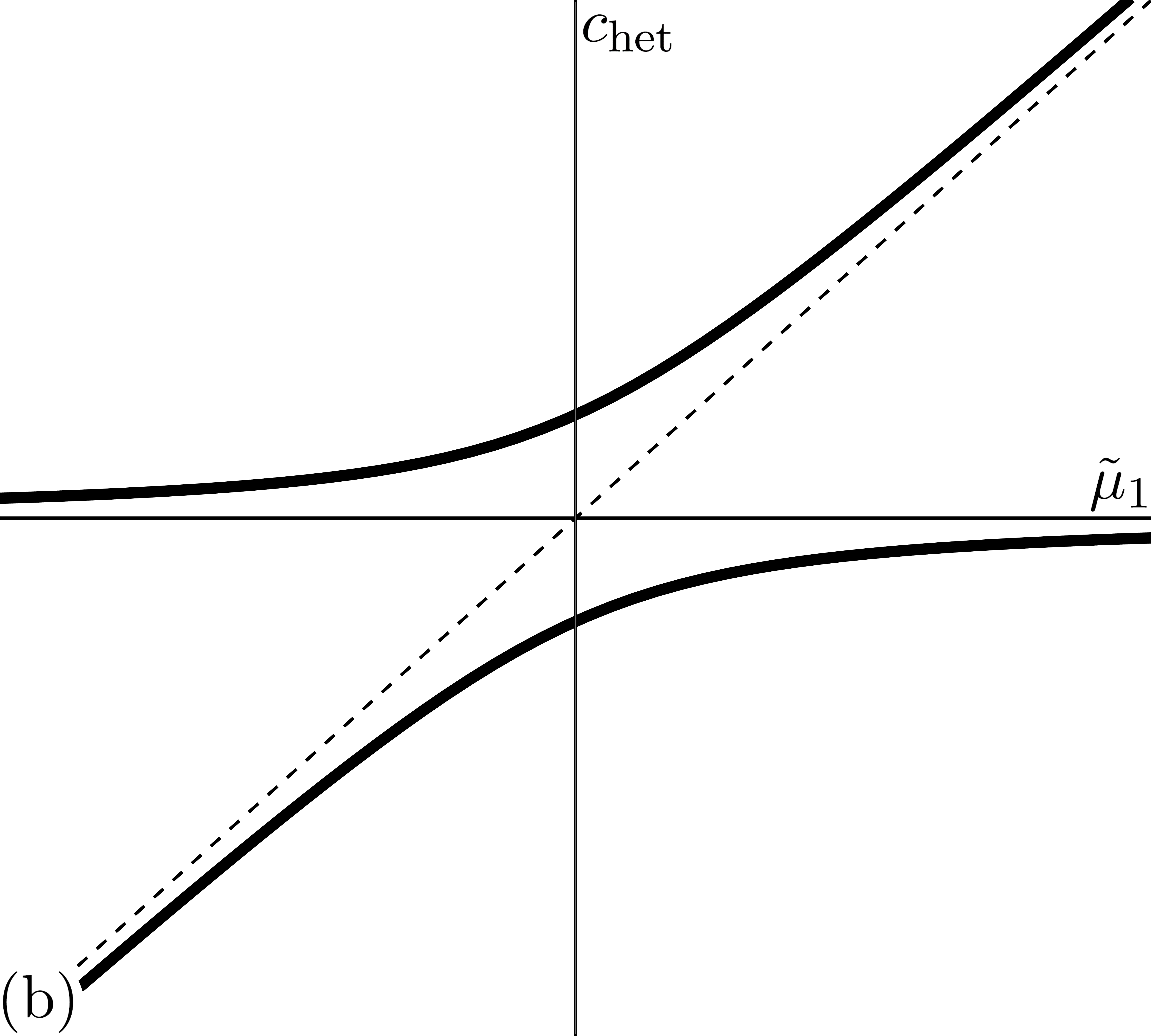}
	\end{minipage}
    \hspace{.5cm}
	\begin{minipage}{0.30\textwidth}
		\includegraphics[width=\linewidth]{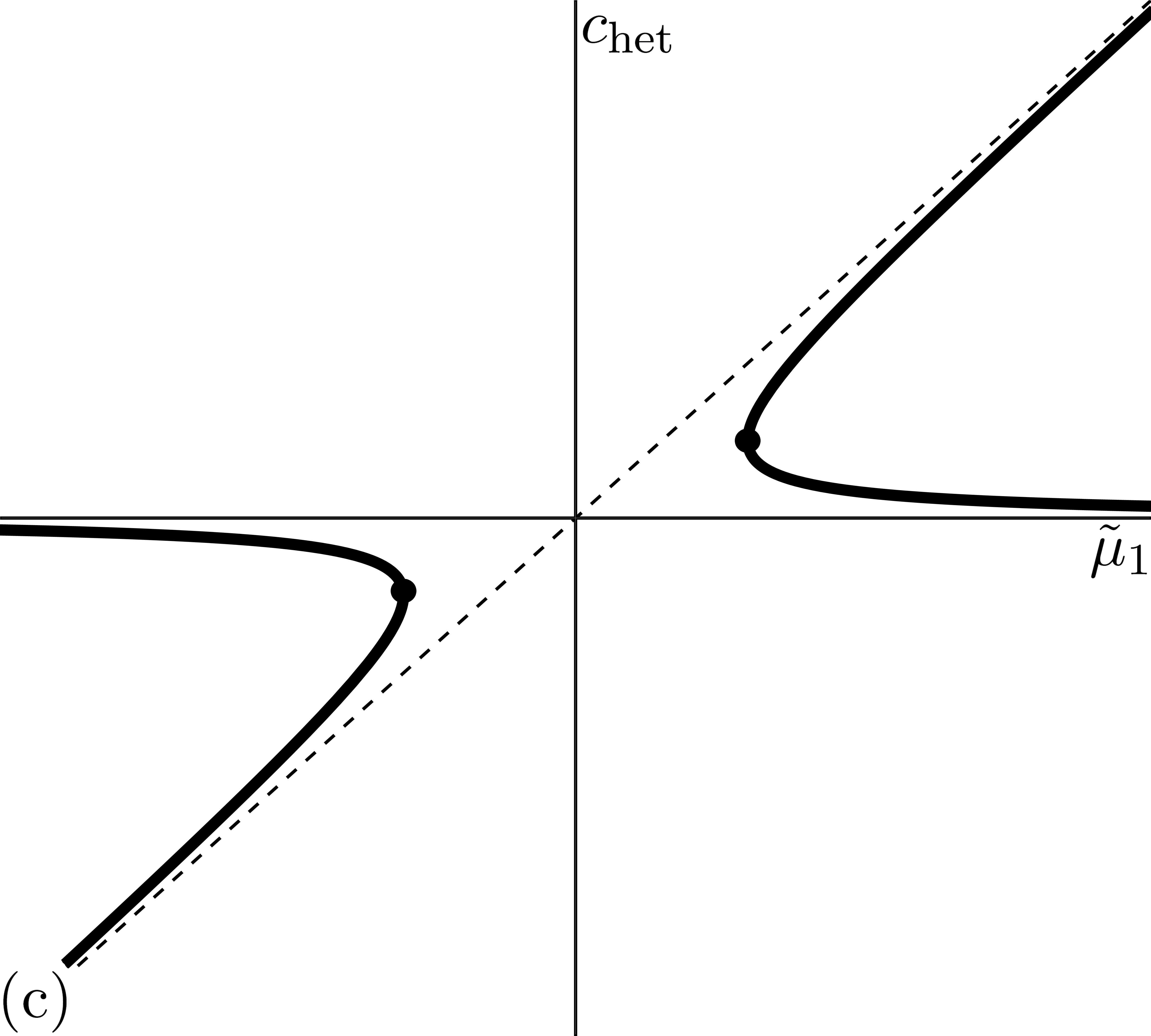}
	\end{minipage}
\caption{\small{The bifurcation diagrams described by Corollary \ref{c:hetbifs} for $\tM_{2cc}(\vmuta) < 0$. (a) Corollary \ref{c:hetbifs}(i): $\tM_{2}(\vmuta) = 0$. (b) Corollary \ref{c:hetbifs}(ii): $\tM_{2}(\vmuta) > 0$. (c) Corollary \ref{c:hetbifs}(iii): $\tM_{2}(\vmuta) < 0$}.}
\label{f:hetbifs}
\end{figure}
\noindent
This corollary is a direct consequence of Theorem \ref{t:exhets}(ii). Note that the stationary fronts of Theorem \ref{t:exhets}(i) are recovered as the $c_{\rm het}^+(\tilde{\mu}_1) \equiv 0$ branch of the case $\tM_{2}(\vmu_{\rm t}^\ast) = 0$ (Corollary \ref{c:hetbifs}(i)), which is the critical case that corresponds to the (well-known) transition as $\tM_{2}$ passes through $0$ of two saddle-node curves (Corollary \ref{c:hetbifs}(iii)) into two separate curves without bifurcations (Corollary \ref{c:hetbifs}(ii)) -- see Fig. \ref{f:hetbifs}.
\\ \\
{\bf Proof of Theorem \ref{t:exhets}.} We continue with our analysis of the expansion (\ref{d:v12q12}) -- with $\vs(X) = \vhe(X)$ -- and note that we know that the intervals of validity for the application of Poincar\'e's Expansion Theorem can be extended to the full real line $\RR$ since we are searching for solutions $\subset W^s((\oV_+,0))$. Since we already explicitly derived the first term $v_1$ (\ref{e:tv1}), we immediately continue with $v_2$, that is determined by
\beq
\label{e:Lav2}
\La v_2 = - c_0^2 \left[\tG_{1c,v}^\ast \vax \tv_1 - \tG_{1c}^\ast \tv_{1,X} -  \frac12 ((f'_\ast)^2 G^\ast_{uu}  + G^\ast_{vv} + 2 f'_\ast G^\ast_{uv} + f''_\ast G^\ast_{u})\tv_1^2\right] - (c_1 \vax \tG_{1c}^\ast + c_0^2 \tau^2 \tG_{2cc}^\ast + \tG_{2}^\ast)
\eeq
(\ref{e:SF}), (\ref{d:tGs}), (\ref{d:fFGast}). Assuming that $v_2(0) = 0$, Lemma \ref{l:Lah} again controls the leading order behavior of $v_2(X)$ for $X \gg 1$. Obviously $v_1(X) \equiv 0$ in the stationary case (\ref{e:tv1}), while the inhomogeneous righthand side of (\ref{e:Lav2}) reduces to only $-\tG_{2}^\ast$, thus it follows by (\ref{e:vgg1}) that indeed $\tM_{2}(\vmu)$ has to be set $0$ for a (potential) heteroclinic connection $\subset W^u((\oV_-,0))\cap W^s((\oV_+,0))$.
\\ \\
In the traveling case, i.e. with $\vmu = \vmu_{\rm t}^\ast + \eps \vmu_1$, it follows by (\ref{e:tv1gg1}) that $\tv_1(X)$ contributes to the -- potentially -- exponentially growing character of $\vhe(X)$ at $\O(\eps^2)$. Hence, we deduce by Lemma \ref{l:Lah} that for $X \gg1$,
\beq
\label{e:vhe-va}
\vhe(X) - \va(X) = - \frac{\eps^2 }{2\be_+ \al_+} \left( c_0 \vec{\nabla} \Ma(\vmu_{\rm t}^\ast) \cdot \vmu_1 + c_0^2 \tM_{2cc}(\vmu_{\rm t}^\ast) + \tM_{2}(\vmu_{\rm t}^\ast) \right) e^{\sqrt{\al_+} X} \left(1 + \O(E_+(X)) \right)
\eeq
at leading order (\ref{d:tM2c2}). Since we are constructing $(\vhe(X), \qhe(X)) \subset W^u((\oV_+,0))$, we know that $\lim_{X \to \infty} |\vhe(X) - \va(X)|$ must be $0$, which indeed necessarily implies (\ref{e:exhets}).
\\ \\
For both cases (i) and (ii), the proof of the Lemma now follows by the (standard, Melnikov-type) observation that $W^u((\oV_-,0))$ as spanned by $(\vhe(X), \qhe(X))$ evolves from a manifold that changes direction asymptotically close to the saddle $(\oV_+,0)$ and returns back in the direction of $(\oV_-,0)$ to a manifold that passes along the saddle $(\oV_+,0)$ (with increasing $v$-coordinate) as a parameter is changed, so that there is a unique value of that parameter for which $(\vhe(X), \qhe(X)) \subset W^u((\oV_-,0)) \cap W^s((\oV_+,0))$. In the stationary case (i), this (family of) parameter(s) is $\vmu$ and the geometric argument determines $\vmu_s(\eps)$ uniquely, in traveling case (ii) the parameter is $c$ and $c_{\rm het}(\eps)$ is determined uniquely. \hfill $\Box$
\\ \\
Like Theorems \ref{t:E-Pers} and \ref{t:exhets}, the upcoming theorem on the existence of nearly (double front) heteroclinic homoclinic orbits will be split into a stationary case and a traveling case. Before we formulate this result, we first motivate why it will now be necessary to consider $\tsi = 2$ in expansion (\ref{d:vmutsi}).
\\ \\
By symmetry (\ref{d:symm}) of (\ref{e:DS}) inherited by (\ref{e:SF}), it follows from Theorem \ref{t:exhets} that if $\tc_0$ solves
\beq
\label{d:tc1}
-\tc_0 \vec{\nabla} \Ma(\vmu_{\rm t}^\ast) \cdot \vmu_1 + \tc_0^2 \tM_{2cc}(\vmu_{\rm t}^\ast) + \tM_{2}(\vmu_{\rm t}^\ast) = 0,
\eeq
there exists a heteroclinic solution $(\vhe(X), -\qhe(-X)) \subset W^u((\oV_+,0))\cap W^s((\oV_-,0))$ of (\ref{e:SF}) -- with $-\qhe(-X) < 0$ that corresponds to a front/interface in (\ref{e:RDE}) that travels with speed $\tc_0 + \O(\eps)$ that connects $(\oU_+,\oV_+)$ for $X \to -\infty$ to $(\oU_-,\oV_-)$ and $X \to +\infty$. Now assume that $c_0 \neq 0$ is close to the solution $c_{\rm het}(0)$ of (\ref{e:exhets}), i.e. that $W^u((\oV_-,0))$ just misses $W^s((\oV_+,0))$ and that the orbit $(\vs(X),\qs(X))$ -- that spans $W^u((\oV_-,0))$ -- returns in the direction of $(\oV_-,0)$ (i.e. it enters the $\{q<0\}$ half-plane). This way, $W^u((\oV_-,0))$ can be made to be arbitrarily close to $W^u((\oV_+,0))$. Since $W^u((\oV_+,0))$ can only connect to $W^s((\oV_-,0))$ for solutions $\tc_0$ of (\ref{d:tc1}) -- i.e. (\ref{e:exhets}) with $c_0 \to -\tc_0$ -- we conclude that $W^u((\oV_-,0))$/$(\vs(X),\qs(X))$ can only be made to be homoclinic to $(\oV_-,0)$ with $c_0 \neq 0$ if $\Ma = \O(\eps^2)$: then, the asymmetric terms in (\ref{e:exhets}) and (\ref{d:tc1}) -- i.e. the terms that are not symmetric under $c_0 \to -c_0$ -- disappear from (\ref{e:exhets})/(\ref{d:tc1}) so that $W^u((\oV_-,0))$ may indeed catch up with $W^s((\oV_-,0))$ --  in fact (\ref{e:exhets}) and (\ref{d:tc1}) have become identical. We conclude that we need to set $\tsi = 2$ in (\ref{d:vmutsi}) to find orbits homoclinic to $(\oV_-,0)$.
\begin{theorem}
\label{t:exhoms}
Let double well potential $\W_0(v)$ (\ref{d:HW0}) be such that $H_{0,+} = 0$ (\ref{d:H0pm}) and let $\eps > 0$ be sufficiently small.
\\
{\bf (i) Stationary pulses.} Let $\vmu \in \RR^m$ be such that $\tM_{2}(\vmu) > 0$, then there exists a homoclinic solution $(\vho(X), \qho(X)) \subset W^u((\oV_-,0))\cap W^s((\oV_-,0))$ of (\ref{e:SF}) with $c=0$ (that merges with the heteroclinic cycle established by Theorem \ref{t:exhets}(i) as $\tM_2 \downarrow 0$).
\\
{\bf (ii) Traveling pulses.} Let $\vmuta$ be such that $\Ma(\vmu_{\rm t}^\ast) = 0$, let $\vmu = \vmu_{\rm t}^\ast + \eps^2 \vmu_2$ so that (\ref{d:vmutsi}) holds with $\tsi = 2$ -- assuming $\vec{\nabla} \Ma(\vmu_{\rm t}^\ast) \neq \vec{0}$ -- and let $c$ be such that $c^2 \tM_{2cc}(\vmu_{\rm t}^\ast) + \tM_{2}(\vmu_{\rm t}^\ast) > 0$ (\ref{d:tM2c2}). Then there is a co-dimension 1 manifold $\V_{\rm t}(c)$ in $\vmu_2$-space, such that there exists a homoclinic orbit $(\vho(X), \qho(X)) \subset W^u((\oV_-,0))\cap W^s((\oV_-,0)) \subset \M_\eps(c)$ in (\ref{e:SF}) for all $\vmu_2 \in \V_{\rm t}(c)$. Moreover, by the symmetry (\ref{d:symm}), the orbit $(\vho(X), \qho(X))$ with $c>0$ is symmetric to the orbit with $\tilde{c} = -c < 0$.
\\
In both cases, the orbits $(\vho(X), \qho(X))$ correspond to slow orbits on $\M_\eps(c)$ homoclinic to the critical point $(\oU_-,0,\oV_-,0)$ in (\ref{e:DS}) and to stationary or traveling localized pulses/stripes $(U(x,y,t),V(x,y,t)) = (\Uho(X),\Vho(X))$ in (\ref{e:RDE})/(\ref{e:RDE-S}) -- with $\Vho(X) = \vho(X)$ and $\Uho(X) = \uho(X) = f(\vho(X)) + \eps f_1(\vho(X),\qho(X)) + \eps^2 f_2(\vho(X),\qho(X)) + \O(\eps^3)$ (\ref{e:Meps}) (\ref{d:fjpj}) -- Fig. \ref{f:frontspulses}(b).
\end{theorem}
\noindent
Note that Theorem \ref{t:exhoms} is not a persistence result, unlike Theorems \ref{t:E-Pers} and \ref{t:exhets}: the slow reduced limit system (\ref{e:RedSF}) for the double well potential with wells of equal depth does not have homoclinic orbits. Moreover, we focused -- for simplicity -- on orbits homoclinic to $(\oV_-,0)$, a completely analogous result on the existence of orbits homoclinic to $(\oV_+,0)$ holds true under the conditions $\tM_{2}(\vmu_{\rm t}^\ast) < 0$ (i) and $c_0^2  \tM_{2cc}(\vmu_{\rm t}^\ast) + \tM_{2}(\vmu_{\rm t}^\ast) < 0$ (ii) (with an equivalent analogous version of Corollary \ref{c:hombifs}). Note also that the formulation of Theorem \ref{t:exhoms}(ii) is somewhat nonnatural: it starts out with a given $c$ and establishes that there is a co-dimension 1 manifold of $\vmu_2$-values for which there is a traveling pulse with this pre-specified speed. It is more natural (and standard) to consider the value of the parameter, i.e. of $\vmu_2$, as given and to ask the question: are there values of $c$ for which there exist nearly heteroclinic homoclinic pulses in (\ref{e:DS})? See Corollary \ref{c:hombifs}.
\\ \\
{\bf Proof of Theorem \ref{t:exhoms}.} This proof consists of two parts. First, we follow the arguments that led to (\ref{e:Xh-0}) and (\ref{e:vaXh}) and establish a leading order approximation of the pulse solutions, for $\vmu = \vmu_{\rm t}^\ast + \eps^2 \vmu_2$ (cf. (\ref{d:vmutsi}) with $\tsi = 2$)). In fact, it is shown that $|\Delta W_{\rm h}(\vmu_2)| \ll \eps$, where $\Delta W_{\rm h}(\vmu_2)$ is the gap between $W^u((\oV_-,0))$ and $W^s((\oV_-,0))$ as they intersect the $\{q=0\}$-axis (in (\ref{e:SF})). In the second part of the proof we study the higher order problems and show that $\Delta W_{\rm h}(\vmu_2) = \O(\eps^2 |\log \eps|)$ and can indeed be closed by varying $\vmu_2$. Many of the results obtained in this part of the proof will return as key ingredients in (the proofs of) the stability results in section \ref{s:Stab} -- especially in Lemma \ref{l:asymptevhom-hot} in section \ref{sss:neartravstab}.
\\ \\
By Lemma \ref{l:Lah}, we know that solution $v_2$ of (\ref{e:Lav2}) is for $X \gg 1$ given by,
\beq
\label{e:v2gg1}
v_2(X) = - \frac{1}{2\al_+ \be_{+}} \left(c_0^2 \tM_{2cc}(\vmu_{\rm t}^\ast) + \tM_{2}(\vmu_{\rm t}^\ast)\right) e^{\sqrt{\al_+} X}(1 + \O(E_{+}))
\eeq
(cf. (\ref{e:vhe-va})), both in case (i) -- in which $c_0 = 0$, so that $v_1(x) \equiv 0$ (\ref{e:tv1}) -- and in case (ii) -- where all terms containing $\tv_1(X)$ in the right hand side of (\ref{e:Lav2}) do not appear in (\ref{e:v2gg1}) by assumption (\ref{d:vmutsi}) with $\tsi = 2$. A solution to (\ref{e:SF}) that is homoclinic to $(\oV_-,0)$ must intersect the $\{q=0\}$-axis, i.e. there must be an $\Xh$ such that $\vhox(\Xh) = 0$. Since $\vho(X) = \va(X)$ at $\O(1)$, we find for $X \gg 1$ (using  (\ref{e:propsvbu})),
\beq
\label{e:vhogg1}
\vhox(X) = \left[\be_{+} \sqrt{\al_+} e^{-\sqrt{\al_+} X} - \frac{\eps^2}{2\be_{+} \sqrt{\al_+}} \left(c_0^2 \tM_{2cc}(\vmu_{\rm t}^\ast) + \tM_{2}(\vmu_{\rm t}^\ast)\right) e^{\sqrt{\al_+} X}\right](1 + \O(E_{+})) + \O(\eps^3 v_3)
\eeq
(cf. (\ref{e:vsgg1})). Clearly, $\vhox(\Xh)$ can be $0$ if $e^{\sqrt{\al_+} \Xh} = \O(1/\eps)$, however this is the singular situation $\si = 1$ in (\ref{e:extPoincare}): it a priori goes beyond the domain of validity of (the extended) Poncar\'e's Expansion Theorem (note also that the leading order correction term $\O(\eps^3 v_3)$ in (\ref{e:vhogg1}) a priori is $\O(1)$ by (\ref{e:extPoincare}) (with $\si = 1$) and thus not asymptotically small). However, the estimates of (\ref{e:extPoincare}) are obtained in the generic situation in which the exponential growth of the inhomogeneous term $h_2(X)$ of the $v_2$ equation is driven by $v_1^2(X))$ so that $h_2(X)$ grows exponentially as $E_+^{-2}(X)$ for $X \gg 1$. Due to assumption (\ref{d:vmutsi}) this is not the case here, as was also already noticed in the proof of Theorem \ref{t:exhets}: for $\tsi = 2$ the growth of $v_1$ has an additional factor $\eps^2$ -- cf. (\ref{e:tv1gg1}) in the $\tsi = 1$ case -- the effect of the $v_1^2$ terms `falls down' to the level of the inhomogeneous term $h_6(X)$ of the $v_6$ equation. In fact, we see from (\ref{e:v2gg1}) that $v_2(X)$ grows as $E_+^{-1}(X)$ for $X \gg 1$ (and thus not as $E_+^{-2}(X)$). Similarly, it follows that $v_3(X)$ will also only grow as $E_+^{-1}(X)$ for $X \gg 1$ -- much slower than the $E_+^{-3}(X)$ growth of the generic setting that yielded (\ref{e:extPoincare}). By Lemma \ref{l:Lah-sharp}, it is straightforward to deduce iteratively that $v_4(X)$ and $v_5(X)$ will grow as $E_+^{-2}(X)$ for $X \gg 1$, $v_6(X)$ and $v_7(X)$ as $E_+^{-3}(X)$, etc.. Thus, we may conclude that we can extend the domain of validity of Poincar\'e's Expansion Theorem to $e^{\sqrt{\al_+} \Xh} = \O(1/\eps)$ under assumption (\ref{d:vmutsi}) and that the correction term $\O(\eps^3 v_3)$ in (\ref{e:vhogg1}) indeed is asymptotically small -- in fact, it is $\O(\eps^2)$.
\\ \\
Hence, we may conclude from (\ref{e:vhogg1}) that
\beq
\label{e:Xh-1}
e^{\sqrt{\al_+}\Xh} = \frac{\be_+ \sqrt{2 \al_+}}{\sqrt{c_0^2 \tM_{2cc}(\vmu_{\rm t}^\ast) + \tM_{2}(\vmu_{\rm t}^\ast)}} \, \frac{1}{\eps} \left(1 + \O(\eps) \right)
\eeq
(cf. (\ref{e:Xh-0})) which yields for case (i) that indeed $\tM_{2}(\vmu_{\rm t}^\ast)$ must positive and for case (ii) that $c_0$ must be so that $c_0^2 \tM_{2cc}(\vmu_{\rm t}^\ast) + \tM_{2}(\vmu_{\rm t}^\ast) > 0$. As a brief side remark we note that,
\beq
\label{e:Xh-2}
X_h = \frac{|\log \eps|}{\sqrt{\al_+}}  + \frac{1}{\sqrt{\al_+}} \log \left( \frac{\be_+ \sqrt{2 \al_+}}{\sqrt{c_0^2 \tM_{2cc}(\vmu_{\rm t}^\ast) + \tM_{2}(\vmu_{\rm t}^\ast)}} \right) + \O(\eps),
\eeq
which will be the source of many $|\log \eps|$ terms in the upcoming analysis. We continue by observing that the $v$-coordinate of $(\vho(X),\qho(X)) \cap \{q=0\}$, i.e. of $W^u((\oV_-,0)) \cap \{q=0\}$, is at leading order given by
\beq
\label{e:vaXhnearhet}
\vho(\Xh) =
\oV_{+} - \frac{\sqrt{2\left(c_0^2 \tM_{2cc}(\vmu_{\rm t}^\ast) + \tM_{2}(\vmu_{\rm t}^\ast)\right)}}{\sqrt{\al_+}} \, \eps
\eeq
(cf. (\ref{e:vaXh})). Since (\ref{e:vaXhnearhet}) only contains quadratic terms in $c_0$, it follows by the symmetry (\ref{d:symm}) of (\ref{e:DS}) that $W^s((\oV_-,0)) \cap \{q=0\}$ is -- at leading order -- also given by (\ref{e:vaXhnearhet}): the gap $\Delta W_{\rm h}(\vmu_2)$ between $W^u((\oV_-,0))$ and $W^s((\oV_-,0))$ as they intersect the $\{q=0\}$-axis is smaller than $\O(\eps)$. Next, we need to consider cases (i) and (ii) separately.
\\ \\
Since the (planar) flow on $\M_\eps(0)$ is reversible and integrable (\ref{e:SFc=0H0}), (\ref{d:H2}) in case (i) ($c = 0$), it immediately follows that if $\tM_{2}(\vmu_{\rm t}^\ast) > 0$, the family of periodic orbits of the integrable flow around the center $(\oV_c, 0)$ must be encircled by a homoclinic orbit $(\vho(X), \qho(X)) \subset W^u((\oV_-,0))\cap W^s((\oV_-,0))$.
\\ \\
For case (ii), we a need to go one step deeper into the perturbation analysis and determine the next order correction to (\ref{e:vaXhnearhet}). This term is made up of several contributing terms. A priori, we need to solve both the equations for $v_3(X)$ and $v_4(X)$ since they both contribute to the $O(\eps^2)$ level -- see the discussion above. Moreover, we need to evaluate the $O(\eps^2)$ effect of the $v_1(X)$ term on the intersection of $W^u((\oV_-,0))$ with the $\{q=0\}$-axis (cf. (\ref{e:vhe-va})).
\\ \\
Note that although the gap $\Delta W_{\rm h}(\vmu_2)$ between $W^u((\oV_-,0)) \cap \{q=0\}$ and $W^s((\oV_-,0)) \cap \{q=0\}$ can only be tuned -- and thus closed -- by terms that are odd in $c$ (the even powers of $c$ do not have an impact by symmetry (\ref{d:symm})), this does not imply that the orbits we are constructing will be symmetric under (\ref{d:symm}): they clearly cannot be so (for $c \neq 0$), since the approximation of $W^s((\oV_-,0))$ is obtained from that of $W^u((\oV_-,0))$ by first applying (\ref{d:symm}) followed by changing $-c$ back into $c$. Nevertheless, a pulse traveling with speed $c \neq 0$ must have its counterpropagating symmetrical counterpart under (\ref{d:symm}). Since the upcoming arguments establish that for given $c \neq 0$ the intersection $W^u((\oV_-,0)) \cap W^s((\oV_-,0))$ determines a unique homoclinic orbit, it follows that the orbit with speed $-c$ must be symmetric to it (under (\ref{d:symm})).
\\ \\
We conclude from (\ref{e:Lav1}) and (\ref{e:Lav2}) that the $c$-dependence of $v_1(X)$ and $v_2(X)$ can be factored out explicitly
\[
v_1(X) = c \tv_1, \; \; v_2 = \tv_2 + c^2 \tv_{2cc}
\]
with $c = c_0 + \eps c_1 + \eps^2 c_2 + $ etc., cf. (\ref{e:tv1}), (\ref{d:tv2s}). Similarly, it can be deduced that
\beq
\label{d:expv34}
v_3 = c \tv_{3c} + c^3 \tv_{3ccc}, \; \; v_4 = \tv_4 + c^2 \tv_{4cc} + c^4 \tv_{4cccc}, \; \; {\rm etc.},
\eeq
which implies that of these two, only $v_3(X)$ contributes to the gap $\Delta W_{\rm h}$ at its $\O(\eps^2)$ level. Naturally also  $v_1$ contributes to $\Delta W_{\rm h}$ (\ref{e:tv1}). In fact, we will find that due to the $|\log \eps|$-term introduced by (\ref{e:Xh-2}), the contribution of $v_1$ to $\Delta W_{\rm h}$ is dominant and that the contribution by $v_3$ reduces to a higher order effect (see however the discussion in section \ref{ss:Obs}). Lemma \ref{l:Lah} does not provide a sufficiently precise approximation of $v_1(\Xh)$, therefore we apply its refinement Lemma \ref{l:Lah-sharp}. First, we note that (\ref{d:tv1}) is of the form (\ref{e:Lah}) with inhomogeneous term $h(x)$ satisfying assumptions (\ref{d:oh01j}) with $j = -1$ (\ref{e:propsvbu}). In fact $M_h = -\Ma(\vmu)$ and $\oh_{0,-1} = - \be_+ \sqrt{\al_+} \oG_{1c}^+ $ with,
\beq
\label{d:oG1+}
\oG_{1c}^+(\vmuta)= \lim_{X \to \infty} \tG_{1c}(\va(X)) = 1- \tau \frac{f'(\oV_+;(\vmuta))\oG^+_u(\vmuta)}{\oF^+_u(\vmuta)}
\eeq
(\ref{d:tG1c-I}), (\ref{d:oFpmetc}). Moreover, $M^+_{{\rm u},-1} = - \tG^+_{{\rm u},-1}$ (\ref{d:Mpubpm1}) with
\[
\tG^+_{{\rm u},-1} = \int_0^{\infty} \left[\tG_{1c}(\va) \vax \vu - \frac{1}{2 \sqrt{\al_+}} \oG_{1c}^+ \right] d\tX.
\]
Hence, it follows by (\ref{d:Ma}), (\ref{d:vmutsi}) from (\ref{d:approxsvX-sharp}) for $j=-1$ that for $X \gg 1$,
\beq
\label{e:tv1Xgg1-next}
\tv_1(X) = \frac{\be_+ \oG_{1c}^+}{2} X e^{-\sqrt{\al_+}X} - \frac{\vec{\nabla} \Ma \cdot \vmu_{2}}{2 \al_+ \be_+} \eps^2 e^{\sqrt{\al_+}X} + \frac{\oG_{1c}^+ + 4 \al_+ \tG^+_{{\rm u},-1}}{4 \sqrt{\al_+}} e^{-\sqrt{\al_+}X} + \O(X E_+^2),
\eeq
so that by (\ref{e:Xh-1}),
\[
\eps v_1(\Xh) =  \eps c_0 \tv_1(\Xh) = \frac{c_0 \sqrt{\tM_{2} + c_0^2 \tM_{2cc}}}{2 \sqrt{2 \al_+}}
\left[
\oG_{1c}^+ \Xh - \frac{2 \vec{\nabla} \Ma \cdot \vmu_{2}}{\tM_{2} + c_0^2 \tM_{2cc}} +
\frac{\oG_{1c}^+ + 4 \al_+ \tG^+_{{\rm u},-1}}{2 \sqrt{\al_+}}
\right] \eps^2 + \O(\eps^3|\log \eps|),
\]
where we note that the $|\log \eps|$ factor in the leading order correction term originates from the $X = \O(|\log \eps|)$ factor in the $\O(X E_+^2)$ term of (\ref{d:approxsvX-sharp}) -- see (\ref{e:Xh-2}). In fact, we may simplify the present -- and upcoming -- analysis considerably by only considering the leading order approximation (\ref{e:Xh-2}) of $\Xh$,
\beq
\label{e:epsv1-leading}
\eps v_1(\Xh) = \eps c_0 \tv_1(\Xh) =
\frac{c_0 \sqrt{\tM_{2} + c_0^2 \tM_{2cc}}}{2 \sqrt{2 \al_+}}
\left[
\frac{\oG_{1c}^+}{\sqrt{\al_+}} - \frac{2 \tmu_{2}}{\tM_{2} + c_0^2 \tM_{2cc}} \right] \, \eps^2 |\log \eps| + \O(\eps^2),
\eeq
where we have extended the range of $\vmu_2$ by introducing $\tilde{\vmu}_2 \in \RR^m$ and $\tmu_2 \in \RR$ by
\beq
\label{d:tmu2}
\vmu = \vmuta + \eps^2 \vmu_2 = \vmuta + \eps^2 |\log \eps| \, \tilde{\vmu}_2, \; \; \tmu_2 = \vnab M_\ast(\vmuta) \cdot \tilde{\vmu}_2
\eeq
(cf. (\ref{d:vmutsi}) with $\tsi = 2$ and recall that we have assumed that $\vec{\nabla} \Ma(\vmu_{\rm t}^\ast) \neq \vec{0}$); note that this -- logarithmically extending the range of $\vmu_2$ -- has no impact on the preceding analysis.
\\ \\
To determine the leading order correction to $\vho(\Xh)$ (\ref{e:vhogg1}) -- and thus to the gap $\Delta W_{\rm h}$ -- we a priori need to take the next order correction of $\Xh$ into account. Therefore, we define $\Yh$ and its expansion in $Y_j$'s, $j \geq 0$, by
\beq
\label{d:Yhj}
e^{\sqrt{\al_+}\Xh} = \be_+ \frac{\Yh}{\eps} = \be_+ \left(Y_0 + \eps Y_1 + \O(\eps^2) \right) \frac{1}{\eps},
\eeq
so that,
\beq
\label{e:tY0}
Y_0 = \frac{\sqrt{2 \al_+}}{\sqrt{\tM_{2}(\vmu_{\rm t}^\ast) + c_0^2 \tM_{2cc}(\vmu_{\rm t}^\ast)}}
\eeq
(\ref{e:Xh-1}). By (\ref{e:vaXgg1}), (\ref{e:v2gg1}), (\ref{e:tv1Xgg1-next}) and (\ref{d:tmu2}) we have for $X \gg 1$,
\beq
\label{e:vhogg1-next}
\vho(X) = \oV_+ - \be_+ e^{-\sqrt{\al_+} X} - \frac{\tM_{2} + c_0^2 \tM_{2cc}}{2 \al_+ \be_{+}} \eps^2 e^{\sqrt{\al_+} X} + \frac{c_0 \be_+ \oG_{1c}^+}{2} \eps X e^{-\sqrt{\al_+}X} - \frac{c_0 \tmu_{2}|\log \eps|}{2 \al_+ \be_+} \eps^3 e^{\sqrt{\al_+}X},
\eeq
with correction terms of $\O(\eps E_+)$, $\O(E_+^2)$, $\O(\eps^2)$ and $\O(\eps^3 E_+^{-1})$ -- that are all $\O(\eps^2)$ for $X = \Xh$ (\ref{e:Xh-1}). Thus, by (\ref{e:Xh-2}), (\ref{d:Yhj}) and (\ref{e:tY0}),
\[
\vhox(\Xh) = - \left[ \frac{\tM_{2} + c_0^2 \tM_{2cc}}{\sqrt{\al_+}} \frac{Y_1}{|\log \eps|} + \frac12 c_0 \left( \frac{\oG_{1c}^+}{Y_0} + \frac{\tmu_2 Y_0}{\sqrt{\al_+}}\right)\right] \eps^2|\log \eps| + \O(\eps^2),
\]
which yields,
\beq
\label{e:tY1}
Y_1 = - \frac12 c_0 \sqrt{2} \, \frac{\tmu_2 \sqrt{\al_+} + \frac12 \oG_{1c}^+ (\tM_{2} + c_0^2 \tM_{2cc})}{(\tM_{2} + c_0^2 \tM_{2cc})\sqrt{\tM_{2} + c_0^2 \tM_{2cc}}} \, |\log \eps| + \O(1) \stackdef c_0 \tY_1 |\log \eps| + \O(1).
\eeq
Note that $Y_1 \to - Y_1$ as $c_0 \to -c_0$ (at leading order): the correction to $\Xh$ is different for the part of the orbit $(\vho(X), \vhox(X))$ with $\vhox(X)> 0$ from that with $\vhox(X)< 0$ -- as expected. Substitution of all of this into (\ref{e:vhogg1-next}) yields,
\[
\vho(\Xh) = \oV_+ - \frac{\sqrt{2(c_0^2 \tM_{2cc} + \tM_{2})}}{\sqrt{\al_+}} \, \eps
- \frac12 c_0 \sqrt{2} \, \frac{\tmu_2 \sqrt{\al_+} - \frac12 \oG_{1c}^+ (\tM_{2} + c_0^2 \tM_{2cc})}{\al_+ \sqrt{\tM_{2} + c_0^2 \tM_{2cc}}} \, \eps^2|\log \eps| + \O(\eps^2)
\]
(cf. (\ref{e:vaXhnearhet})), where we notice that the contributions by $Y_1$ dropped out completely (at leading order). Thus,
the gap $\Delta W_{\rm h}$ between $W^u((\oV_-,0))$ and $W^s((\oV_-,0))$ as they intersect the $\{q=0\}$-axis is at leading order given by
\beq
\label{e:gapOeps2inhomasymm}
\Delta W_{\rm h}(\tmu_2) =
c_0 \sqrt{2} \, \frac{\tmu_2 \sqrt{\al_+} - \frac12 \oG_{1c}^+ (\tM_{2} + c_0^2 \tM_{2cc})}{\al_+ \sqrt{\tM_{2} + c_0^2 \tM_{2cc}}} \, \eps^2|\log \eps| + \O(\eps^2),
\eeq
which is in fact identical to (\ref{e:epsv1-leading}), up to a factor $-2$: the leading order contribution to the gap $\Delta W_{\rm h}(\tmu_2)$ is completely determined by $v_1(\Xh)$. For any given $c_0$ we can clearly tune $\tmu_2$ such that $\Delta W_{\rm h}(\tmu_2) = 0$ -- see also (\ref{d:Vtc}) below. Moreover, it follows from the geometry of the set-up that this quantity is determined uniquely (see however section \ref{ss:Obs}). Hence, we may conclude that for any given $c$ such that $c_0^2 \tM_{2cc}(\vmu_{\rm t}^\ast) + \tM_{2}(\vmu_{\rm t}^\ast) > 0$, there (locally) is a co-dimension 1 manifold in parameter space for which the gap $\Delta W_{\rm h}$ is closed, and thus that there is a homoclinic orbit $(\vho(X), \qho(X)) \subset W^u((\oV_-,0))\cap W^s((\oV_-,0))$ on $\M_\eps(c)$ for these values of $\vmu = \vmu_{\rm t}^\ast + \eps^2 \tilde{\vmu}_2 |\log \eps|$.  \hfill $\Box$
\\ \\ 
Note that we may conclude from the leading order approximation (\ref{e:gapOeps2inhomasymm}) of the gap $\Delta W_{\rm h}(\tmu_2)$ that the co-dimension 1 manifold $\V_{\rm t}(c)$ in the statement of Theorem \ref{t:exhoms}(ii) is at leading order (in $\eps$) given by the relatively simple hyperplane
\beq
\label{d:Vtc}
\V_{\rm t}(c) = \{\vmu_2 \in \RR^m: \vec{\nabla} \Ma(\vmu_{\rm t}^\ast) \cdot \vmu_2 =
\frac{\oG_{1c}^+(\tM_{2}(\vmu_{\rm t}^\ast) + c_0^2 \tM_{2cc}(\vmu_{\rm t}^\ast))}{2 \sqrt{\al_+}}|\log \eps| + \O(1)\}
\eeq
(\ref{d:tM2c2}), (\ref{d:oG1+}), (\ref{d:tmu2}). This result also provides the setting to consider the (natural) question whether for a give choice of the parameter $\vmu_2$ there are values of $c$ for which nearly heteroclinic traveling pulse patterns exist in (\ref{e:DS})/(\ref{e:RDE}).
\begin{corollary}
\label{c:hombifs}
Consider the setting of Theorem \ref{t:exhoms}, assume that $\|\vmu_2\| = \O(|\log \eps|)$, i.e. that $\|\vmu - \vmu_{\rm t}^\ast\| = \O(\eps^2 |\log \eps|)$ (cf. (\ref{d:vmutsi}) with $\tsi = 2$) and consider $\tmu_2$ as introduced in (\ref{d:tmu2}). Impose the additional non-degeneracy conditions $\tM_{2cc}(\vmuta) \neq 0$ (\ref{d:tM2c2}) and $\oG_{1c}^+(\vmuta) \neq 0$ (\ref{d:oG1+}) and define $C_{\rm hom}$ and $\tmu_{\rm hom}^{\rm TW}$ by
\beq
\label{d:Ctmuhet}
C_{\rm hom}(\tmu_2) = \frac{2 \sqrt{\al_+} \, \tmu_2 - \oG_{1c}^+ \tM_{2}}{\oG_{1c}^+ \tM_{2cc}}, \; \;
\tmu_{\rm hom}^{\rm TW} = \frac{\oG_{1c}^+\tM_{2}}{2 \sqrt{\al_+}}.
\eeq
{\bf (i) $\tM_{2} > 0$ and $\tM_{2cc} > 0$:} for any $\tmu_2 = \O(1)$ there is a {\rm stationary} homoclinic pulse solution $(\Uho(X),\Vho(X))$ (as established by Theorem \ref{t:exhoms}(i)); as $\tmu_2$ increases for $\oG_{1c} >0$, resp. decreases for $\oG_{1c} <0$, through $\tmu_{\rm hom}^{\rm TW}$ -- at leading order given by (\ref{d:Ctmuhet}) -- a {\rm bifurcation into traveling waves} takes place from which two counterpropagating pulses $(\Uho(X),\Vho(X))$ appear; these pulses are symmetric under (\ref{d:symm}) as homoclinic orbits in (\ref{e:DS}) with $c = c_{\rm hom}^\pm(\tmu_2)$ at leading order given by $c_0 = \pm \sqrt{C_{\rm hom}(\tmu_2)}$.
\\
{\bf (ii) $\tM_{2} > 0$ and $\tM_{2cc} < 0$:} for any $\tmu_2 = \O(1)$ there is a stationary homoclinic pulse solution $(\Uho(X),\Vho(X))$. For $0 < \tmu_2 < \tmu_{\rm hom}^{\rm TW}$ and $\oG_{1c} >0$, resp. $\tmu_{\rm hom}^{\rm TW} < \tmu_2 < 0$ and $\oG_{1c} <0$, there exist two counterpropagating pulses $(\Uho(X),\Vho(X))$ in (\ref{e:RDE})/(\ref{e:RDE-S}) traveling with leading order speed $c_0 = \pm \sqrt{C_{\rm hom}(\tmu_2)}$ (and symmetric under (\ref{d:symm})); as $\tmu_2$ decreases, resp. increases, through $\tmu_{\rm hom}^{\rm TW}$, these pulses appear from a bifurcation into traveling waves, as $\tmu_2$ decreases, resp. increases, further and $\tmu_2 \to 0$, the homoclinic pulses merge into a heteroclinic cycle and split up into two pairs of heteroclinic fronts traveling with distinct speeds $c_{\rm het}$, close to the bifurcation each at leading order given by $c_0 = \pm \sqrt{C_{\rm hom}(0)} = \pm \sqrt{-\tM_{2}/\tM_{2cc}} \stackdef \pm c_0^{\rm h-c}$ (cf. Corollary \ref{c:hetbifs}(ii)).
\\
{\bf (iii) $\tM_{2} < 0$ and $\tM_{2cc} > 0$:} for any $\tmu_2 > 0$ and $\oG_{1c} >0$, resp. $\tmu_2 < 0$ and $\oG_{1c} <0$, there exist two counterpropagating pulses $(\Uho(X),\Vho(X))$ in (\ref{e:RDE})/(\ref{e:RDE-S}) traveling with leading order speed $c_0 = \pm \sqrt{C_{\rm hom}(\tmu_2)}$; as $\tmu_2 \to 0$, the homoclinic pulses split up in a two pairs of heteroclinic fronts traveling with (distinct) speeds (near the bifurcaton) at leading order given by $c_0 = \pm c_0^{\rm h-c}$.
\\
{\bf (iv) $\tM_{2} < 0$ and $\tM_{2cc} < 0$:} there are no traveling pulse solutions of the type established by Theorem \ref{t:exhoms}(ii).
\end{corollary}
\begin{figure}[t]
\centering
	\begin{minipage}{.3\textwidth}
		\centering
		\includegraphics[width =\linewidth]{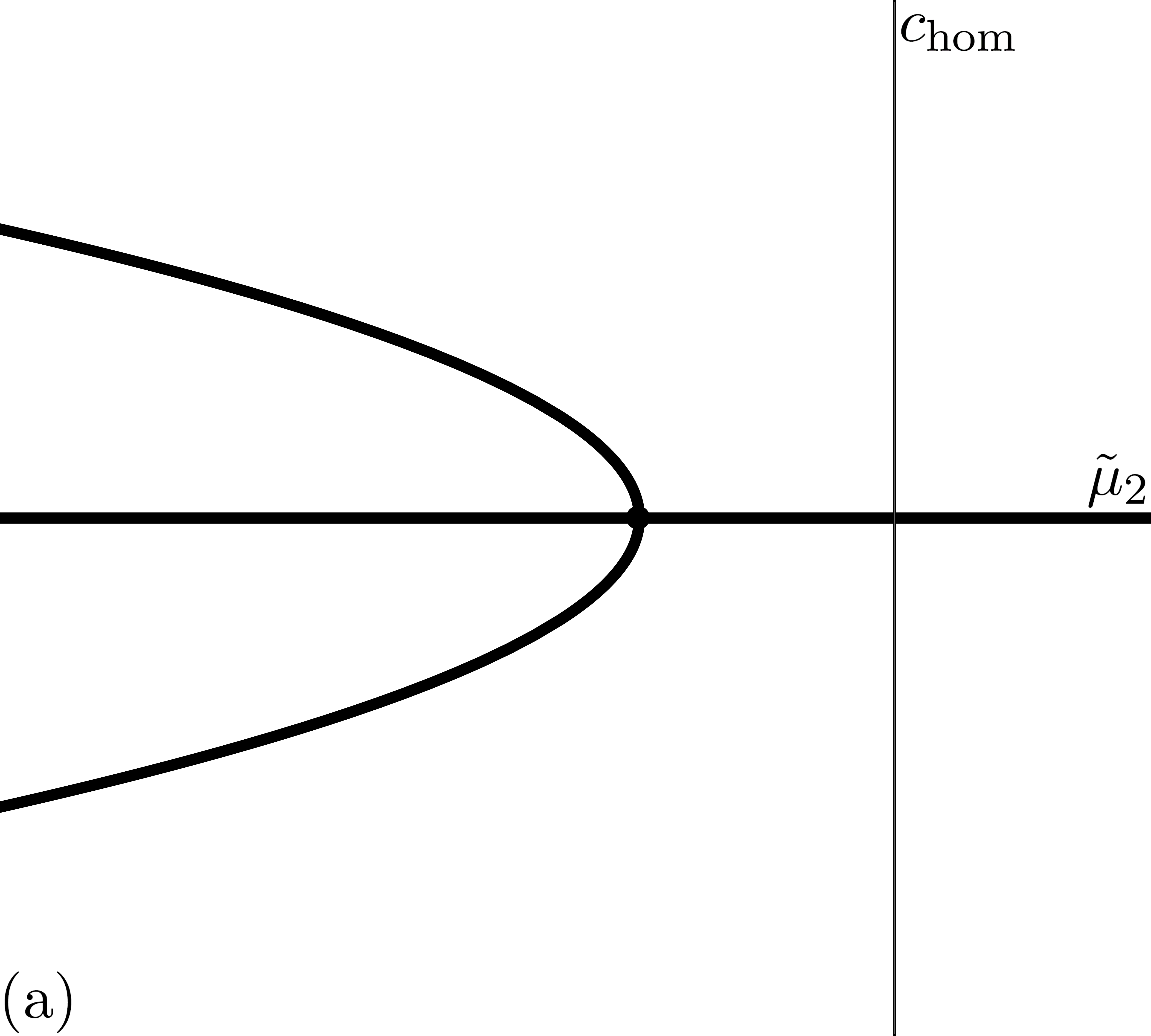}
	\end{minipage}%
	\hspace{.5cm}
	\begin{minipage}{0.3\textwidth}
		\includegraphics[width=\linewidth]{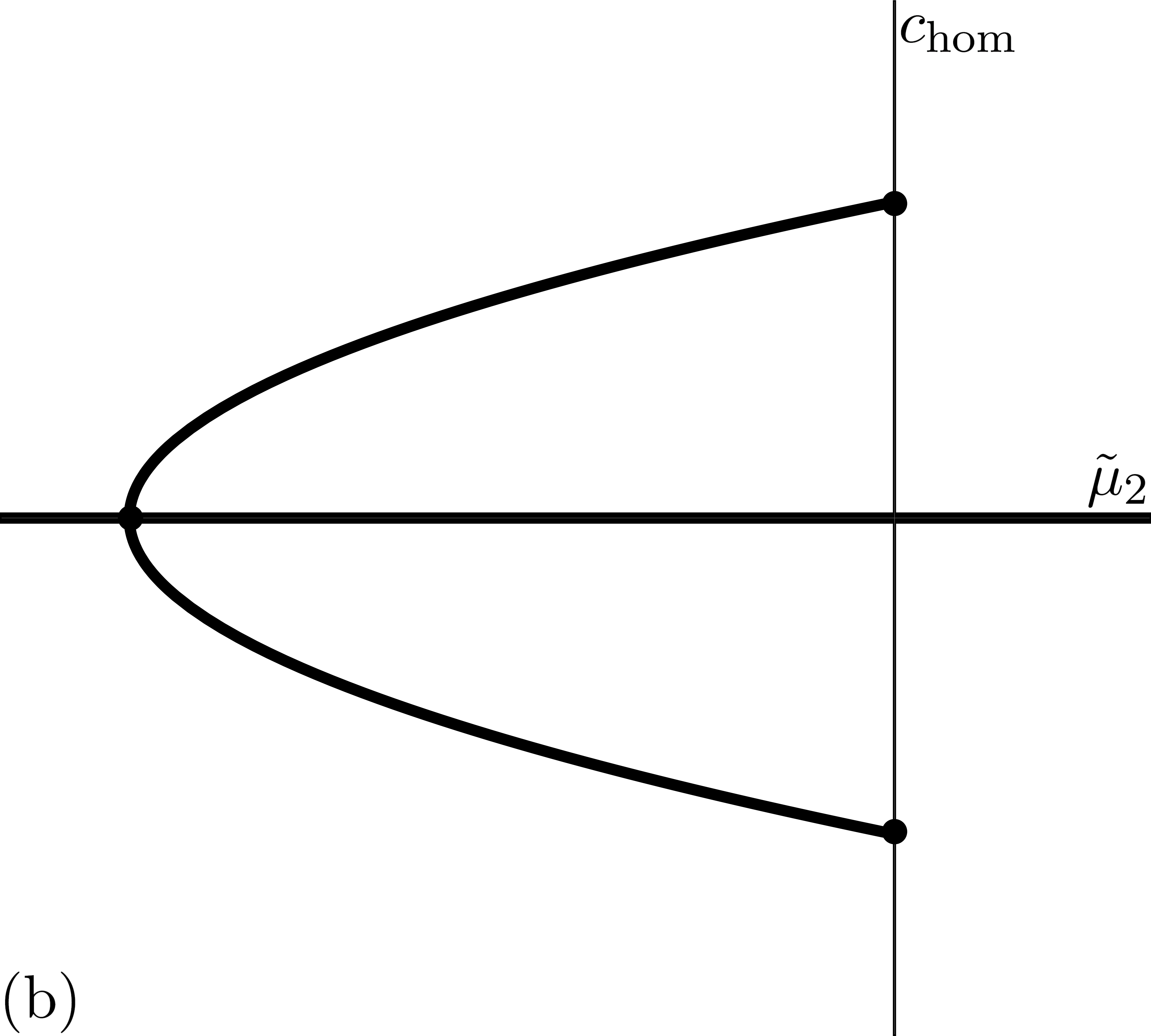}
	\end{minipage}
    \hspace{.5cm}
	\begin{minipage}{0.3\textwidth}
		\includegraphics[width=\linewidth]{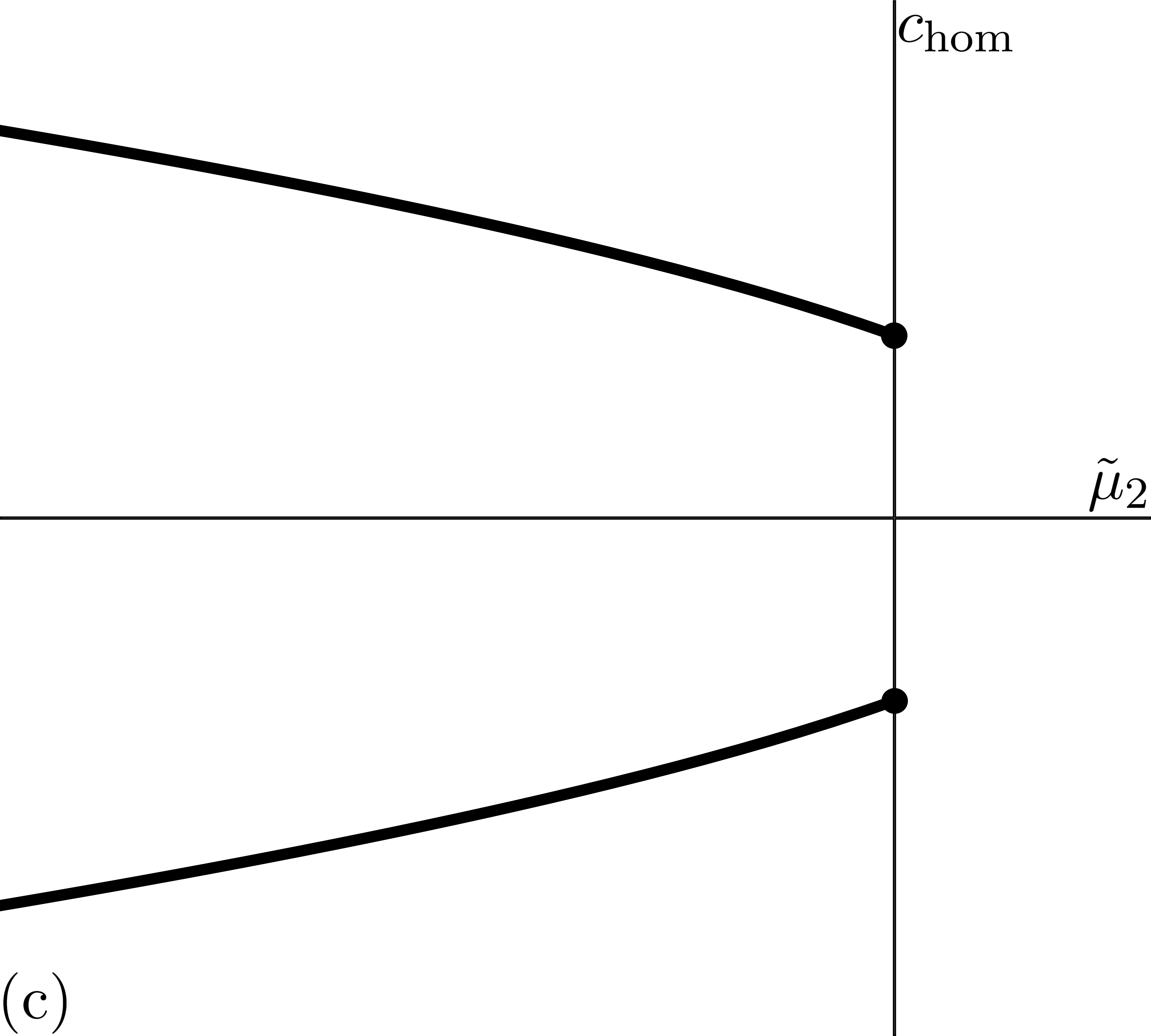}
	\end{minipage}
\caption{\small{The bifurcation diagrams described by Corollary \ref{c:hombifs} for $\oG_{1c}(\vmuta) < 0$. (a) Corollary \ref{c:hombifs}(i): $\tM_{2}(\vmuta) > 0, \tM_{2cc}(\vmuta) > 0$. (b) Corollary \ref{c:hombifs}(ii): $\tM_{2}(\vmuta) > 0, \tM_{2cc}(\vmuta) < 0$. (c) Corollary \ref{c:hombifs}(iii): $\tM_{2}(\vmuta) < 0, \tM_{2cc}(\vmuta) > 0$.}}
\label{f:hombifs}
\end{figure}
\noindent
{\bf Proof of Corollary \ref{c:hombifs}.} This corollary is a direct consequence of Theorem \ref{t:exhoms} -- and especially the existence condition $\tM_{2} + c_0^2 \tM_{2cc} > 0$ of Theorem \ref{t:exhoms}(ii) -- and the leading order existence relation (\ref{d:Vtc}) between $\vmu_2$ and $c_0$. However, to understand the split up of the homoclinic pulses into two pairs of fronts -- Corollary \ref{c:hombifs}(ii) and (iii) -- we need to add one additional ingredient.
\\ \\
The boundary of the existence condition $\tM_{2} + c_0^2 \tM_{2cc} > 0$ of Theorem \ref{t:exhoms}(ii) -- i.e $\tM_{2} + c_0^2 \tM_{2cc} = 0$ -- corresponds to $\vmu_1 \to 0$ in condition (\ref{e:exhets}) of Theorem \ref{t:exhets}(ii) (\ref{d:vmutsi}) -- i.e the case that limits on the present $\Ma(\vmu) = O(\eps^2|\log \eps|)$ setting. In this limit, two traveling fronts connecting $(\bU_-,\bV_-)$ for $X \to -\infty$ to $(\bU_+,\bV_+)$ for $X \to \infty$  exist (Corollary \ref{c:hetbifs}(ii)). By symmetry (\ref{d:symm}), these fronts have two symmetrical counterparts -- i.e. fronts connecting $(\bU_+,\bV_+)$ for $X \to -\infty$ to $(\bU_-,\bV_-)$ for $X \to \infty$ traveling with the same (opposite) speeds. Thus, as $\tmu_2 \downarrow 0$ (both in case (i) and (ii)), the two traveling pulses -- with speeds at leading order given by $c_0 = \pm \sqrt{-\tM_{2}/\tM_{2cc}} = \pm c_0^{\rm h-c}$ (\ref{d:Ctmuhet}) -- merge with the two heteroclinic cycles (at leading order), naturally also traveling with speeds $c_0 = \pm c_0^{\rm h-c}$ (Corollary \ref{c:hetbifs}(ii)). \hfill $\Box$
\begin{remark}
\label{r:verttopar}
\rm
Traveling patterns only exist in (\ref{e:DS}) under a condition on the parameter(s) $\vmu$ -- see (\ref{e:DHmu=0}) in Theorem \ref{t:E-Pers}(ii) and the condition $\vmu = \vmu_{\rm t}^\ast$ (such that $\Ma(\vmu_{\rm t}^\ast) = 0$ (\ref{d:Ma}), (\ref{d:vmutsi})) in Theorems \ref{t:exhets}(ii) and \ref{t:exhoms}(ii). In Theorem \ref{t:E-Pers}(ii) it is claimed that if condition (\ref{e:DHmu=0}) holds, the $\eps \to 0$ solutions persist for any $c \in \RR$ (of $\O(1)$), without zooming in near the critical value $\vmu_{\rm p}$ of $\vmu$ (at which (\ref{e:DHmu=0}) holds), while Corollaries \ref{c:hetbifs} and \ref{c:hombifs} show that zooming in near $\vmu = \vmu_{\rm t}^\ast$ -- through (\ref{d:vmutsi}) -- yields conditions on $c$ in both cases. Nevertheless, the statement of Theorem \ref{t:E-Pers}(ii) is correct. To establish this, we need to set up a higher order analysis along the lines of the above proof of Theorem \ref{t:exhoms}.
\\ \\
The starting point of the analysis is zooming in to $\O(\eps^2)$ near $\vmu_{\rm p}$, i.e. to introduce $M_{\rm p}(\vmu)$ by setting $\tilde{\Delta} \H_0(H_0, \vmu) = 2c M_{\rm p}(\vmu)$ (\ref{e:DHmu=0}) -- with $M_{\rm p} = M_{\rm hom}$ (\ref{d:M0vmuW0}) in the case of (unperturbed) homoclinic orbits -- so that $M_{\rm p}(\vmu_{\rm p}) = 0$, to introduce $\vmu_2$ by $\vmu = \vmu_{\rm p} + \eps^2 \vmu_2$, and to expand $M_{\rm p}(\vmu)$: $M_{\rm p}(\vmu) = \vec{\nabla} M_{\rm p}(\vmu_{\rm p}) \cdot \vmu_2 \, \eps^{2} + \O(\eps^{3})$ as in (\ref{d:vmutsi}). Without going into the details, we note that only even powers $c$ will appear in the analysis at $\O(\eps^2)$, thus, by symmetry (\ref{d:symm}) -- and as in the proof of Theorem \ref{t:exhoms} -- we cannot `close' the perturbed orbit by tuning $c$. On the other hand, the situation is different from that of Theorem \ref{t:exhoms}, since {\rm we do not get a sign condition} like $c^2 \tM_{2cc}(\vmu_{\rm t}^\ast) + \tM_{2}(\vmu_{\rm t}^\ast) > 0$: this condition distinguishes between the 2 ways the nearly (double front) heteroclinic orbit may pass along the saddle $(\oV_+,0)$ in the case of Theorem \ref{t:exhoms} -- there is no similar separation for the orbits considered in Theorem \ref{t:E-Pers}. At $\O(\eps^3)$, we can define the gap $\Delta W_{\rm p}(\vmu_2)$, similar to $\Delta W_{\rm h}(\vmu_2)$ introduced in the proof of Theorem \ref{t:exhoms}, and determine its leading order approximation,
\[
\Delta W_{\rm p}(\vmu_2) = C_{3, \rm p} \left[\oG_{1, {\rm p}} + \tM_{3c, {\rm p}} + c_0^2 \tM_{3ccc, {\rm p}} - \vec{\nabla} M_{\rm p} \cdot \vmu_{2} \right] c_0 \eps^3
\]
(cf. (\ref{e:gapOeps2inhomasymm})) -- where the terms $C_{3, \rm p}$, $\oG_{1, {\rm p}}$, $\tM_{3c, {\rm p}}$ and $\tM_{3ccc, {\rm p}}$ can in principle all be expressed explicitly (and are all of $\O(1)$). Thus, we conclude that indeed for any $c_0$ given, there is a co-dimension 1 manifold of values $\vmu_2$ for which the gap can be closed (and the original orbit persists). Moreover, we also observe that the condition $\Delta W_{\rm p}(\vmu_2) = 0$ determines a parabolic relation between $c_0$ and $\vmu_2$: the vertical line describing the bifurcation into traveling waves as formulated in Theorem \ref{t:E-Pers}(ii) is (locally) indeed replaced by the expected parabola.
\end{remark}

\section{Stability}
\label{s:Stab}

\subsection{The instability of the persisting homoclinic pulses}
\label{ss:S-RegPuls}

Before we consider the stability of the homoclinic patterns as established by Theorem \ref{t:E-Pers}, we derive the linearized stability problem associated to any (traveling) pattern in (\ref{e:RDE})/(\ref{e:RDE-S}) that corresponds to any (slow) bounded solution of (\ref{e:ODE})/(\ref{e:DS}) on $\M_\eps(c)$.

\subsubsection{The linearized stability problem and its expansion}
\label{sss:S-LinStab}
Let $\gas(X) = (\us(X),\ps(X),\vs(X),\qs(X)) \subset \M_\eps(c)$ be a -- homoclinic, heteroclinic, periodic -- solution of (\ref{e:DS}) that corresponds to the slow (stationary or traveling) pattern $(U(X,Y,t),V(X,Y,t)) \equiv (\Us(X),\Vs(X)) = (\us(X),\vs(X))$ in (\ref{e:RDE-S}) -- the slow version of (\ref{e:RDE}). Perturbing $(\Us(X),\Vs(X))$ by,
\beq
\label{d:L}
(U(X,Y,t),V(X,Y,t)) = (\Us(X) + \bu(X) e^{\la t + i LY}, \Vs(X) + \bv(X) e^{\la t + i LY}), \; \; \la \in \CC, L \in \RR,
\eeq
yields the linearized system,
\begin{equation}
\label{e:ODE-lin}
\left\{	
\begin{array}{rcrclcccccc}
\tau \la \bu &=& \eps^2(\bu_{XX} & - & L^2 \bu) & + & \eps c \tau \bu_X & + & F_u(\us(X),\vs(X)) \bu & + & F_v(\us(X),\vs(X)) \bv\\
\la \bv &=& \bv_{XX} & - & L^2 \bv & + & \eps c \bv_X & + & G_u(\us(X),\vs(X)) \bu & + & G_v(\us(X),\vs(X)) \bv
\end{array}
\right.
\end{equation}
The pattern $(\Us(X),\Vs(X))$ is a spectrally stable solution of (\ref{e:RDE-S}), if for all $L \in \RR$ (\ref{e:ODE-lin}) does only have bounded nontrivial solutions $(\bu(X), \bv(X))$ -- i.e. sup$_{x \in \RR} |(\bu(X), \bv(X))| < C$ for some $C > 0$ and $(\bu(X), \bv(X)) \not\equiv (0,0)$ -- for Re$(\la) =$ Re$(\la(L)) \leq 0$.
\\ \\
As in the existence analysis, we will study this problem by an expansion in $\eps$. Therefore, we write,
\beq
\label{d:usvsexp}
\us(X) = u_0(X) + \eps u_1(X) + \eps^2 u_2(X) + \O(\eps^3), \;
\vs(X) = v_0(X) + \eps v_1(X) + \eps^2 v_2(X) + \O(\eps^3),
\eeq
with $v_0(X)$ a solution of the reduced slow flow (\ref{e:RedSF}), etc. (cf. (\ref{d:v12q12}) where $v_0 = \va$). Since $\gas(X) = (\us(X),\ps(X),\vs(X),\qs(X)) \subset \M_\eps(c)$, it follows by (\ref{e:Meps}), (\ref{d:fjpj}) that,
\beq
\label{e:u12}
\begin{array}{l}
u_0  = f(v_0)
\; \; \; \; \; \; \; \; \; \;
u_1 = f'(v_0) v_1 - c_0 \tau q_0 \tf_1(v_0)
\\
u_2 = f'(v_0) v_2 + \frac12 f''(v_0) v_1^2 - c_0 \tau [v_{1,X}\tf_1(v_0) + q_0 \tf_1'(v_0)v_1 ] -c_1 \tau q_0 \tf_1(v_0) - c_0^2 \tau^2 \tF_{2cc} - \tF_2
\end{array}
\eeq
where we have dropped the $(v,q^2)$-dependence of $\tF_{2cc}$ and $\tF_2$, $q_0(X) = v_{0,X}(X)$ and we have reintroduced the expansion of $c$ ($=c_0 + \eps c_1 + \O(\eps^2)$). Moreover, we expand $\la$,
\beq
\label{d:laj}
\la = \la_0 + \eps \la_1 + \eps^2 \la_2 + \O(\eps^3),
\eeq
and $\bu(X)$ and $\bv(X)$,
\beq
\label{d:busbvsexp}
\bu(X) = \bu_0(X) + \eps \bu_1(X) + \eps^2 \bu_2(X) + \O(\eps^3), \; \;
\bv(X) = \bv_0(X) + \eps \bv_1(X) + \eps^2 \bv_2(X) + \O(\eps^3)
\eeq
(cf. (\ref{e:specstabexpO1})). Naturally, this yields leading order system (\ref{e:ODE-lin-0}) for $(\bu_0(X),\bv_0(X))$ that was already derived in the Introduction and thus also the nonlinear eigenvalue problem (\ref{e:O1EigPb}) -- where we now note that the denominator in the intermediate step (\ref{e:bu0inbv0}) is non-degenerate expression as long as Re$(\la_0) > - \kappa/\tau$ (\ref{e:NormHyp}): to establish the (in)stability of a pattern $(\Us(X), \Vs(X))$, we may use that expressions (\ref{e:bu0inbv0}) and (\ref{e:O1EigPb}) vary smoothly as function of $X$.
\\ \\
In this and the next (sub)section we will consider the stability of $(\Us(X),\Vs(X))$ as 1-dimensional pattern, which means that we will set $L=0$ in (\ref{e:ODE-lin}) and the equations/expressions deduced from it. We consider the stability of the patterns constructed in section \ref{s:Ex} as planar interface/stripe patterns in 2 space dimensions in section \ref{ss:IntStr}. Since we only need the leading order equations to establish the instability of the regular pulses of Theorem \ref{t:E-Pers} -- more precisely (\ref{e:ODE-lin-0}) and (\ref{e:O1EigPb}) (with $L=0$) -- we do not consider the higher order spectral problems here. However, these will be crucial for our analysis of the stability of the (nearly) heteroclinic fronts and pulses in section \ref{ss:S-NearHet}.
\\ \\
Throughout this section, we also assume that the localized patterns we study are (bi-)asymptotic to stable background states, i.e. we assume that the trivial solutions $(U(X,t),V(X,t)) \equiv (\oU_\pm,\oV_\pm)$ of (\ref{e:RDE-S}) are stable. This is the case under conditions (\ref{e:normhypsaddle}) and (\ref{e:stabtriv}) formulated in Lemma \ref{l:sigmaess} below. Hence, under these conditions, the essential spectrum associated to any of the localized patterns $(U(X,t),V(X,t)) = (\us(X),\vs(X))$ of (\ref{e:RDE-S}) constructed in Theorems \ref{t:E-Pers}, \ref{t:exhets} and \ref{t:exhoms} must be in the stable half-plane, which yields that their (spectral) stability is determined by their discrete spectrum (and thus that we may restrict our search for nontrivial solutions of (\ref{e:ODE-lin}) to integrable -- in fact exponentially converging -- (eigen)functions \cite{KP13}).
\begin{lemma}
\label{l:sigmaess}
The critical points $(\oU_\pm,0,\oV_\pm,0)$ of (\ref{e:DS}) are saddles for the slow flow on the normally hyperbolic slow manifold $\M_\eps$ if
\beq
\label{e:normhypsaddle}
\oF^\pm_u < 0, \; \; \oF^\pm_u \oG^\pm_v - \oF_v^\pm \oG^\pm_u > 0
\eeq
(\ref{d:oFpmetc}). These critical points correspond to stable trivial background states $(U(X,t),V(X,t)) \equiv (\oU_\pm,\oV_\pm)$ of (\ref{e:RDE-S}) -- in 1 space dimension -- if additionally
\beq
\label{e:stabtriv}
\oF^\pm_u + \tau \oG^\pm_v < 0.
\eeq
\end{lemma}
\noindent
{\bf Proof.}
First, we note that the condition on the critical points  $(\oU_\pm,0,\oU_\pm,0)$ of (\ref{e:DS}) to be saddles (\ref{e:Gfvvlin}) on a normally hyperbolic slow manifold $\M_\eps$ (\ref{e:f'}), (\ref{e:NormHyp}) is indeed equivalent to (\ref{e:normhypsaddle}). Next, we reduce (\ref{d:L}) to
\[
(U(X,t),V(X,t)) = (\oU + \bu e^{(\la + i\eps ck) t + i kX}, \oV + \bv e^{(\la + i\eps ck)t + i kX}), \; \; \la = \la(k) \in \CC, k \in \RR,
\]
with $(\oU, \oV) = (\oU_\pm, \oV_\pm)$, so that (\ref{e:ODE-lin}) simplifies to
\beq
\label{d:oFoG}
\left(
\begin{array}{cc}
\oF_u - \eps^2 k^2 - \tau \la & \oF_v\\
\oG_u & \oG_v - k^2 - \la
\end{array}
\right)
\left(
\begin{array}{c}
\bu \\ \bv
\end{array}
\right)
=
\left(
\begin{array}{c}
0 \\ 0
\end{array}
\right)
\eeq
We know \cite{Doe19} -- or can directly check -- that background states $(\oU_\pm,\oV_\pm)$ that correspond to saddles on a normally hyperbolic slow manifold are stable in singularly perturbed 2-component reaction-diffusion system (\ref{e:RDE-S}) if they are stable in the associated reaction ODE. Thus, we set $k=0$ in (\ref{d:oFoG}) and assume that the solutions $\la_\pm(0)$ of the characteristic polynomial
\beq
\label{e:charpol}
\tau \la^2 - \left[\oF_u + \tau \oG_v \right] \la + \left[\oF_u \oG_v - \oF_v \oG_u\right] = 0
\eeq
satisfy Re$(\la_\pm(0)) < 0$. By (\ref{e:normhypsaddle}) we conclude that the stability of $(\oU_\pm,\oV_\pm)$ in (\ref{e:RDE}) is indeed settled if (\ref{e:stabtriv}) holds (additional to (\ref{e:normhypsaddle})). \hfill $\Box$

\subsubsection{The 1-parameter family of Sturm-Liouville operators $\L_{\rho}(X)$}
\label{sss:famSL}

We set up the (spectral) stability analysis of the full pulse $(\Us(X),\Vs(X))$ by introducing $\rho = \tau \la_0$ in (\ref{d:Ls}) and defining the family of (smooth) operators $\L_\rho(X)$ associated to (\ref{e:ODE-lin-0})
\beq
\label{d:Lrho}
\L_\rho(X) = \L_s(X) - \frac{\rho f'(v_0(X)) G_u(f(v_0(X)),v_0(X))}{\rho - F_u(f(v_0(X)),v_0(X))}
\eeq
parameterized by $\rho > -\kappa$ (\ref{e:NormHyp}). Since $v_0(X)$ is a homoclinic orbit that is even in $X$, the operator $\L_\rho(X)$ is a Sturm-Liouville operator (for $X \in \RR$) that is also even in $X$. Thus, for any $\rho \; (> -\kappa)$ its point spectrum $\si_{\rm pt}(\L_\rho)$ consists of $J = J(\rho) + 1$ simple, real eigenvalues $\la_{J}(\rho) < ... < \la_{1}(\rho) < \la_{0}(\rho)$ -- with possibly $J(\rho) = - 1$, i.e. $\si_{\rm pt}(\L_\rho) = \emptyset $. The associated eigenfunctions -- with a slight abuse of notation denoted by $\bv_{j}(X; \rho)$, $j=0,1,...,J(\rho)$ -- are orthogonal, even/odd as function of $X$ for $j$ even/odd (since $v_0(X)$ is even as function of $X$) and have exactly $j$ zeroes; its essential spectrum $\si_{\rm ess}(\L_\rho)$ is given by $(-\infty, \partial \si_{\rm ess}(\L_\rho)]$ with $\partial \si_{\rm ess}(\L_\rho) < \la_{J}(\rho)$ -- see for instance \cite{KP13}. Moreover, all $\la_{j}(\rho)$ and $\bv_{j}(X; \rho)$, $j=0, 1, ..., J(\rho)$ vary smoothly in $\rho$.
\\
\begin{figure}[t]
\centering
\includegraphics[width=15cm]{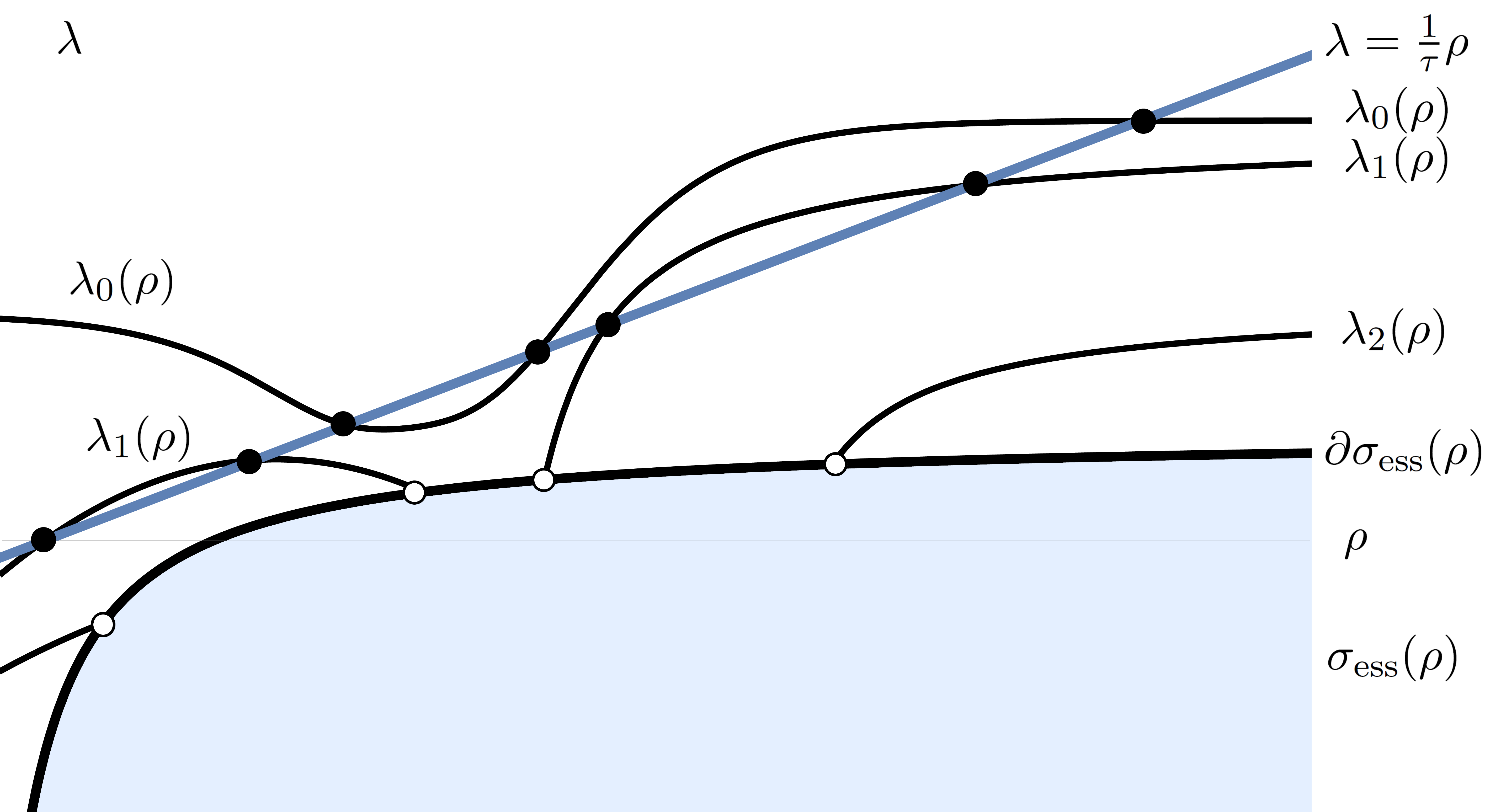}
\caption{\small{Sketches of the eigenvalue branches $\la_j(\rho)$ associated to the Sturm-Liouville operator $\L_\rho(X)$ (\ref{d:Lrho}), (the boundary of) its essential spectrum $\si_{\rm ess}(\rho)$ (\ref{e:psiessLrho}), and the intersections $\la_j(\rho) \cap \{\la = \rho/\tau\}$ -- that determine the positive eigenvalues of (leading order) nonlinear eigenvalue problem (\ref{e:O1EigPb}) and thus the stability of the homoclinic pulse patterns of Theorem \ref{t:E-Pers}. Note that the sketch considers the case of Corollary \ref{c:1laast}(ii), i.e. that $\la_1'(0) > \frac{1}{\tau}$, so that there must be at least 2 intersections/positive eigenvalues. Note also that the open bullets indicate the points $(\rho_{{\rm e},j}, \la_j(\rho_{{\rm e},j}))$, respectively $(\tilde{\rho}_{{\rm e},i}, \la_j(\tilde{\rho}_{{\rm e},i}))$, at which eigenvalues of $\L_\rho(X)$ disappear into, resp. appear from, $\si_{\rm ess}(\rho)$ as $\rho$ is varied (Lemma \ref{l:Lrho}). See Remark \ref{r:larhosin} for a class of models for which the curves $\la_j(\rho)$ indeed may fluctuate and for which a set-up as sketched here -- with 6 (or even more) positive eigenvalues -- may be expected to occur.}}
\label{f:Lambdajrho}
\end{figure}
\\
For $\rho = 0$, $\L_\rho(X) = \L_0(X) = \Ls(X)$ (\ref{d:Ls}), the operator associated to the linearization about the homoclinic pulse solution $V(X,t) = \Vho(X) = v_0(X)$ of slow reduced PDE (\ref{e:RedPDE}), hence we know that $\la_{0}(0) > 0$, $\la_{1}(0) = 0$ and $\bv_{1}(X; 0) = v_{0,X}(X)$. However, it should be noted that only for $\rho = 0$ there is a direct relation between $\L_\rho(X)$ and a linearization in a scalar PDE: $\la = 0$ only occurs in its special role associated to the translational invariance of an underlying PDE for $\rho \, (= \tau \la) = 0$. Moreover, since $v_0(X) \to \oV_-$ as $X \to \pm \infty$ it follows directly that,
\beq
\label{e:psiessLrho}
\partial \si_{\rm ess}(\L_\rho) = \partial \si_{\rm ess}(\rho) = \frac{\oF^-_u \oG^-_v - \oF_v^- \oG^-_u}{\oF_u^-} + \frac{\oF_v^- \oG^-_u}{\oF_u^-} \frac{\rho}{\rho -\oF_u^-},
\eeq
(\ref{d:oFpmetc}) so that,
\beq
\label{e:psiessLrho0infty}
\partial \si_{\rm ess}(0) = \frac{\oF^-_u \oG^-_v - \oF_v^- \oG^-_u}{\oF_u^-} < 0, \; \;
\lim_{\rho \to \infty} \partial \si_{\rm ess}(\rho) = \oG^-_v < -\frac{\oF_u^-}{\tau}
\eeq
(\ref{e:normhypsaddle}), (\ref{e:stabtriv}), where we note that $\si_{\rm ess}(\rho)$ may thus enter into the $\{\la > 0\}$ half-plane for increasing $\rho$. We have the following result on the eigenvalues $\la_{\rho, j}$ (which will be proven after we have established the instability result of Theorem \ref{t:InstabRegHom}).
\begin{lemma}
\label{l:Lrho}
Consider for $\rho \geq 0$ the $J(\rho) + 1$ eigenvalues $\la_{j}(\rho)$ associated to operator $\L_\rho(X)$ (\ref{d:Lrho}). For each $j$, $j=0,1,...,J(0)$, there either is a critical value $\rho_{{\rm e},j} > 0$ such that $\la_{j}(\rho)$ exists for $\rho \in [0, \rho_{{\rm e},j})$ and $\lim_{\rho \to \rho_{{\rm e},j}} \la_{j}(\rho) = \partial \si_{\rm ess}(\rho_{{\rm e},j})$ (\ref{e:psiessLrho}), or $\la_{j}(\rho)$ exists for all $\rho \in [0, \infty)$ and
$\lim_{\rho \to \infty} \la_{j}(\rho) = \la_{\infty, j} < \infty $ exists. Moreover, for all $\rho$ such that $\la_{j}(\rho)$ exists, all $\la_{i}(\rho)$ with $i = 0, 1, ..., j-1$ also exist and remain ordered: $\la_{j}(\rho) < ... < \la_{1}(\rho) < \la_{0}(\rho)$ (which also infers that either $\rho_{{\rm e},i} > \rho_{{\rm e},j}$ or that $\la_{i}(\rho)$ exists for all $\rho \geq 0$). For $\rho$ such that $\la_{j}(\rho)$ exists, the evolution of $\la_{j}(\rho)$ is explicitly given by
\beq
\label{d:larho'}
\la_{j}'(\rho) = \frac{d}{d \rho} \la_{j}(\rho) = -\frac{1}{\|\bv_{j}(\rho)\|_2^2} \int_{-\infty}^{\infty} \frac{F_v(f(v_0),v_0) G_u(f(v_0),v_0)}{(\rho - F_u(f(v_0),v_0))^2} \bv_{j}^2(\rho) \, dX.
\eeq
\end{lemma}
\noindent
Note that this Lemma only considers the branches $\la_j(\rho)$, $j = 0,1,...,J(0)$, that connect to $\rho=0$. However, just as these branches may end in the essential spectrum at $\rho_{{\rm e},j}$, new branches $\la_i(\rho)$ may appear from $\si_{\rm ess}(\rho)$ at values $\tilde{\rho}_{{\rm e},i}$, etc. -- see the open bullets in Fig. \ref{f:Lambdajrho}. Moreover, one could also consider eigenvalue branches $\la_j(\rho)$ for $\rho < 0$. However, we do not pursue this here since we only need Lemma \ref{l:Lrho} to establish the following instability result.
\begin{theorem}
\label{t:InstabRegHom}
Every homoclinic pulse pattern $(U(x,y,t),V(x,y,t)) = (\Us(X),\Vs(X))$ of (\ref{e:RDE}) as established by Theorem \ref{t:E-Pers} is unstable.
\end{theorem}
\noindent
The proof of this theorem is completely based on the observation that eigenvalues of the full (leading order) nonlinear eigenvalue problem (\ref{e:O1EigPb}) associated to the stability of the pulse pattern $(\Us(X),\Vs(X))$ (in 1 space dimension) correspond to solutions of
\beq
\label{d:lajbyLrho}
\la_{j}(\rho) = \frac{\rho}{\tau}, \; \; j=0, 1, ..., J(\rho),
\eeq
see also Fig. \ref{f:Lambdajrho}. By Lemma \ref{l:Lrho} we can establish that (\ref{d:lajbyLrho}) always has at least one solution $\rho_\ast = \tau \la_\ast > 0$, i.e. that there must be positive eigenvalue $\la_\ast$. The situation thus is similar to the instability result of pulse solutions to scalar RDEs \cite{KP13}. However, we know that for scalar RDEs the spectral stability problem associated to a homoclinic pulse has exactly one unstable eigenvalue \cite{KP13}, in the present case there can be more than 1 unstable eigenvalue -- see Fig. \ref{f:Lambdajrho} and Remark \ref{r:larhosin}. Lemma \ref{l:Lrho} also provides explicit conditions that ensure that there either is exactly one unstable eigenvalue or at least 2. In fact, the present set-up does not even exclude the possibility of complex eigenvalues: operator $\L_\rho(X)$ no longer is self-adjoint (and thus Sturm-Liouville) for $\rho \in \CC$, so that no longer necessarily $\la_{j}(\rho) \in \RR$ and (\ref{d:lajbyLrho}) thus may have solutions $\la \notin \RR$ -- see sections \ref{sss:nearstatstab} and \ref{sss:neartravstab}.
\begin{corollary}
\label{c:1laast}
Let $v_0^{\rm max}$ be the maximal value of the homoclinic orbit $v_0(X)$, i.e. $v_0^{\rm max} \in (\oV_c, \oV_+)$ satisfies $\W_0(v) = 0$ with $\W_0(v)$ as defined in (\ref{d:HW0}) with $\W_0(\oV_+) < \W_0(\oV_-) = 0$ (cf. (\ref{d:M0vmuW0})), and let $(U(x,t),V(x,t)) = (\Us(X),\Vs(X))$ of (\ref{e:RDE}) be an associated (stationary or traveling) homoclinic pulse pattern (in 1-space dimension) as established by Theorem \ref{t:E-Pers}.
\\
{\bf (i)} If
\beq
\label{d:lajrhodecreases}
F_v(f(v),v) G_u(f(v),v) \geq -\frac{\kappa^2}{\tau} \; \; {\rm for~all} \; \; v \in [\oV_-,v_0^{\rm max}]
\eeq
with $\kappa$ as defined in (\ref{e:NormHyp}), then the spectral problem associated to the stability $(\Us(X),\Vs(X))$ has exactly 1 eigenvalue $\la_\ast > 0$.
\\
{\bf (ii)} If $(\Us(X),\Vs(X))$ is a stationary pattern (as in Theorem \ref{t:E-Pers}(i)) with $\la_{1}'(0) > \frac{1}{\tau}$ -- or equivalently $M_{{\rm hom}} < 0$ (\ref{d:M0vmuW0}) -- then its spectral stability problem has at least 2 unstable eigenvalues.
\end{corollary}
\noindent
This Corollary is of course not relevant for the stability of the pulse, it is formulated here as an example of the usefulness of conditions like (\ref{d:lajrhodecreases}) -- obtained from (\ref{d:larho'}) -- and the direct link between (\ref{d:larho'}) -- with $j=1$ -- and the Melnikov conditions of the existence analysis in section \ref{ss:E-Pers} -- see also Remark \ref{r:larhosin} and the upcoming sections. Note that part (ii) of the Corollary cannot be applied to the traveling pulse patterns of Theorem \ref{t:E-Pers}(ii) since these necessarily have $M_{{\rm hom}} = 0$.
\\ \\
{\bf Proof of Theorem \ref{t:InstabRegHom} and Corollary \ref{c:1laast}.}
As was already noted, the pulse $(\Us(X),\Vs(X))$ can only be stable if it is bi-asymptotic to a stable background state $(\oU_-,\oV_-)$, see Lemma \ref{l:sigmaess} and especially assumptions (\ref{e:normhypsaddle}) and (\ref{e:stabtriv}) -- that we impose throughout this paper. These assumption also imply that the line $\la = \rho/\tau$ cannot intersect the boundary of the essential spectrum $\si_{\rm ess}(\rho)$: intersecting $\la = \rho/\tau$ with the curve (\ref{e:psiessLrho}) yields the following equation for $\rho$,
\[
\frac{1}{\tau} \rho^2 - \left[ \frac{\oF_u^-}{\tau} + \oG_v^- \right] \rho + \left[\oF_u^- \oG_v^- - \oF_v^- \oG_u^- \right] = 0,
\]
which is identical to (\ref{e:charpol}) -- after substitution of $\rho = \tau \la$ -- and which has, by (\ref{e:normhypsaddle}) and (\ref{e:stabtriv}), no positive solutions $\la$/$\rho$.
\\ \\
Next, we note that it follows from Lemma \ref{l:Lrho} that the line $\la = \rho/\tau$ {\it must} intersect the curve $\la_0(\rho)$. For $\rho = 0$, $\la_0(0) >0$, i.e. it is above the intersection of $\la = \rho/\tau$ with the $\{\rho=0\}$-axis. For increasing $\rho$, $\la = \rho/\tau$ increases linearly, while $\la_0(\rho)$ either ends at $\la_0(\rho_{{\rm e},0}) \in \partial \si_{\rm ess}(\rho)$ or approaches $\la_{\infty, 0}$ as $\rho \to \infty$ (Lemma \ref{l:Lrho}): since $\la = \rho/\tau$ cannot intersect $\partial \si_{\rm ess}(\rho)$, it necessarily passes through $\la_0(\rho)$ (at least once) at a certain value $\rho_{\ast} > 0$. Thus, stability problem (\ref{e:ODE-lin}) has a leading order eigenvalue $\la_{0,\ast} = \rho_{\ast}/\tau > 0$ with associated leading order eigenfunction $(\bv_0(X; \rho_{\ast}), \bu_0(X; \rho_{\ast}))$ -- with $\bu_0(X; \rho_{\ast})$ determined in terms of $\la_{0,\ast}$ and $\bv_0(X; \rho_{\ast})$ by (\ref{e:bu0inbv0}). It follows by standard methods that (\ref{e:ODE-lin}) must have an eigenvalue $\la_\ast = \la_\ast(\eps)$  with $|\la_\ast(\eps) - \la_{0,\ast}| < C \eps$ and eigenfunction $(\bv_\ast(X), \bu_\ast(X))$ with $\|(\bv_\ast(X), \bu_\ast(X)) - (\bv_0(X; \rho_{\ast}), \bu_0(X; \rho_{\ast}))\| < C \eps$ (pointwise) for some $C > 0$ and $\eps$ sufficiently small.
\\ \\
Both Corollary \ref{c:1laast} (i) and (ii) are directly induced by explicit expression (\ref{d:larho'}). Assumption (\ref{d:lajrhodecreases}) implies by (\ref{e:NormHyp}) that
\beq
\label{e:laj's1tau}
\la_{j}'(\rho)
\leq \frac{1}{\|\bv_{j}(\rho)\|_2^2} \int_{-\infty}^{\infty} \frac{\kappa^2}{\tau(\rho - F_u(f(v_0),v_0))^2} \bv_{j}^2(\rho) \, dX < \frac{1}{\|\bv_{j}(\rho)\|_2^2} \int_{-\infty}^{\infty} \frac{\kappa^2}{\tau(\rho + \kappa)^2} \bv_{j}^2(\rho) \, dX
\leq \frac{1}{\tau}
\eeq
for all (relevant) $j$ and $\rho \geq 0$. Since $\la_{j}(0) < \la_{1}(0) = 0$ for all $j=2,...,J(0)$ (and thus also $\partial \si_{\rm ess}(0) < 0$) and $\la_{0}(0) > 0$, $\la_{0}(\rho)$ is (for $\rho > 0$) the only $\la_{j}(\rho)$-branche attached to the $\{\rho=0\}$-axis that can have an intersection with the line $\la = \rho/\tau$. Naturally, new $\la_{j}(\rho)$-branches may appear from $\si_{\rm ess}(\rho)$, but assumption (\ref{d:lajrhodecreases}) also implies that $\frac{d}{d\rho} \partial \si_{\rm ess}(\rho) < 1/\tau$ for all $\rho \geq 0$ (\ref{e:psiessLrho}): new branches possibly appearing from $\si_{\rm ess}(\rho)$ -- of which the growth is also bounded by (\ref{e:laj's1tau}) -- cannot catch up with $\la = \rho/\tau$. Thus, condition (\ref{d:lajrhodecreases}) indeed ensures that the only possible eigenvalue $\la_\ast$ is determined by $\{\la = \rho/\tau\} \cap \{\la = \la_{0}(\rho)\}$ and is necessarily unique.
\\ \\
Since $\la_1(0) = 0$, the condition $\la_1'(0) > \frac{1}{\tau}$ implies that $\la_1(\rho) > \frac{\rho}{\tau}$ for $\rho > 0$ sufficiently close to $0$ -- see Fig. \ref{f:Lambdajrho}. Thus it follows by the same arguments as applied above to the $\la_0(\rho)$-branch, that the $\la_1(\rho)$-branch must also intersect the line $\frac{\rho}{\tau}$, which indeed yields the second unstable eigenvalue $\la_{\ast, 2} > 0$. Moreover, the condition $\tau \la_1'(0) - 1 > 0$ is equivalent to
\[
\frac{-\int_{-\infty}^{\infty} \left[\tau \frac{F_v(f(v_0),v_0) G_u(f(v_0),v_0)}{(F_u(f(v_0),v_0))^2} + 1 \right] v_{0,X}^2 dX}{\|v_{0,X}\|_2^2} =
\frac{2 \int_{\oV_-}^{v_0^{\rm max}}  \left[\tau \frac{f'(v) G_u(f(v),v)}{F_u(f(v),v)} -1 \right] \sqrt{2 \W_0(v)} dv}{\|v_{0,X}\|_2^2} =
\frac{- 2M_{{\rm hom}}}{\|v_{0,X}\|_2^2} > 0
\]
(\ref{d:larho'}), (\ref{e:f'}), (\ref{d:HW0}), (\ref{d:M0vmuW0}). \hfill $\Box$
\\ \\
Note that condition (\ref{d:lajrhodecreases}) guarantees that $\la_0(\rho)$, i.e. the curve connected to $\la_0(0) > 0$, must have a unique intersection with the line $\la = \rho/\tau$. Since $\la_1(0) \downarrow 0$ as the homoclinic orbit $v_0(X)$ approaches the heteroclinic cycle associated to the limit potential $\W_0(v)$ with wells of equal depth -- and since (\ref{d:lajbyLrho}) only gives leading order approximations -- this opens up the possibility of having stable homoclinic pulse patterns of the type established by Theorem \ref{t:exhoms} -- see Theorem \ref{t:nearstabtravfrontspulses} and Corollary \ref{c:stabstandfrontpulseMsmall}.
\\ \\
{\bf Proof of Lemma \ref{l:Lrho}.}
Since $\L_\rho(\rho)$ is a Sturm-Liouville operator for all $\rho \geq 0$, we know that for any $\rho$ fixed, the associated eigenvalue problem has point spectrum $\si_{\rm pt}(\L_\rho) = \{\la_{j}(\rho)\}_{j=0}^J$ with  $\partial \si_{\rm ess}(\rho) < \la_{J}(\rho) < ... < \la_{1}(\rho) < \la_{0}(\rho)$. The number of eigenvalues may change as function of $\rho$ (i.e. $J = J(\rho)$): $J$ decreases to $J-1$ as the smallest eigenvalue $\la_{J}(\rho)$ merges with the essential spectrum $\si_{\rm ess}(\rho)$ -- which defines the value $\rho_{{\rm e},J}$ ($J$ increases to $J+1$ as a new eigenvalue $\la_{J+1}(\rho) < \la_{J}(\rho)$ appears from $\si_{\rm ess}(\rho)$). Except for the behavior of $\la_{j}(\rho)$ as $\rho \to \infty$ and expression (\ref{d:larho'}) for $\la_{j}'(\rho)$, all statements in the lemma follow from these classical observations.
\\ \\
Clearly, the $\rho \to \infty$ limit of $\L_{\rho}(X)$ is well-defined,
\[
\lim_{\rho \to \infty} \L_\rho(X) \stackdef \L_\infty(X) = \frac{d^2}{dX^2} + G_v(f(v_0(X)),v_0(X))
\]
(\ref{d:Ls}), (\ref{d:Lrho}). If $G_v(f(v_0(X)),v_0(X)) \not\equiv$ a constant, then $\L_\infty(X)$ is again a Sturm-Liouville operator with $\si_{\rm pt}(\L_\infty) = \{\la_{j}(\infty)\}_{j=0}^{J(\infty)}$ and the statement(s) of the lemma follow(s). It should be noted though that, for instance, the fact that $\la_{0}(\infty)$ exists does not necessarily imply that $\lim_{\rho \to \infty} \la_{0}(\rho) = \la_{0}(\infty)$ ($ = \la_{\infty, 0}$): it can in general not be excluded that the $\la_{0}(\rho)$-branch connected to the $\{\rho=0\}$-axis first disappears into $\si_{\rm ess}(\rho)$ at $\rho_{{\rm e},0}$, reappears again from $\si_{\rm ess}(\rho)$ at $\tilde{\rho}_{{\rm e},0} > \rho_{{\rm e},0}$ and then continues to exist for all $\rho > \tilde{\rho}_{{\rm e},0}$ (with $\lim_{\rho \to \infty} \la_{0}(\rho) = \la_{0}(\infty) = \la_{\infty, 0}$). However, this does not contradict the statement of the lemma.
\\ \\
The alternative case in which $G_v(f(v_0(X)),v_0(X)) \equiv \alpha$, a constant, is more subtle. Note that this situation for instance occurs for $G(u,v) = \alpha v + \tilde{g}(u)$, which is not especially atypical (in fact, this may happen in  Gierer-Meinhardt-type models \cite{GM72}). In this case, $\L_\rho(X)$ is Sturm-Liouville for all $\rho < \infty$, but all $X$-dependent terms disappear from $\L_\rho(X)$ in the limit $\rho \to \infty$. Assuming that $J(\rho) \geq 0$ for large $\rho$, i.e. that there are $J(\rho) + 1 \geq 1$ eigenvalues $\la_{j}(\rho)$, then the question is: what happens to $\la_{j}(\rho)$ as $\rho \to \infty$? A priori it cannot be excluded that $\la_{j}(\rho) \to \infty$, which would not only invalidate the lemma, but would also undermine the proof of Theorem \ref{t:InstabRegHom}.
\\ \\
Therefore, we introduce the (artificial) small parameter $0 < \de \ll 1$ by $\rho = 1/\de \gg 1$, so that
\beq
\label{d:Q}
\begin{array}{rcl}
\L_\rho(X) & = & \frac{d^2}{dX^2} + f'(v_0) G_u(f(v_0),v_0) + G_v(f(v_0),v_0) - f'(v_0) G_u(f(v_0),v_0) \frac{1}{1 - \de F_u(f(v_0),v_0)}
\\
 & \stackdef & \frac{d^2}{dX^2} + \alpha + \delta Q(X)
\end{array}
\eeq
(\ref{d:Ls}), (\ref{d:Lrho}). It can be deduced in a straightforward fashion that $\la_{j}(\rho=1/\de)$ must be bounded (if it exists) -- and thus that $\lim_{\rho \to \infty} \la_{j}(\rho) < \infty$ -- which would prove the lemma. However, we can be more explicit: we define $Q_\infty$, $\tilde{Q}(X)$ and $\tQ_\RR$ by
\beq
\label{d:otherQs}
\lim_{X \to \pm \infty} Q(X) = Q_\infty, \; \tilde{Q}(X) = Q(X) - Q_\infty, \; \tQ_\RR = \int_{-\infty}^\infty \tilde{Q}(X) \, dX,
\eeq
where we note that both $|Q_\infty|, |\tQ_\RR| < \infty$ by the exponential decay of $v_0(X)$. By introducing
\beq
\label{d:tXtla}
\la = \alpha + \de Q_\infty + \de^2 \tla, \; \tX = \de X,
\eeq
the eigenvalue problem for $\L_\rho(X)$ with $\rho = 1/\de$ (\ref{d:Q}) can be written as a `locally kicked' system
\beq
\label{e:evde}
\bv_{\tX\tX} - \left[ \tla - \frac{1}{\de} \tQ \left(\frac{\tX}{\de} \right) \right] \bv = 0,
\eeq
from which it follows that,
\[
\begin{array}{ll}
\tQ_\RR > 0: & \si_{\rm pt}(\L_{1/\de}) = \{\la_0(\frac{1}{\de})\} = \{\alpha + \de Q_\infty +  \frac{1}{4} \de^2 \tQ_\RR^2 + \O(\de^3)\}\\
\tQ_\RR < 0: & \si_{\rm pt}(\L_{1/\de}) = \emptyset
\end{array}
\]
which can be deduced by a simple phase plane analysis of (\ref{e:evde}) -- see also \cite{DvHS18}. Thus, we conclude that also in the degenerate case that $G_v(f(v_0(X)),v_0(X)) \equiv \alpha$ the lemma holds, and that in that case $J(\frac{1}{\de}) = -1$ or $0$: all other (potential) eigenvalues $\la_j(\rho)$ ($j \geq 1$) must already have merged with $\si_{\rm ess}(\rho)$.
\\ \\
Finally, we note that $\L_{\rho}\bv_j(\rho) = \la_j(\rho)\bv_j(\rho)$ implies that
\[
0 = \frac{d}{d \rho} \left[\L_{\rho}\bv_j(\rho) - \la_j(\rho)\bv_j(\rho) \right] = \left[\L_{\rho} - \la_j(\rho)\right] \frac{d}{d \rho}\bv_j(\rho) + \left[\frac{d}{d \rho}\L_{\rho} - \la_j'(\rho)\right] \bv_j(\rho).
\]
Thus, $\tilde{v}(X) = \frac{d}{d \rho}\bv_j(X;\rho)$ solves
\[
\left[\L_{\rho} - \la_j(\rho)\right] \tilde{v} = \left[\frac{F_v(f(v_0),v_0) G_u(f(v_0),v_0)}{(\rho - F_u(f(v_0,v_0))^2} + \la_j'(\rho)\right] \bv_j(\rho)
\]
(\ref{e:f'}), (\ref{d:Lrho}). Since $\bv_j(\rho)$ spans the kernel of $\L_{\rho} - \la_j(\rho)$ this yields (\ref{d:larho'}).
\hfill $\Box$
\begin{remark}
\label{r:f'=0'stab}
\rm
It follows from (\ref{d:larho'}), or directly from (\ref{d:Lrho}, that $\la_j(\rho) \equiv \la_j(0)$ for all $j$ in the vertical case $f'(v) \equiv 0$. Thus, in the vertical case, the stability of the slow patterns indeed is determined by slow reduced PDE (\ref{e:RedPDE}) -- with $f(v) \equiv \oU (= \oU_\pm)$. Note however that this is only the case for spectrum with (Re)$\la > - \kappa/\tau$: for $f'(v) \equiv 0$, the equation for $\bu$ decouples from the system in (\ref{e:ODE-lin}),
\[
\bu_{XX} + \frac{1}{\eps^2} \left[ F_u(\oU,v_0(X)) - \tau \la \right] \bu = 0
\]
(with $L=0$ and $c=0$ (Remark \ref{r:f'0ex})). This also is a Sturm-Liouville problem, but since $F_u(\oU,v_0(X)) - \tau \la < 0$ for $\la > - \kappa/\tau$ (\ref{e:NormHyp}), it follows that $\bu \equiv 0$ for $\la > - \kappa/\tau$, so that indeed the stability is determined by the Sturm-Liouville problem associated to (\ref{e:RedPDE}) (i.e. the second line of (\ref{e:ODE-lin}) with $\bu = 0$).
\end{remark}
\begin{remark}
\label{r:larhosin}
\rm
As in Corollary \ref{c:1laast}, it will be necessary in the upcoming sections to know whether it is possible to have (some) control over the evolution of the eigenvalues $\la_j(\rho)$ of $\L_\rho(X)$. In the most general setting, it is not directly obvious how to obtain explicit and sufficiently sharp estimates. Nevertheless, condition (\ref{d:lajrhodecreases}) clearly provides such control. However, it is a rather strong condition: Corollary \ref{c:1laast}(i) also holds under the weaker -- but also less explicit and still stronger than strictly necessary -- condition $\la_j'(\rho) \leq 1/\tau$ for all $j$ -- see Theorem \ref{t:nearstabtravfrontspulses}. On the other hand, we may consider briefly the following (conceptual) class of explicit example systems (\ref{e:RDE}) with,
\[
F(U,V) = - (U - f(V)), \; \; G(U,V) = -\tilde{\W}'_0(V) + H_1(U -f(V))H_2(V),
\]
with $\tilde{\W}_0(V)$ a double well potential as introduced in section \ref{ss:E-Pers} and $H_1(0) = 0$ as only condition on the (sufficiently smooth) functions $f$, $H_1$ and $H_2$. Clearly, there is a unique slow manifold $\M_\eps(c) = \{u = f(V) + \O(\eps), p = \O(\eps) \}$ (\ref{e:Meps}) with reduced slow flow $v_{XX} - \tilde{\W}'_0(V) = 0$ (\ref{e:RedSF}), that is uniformly normally hyperbolic (since $F_U(f(v),v) \equiv -1$ (\ref{e:NormHyp})). Thus, by choosing $f(V)$ we have complete freedom in `tuning' the (leading order) `shape' of $\M_\eps(c)$. Moreover,
\[
\la_{j}'(\rho) = - \frac{1}{\|\bv_{j}(\rho)\|_2^2} \frac{H_1'(0)}{(1+\rho)^2} \int_{-\infty}^{\infty} f'(v_0(X))H_2(v_0(X)) \bv_{j}^2(X; \rho) \, dX
\]
(\ref{d:larho'}), so that even within this restricted class we have freedom to tune $\la_{j}'(\rho)$ -- without effecting the (leading order) slow flow on $\M_\eps(c)$ (since $H_1(0) = 0$). In fact, we may choose $f(V) = \sin (\mu_1 V + \mu_2)$, $H_1(w) = \mu_3 w$ and $H_2(w) \equiv 1$ so that,
\[
\la_{j}'(\rho) = - \frac{1}{\|\bv_{j}(\rho)\|_2^2} \frac{\mu_1 \mu_3}{(1+\rho)^2} \int_{-\infty}^{\infty} \cos \left(\mu_1 v_0(X) + \mu_2 \right) \bv_{j}^2(X; \rho) \, dX.
\]
A priori, this expression suggests that $\la_{j}'(\rho)$ may change sign a number of times by tuning $\mu_1$ and $\mu_2$ appropriately -- although this of course strongly depends on the nature of eigenfunctions $\bv_{j}(\rho)$ (vs. $\cos (\mu_1 v_0 + \mu_2)$). Note by (\ref{e:psiessLrho0infty}) that $\lim_{\rho \to \infty} \partial \si_{\rm ess}(\rho) = -\al_- - \mu_1 \mu_3 \cos (\mu_1 \oV_- + \mu_2)$ (\ref{d:muEnpm}), which also may be made either positive or negative by tuning $\mu_3$ (for given $\mu_{1,2}$). Finally, we note that, since the fast reduced flow is linear, $\M_\eps(c)$ is the only slow manifold and there cannot be a `return mechanism' for this class of systems: orbits that `jump off' from $\M_\eps(c)$ cannot `touch down' again on $\M_\eps(c)$. As consequence, all possible bounded orbits of (\ref{e:DS}) must be on $\M_\eps(c)$ for all $X$ and are thus necessarily slow, i.e. the orbits considered in this paper are especially relevant for this class of systems.
\end{remark}

\subsection{Near the heteroclinic limit: the stability of fronts and pulses}
\label{ss:S-NearHet}
In this subsection we consider the stability against perturbations of the nearly heteroclinic fronts and pulses as established by Theorems \ref{t:exhets} and \ref{t:exhoms} as (stationary or traveling) solutions of the 1-dimensional version of (\ref{e:RDE}), or equivalently, its rescaled co-moving version (\ref{e:RDE-S}),
\begin{equation}
\label{e:RDE-S1}
\left\{	
\begin{array}{rcrcrcr}
\tau U_t &=& \eps^2U_{XX} & + & \eps c \tau U_X & + & F(U,V)\\
V_t &=& V_{XX} & + & \eps c V_X & + & G(U,V)
\end{array}
\right.
\end{equation}
where we note that we come back to the stability of these patterns as interfaces or stripes in the full 2-dimensional systems (\ref{e:RDE})/(\ref{e:RDE-S}) in section \ref{ss:IntStr}. As a necessary preparatory step we first derive the next order terms -- compared to section \ref{sss:S-LinStab} --  of the asymptotic spectral stability problems (postponing most details to Appendix \ref{a:explicit}). Then we consider the stability of the standing fronts and pulses of Theorems \ref{t:exhets}(i) and \ref{t:exhoms}(i) under the condition that $\Ma(\vmu) = \O(1)$ and non-zero. In the final subsection of this section, the (spectral) stability of the traveling fronts and pulses of Theorems \ref{t:exhets}(ii) and \ref{t:exhoms}(ii) -- that have $\Ma(\vmu) = \O(\eps)$ or  $\Ma(\vmu) = \O(\eps^2|\log \eps|)$ -- is investigated. The missing case of (the stability of) stationary fronts and pulses with $|\Ma(\vmu)|$ asymptotically small is considered at the end of this section.

\subsubsection{The higher order spectral problems}
\label{sss:higherordereqs}
To establish the stability of localized patterns in the nearly heteroclinic limit, we need to consider the case of asymptotically small eigenvalues -- and especially $\la = \O(\eps)$ --  Therefore, we set $L = 0$ in (\ref{e:ODE-lin}) -- since we consider stability in (\ref{e:RDE-S1}) in this section -- and continue the expansion of spectral stability problem only for the special case that $\la_0 = 0$. Moreover, we replace the leading order component $v_0(X)$ (\ref{d:usvsexp}) of the $v$-component of the general homoclinic/heteroclinic/periodic orbit $\gas(X) \subset \M_\eps(c)$ of section \ref{sss:S-LinStab} by the unperturbed heteroclinic connection $\va(X)$ between the background states $\oV_-$ and $\oV_+$ as introduced in section \ref{ss:E-NearHet}. Naturally it follows by (\ref{e:u12}) that the leading order $u$-component $u_0(X)$ of $\gas(X)$ is replaced by $\ua(X) = f(\va(X))$. By (\ref{d:usvsexp}), (\ref{d:busbvsexp}) we thus find at the $\O(\eps)$ level of (\ref{e:ODE-lin}),
\begin{equation}
\label{e:ODE-lin-1}
\left\{	
\begin{array}{rcl}
\tau \la_1 \bu_0 &=& c_0 \tau \bu_{0,X} + \Fa_u \bu_1 + \Fa_v \bv_1 + (u_1 \Fa_{uu} + v_1 \Fa_{uv})\bu_0 + (u_1 \Fa_{uv} + v_1 \Fa_{vv})\bv_0 \\
\la_1 \bv_0 &=& \bv_{1,XX} + c_0 \bv_{0,X} + \Ga_u \bu_1 + \Ga_v\bv_1 + (u_1 \Ga_{uu} + v_1 \Ga_{uv})\bu_0 + (u_1 \Ga_{uv} + v_1 \Ga_{vv})\bv_0
\end{array}
\right.
\end{equation}
where we have reintroduced the notation of (\ref{d:fFGast}). In the present setting, we can determine solutions $\bu_0$ and $\bv_0$ of the leading order problem (\ref{e:ODE-lin-0}) explicitly by (\ref{e:bu0inbv0}), (\ref{e:O1EigPb}) with $\la_0 = L = 0$ and $\Ls(X)$ replaced by $\La(X)$ (\ref{d:La}),
\beq
\label{e:bv0bu0}
\bv_0(X) = \vax(X), \; \; \bu_0(X) = f'(\va(X)) \vax(X).
\eeq
Moreover, if $\vs(X) = \va(X) + \eps v_1(X) + \O(\eps^2)$ is the $v$-component of either a heteroclinic orbit of Theorem \ref{t:exhets} or of a nearly heteroclinic homoclinic orbit of Theorem \ref{t:exhets}, then $v_1(X) = c_0 \tv_1(X)$ (\ref{e:Lav1}), with $\tv_1(X)$ explicitly determined by (\ref{e:tv1}). It thus follows by (\ref{e:u12}) that it is also possible to split off a factor $c_0$ from $u_1$,
\beq
\label{d:tu1}
u_1(X) = c_0 \tu_1(X) = c_0 \left[f'(\va(X)) \tv_1(X) - \tau \qa(X) \tf_1(\va(X)) \right].
\eeq
Defining
\beq
\label{d:tFauuvv}
\tF^\ast_{uu} = \tu_1 \Fa_{uu} + \tv_1 \Fa_{uv}, \; \;
\tF^\ast_{vv} = \tu_1 \Fa_{uv} + \tv_1 \Fa_{vv}
\eeq
(\ref{d:fFGast}), we thus find from the first line of (\ref{e:ODE-lin-1}) that
\beq
\label{e:bu1inbv1}
\bu_1 = \fa' \bv_1 + \la_1 \frac{\tau\fa'\vax}{\Fa_u} - c_0 \frac{\tau (\fa'\vax)_X + (\tF^\ast_{uu} \fa' + \tF^\ast_{vv})\vax}{\Fa_u},
\eeq
so that we arrive at the following inhomogeneous equation for $\bv_1$,
\beq
\label{e:OepsEigPb}
\La \bv_1 = \la_1 \left[1-\tau \frac{\fa' \Ga_u}{\Fa_u}\right]\vax +  c_0 \hJ_{1c},
\eeq
with
\beq
\label{d:hJ1c}
\hJ_{1c} = \left[\frac{\tau \Ga_u (\fa'\vax)_X}{\Fa_u} - \vaxx \right] +
\frac{(\tF^\ast_{uu} \fa' + \tF^\ast_{uu}) \Ga_u - (\tG^\ast_{uu} \fa' + \tG^\ast_{uu})\Fa_u}{\Fa_u} \vax
\eeq
and $\tG^\ast_{uu}$, $\tG^\ast_{vv}$ as in (\ref{d:tFauuvv}). Recall that we are looking for an (integrable) eigenfunction, i.e. $\bv_1(X) \to 0$ as $X \to -\infty$, which implies that $\bv_1(X)$ is determined uniquely by setting $\bv_1(0) = 0$ (Lemma \ref{l:Lah}). The solution of (\ref{e:OepsEigPb}) can thus be written as
\beq
\label{e:bv1lc}
\bv_1 = - \la_1 \tv_1 + c_0 \hv_{1c},
\eeq
with $\tv_1(X)$ as already defined (\ref{d:tv1}), and solved (\ref{e:tv1}), in the existence analysis and solution $\hv_{1c}$ of $\La \hv_{1c} = \hJ_{1c}$ uniquely determined  by $\lim_{X \to -\infty} \hv_{1c}(X) =0$ and $\hv_{1c}(0)=0$ -- where we note that neither $\tv_1$ nor $\hv_{1c}$ depend on $\la_1$ or $c_0$ (and that they do vary with $\tau$ and $\vmu$). By (\ref{e:bu1inbv1}), this gives a similar expression for $\bu_1(X)$,
\beq
\label{e:bu1lc}
\bu_1 = \la_1 \hu_{1\la} + c_0 \hu_{1c},
\eeq
with
\beq
\label{d:hu1cl}
\hu_{1\la} = -\fa' \tv_1 + \frac{\tau\fa'\vax}{\Fa_u}, \; \;
\hu_{1c} =  \fa'\hv_{1c} - \frac{\tau (\fa'\vax)_X + (\tF^\ast_{uu} \fa' + \tF^\ast_{vv})\vax}{\Fa_u}.
\eeq
We use this as input in determining the spectral equations at $\O(\eps^2)$ -- noticing that the $\la_1$ pre-factors of $\bu_1$ and $\bv_1$ will give rise to terms containing $\la^2_1$ at this next level (through the expansions of the left hand sides of (\ref{e:ODE-lin})). The details are given in Appendix \ref{a:explicit}, here we only provide the endproduct, the inhomogeneous equation for $\bv_2$,
\beq
\label{e:Oeps2EigPb}
\La \bv_2 = \la_2 \left[1-\tau \frac{\fa' \Ga_u}{\Fa_u}\right]\vax +  c_1 \hJ_{1c} - \la_1^2 \hJ_{2\la \la} + c_0 \la_1 \hJ_{2c\la} + c_0^2 \hJ_{2cc} + \hJ_{2},
\eeq
with $\hJ_{1c}$ as in (\ref{d:hJ1c}),
\beq
\label{d:hJ2ll}
\hJ_{2\la \la} = (1-\tau \frac{\fa' \Ga_u}{\Fa_u})\tv_1 + \frac{\tau^2\Ga_u\fa'\vax}{(\Fa_u)^2},
\eeq
and $\hJ_{2c\la}(X)$, $\hJ_{2cc}(X)$ and $\hJ_{2}(X)$ given in (\ref{d:hJ2s}). Note that all terms $c_0$, $c_1$, $\la_1$ or $\la_2$ have been factored out explicitly in the right hand side of (\ref{e:Oeps2EigPb}), i.e. neither of remaining expressions depends on either of these 4 terms.

\subsubsection{$\Ma(\vmu) \neq 0$: the stability of standing patterns}
\label{sss:nearstatstab}

Unlike in the existence analysis, we combine the results on the stability of both types of standing localized patterns into one theorem.
\begin{theorem}
\label{t:nearstabstat}
Consider the stationary localized patterns $(\Uhe(X),\Vhe(X))$ and $(\Uho(X),\Vho(X))$ of (\ref{e:DS}) with $c=0$ as constructed in Theorems \ref{t:exhets}(i) and \ref{t:exhoms}(i). Assume that conditions (\ref{e:normhypsaddle}) and (\ref{e:stabtriv}) hold and let $\Ma(\vmu) \neq 0$ be $\O(1)$ with respect to $\eps$ (\ref{d:Ma}) and let $\eps$ be sufficiently small.
\\
{\bf (i-a)} If $\Ma(\vmu) < 0$, then the standing front $(\Uhe(X),\Vhe(X))$ is unstable.
\\
{\bf (i-b)} If $F_v(f(v),v) G_u(f(v),v) \geq -\frac{\kappa^2}{\tau}$ for all $v \in [\oV_-,\oV_+]$, then $(\Uhe(X),\Vhe(X))$ is (spectrally) stable as solution of (\ref{e:RDE-S1}).
\\
{\bf (ii)} All stationary pulses $(\Uho(X),\Vho(X))$ are unstable.
\end{theorem}
\noindent
Thus, unlike in scalar reaction-diffusion equations, slow fronts may be unstable. In fact, the situation is very similar to that of Corollary \ref{c:1laast} that considers essentially the same issue of the (possibly) non-Sturm-Liouville character of the spectral problem associated to slow localized patterns. Indeed, the condition on $F_v(f(v),v) G_u(f(v),v)$ in Theorem \ref{t:nearstabstat}(i-b) is the heteroclinic counterpart of (\ref{d:lajrhodecreases}) in Corollary \ref{c:1laast}. As was already noted about condition (\ref{d:lajrhodecreases}) in Remark \ref{r:larhosin}: the condition in Theorem \ref{t:nearstabstat}(i-b) can be replaced by the weaker but less explicit (and still stronger than necessary) condition $\la_j'(\rho) \leq 1/\tau$ for all $j$. We will opt for this formulation in Theorem \ref{t:nearstabtravfrontspulses}.
\\ \\
{\bf Proof of Theorem \ref{t:nearstabstat}}. We note that it follows from Lemma \ref{l:sigmaess} -- and conditions (\ref{e:normhypsaddle}) and (\ref{e:stabtriv}) -- that we only need to consider the discrete spectrum associated to the patterns $(\Uhe(X),\Vhe(X))$ and $(\Uho(X),\Vho(X))$.  We will first explicitly construct all asymptotically small eigenvalues (for both cases) and next consider the (potential) $\O(1)$ eigenvalues.
\\ \\
For asymptotically small $\la$ we have $\la_0 = 0$ and, by (\ref{e:tv1}), (\ref{d:busbvsexp}), (\ref{e:bv0bu0}), (\ref{e:bv1lc}),
\beq
\label{e:bvleading}
\bv = \vax -\eps \la_1\left\{\left[\int_0^X \tG_{1c}(\va) \vax \vu d\tX\right] \vax - \left[\int_{-\infty}^X \tG_{1c}(\va) \vax^2 d\tX\right] \vu \right\} + \O(\eps^2 \bv_2)
\eeq
(since $c_0 = 0$). To have an eigenfunction, the $X \to \infty$ limit of $\bv(X)$ needs to exist (and be $0$). Clearly, this implies by (\ref{e:propsvbu}) that $\la_1 \Ma = 0$, hence it follows -- by the assumption on $\Ma$ -- that $\la_1 = 0$. More importantly, it follows that there is only one possible asymptotically small eigenvalue. Since we know by translation invariance that $\la = 0$ must be an eigenvalue, we conclude that the present approximation procedure only leads to recovering the trivial $\la = 0$ eigenvalue.
\\ \\
This is different in the homoclinic case. Since $c=0$, the homoclinic pattern $(\Uho(X),\Vho(X)) = (\uho(X),\vho(X))$ (Theorem \ref{t:exhoms}) is symmetric -- even -- around its center/maximum, so that (also) the linearized stability problem (\ref{e:ODE-lin}) is symmetric under $X \to -X$, with $X = 0$ by definition pinned at the maximum of $\vho(X)$. Thus, solutions of (\ref{e:ODE-lin}) must be either odd or even in $X$: at the point of symmetry, either $\bv = 0$ -- for an odd eigenfunction -- or $\bv_X = 0$ -- for an even eigenfunction. Again, we also must recover the trivial eigenvalue $\la =0$ and its associated eigenfunction must be odd (since it is given by $(\uhox(X),\vhox(X))$).
\\ \\
By translation, the point of symmetry corresponds to $X= \Xh$ (\ref{e:Xh-1}), both in the approximation procedure of the proof of Theorem \ref{t:exhoms} that established the existence of the homoclinic solution on $\M_\eps$ of (\ref{e:DS}), as well as in the expansion of the spectral problem initiated in sections \ref{sss:S-LinStab} and \ref{sss:higherordereqs}. Thus, in order to construct a potential odd, respectively even, eigenfunction, we need to evaluate $\bv(\Xh)$, resp. $\bv_X(\Xh)$. It follows by (\ref{e:propsvbu}), (\ref{e:Xh-1}) from (\ref{e:bvleading}) that
\[
\bv(\Xh) = \frac{\la_1 \Ma}{\sqrt{2 \al_+ \tM_2}} + \O(\eps), \; \;
\bv_X(\Xh) = \frac{\la_1 \Ma}{\sqrt{2 \tM_2}} + \O(\eps)
\]
(where we have already used that $\bv_2(\Xh), \bv_{2,X}(\Xh)  = \O(1/\eps)$ -- see below). Thus, it follows that $\la_1 = 0$, both for the odd and the even case. Since $c_0=c_1=\la_1=0$, equation (\ref{e:Oeps2EigPb}) for $\bv_2$ simplifies drastically,
\[
\La \bv_2 = \la_2 \left[1-\tau \frac{\fa' \Ga_u}{\Fa_u}\right]\vax + \hJ_{2},
\]
with $\hJ_{2}(X)$ as in (\ref{d:hJ2s}). It is straightforward to check that $\hJ_2 (X) \vax(X)$ is integrable, thus we define
\beq
\label{d:hN2}
\hN_2(\vmu) = \int_{-\infty}^{\infty} \hJ_2 (X) \vax(X) \, dX,
\eeq
and conclude by Lemma \ref{l:Lah} that for $X \gg 1$,
\beq
\label{e:bvnext}
\bv(X) = \left[ \be_+ \sqrt{\al_+} e^{-\sqrt{\al_+} X} + \eps^2 \frac{\la_2 \Ma + \hN_2}{2 \al_+ \be_+} e^{+\sqrt{\al_+} X} \right] (1 + \O(E_+)),
\eeq
where we note that we again used the extension of Poincar\'e's Expansion Theorem as introduced in the proof of Theorem \ref{t:exhoms} (which is possible since $\la_1 = 0$ so that $\bv_1(X) \equiv 0$ (\ref{e:bv1lc})). Thus, by (\ref{e:Xh-1}),
\[
\bv(\Xh) = \frac{\eps}{\sqrt{2 \tM_2}} \left[ \tM_2 + \frac{\la_2 \Ma + \hN_2}{\sqrt{\al_+}} \right] +\O(\eps^2).
\]
Since $\la_2 = 0$ must be an eigenvalue for the odd case this necessarily implies that $\hN_2 = - \sqrt{\al_+} \tM_2$. (In fact, in the odd case, the present approximation scheme is reconstructing $\vhox$.) It also follows from (\ref{e:bvnext}) that
\beq
\label{e:bvXXhstat}
\bv_X(\Xh) = \frac{\eps}{\sqrt{2 \tM_2}} \left[-\sqrt{\al_+} \tM_2 + (\la_2 \Ma + \hN_2) \right] +\O(\eps^2).
\eeq
Since $\hN_2 = - \sqrt{\al_+} \tM_2$, we conclude that an even eigenfunction may exist for
\beq
\label{e:la2stat}
\la_2 = \frac{2 \sqrt{\al_+} \tM_2}{\Ma}.
\eeq
It follows by a standard Melnikov-type argument that there indeed must be an eigenvalue associated to an even eigenfunction that is at leading order given by (\ref{e:la2stat}). (The sign of $\bv_X(\Xh)$ can be changed by varying $\la_2$ (\ref{e:bvXXhstat}): there must be (a unique) value of $\la_2$ -- at leading order given by (\ref{e:la2stat}) -- for which $\bv_X(\Xh)$ is exactly $0$; by the application of the reversibility symmetry this yields the even eigenfunction.)
\\ \\
Thus we may conclude from (\ref{e:la2stat}) that $(\Uho(X),\Vho(X))$ is unstable for $\Ma > 0$ (recall that $\tM_2 > 0$, Theorem \ref{t:exhoms}(i)). There is no unstable asymptotically small spectrum for $\Ma < 0$, by an intuition based on scalar reaction-diffusion equations, we might expect that there cannot be non-asymptotically small unstable eigenvalues near a heteroclinic limit as in the present case: $(\Uho(X),\Vho(X))$ could thus be stable? However, this is not the case, if $\Ma < 0$ there must be $\O(1)$ unstable eigenvalues  -- as we shall show by an argument that is essentially the same as that of the proof of Corollary \ref{c:1laast}(ii).
\\ \\
But first we consider the possibility of $\O(1)$ (unstable) eigenvalues for the standing heteroclinic fronts of Theorem \ref{t:exhets}(i). Following the approach of section \ref{sss:famSL}, we define the family of Sturm-Liouville operators $\L_{\rho}^{\ast}(X)$ -- for $\rho \in \RR$ -- as in (\ref{d:Lrho}) over the leading order (unperturbed) heteroclinic connection $\va(X)$ (thus replacing the unperturbed homoclinic connection $v_0(X)$ in (\ref{d:Lrho}) by $\va(X)$). We denote the eigenvalues of $\L_{\rho}^{\ast}(X)$ by $\la_j^\ast(\rho)$ and conclude that Lemma \ref{l:Lrho} also holds in this case (with a slight but straightforward modification in the explicit analysis for the special case $G_v(f(\va(X)),\va(X)) \equiv \alpha$). Thus, we may follow the argument of the proof of Theorem \ref{t:InstabRegHom}: all $\O(1)$ eigenvalues of the spectral problem associated to $(\Uhe(X),\Vhe(X))$ are determined by the solutions of $\la_j^\ast(\rho) = \rho/\tau$ (cf. (\ref{d:lajbyLrho})).
\\ \\
Identical to its derivation in the proof of Corollary \ref{c:1laast}, we conclude that $\tau \la_0^\ast(0) - 1 > 1$ if $\Ma < 0$ -- where we notice that in the present heteroclinic case, the trivial translational eigenvalue at $\rho = 0$ (associated to the eigenfunction $\vax(X)$) is the critical (largest) eigenvalue, i.e. $\la^\ast_0(0) = 0$, while $0=\la^\ast_1(0) < \la^\ast_0(0)$ in the homoclinic case of Theorem \ref{t:InstabRegHom} and Corollary \ref{c:1laast}. Thus, the branch $\la^\ast_0(\rho)$ intersects the line $\rho/\tau$ at $\rho = 0$ and lies above it for $\rho > 0$ (and sufficiently small). This implies by (the equivalent of) Lemma \ref{l:Lrho} that there must at least be 1 next intersection of $\la^\ast_0(\rho)$ and $\rho/\tau$: if $\Ma < 0$, the spectral problem associated to $(\Uhe(X),\Vhe(X))$ has at least 1 $\O(1)$ positive eigenvalue -- which settles the instability result of Theorem \ref{t:nearstabstat}(i-a).
\\ \\
As in Corollary \ref{c:1laast}(i) (and its proof), the condition of Theorem \ref{t:nearstabstat}(i-b) implies that $\frac{d}{d \rho} \la_j^\ast(\rho) < \frac{1}{\tau}$ for all (relevant) $j$ and $\rho \geq 0$ (cf. (\ref{e:laj's1tau})). Since $\la_j^\ast(\rho) < \la_0^\ast(\rho)$ and $\la^\ast_0(\rho)$ intersects the line $\rho/\tau$ at $\rho = 0$, it follows that there cannot be any positive eigenvalues in the case of Theorem \ref{t:nearstabstat}(i-b). Nevertheless, we cannot yet conclude that $(\Uhe(X),\Vhe(X))$ is stable: we need to exclude the possibility of pairs of unstable {\it complex} conjugate eigenvalues $\la_\pm$ (i.e. with Re$\la_+ =$ Re$\la_+ > 0$). The set-up of our approach by the Sturm-Liouville operators $\L_\rho^\ast(X)$ strongly suggests that the spectral problem associated to $(\Uhe(X),\Vhe(X))$ cannot have eigenvalues $\la \in \CC\backslash\RR$, but $\L_\rho^\ast(X)$ only is Sturm-Liouville for $\rho \in \RR$ (it is not self-adjoint for $\rho \notin \RR$). Moreover, there is a priori no reason why the underlying {\it nonlinear} eigenvalue problem (\ref{e:O1EigPb}) -- with $L = 0$ and $v_0(X)$ replaced by $\va(X)$ -- cannot have nonreal eigenvalues.
\\ \\
In fact, our approach to finding eigenvalues as intersections of the branches $\la^\ast_j(\rho)$ with the line $\rho/\tau$ shows that (generically) pairs of complex conjugate eigenvalues will appear in nonlinear eigenvalue problem (\ref{e:O1EigPb}) at tangencies of $\la^\ast_j(\rho)$ and $\rho/\tau$: as parameter $\vmu$ is varied through such a tangency, 2 real eigenvalues merge and come out as a complex conjugate pair $\la_\pm$. Thus, one must in general expect (\ref{e:O1EigPb}) to have eigenvalues $\la \notin \RR$. Nevertheless, this cannot be the case under the condition of Theorem \ref{t:nearstabstat}(i-b). No tangencies between $\la^\ast_j(\rho)$ and $\rho/\tau$ can occur within the full $(F,G)$ family of `reaction-terms' $F(U,V)$ and $G(U,V)$ determined by this condition. Moreover, this family contains the special case $F_v(f(v),v) \equiv 0$ which determines the vertical case $f'(v) \equiv 0$ (\ref{e:f'}). We know by Remark \ref{r:f'=0'stab} that $\L_\rho^\ast(X) \equiv \L_s^\ast(X)$ in this case. Thus, within the $(F,G)$ family determined by Theorem \ref{t:nearstabstat}(i-b), there is a subfamily of spectral problems for which the full $\O(1)$ nonlinear eigenvalue problem (\ref{e:O1EigPb}) is a simple (scalar) Sturm-Liouville problem and thus cannot have spectrum in $\CC\backslash\RR$. By `straightening out' the $\la_j^\ast(\rho)$ branches to horizontal lines, any element of the full $(F,G)$ family of Theorem \ref{t:nearstabstat}(i-b) can be smoothly homotopied to the vertical subfamily without leaving the original family, and thus without going through situations in which there are tangencies between $\la^\ast_j(\rho)$ and $\rho/\tau$ (since $\frac{d}{d \rho} \la_j^\ast(\rho) < \frac{1}{\tau}$). Hence there cannot be (unstable) complex eigenvalues under the condition of Theorem \ref{t:nearstabstat}(i-b): $(\Uhe(X),\Vhe(X))$ must be spectrally stable.
\\ \\
Finally, we return to the (in)stability of the stationary homoclinic pulses $(\Uho(X),\Vho(X))$, i.e. to Theorem \ref{t:nearstabstat}(ii). By their construction in Theorem \ref{t:exhoms}, we know that the associated homoclinic orbits $(\uho(X),\pho(X),\vho(X),\qho(X)) \subset \M_\eps(0)$ are asymptotically close to the orbit on $\M_0$ spanned by $(\va(X),\qa(X))$ for $X$ such that $\qa(X) \geq 0$, and to its symmetrical counterpart $(\va(-X),-\qa(-X))$ for $X$ such that $\qa(X) \leq 0$, and that $\vho(X)$ is asymptotically close to $\oV_+$ for a range of $X$ values that is asymptotically large (in fact of $\O(|\log \eps|)$ -- see the proof of Theorem \ref{t:exhoms}. This implies that the leading order spectral stability problem associated to $(\Uho(X),\Vho(X))$ parameterized by $\rho$ -- i.e. the nearly heteroclinic limit of the homoclinic eigenvalue problem for $\L_\rho(X)$ of (\ref{d:Lrho}) -- is also determined by the above heteroclinic operator $\L_\rho^\ast(X)$ and that for any $\la_e = \la_e(\rho)$ for which the spectral problem has an even eigenfunction, there must be a $\la_o(\rho)$ asymptotically close to $\la_e(\rho)$ for which there is an odd eigenfunction (and vice versa). More specifically, for the stability problem associated to $(\Uho(X),\Vho(X))$, each intersection of a branch $\la^\ast_j(\rho)$ with the line $\rho/\tau$ provides the leading order approximation of 2 eigenvalues of the full eigenvalue problem (\ref{e:ODE-lin}) (with $c = L = 0$), $\la_e$ (associated to an even eigenfunction) and $\la_o$ (with an odd eigenfunction). (Note that this is confirmed by the above analysis: the intersection at $\rho = 0$ of $\la^\ast_0(\rho)$ and $\rho/\tau$ yields the trivial odd eigenvalue $\la = 0$ and the even eigenvalue given by (\ref{e:la2stat}) -- both obviously asymptotically close to $0$.)
\\ \\
Naturally, the instability of $(\Uho(X),\Vho(X))$ for $\Ma < 0$ follows by the same argument by which Theorem \ref{t:nearstabstat}(i-a) was established: for $\Ma < 0$, $\la^\ast_0(\rho)$ and $\rho/\tau$ must intersect at (at least) 1 $\O(1)$ value of $\rho > 0$, so that there must be at least 2 $\O(1)$ unstable eigenvalues. \hfill $\Box$

\subsubsection{$|\Ma(\vmu)|$ asymptotically small: the stability of fronts and pulses}
\label{sss:neartravstab}
\begin{figure}[t]
\begin{minipage}{.30\textwidth}
		\centering
		\includegraphics[width =\linewidth]{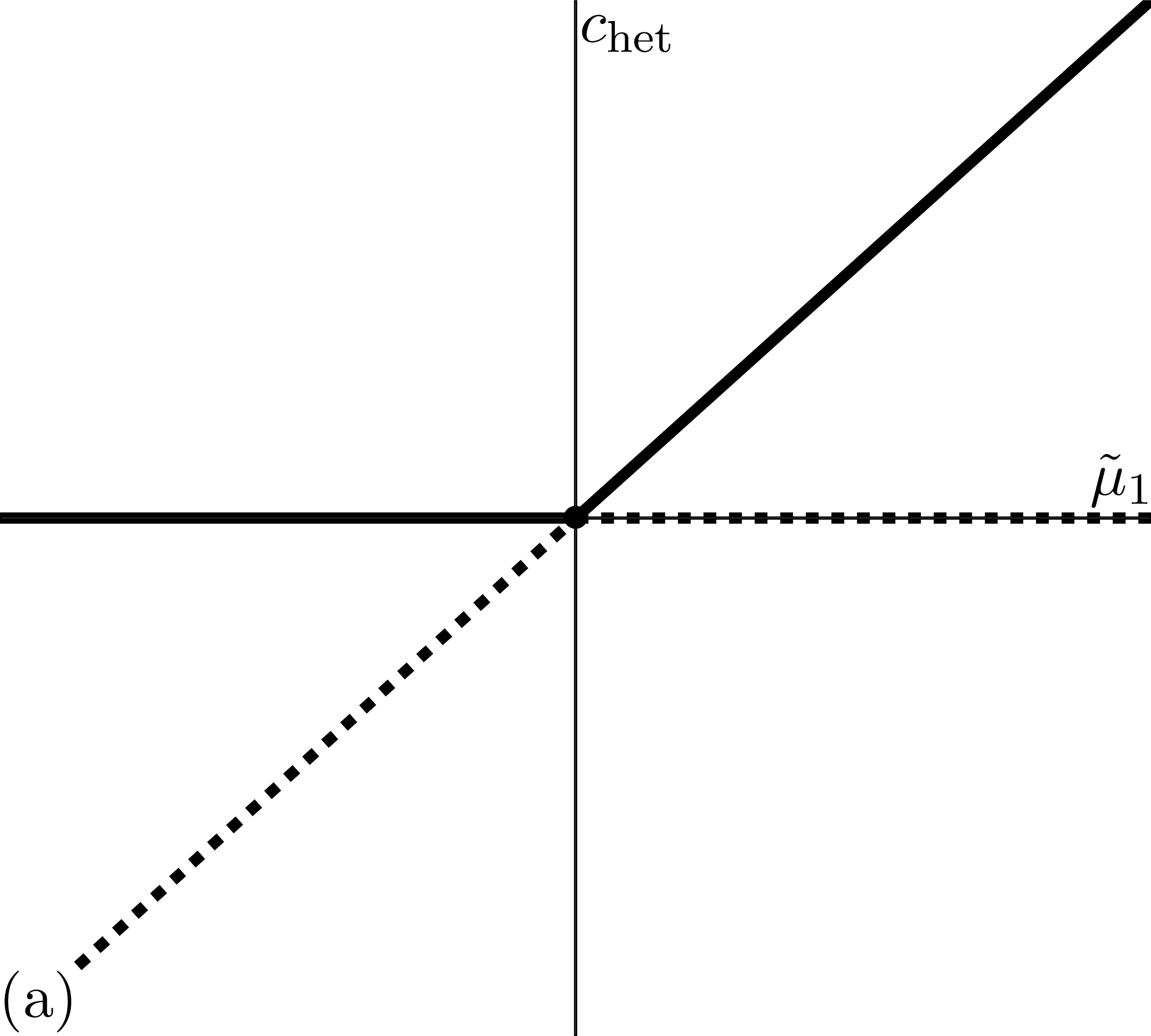}
\end{minipage}%
\hspace{.5cm}
\begin{minipage}{0.30\textwidth}
		\includegraphics[width=\linewidth]{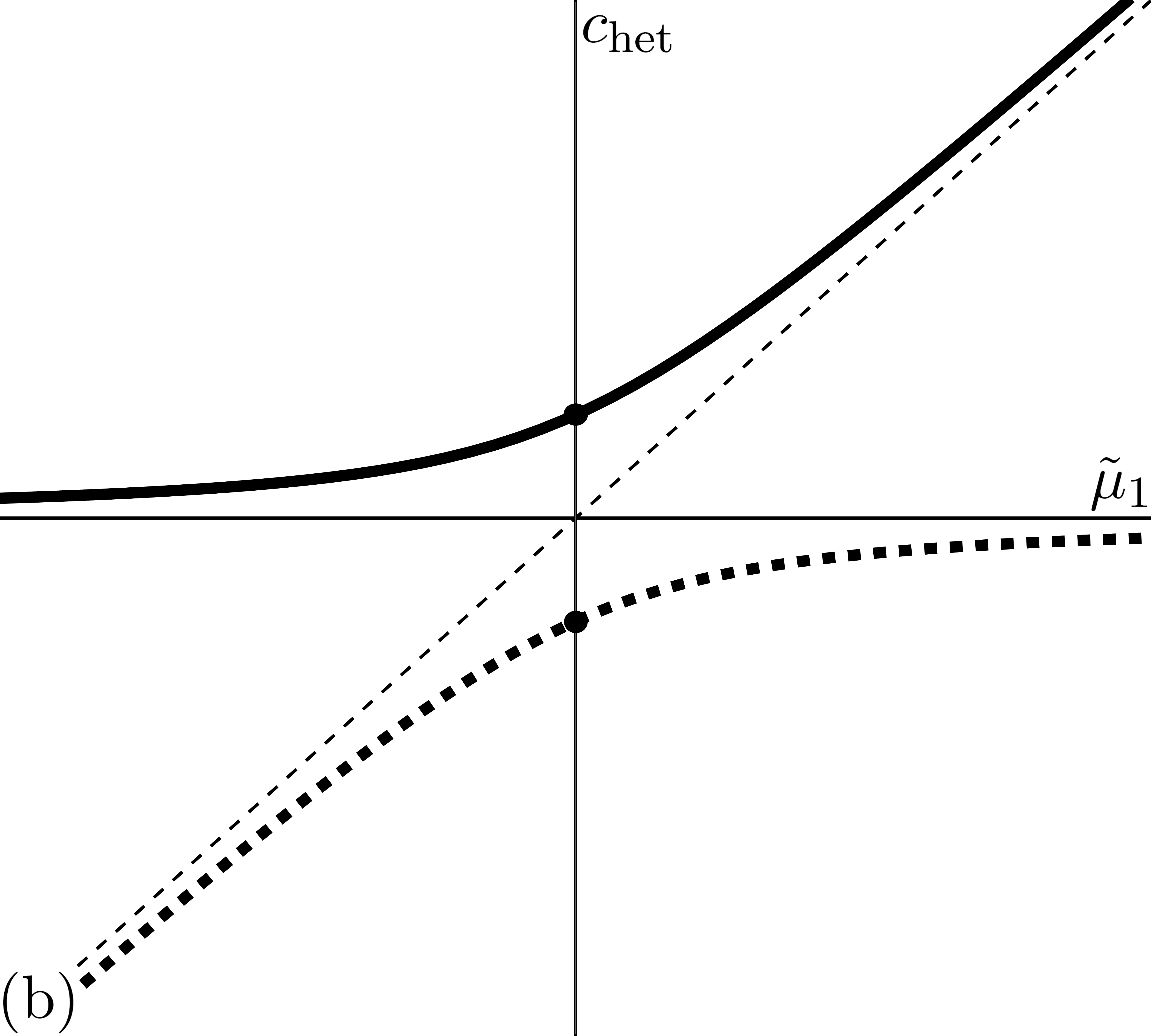}
\end{minipage}
    \hspace{.5cm}
\begin{minipage}{0.30\textwidth}
		\includegraphics[width=\linewidth]{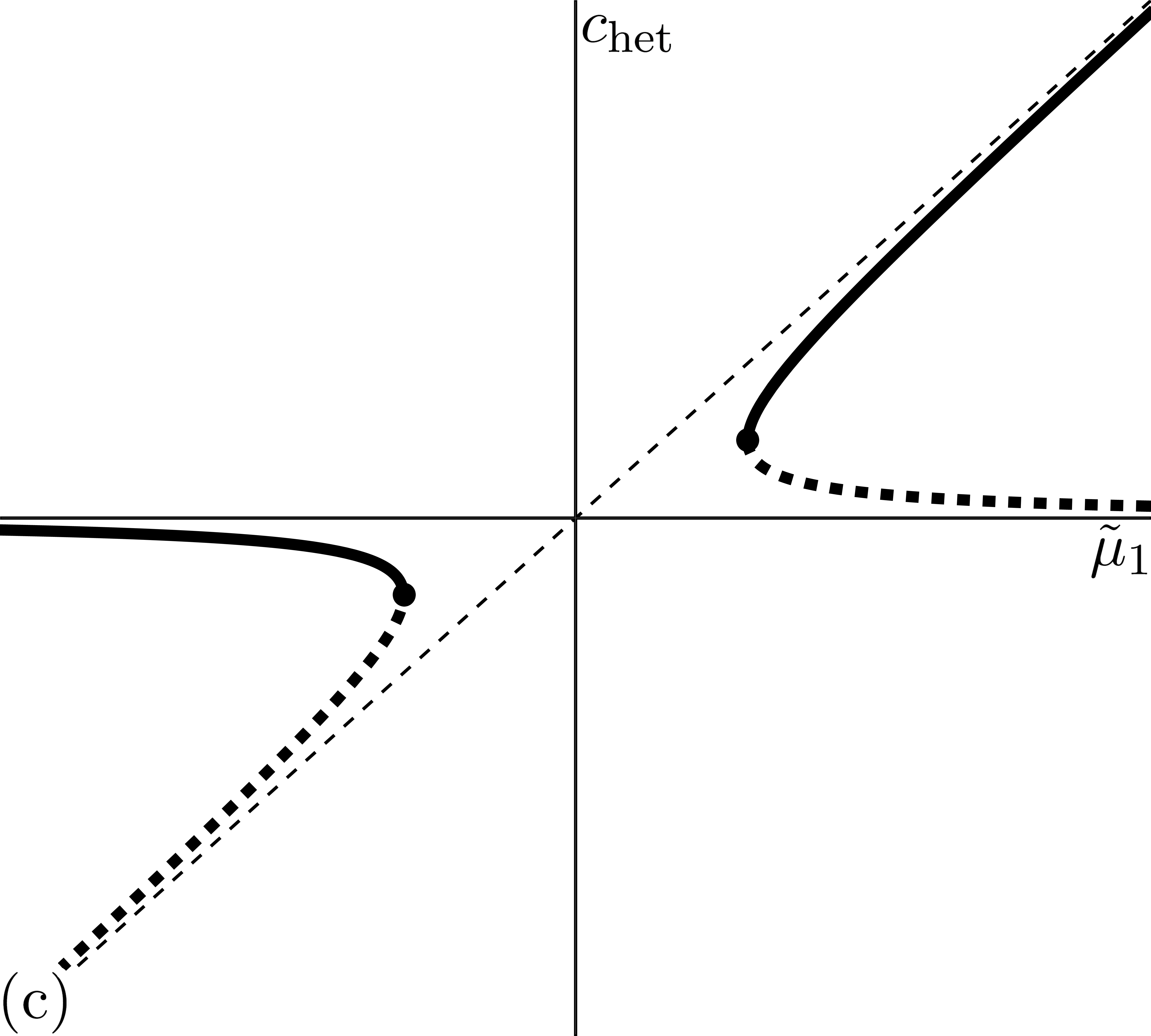}
\end{minipage}
\\
\vspace{.25cm}
\\
\begin{minipage}{.30\textwidth}
		\centering
		\includegraphics[width =\linewidth]{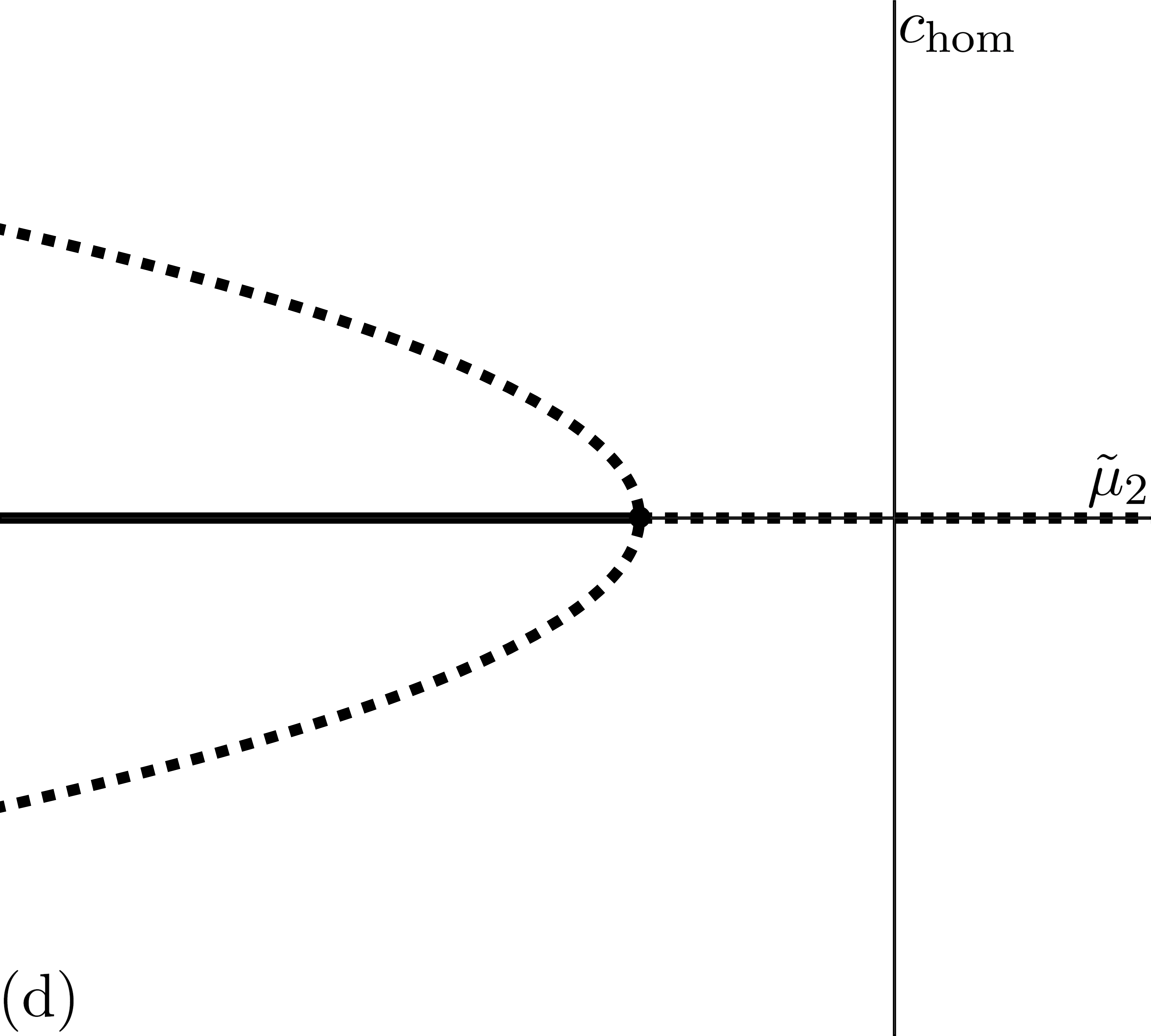}
\end{minipage}
\hspace{.5cm}
\begin{minipage}{0.30\textwidth}
		\includegraphics[width=\linewidth]{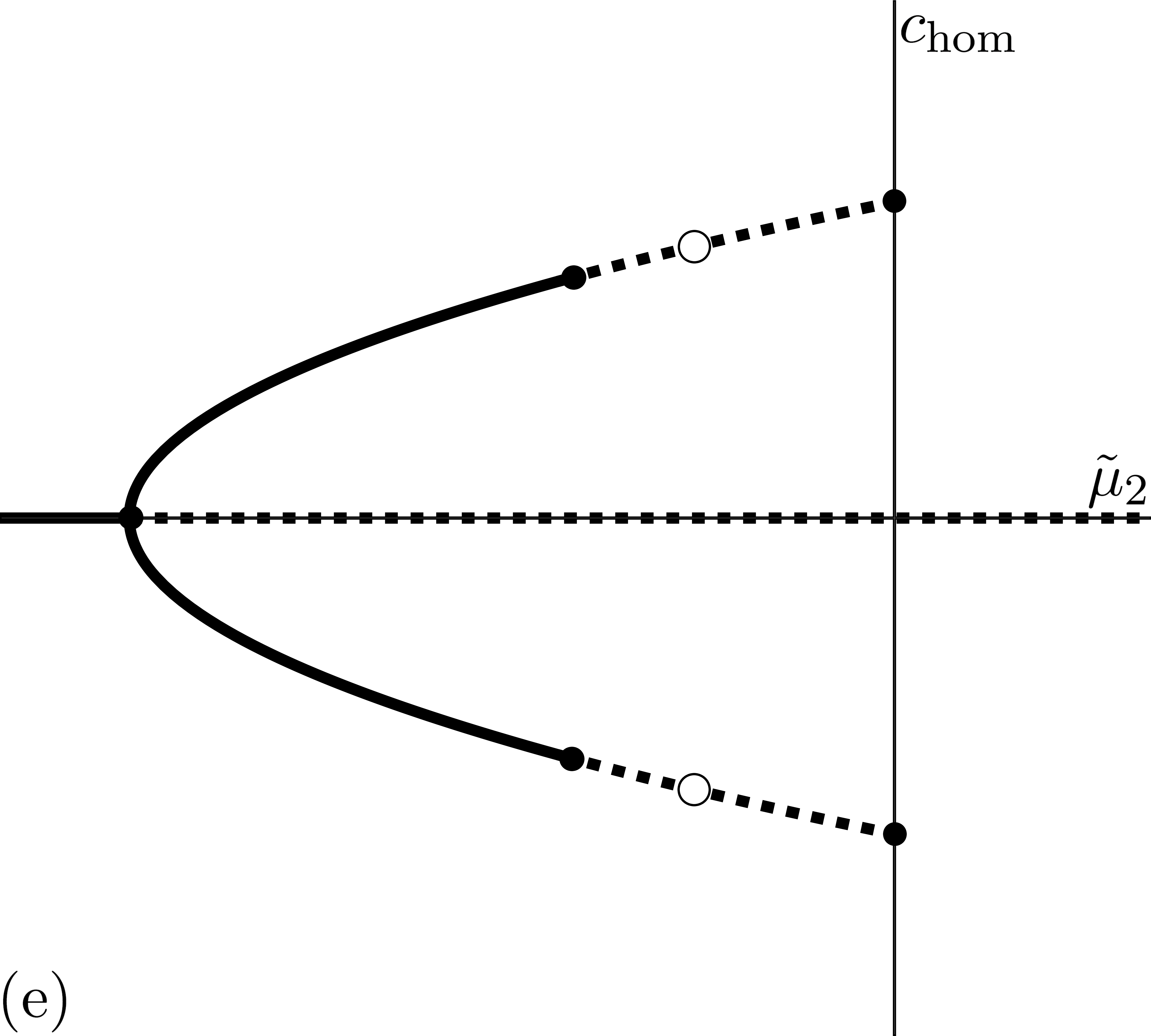}
\end{minipage}
    \hspace{.5cm}
\begin{minipage}{0.30\textwidth}
		\includegraphics[width=\linewidth]{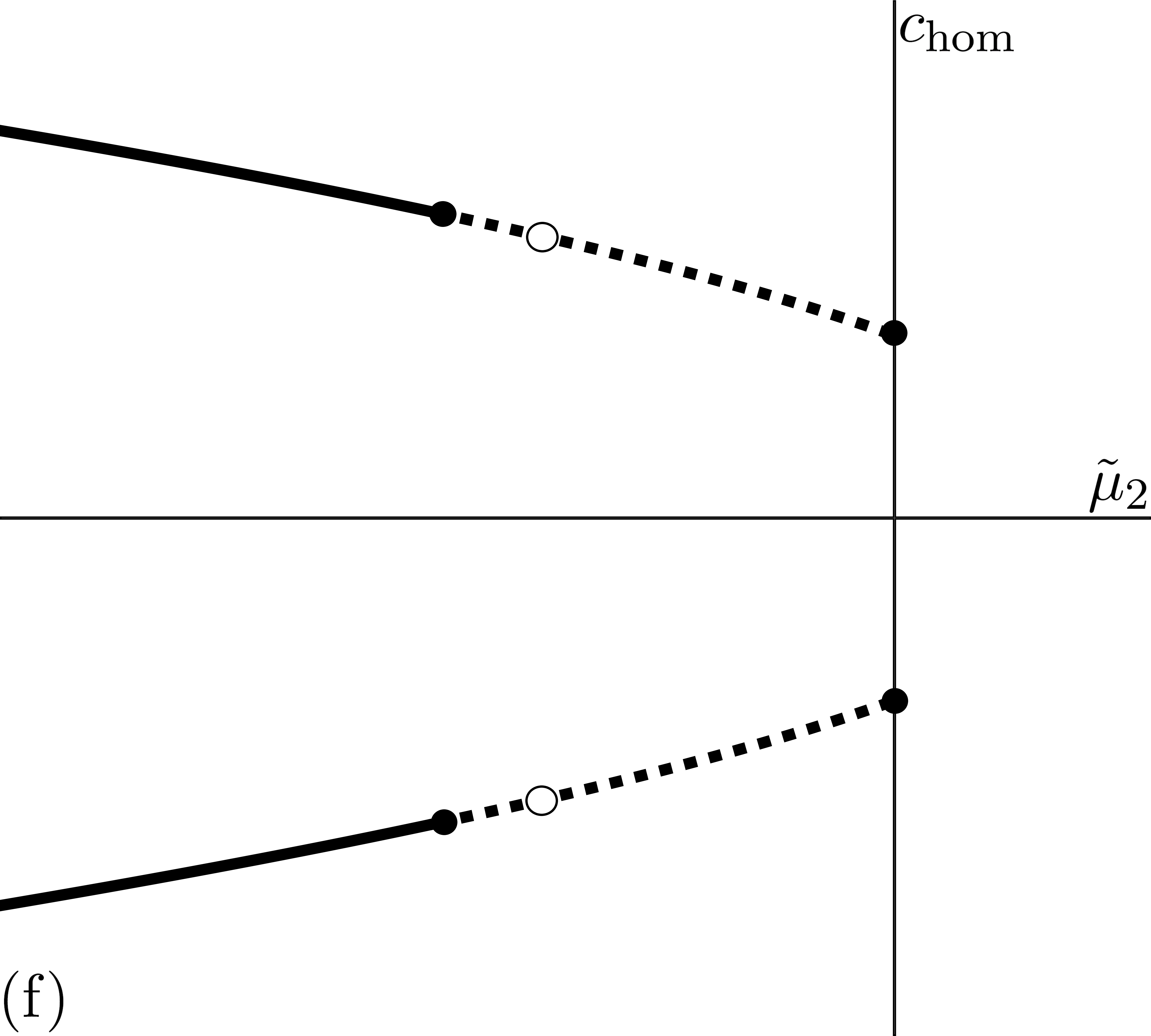}
\end{minipage}
\caption{\small{(a, b, c) The bifurcation diagrams of the heteroclinic fronts described by Corollary \ref{c:hetbifs} for $\tM_{2cc}(\vmuta) < 0$ of Fig. \ref{f:hetbifs}, now including the stability results of the present section (under the additional assumption that $\hN_{2 \la \la}(\vmuta) > 0$). (a) $\tM_{2}(\vmuta) = 0$. (b) $\tM_{2}(\vmuta) > 0$. Note that the bullets on the $\{\tmu_1=0\}$-axis indicate the bifurcation at which the traveling pulses split up in pairs of traveling fronts (Corollary \ref{c:hombifs}(ii), (iii) and Figs. \ref{f:hombifs}(b, c) and (e, f) below). (c) Corollary \ref{c:hetbifs}(iii): $\tM_{2}(\vmuta) < 0$. (d, e, f): The bifurcation diagrams of the homoclinic pulses described by Corollary \ref{c:hombifs} for $\oG_{1c}(\vmuta) < 0$ of Fig. \ref{f:hombifs}, now including stability (under the additional assumptions that $\hN_{2 \la \la}(\vmuta) > 0$ and that $\Q(c_0^2, c_0^2 \mu_N)$, i.e. $\laho^{1,2}(c_0)$ (\ref{d:la12hom-next}), does not change sign along the depicted branches of traveling pulses, cf. Corollary \ref{c:Q00}). (d) $\tM_{2}(\vmuta) > 0, \tM_{2cc}(\vmuta) > 0$. (e) $\tM_{2}(\vmuta) > 0, \tM_{2cc}(\vmuta) < 0$ and (f) $\tM_{2}(\vmuta) < 0, \tM_{2cc}(\vmuta) > 0$, where we note that the open dots represent the points $(\tmu_{\rm hom}^{\rm m}, \pm c_0^{\rm m})$ (\ref{d:c0etmue}), i.e. the $c_0= c_0(\tmu_2)$-values at which leading order eigenvalue $\la_{1, {\rm h}}(c_0)$ (\ref{d:la2hom}) moves from $\RR$ to $i\RR$ or vice versa (\ref{e:laho2pm=Im}), and that the order on the branch of traveling pulses between this bullet and the bullet that indicates the Hopf bifurcation associated to the eigenvalue pair $\laho^{2,\pm}(c_0)$ (\ref{d:la2hom-next}) may change by varying $\mu_N$ -- so that the Hopf bifurcation may disappear.}}
\label{f:hethombifs-stab}
\end{figure}

In the upcoming stability analysis, a crucial role is played by the sign of the (Melnikov-type) expression $\hN_{2\la\la}(\vmu)$ at $\vmu = \vmu_{\rm t}^\ast$,
\beq
\label{d:hN2ll}
\hN_{2\la\la}(\vmu_{\rm t}^\ast) = \int_{-\infty}^{\infty} \hJ_{2\la\la}(X;\vmu_{\rm t}^\ast) \vax(X) \, dX =
\int_{-\infty}^{\infty}\left[(1-\tau \frac{\fa' \Ga_u}{\Fa_u})\tv_1 + \frac{\tau^2\Ga_u\fa'\vax}{(\Fa_u)^2}\right] \vax \, dX
\eeq
(\ref{d:hJ2ll}). The situation is similar to that of $\tM_{2cc}$ in (\ref{d:tM2c2}): for general $\mu$, the term $(1-\tau \frac{\fa' \Ga_u}{\Fa_u})\tv_1 \vax$ does not decay as $X \to \infty$ (\ref{e:propsvbu}), (\ref{e:tv1gg1-0}), so that this integral does not converge. However, $\Ma(\vmu_{\rm t}^\ast) = 0$, so that this term does decay at $\vmu=\vmu_{\rm t}^\ast$ (see also $\vmu_1 = 0$ in (\ref{e:tv1gg1})): $|\hN_{2\la\la}(\vmu_{\rm t}^\ast)| < \infty$. Throughout the rest of this paper we impose the non-degeneracy condition $\hN_{2\la\la}(\vmu_{\rm t}^\ast) \neq 0$ (see however section  \ref{ss:Bifs}). Moreover, we also define the two similar expressions,
\beq
\label{d:hN2s}
\hN_{2c\la}(\vmu_{\rm t}^\ast) =  \int_{-\infty}^{\infty} \hJ_{2c\la}(X; \vmu_{\rm t}^\ast) \, \vax(X) \, dX, \;
\hN_{2cc}(\vmu_{\rm t}^\ast)  =  \int_{-\infty}^{\infty} \hJ_{2cc}(X; \vmu_{\rm t}^\ast) \, \vax(X) \, dX
\eeq
(\ref{d:hJ2s}) and conclude after a careful study of (\ref{d:hI2s}), (\ref{d:hJ2s}) that also $|\hN_{2c\la}(\vmu_{\rm t}^\ast)|, |\hN_{2cc}(\vmu_{\rm t}^\ast)| < \infty$.
\\ \\
We first consider the asymptotically small eigenvalues for the traveling patterns in three separate lemmas -- Lemma \ref{l:asymptevhet} for the heteroclinic fronts of Theorem \ref{t:exhets}(ii) and Lemmas \ref{l:asymptevhom} and \ref{l:asymptevhom-hot} for the nearly heteroclinic pulses of Theorem \ref{t:exhoms}(ii) -- before we formulate the main result(s) of this section, Theorem \ref{t:nearstabtravfrontspulses}. As a corollary to Theorem \ref{t:nearstabtravfrontspulses} -- and its preceding lemmas -- we consider the so far missing case of standing fronts pulses for $|\Ma(\vmu)|$ asymptotically small (Corollary \ref{c:stabstandfrontpulseMsmall}).
\begin{lemma}
\label{l:asymptevhet}
Consider the setting of Theorem \ref{t:exhets}(ii) and let $(\Uhe(X),\Vhe(X))$ be the traveling front established by Theorem \ref{t:exhets}(ii) that travels with speed $\che(\vmu_1)$ at leading order given by the solution $c_0 = c_0(\vmu_1)$ of (\ref{e:exhets}) (with $\vmu_1$ as in (\ref{d:vmutsi}) with $\tsi = 1$). Then, the spectral problem associated to the stability of $(\Uhe(X),\Vhe(X))$ as solution of (\ref{e:RDE-S1}) has 2 asymptotically small eigenvalues, the translational $\lahe^1(\vmu_1) \equiv 0$ and
\beq
\label{d:la2het}
\lahe^2(\vmu_1) = \frac{\vec{\nabla} \Ma \cdot \vmu_1 + c_0(\vmu_1) \hN_{2c\la}}{\hN_{2\la\la}} \eps + \O(\eps^2) = \pm \frac{\sqrt{\tmu_1^2 - 4\tM_{2cc} \tM_2}}{\hN_{2\la\la}} \eps + \O(\eps^2)
\eeq
(\ref{d:tmu1}), under the non-degeneracy condition that $\hN_{2\la\la}(\vmu_{\rm t}^\ast) \neq 0$ and where the $\pm$ is decided by solution $c_0(\vmu_1)$ of (\ref{e:exhets}).
\end{lemma}
\noindent
Note that $\lahe^2(\vmu_1) \in \RR$ by Corollary \ref{c:hetbifs} and that its sign is (also) determined by the sign of $\hN_{2\la\la}(\vmuta)$. Note furthermore that in the limit $\vmu_1 \to 0$, i.e the case that coincides with the $\Ma(\vmu) = O(\eps^2)$ setting of Theorem \ref{t:exhoms} in which (\ref{e:exhets}) reduces to $c_0^2 \tM_{2cc}(\vmu_{\rm t}^\ast) + \tM_{2}(\vmu_{\rm t}^\ast) = 0$, the stability of the front is determined by (the sign of) its speed. Thus, of the two traveling fronts connecting $(\bU_-,\bV_-)$ for $X \to -\infty$ to $(\bU_+,\bV_+)$ for $X \to \infty$  that exist in this limit (Corollary \ref{c:hetbifs}(ii) and Figs. \ref{f:hetbifs}(b) and \ref{f:hethombifs-stab}(b)), only one may be stable. Naturally, their two symmetrical counterparts -- i.e. the fronts connecting $(\bU_+,\bV_+)$ for $X \to -\infty$ to $(\bU_-,\bV_-)$ for $X \to -\infty$ induced by symmetry (\ref{d:symm}) -- have the same stability characteristics. Together, these pairs of fronts determine a bifurcation of the traveling homoclinic pulses (if $\tM_{2} \tM_{2cc} < 0$, Corollary \ref{c:hombifs}(ii), (iii) and Figs. \ref{f:hombifs}(b,c) and \ref{f:hethombifs-stab}(e,f)). Thus, it is natural to expect that the limiting traveling homoclinic pulses of Theorem \ref{t:exhoms}(ii) with (limiting) speed $c_0 \to c_0^{\rm h-c} = \sqrt{-\tM_{2}/\tM_{2cc}} = \sqrt{-2\tM_{2}/\hN_{2c\la}}$ (\ref{e:hN2tM}) have two eigenvalues given by (\ref{d:la2het}) with $\tmu_1 = 0$ and $c_0$ replaced by $\pm c_0^{\rm h-c}$. This is indeed confirmed by Lemma \ref{l:asymptevhom}.
\\ \\
{\bf Proof of Lemma \ref{l:asymptevhet}.} This proof follows the by now standard approach of this paper, therefore we keep the analysis compact. First, we conclude from the $\O(\eps)$ problem that
\beq
\label{e:bvOeps}
\bv = \vax + \eps(-\la_1 \tv_1 + c_0 \hv_{1c}) + \O(\eps^2 \bv_2)
\eeq
(\ref{e:OepsEigPb}), (\ref{e:bv1lc}). We define $\hN_{1c}$ by
\beq
\label{d:hN1c}
\hN_{1c}(\vmu) = \int_{-\infty}^{\infty} \hJ_{1c}(X; \vmu_{\rm t}^\ast) \, \vax(X) \, dX,
\eeq
and conclude from (\ref{d:hJ1c}) that $|\hN_{1c}(\vmu)| < \infty$. We know by Lemma \ref{l:Lah} and (\ref{d:vmutsi}), with $\tsi = 1$, that $\bv(X)$ grows exponentially if $c_0 \hN_{1c}(\vmu_{\rm t}^\ast) \neq 0$. Since we know that the translational eigenvalue $\lahe^1 = 0$ must exist, we conclude that $\hN_{1c}(\vmu_{\rm t}^\ast) = 0$ for non-stationary fronts. Since $\hN_{1c}(\vmu)$ only depends on $\eps$ through $\vmu$ (\ref{d:hJ1c}), (\ref{d:hN1c}) it follows that
\beq
\label{d:expandhN1c}
\hN_{1c}(\vmu_{\rm t}^\ast + \eps^{\tsi} \vmu_{\tsi}) = \vec{\nabla} \hN_{1c}(\vmu_{\rm t}^\ast) \cdot \vmu_{\tsi} \, \eps^{\tsi} + \O(\eps^{\tsi + 1})
\eeq
(cf. (\ref{d:vmutsi})), with here $\tsi = 1$. Approximation $\bv_2(X)$ is determined as solution of (\ref{e:Oeps2EigPb}) and its leading order behavior for $X \gg 1$ is once again controlled by Lemma \ref{l:Lah}. Naturally, we impose that also at $\O(\eps^2)$ $\la_1 = 0$ must yield a bounded -- in fact converging -- solution (that corresponds to the $\O(\eps^2)$ approximation of $v_{{\rm het}, X}(X)$). Thus, we conclude by (\ref{e:Oeps2EigPb}), (\ref{d:hN2s}), (\ref{d:expandhN1c}) that
\beq
\label{d:condder2}
c_0^2 \hN_{2cc}(\vmu_{\rm t}^\ast) + \hN_{2}(\vmu_{\rm t}^\ast) + c_0 \vec{\nabla} \hN_{1c}(\vmu_{\rm t}^\ast) \cdot \vmu_1 = 0.
\eeq
For $\la_1 \neq 0$, $\bv_2(X)$ grows exponentially for $X \gg 1$, unless
\beq
\label{e:la1stravfront}
\la_1 \vec{\nabla} \Ma(\vmu_{\rm t}^\ast) \cdot \vmu_1 - \la_1^2 \hN_{2\la\la} (\vmu_{\rm t}^\ast) + c_0 \la_1 \hN_{2c\la} (\vmu_{\rm t}^\ast) =0
\eeq
(\ref{d:hN2ll}), (\ref{d:hN2s}) -- from which the first part of (\ref{d:la2het}) follows. Next, we note that $\lahe^2(\vmu_1)$ necessarily has to be $0$ at the saddle-node bifurcations of Corollary \ref{c:hetbifs}(iii) -- see Fig. \ref{f:hetbifs}(c) (and we observe that $\lahe^2 = 0$ at the transcritical bifurcation at $\tmu_1 = c_0 = 0$ for $\tM_2 = 0$ of Corollary \ref{c:hetbifs}(i), see Fig. \ref{f:hetbifs}(a)). Re-introducing $\tmu_1$ (\ref{d:tmu1}) and using that $c_0(\tilde{\mu}_{{\rm het}-SN}) = -\tmu_{{\rm het}-SN}/(2 \tM_{2cc})$ (Corollary \ref{c:hetbifs} (iii)), we find
\[
\lahe^2(\tmu_{{\rm het}-SN})= \frac{\tmu_{{\rm het}-SN}}{\hN_{2\la\la}} \left(1 - \frac{\hN_{2c\la}}{2\tM_{2cc}}\right),
\]
from which we conclude that
\beq
\label{e:hN2tM}
\hN_{2c\la}(\vmu_{\rm t}^\ast) = 2\tM_{2cc}(\vmu_{\rm t}^\ast),
\eeq
see also Remark \ref{r:hN2cl2tM2cc}. The second part of (\ref{d:la2het}) now follows by solving $c_0$ from (\ref{e:exhets}). \hfill $\Box$
\begin{lemma}
\label{l:asymptevhom}
Consider the setting of Theorem \ref{t:exhoms}(ii) and let $(\Uho(X),\Vho(X))$ be a traveling pulse established by Theorem \ref{t:exhoms}(ii) that travels with speed $\cho(\tmu_2) = \pm \sqrt{C_{\rm hom}(\tmu_2)} = c_0 + \O(\eps)$ (Corollary \ref{c:hombifs}, (\ref{d:Ctmuhet})), so that $\cho^2 \tM_{2cc}(\vmu_{\rm t}^\ast) + \tM_{2}(\vmu_{\rm t}^\ast) > 0$ -- with $\tmu_2$ as defined in (\ref{d:tmu2}). Then, the spectral problem associated to the stability of $(\Uho(X),\Vho(X))$ as solution of (\ref{e:RDE-S1}) has 4 asymptotically small eigenvalues, 2 that are associated to the eigenvalue $\lahe^1$ of the heteroclinic stability problem considered in Lemma \ref{l:asymptevhet},
\beq
\label{d:la1hom}
\laho^{1,1}(\cho(\tmu_2)) \equiv 0, \; \; \laho^{1,2}(\cho(\tmu_2)) = 0 +\O(\eps^2 |\log \eps|),
\eeq
and 2 associated to the eigenvalue $\lahe^2$ (\ref{d:la2het}),
\beq
\label{d:la2hom}
\laho^{2,\pm}(\cho(\tmu_2)) = \pm \eps \frac{\sqrt{c^2_0 \hN_{2c \la}^2 - 2 \sqrt \al_+ \left(c_0^2 \tM_{2cc} + \tM_{2}\right)\hN_{2 \la \la}}}{\hN_{2 \la \la}} + \O(\eps^2|\log \eps|) \stackdef \pm \eps \la_{1, {\rm h}}(c_0) + \O(\eps^2|\log \eps|)
\eeq
with $\al_+$, $\tM_{2cc}(\vmu_{\rm t}^\ast)$, $\tM_{2}(\vmu_{\rm t}^\ast)$, $\hN_{2c \la}(\vmu_{\rm t}^\ast)$ and $\hN_{2 \la \la}(\vmu_{\rm t}^\ast) \neq 0$ as defined in (\ref{d:muEnpm}), (\ref{d:tM2c2}), (\ref{d:hN2ll}), (\ref{d:hN2s}).
\end{lemma}
\noindent
Note that the fact $\la = 0$ is a double eigenvalue (at leading order in $\eps$) was to be expected: it is caused by the (leading order) vertical character of the bifurcation into traveling waves -- see (the proof of) Theorem \ref{t:exhoms} and Corollary \ref{c:hombifs}. Note also that the second pair $\laho^{2,\pm}$ indeed merges with the (non-trivial) eigenvalue (\ref{d:la2het}) of the stability problem associated to the traveling front at the bifurcation at which the homoclinic pulse merges with a heteroclinic cycle that subsequently splits up into a 2 traveling fronts, i.e. for $\vmu_1 \to 0$ in (\ref{d:la2het}) and $c_0$ such that $c_0^2 \tM_{2cc} + \tM_{2} \downarrow 0$, i.e $c_0 \to c_0^{\rm h-c} = \sqrt{-\tM_{2}/\tM_{2cc}}$ -- see again Corollary \ref{c:hombifs} and Fig. \ref{f:hombifs}.
\\ \\
More importantly, it should be observed that unlike Lemma \ref{l:asymptevhet} for $(\Uhe(X),\Vhe(X))$, Lemma \ref{l:asymptevhom} does not provide decisive insight in the (potential) stability properties of $(\Uho(X),\Vho(X))$. At least, if $\hN_{2 \la \la}(\vmu_{\rm t}^\ast) < 0$ then obviously $\laho^{2,+}(c_0) > 0$ (\ref{d:la2hom}), so that $(\Uho(X),\Vho(X))$ is unstable (recall that $c_0^2 \tM_{2cc}(\vmu_{\rm t}^\ast) + \tM_{2}(\vmu_{\rm t}^\ast) > 0$). However, if $\hN_{2 \la \la}(\vmu_{\rm t}^\ast) > 0$, a pair real eigenvalues $\laho^{2,\pm}(c_0)$ may merge and form a complex conjugate pair of purely imaginary eigenvalues as $|c_0| = |c_0(\tmu_2)|$ passes through a critical value $c_0^{\rm m}$ (see (\ref{d:c0etmue}) below) -- or vice versa. Although this is against a Sturm-Liouville based intuition, this is possible, for instance if additional to $\hN_{2 \la \la}(\vmu_{\rm t}^\ast) > 0$ also $\tM_{2}(\vmu_{\rm t}^\ast) > 0$, since in that case a bifurcation into traveling waves takes place at $c_0(\tmu_2) = 0$ (Corollary \ref{c:hombifs}, Fig. \ref{f:hombifs}). More specifically, 
\beq
\label{e:laho-TW}
\laho^{2,\pm}(0) = \pm \eps \frac{\sqrt{- 2 \sqrt{\al_+} \tM_{2} \hN_{2 \la \la} }}{\hN_{2 \la \la}} \in i \RR
\eeq
at leading order. Thus, in this case -- and in all cases for which Re$(\laho^{2,\pm}(c_0)) = 0$ at leading order, see (\ref{e:laho2pm=Im}) below -- we do need to perform a higher order analysis to determine the sign of Re$(\laho^{2,\pm}(c_0))$ and establish the possible (in)stability of $(\Uho(X),\Vho(X))$. Moreover, we also need to do a higher order analysis to determine the sign of $\laho^{1,2}(\vmu_1)$ (\ref{d:la1hom}). Nevertheless, at this point we may expect that all 3 nontrivial asymptotically small eigenvalues $\laho^{1,2}(\vmu_1)$ and $\laho^{2,\pm}(c_0)$ could possibly be in the stable (complex) half-plane and that the traveling pulse pattern $(\Uho(X),\Vho(X))$ may be stable as solution of (\ref{e:RDE-S1}). See Lemma \ref{l:asymptevhom-hot} and Theorem \ref{t:nearstabtravfrontspulses}(ii).
\\ \\
{\bf Proof of Lemma \ref{l:asymptevhom}.}
Unlike the proof of Lemma \ref{l:asymptevhom}, establishing expressions (\ref{d:la1hom}) and (\ref{d:la2hom}) does not go exactly along the lines of the approach so far developed. The main difference is that even for the homoclinic pulses we so far could focus on the first part of the solution that is at leading order determined by $\va(X)$ and approaches $\bV_-$ as $X \to -\infty$ -- except for the final part of the proof of Theorem \ref{t:exhoms}(ii) in which we unfolded the leading order vertical structure of the bifurcation diagrams by zooming in $\O(\eps^2 |\log \eps|)$ close to $\vmuta$ -- see Corollary \ref{c:hombifs} and Figs. \ref{f:HetHomBifsIntro}(b), \ref{f:hombifs}. Like in the proof of Theorem \ref{t:exhoms}(ii), we here need to explicitly consider the second part of the solution that is at leading order determined by $\va(-X)$ and approaches $\bV_-$ as $X \to \infty$.
\\
\begin{figure}[t]
\centering
	\begin{minipage}{.2425\textwidth}
		\includegraphics[width =\linewidth]{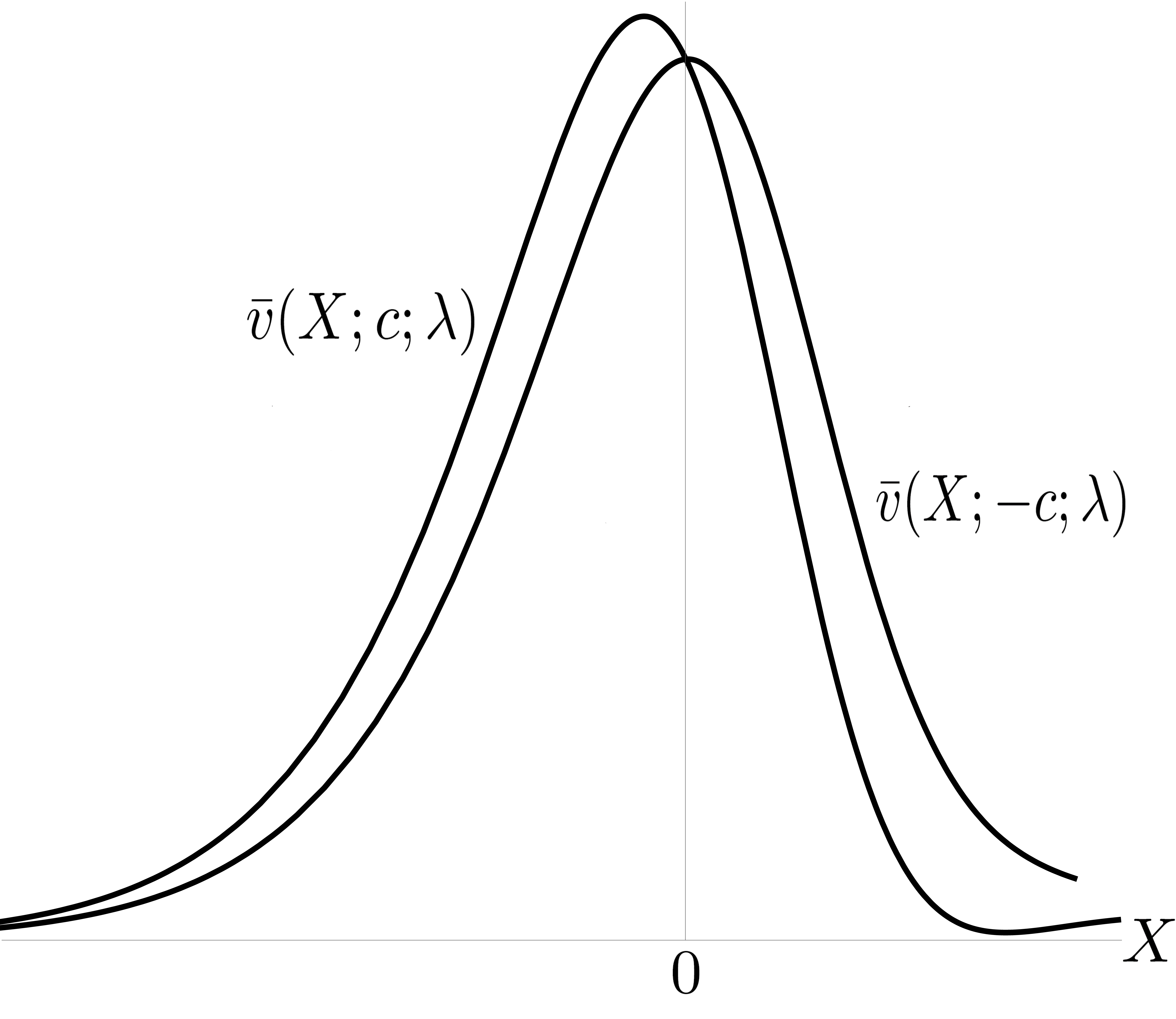}
	\end{minipage}%
	\hspace{.01cm}
	\begin{minipage}{0.2425\textwidth}
		\includegraphics[width=\linewidth]{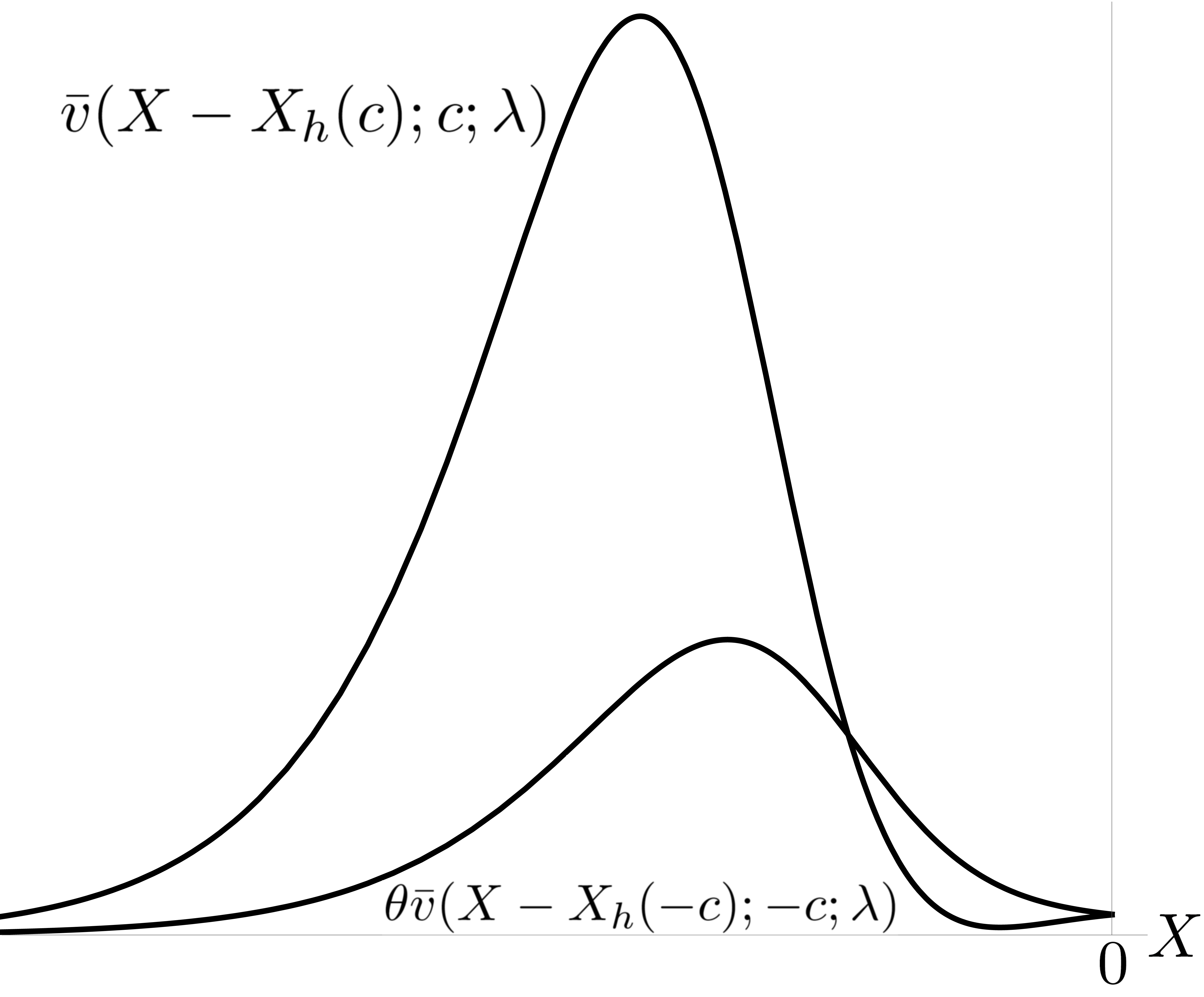}
	\end{minipage}
    \hspace{.01cm}
	\begin{minipage}{0.49\textwidth}
		\includegraphics[width=\linewidth]{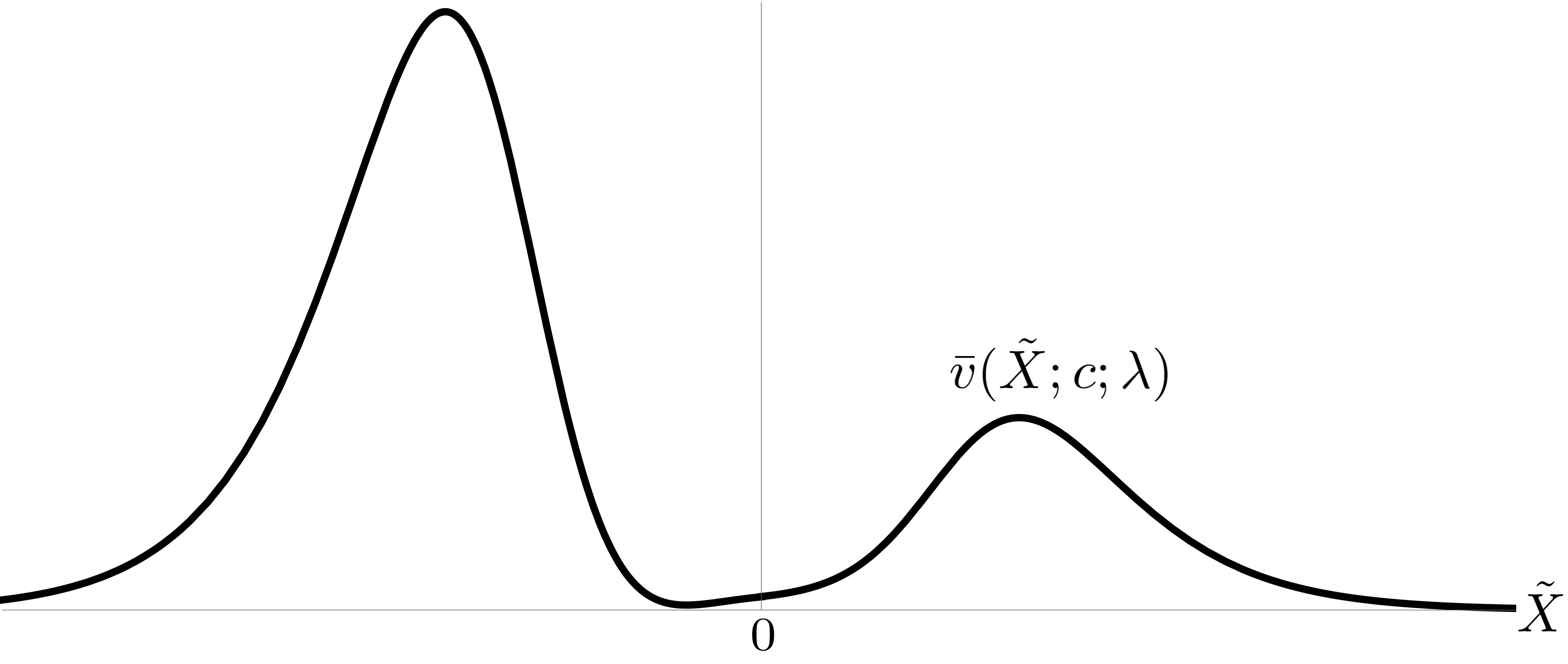}
	\end{minipage}
\caption{\small{The construction of the eigenfunctions associated to the stability of the nearly heteroclinic homoclinic pulse $(\Uho(X), \Vho(X))$. The first sketch shows $\bv(X; c; \la)$ for $X \in (-\infty, X_h(c))$ and $\bv(X; -c; \la)$ for $X \in (-\infty, X_h(-c))$, where we note that in general $X_h(c)) \neq X_h(-c)$. Next, the translated and scaled functions $\bv(X + X_h(c); c; \la)$ and $\theta \bv(X + X_h(-c); -c; \la)$ for $X \in (-\infty, 0)$ with $\theta$ and $\la$ such that the conditions in (\ref{e:bvsmoothXh}) are satisfied. The final sketch is obtained by the application of the $X \to -X$, $-c \to c$ symmetry (\ref{d:symm}) to $\theta \bv(X + X_h(-c); -c; \la)$, resulting in a smooth eigenfunction  $\bv(\tX; c; \la)$ defined on $\RR$.}}
\label{f:EigenFctsHom}
\end{figure}
\\
To explicitly construct an eigenfunction, and thus to determine the (asymptotically small) eigenvalues, we first determine a sufficiently accurate approximation of $\bv(X)$, at leading order given by $\vax(X)$ (\ref{e:bv0bu0}), for its first part, i.e. for $X \in (-\infty, X_h(c)]$, where we note that $\bv(X)$ of course must converge (exponentially fast) for $X \to \pm \infty$ and recall that $X_h(c)$ determines the center of the underlying pulse $\vho(X)$, since $\vho(X)$  attains its maximal value (by definition) at $X = X_h$ (\ref{e:Xh-1}), (\ref{e:vaXhnearhet}). Next, we set up the second part of the construction by first changing $c$ to $-c$ in the first approximation -- thus obtaining an approximating function on $(-\infty, X_h(-c)]$ (recall that $X_h(c) \neq X_h(-c)$ (\ref{d:Yhj}), (\ref{e:tY0}), (\ref{e:tY1})) -- followed by the application of symmetry (\ref{d:symm}) -- see Fig. \ref{f:EigenFctsHom}. Note that the idea behind this approach has already been developed at several places in the existence analysis and that it has been made explicit in the proof of Theorem \ref{t:exhoms}. However, there is a difference between the present setting and the existence analysis: (\ref{e:ODE-lin}) is a linear system, its solutions are defined up to a factor $\theta \in \CC$ (in general, $\theta \in \RR$ in the case of a real eigenvalue). Naturally, we may scale $\bv(X)$, i.e. we may choose $\theta$ to be 1 for the first part of $\bv(X)$, however, $\theta$ is not necessarily equal to $1$ as factor of the second part -- see Fig. \ref{f:EigenFctsHom} and (\ref{e:thetaO1}) below. Nevertheless, after the application of symmetry (\ref{d:symm}), the connection between the first part of $\bv(X) = \bv(X;c;\la)$ for $X \in (-\infty, X_h(c)]$ to its second part $\theta \bv(X; -c;\la)$ that is constructed for $X \in (-\infty, X_h(-c)]$ must be smooth as $X$ passes over the maximum of $\vho(X)$ -- which corresponds to $X = X_h(c)$ for the first part and $X = X_h(-c)$ for the second. This yields 2 conditions at the maximum of $\vho(X)$, on $\bv(X)$ and on $\bv_X(X)$,
\beq
\label{e:bvsmoothXh}
\begin{array}{ccc}
\bv(X_h(c); c; \la)
& = & \theta \bv(X_h(-c); -c; \la)
\\
\bv_X(X_h(c); c; \la)& = & - \theta \bv_X(X_h(-c); -c; \la)
\end{array}
\eeq
Note that the apparent non-smoothness associated to (\ref{e:bvsmoothXh}) is caused by the fact that we split up the approximation of $\bv$ into 2 distinct half-line approximations defined on different half-lines ($(-\infty, X_h(c)]$ and $(-\infty, X_h(-c)]$) and that the approximation over the full domain $\RR = (-\infty, 0] \cup [0, \infty)$ only followed after shifting both half-line approximations in $X$ and applying symmetry (\ref{d:symm}) to the second part (Fig. \ref{f:EigenFctsHom}). A priori, (\ref{e:bvsmoothXh}) only provides a $C^1$ smooth connection between both sides, however, we may conclude that the combined function $\bv(X)$ is smooth over $X=0$ (and thus on $\RR$) from the uniqueness of solutions to (\ref{e:ODE-lin}). Note also that for a given $c$, conditions (\ref{e:bvsmoothXh}) indeed determine discrete pairs $(\theta_j(c),\la_j(c))$, that $\la_j(c)$ thus are the asymptotically small eigenvalues associated to the stability of the pulse $(\Uho(X), \Vho(X))$ traveling with speed $c$, and that the full solution $(\bu(X),\bv(X))$ of (\ref{e:ODE-lin}) indeed must be an eigenfunction. Moreover, the geometry of the set-up guarantees the existence and local uniqueness of these (eigenvalue, eigenfunction) pairs -- in fact, the present set-up can be recast directly in terms of an Evans function approach \cite{AGJ90,DGK01}.
\\ \\
Now that we have deduced the eigenvalue conditions (\ref{e:bvsmoothXh}), we may return to the standard approach set up in this paper: we need to determine the leading order approximations of $\bv(X_h)$ and $\bv_X(X_h)$. As usual, and as in the proof of Lemma \ref{l:asymptevhet}, we start at the $\O(\eps)$ level and rederive (\ref{e:bvOeps}). This implies by Lemma \ref{l:Lah}, (\ref{d:Ma}), (\ref{e:Xh-1}), and (\ref{d:hN1c}) that (at leading order)
\[
\bv(X_h) = \frac{1}{\sqrt{2 \al_+}}\frac{c_0 \hN_{1c}(\vmu_{\rm t}^\ast)}{\sqrt{c_0^2 \tM_{2cc}(\vmu_{\rm t}^\ast) + \tM_{2}(\vmu_{\rm t}^\ast)}},
\]
since $\Ma(\vmu_{\rm t}^\ast) = 0$. Clearly, the derivative $v_{{\rm hom},X}(X)$ must be a eigenfunction (with trivial eigenvalue $\la_1 = 0$), hence its approximation must have $\bv(X_h) = 0$, which implies that $c_0 \hN_{1c}(\vmuta) = 0$ -- as in the heteroclinic case. Moreover, $\hN_{1c}(\vmu)$ must again be expanded -- as in the proof of Lemma \ref{l:asymptevhet} -- but now  with $\tsi = 2$ in (\ref{d:expandhN1c}). As a consequence -- and by (\ref{d:vmutsi}) with $\tsi = 2$ -- we may again conclude that we can employ the extension of Poincar\'e's Expansion Theorem as introduced in the proof of Theorem \ref{t:exhoms}. Hence, it follows by (\ref{e:Oeps2EigPb}), Lemma \ref{l:Lah}, (\ref{e:Xh-1}), (\ref{d:hN2ll}), (\ref{d:hN2s}), (\ref{d:expandhN1c}) that
\beq
\label{e:bvbvXOeps}
\begin{array}{ccl}
\bv(X_h; c_0; \eps \la_1) & = & \frac{\left[- \la_1^2 \hN_{2 \la \la} + c_0^2 \hN_{2cc} + \hN_{2} + c_0 \la_1 \hN_{2c \la}\right] + \sqrt \al_+ \left(c_0^2 \tM_{2cc} + \tM_{2}\right)}{\sqrt{2 \al_+(c_0^2 \tM_{2cc} + \tM_{2})}}  \eps
\\[3mm]
\bv_X(X_h; c_0; \eps \la_1) & = & \frac{\left[-\la_1^2 \hN_{2 \la \la} + c_0^2 \hN_{2cc} + \hN_{2} + c_0 \la_1 \hN_{2c \la}\right] - \sqrt \al_+ \left(c_0^2 \tM_{2cc} + \tM_{2}\right)}{\sqrt{2(c_0^2 \tM_{2cc} + \tM_{2})}}  \eps
\end{array}
\eeq
at leading order in $\eps$. We again use that $\la_1 = 0$ must be an eigenvalue with eigenfunction $v_{{\rm hom},X}(X)$ that has $\bv(X_h; c_0; 0) = 0$ and obtain the equivalent of (\ref{d:condder2}),
\[
c_0^2 \hN_{2cc}(\vmu_{\rm t}^\ast) + \hN_{2}(\vmu_{\rm t}^\ast) + \sqrt \al_+ \left(c_0^2 \tM_{2cc}(\vmuta) + \tM_{2}(\vmuta)\right) = 0.
\]
For notational simplicity we can now introduce $N_c = N_c(c_0^2)$
\beq
\label{d:Nc}
N_c (c_0^2) = -\left(c_0^2 \hN_{2cc}(\vmu_{\rm t}^\ast) + \hN_{2}(\vmu_{\rm t}^\ast)\right) = \sqrt \al_+ \left(c_0^2 \tM_{2cc}(\vmu_{\rm t}^\ast) + \tM_{2}(\vmu_{\rm t}^\ast)\right) > 0
\eeq
(cf. Theorem \ref{t:exhoms}(ii)). Combining (\ref{e:bvsmoothXh}) with (\ref{e:bvbvXOeps}) and (\ref{d:Nc}) yields,
\beq
\label{e:thetalaOeps2}
\begin{array}{rcl}
-\la_1^2 \hN_{2 \la \la} + c_0 \la_1 \hN_{2c \la}
& = & \theta \left[-\la_1^2 \hN_{2 \la \la} - c_0 \la_1 \hN_{2c \la}\right]
\\[2mm]
- 2 N_c -\la_1^2 \hN_{2 \la \la} + c_0 \la_1 \hN_{2c \la}
& = & \theta \left[2 N_c + \la_1^2 \hN_{2 \la \la} + c_0 \la_1 \hN_{2c \la}\right]
\end{array}
\eeq
so that it follows by eliminating $\theta$ that
\beq
\label{e:quartla1}
(\hN_{2 \la \la})^2 \la_1^4 - \left[c^2_0 (\hN_{2c \la})^2 - 2 N_c \hN_{2 \la \la} \right] \la_1^2 = 0,
\eeq
which indeed yields the leading order terms of (\ref{d:la1hom}) and (\ref{d:la2hom}). For future reference, i.e. not as necessary ingredient of this proof, we define
\beq
\label{d:theta0cl}
\theta_0(c,\la) = - \frac{c \hN_{2c \la} - \la \hN_{2 \la \la}}{c \hN_{2c \la} + \la \hN_{2 \la \la}},
\eeq
and note that the $\lambda$-solutions of (\ref{e:thetalaOeps2}) are associated to three different $\theta$-solutions,
\beq
\label{e:thetaO1}
\theta_0 = -1, \; \;
\theta_0 = \theta_0(c_0,\pm \la_{1, {\rm h}}) = - \frac{\left(c_0 \hN_{2c \la} \mp \sqrt{c^2_0 (\hN_{2c \la})^2 - 2 N_c \hN_{2 \la \la}} \right)^2}{2 N_c \hN_{2 \la \la}}
\eeq
with $\la_{1, {\rm h}}$ as defined in (\ref{d:la2hom}). Thus, $\theta_0(c_0,\pm \la_{1, {\rm h}}) \to 1$ as $c_0 \to 0$ -- i.e. at the bifurcation into traveling waves (Corollary \ref{c:hombifs}) -- but in general clearly $|\theta_\pm| \neq 1$: the amplitudes of the associated eigenfunctions for $X < 0$ and $X > 0$ indeed differ an $\O(1)$ factor -- see Fig. \ref{f:EigenFctsHom}.
\\ \\
We postpone the validation of the (magnitude of the) next order corrections to Lemma \ref{l:asymptevhom-hot}.  \hfill $\Box$
\begin{lemma}
\label{l:asymptevhom-hot}
Consider the setting of Lemma \ref{l:asymptevhom-hot} and let $ \laho^{1,2}(\cho(\tmu_2))$ and $\laho^{2,\pm}(\cho(\tmu_2))$ be the 3 nontrivial asymptotically small eigenvalues of the spectral problem associated to the stability of the pulse solution $(\Uho(X),\Vho(X))$  of (\ref{e:RDE-S1}) established by Theorem \ref{t:exhoms}(ii). Define $\mu_N$  by
\beq
\label{d:muN}
\mu_N = \vnab \hN_{1c}(\vmuta) \cdot \tilde{\vmu}_2
\eeq
(\ref{d:hN1c}) and assume that dim$($span$\{\vnab \Ma(\vmuta), \vnab \hN_{1c}(\vmuta) \}) = 2$, so that for $\vmu = \vmuta + \eps^2 \tilde{\vmu} |\log \eps| \in \RR^m$ and $\tmu_2 = \vnab M_\ast(\vmuta) \cdot \tilde{\vmu}_2$ fixed (\ref{d:tmu2}), $\mu_N = \O(1)$ can be varied over $\RR$.
\\
{\bf (i)} Define $\A(c,\la)$, $\B(c,\la)$, $\C(c,\la)$ and $\D(c,\la)$ by
\beq
\label{d:ABCD}
\begin{array}{ll}
\A(c,\la) =
\frac{2c(c \la \hN_{2 c \la} + \N(1-\theta_0))}{\D(c,\la)}
&
\B(c,\la) =
\frac{\oG^+_{1c}(c \la \hN_{2 c \la} - \N(1))}{2 \sqrt{\al_+}},
\\
\C(c,\la) =
\frac{\la (1+\theta_0)(c \la \hN_{2 c \la} + \N(2))\oG^+_{1c} N_c}{\al_+ \D(c,\la)}
&
\D(c,\la) = c \hN_{2 c \la} \N(3 + \theta_0) + \la \hN_{2\la \la} \N(2 + 2\theta_0)
\end{array}
\eeq
with $\oG^+_{1c}$, $\theta_0 = \theta_0(c, \la)$ as defined in (\ref{d:oG1+}), (\ref{d:theta0cl}) and,
\beq
\label{d:Nz}
\N(z; c^2) = c^2 \frac{\hN_{2c \la}^2}{\hN_{2\la \la}} - z N_c(c^2) \stackdef c^2 N_r - z N_c(c^2)
\eeq
(\ref{d:hN2ll}), (\ref{d:hN2s}), (\ref{d:Nc}). Assume that $\hN_{2\la \la}(\vmuta) \neq 0$ and $\D(c,\la) \neq 0$. Then, the next order correction to $\laho^{2,\pm}(\cho(\tmu_2))$ is given by
\beq
\label{d:la2hom-next}
\laho^{2,\pm}(\cho(\tmu_2)) =  \pm \eps \la_{1, {\rm h}} +
\eps^2\left[\A(c_0, \pm \la_{1, {\rm h}}) \left(\mu_N - \B(c_0, \pm \la_{1, {\rm h}}) \right) + \C(c_0, \pm \la_{1, {\rm h}}) \right]|\log \eps| + \O(\eps^2),
\eeq
with $\la_{1, {\rm h}}$ as defined in (\ref{d:la2hom}) and $\cho(\tmu_2) = c_0 + \O(\eps)$ ((\ref{d:Ctmuhet}), Lemma \ref{l:asymptevhom}).
\\
{\bf (ii)} The next order correction to $\laho^{1,2}(\cho(\tmu_2))$ is given by
\beq
\label{d:la12hom-next}
\laho^{1,2}(\cho(\tmu_2)) =  \eps^2 c_0^2 \Q(c_0^2, c_0^2 \mu_N) |\log \eps| + \O(\eps^2)
\eeq
in which $\Q(C, C\mu) \in \RR$ is an expression that can be determined explicitly and that varies smoothly in $(C, \mu) \in \RR^2$.
\end{lemma}
\noindent
{\bf Proof of Lemma \ref{l:asymptevhom-hot}.} The approach of this proof is similar to that of the second part of the proof of Theorem \ref{t:exhoms}: we proceed beyond the leading order analysis of Lemma \ref{l:asymptevhom} using Lemma \ref{l:Lah-sharp} to filter out the dominant terms. Also similar to the second part of the proof of Theorem \ref{t:exhoms}, we will find that the leading order corrections to the $\O(\eps)$ approximations of (\ref{e:bvbvXOeps}) in the proof of Lemma \ref{l:asymptevhom} all originate from the $\O(\eps)$ term $\bv_1(X)$ in expansion (\ref{d:busbvsexp}) of $\bv(X)$.
\\ \\
Since $\bv_1 = -\la_1 \tv_1 + c_0 \hv_{1c}$ (\ref{e:bv1lc}), we may use (\ref{e:epsv1-leading}) as leading order approximation of the first part. The second part, i.e. $\hv_{1c}$, is determined by inhomogeneous problem (\ref{e:Lah}) with $h(X) = \hJ_{1c}(X)$ (\ref{e:OepsEigPb}), (\ref{d:hJ1c}). By Lemma \ref{l:Lah}, the growth of $\hv_{1c}(X)$ for $X \gg 1$ is determined by $\hN_{1c}(\vmu)$ (\ref{d:hN1c}) and we know from the proof of Lemma \ref{l:asymptevhom} that $\hN_{1c}(\vmuta) = 0$ (cf. (\ref{d:expandhN1c})). Thus, for $\tilde{\vmu}_2$ and $\mu_N$ as defined in (\ref{d:tmu2}), (\ref{d:muN}),
\beq
\label{e:expN1c}
\hN_{1c}(\vmu) = \hN_{1c}(\vmuta + \eps^2 |\log \eps| \, \tilde{\vmu}_2 + \O(\eps^2)) = \eps^2 \vnab \hN_{1c}(\vmuta) \cdot \tilde{\vmu}_2 |\log \eps| + \O(\eps^2) = \eps^2 \mu_N |\log \eps| + \O(\eps^2).
\eeq
However, this is not the only leading order term in the approximation of $\hv_{1c}(X)$ for $X \gg 1$: we need to apply the refinement of Lemma \ref{l:Lah}, Lemma \ref{l:Lah-sharp}. It follows from (\ref{d:hJ1c}) that the inhomogeneous problem for $\hv_{1c}$ is of the form considered in Lemma \ref{l:Lah-sharp} with $j=-1$ in (\ref{d:oh01j}). In fact, since $\tv_1(X), \tu_1(X) = O(X E_+(X))$ at leading order in $\eps$ for $X \gg1$ (\ref{e:tv1Xgg1-next}), (\ref{d:tu1}), we find by (\ref{d:fFGast}), (\ref{d:tFauuvv}) that,
\[
\begin{array}{ll}
\lim_{X \to \infty} \hJ_{1c}(X) e^{\sqrt{\al_+}X}
&=
\lim_{X \to \infty} \left[\left[\frac{\tau \Ga_u (\fa'\vax)_X}{\Fa_u} - \vaxx \right] +
\frac{(\tF^\ast_{uu} \fa' + \tF^\ast_{uu}) \Ga_u - (\tG^\ast_{uu} \fa' + \tG^\ast_{uu})\Fa_u}{\Fa_u} \vax
\right] e^{\sqrt{\al_+}X}
\\
&=
\lim_{X \to \infty} \left[\frac{\tau \Ga_u (\fa''(\vax)^2 + \fa' \vaxx)}{\Fa_u} - \vaxx \right]  e^{\sqrt{\al_+}X}
\\
&=
\lim_{X \to \infty} \vaxx e^{\sqrt{\al_+}X} \left[\frac{\tau \oG^+_u \fa'(\oV^+)}{\oF^+_u} - 1 \right]
=
\al_+ \be_+ \left[1- \frac{\tau \oG^+_u \fa'(\oV^+)}{\oF^+_u} \right] = \al_+ \be_+ \oG^+_{1c}
\end{array}
\]
(\ref{d:oFpmetc}), (\ref{e:propsvbu}), (\ref{d:oG1+}), at leading order in $\eps$. It thus follows by (\ref{d:approxsvX-sharp}) with $j=-1$ that, at leading order,
\beq
\label{e:bv1Xgg1}
\bv_1(X) = -\frac12 \be_+ (\la_1 + c_0\sqrt{\al_+}) \oG_{1c}^+ X e^{-\sqrt{\al_+}X} + \eps^2 \frac{(\la_1 \vec{\nabla} \Ma + c_0 \vnab \hN_{1c})\cdot \tilde{\vmu}_2}{2 \al_+ \be_+} |\log \eps| e^{\sqrt{\al_+}X}
\eeq
for $X \gg 1$ (cf. (\ref{e:tv1Xgg1-next})). It is straightforward to check that all other correction terms to (\ref{e:bvbvXOeps}) are of $\O(\eps^2)$, thus we find,
\beq
\label{e:bvbvX-next}
\begin{array}{rcl}
\bv(X_h; c_0; \eps \tla) & = & \eps \left[\frac{\sqrt{\al_+}}{\Yh} + \frac{\Yh(- \tla^2 \hN_{2 \la \la} + c_0^2 \hN_{2cc} + \hN_{2} + c_0 \tla \hN_{2c \la})}{2 \al_+}\right]
\\
& & \; \; \; \; \; + \eps^2 \left[ \frac{\Yh(\tla \vec{\nabla} \Ma + c_0 \vnab \hN_{1c})\cdot \tilde{\vmu}_2}{2 \al_+} - \frac{(\tla + c_0\sqrt{\al_+}) \oG_{1c}^+}{2 \sqrt{\al_+} \Yh} \right] |\log \eps| + \O(\eps^2)
\\[3mm]
\bv_X(X_h; c_0; \eps \tla) & = & \eps \left[-\frac{\al_+}{\Yh} + \frac{\Yh(- \tla^2 \hN_{2 \la \la} + c_0^2 \hN_{2cc} + \hN_{2} + c_0 \tla \hN_{2c \la})}{2 \sqrt{\al_+}}\right]
\\
& & \; \; \; \; \; + \eps^2 \left[ \frac{\Yh(\tla \vec{\nabla} \Ma + c_0 \vnab \hN_{1c})\cdot \tilde{\vmu}_2}{2 \sqrt{\al_+}} + \frac{(\tla + c_0\sqrt{\al_+}) \oG_{1c}^+}{2\Yh} \right] |\log \eps| + \O(\eps^2)
\end{array}
\eeq
with $\Yh$ as in (\ref{d:Yhj}) and $\tla = \la_1$ at leading order. Note that we indeed recover leading order approximation (\ref{e:bvbvXOeps}) when we replace $\Yh = Y_0 + \eps c_0 \tY_1 |\log \eps| + \O(\eps)$ (\ref{d:Yhj}), (\ref{e:tY0}), (\ref{e:tY1}) by $Y_0$ and $\tla$ by $\la_1$. Clearly, we must expect that the next order corrections to Lemma \ref{l:asymptevhom}'s $\laho^{1,2}$  and $\laho^{2,\pm}$ are of $\O(\eps^2 |\log \eps|)$, therefore we adapt expansion (\ref{d:laj}) and expand $\theta$ (cf. (\ref{e:bvsmoothXh})) likewise,
\beq
\label{d:tla2ttheta1}
\la = \eps \tla = \eps \la_1 + \eps^2 \tla_2 |\log \eps| + \O(\eps^2), \; \; \theta = \theta_0 + \eps \tth_1 |\log \eps| + \O(\eps).
\eeq
Substitution of all these expansions into (\ref{e:bvsmoothXh}) yields -- after considerable algebra -- a $2 \times 2$ system for $(\tla_2,\tth_1)$,
\beq
\label{d:MR}
\mathbb{M}(c_0,\la_1)
\left(
\begin{array}{c}
\tla_2
\\
\tth_1
\end{array}
\right)
=
\left(
\begin{array}{cc}
m_{11}(c_0,\la_1) & m_{12}(c_0,\la_1)
\\
m_{21}(c_0,\la_1) & m_{22}(c_0,\la_1)
\end{array}
\right)
\left(
\begin{array}{c}
\tla_2
\\
\tth_1
\end{array}
\right)
=
\left(
\begin{array}{c}
r_1(c_0, \la_1, \mu_N)
\\
r_2(c_0, \la_1, \mu_N)
\end{array}
\right)
\eeq
with,
\beq
\label{e:Mcl}
\left(
\begin{array}{cc}
m_{11} & m_{12}
\\
m_{21} & m_{22}
\end{array}
\right)
=
\left(
\begin{array}{cc}
-2 \la_1 \hN_{2\la \la}(1-\theta_0) + c_0 \hN_{2 c \la}(1+\theta_0) & \la_1^2 \hN_{2\la \la} + c_0 \la_1 \hN_{2 c \la}
\\
-2 \la_1 \hN_{2\la \la}(1+\theta_0) + c_0 \hN_{2 c \la}(1-\theta_0) & -\la_1^2 \hN_{2\la \la} - c_0 \la_1 \hN_{2 c \la} -2N_c
\end{array}
\right)
\eeq
(\ref{d:Nc}), and
\beq
\label{e:Rclm}
\left(
\begin{array}{c}
r_1
\\
r_2
\end{array}
\right)
=
- \left(
\begin{array}{c}
c_0 (1+\theta_0)
\left[
\mu_N + \frac{\left( \la_1^2 \hN_{2\la \la} - c_0 \la_1 \hN_{2 c \la} + N_c \right)\oG^+_{1c}}{2 \sqrt{\al_+}}
\right]
\\
\frac{\la_1 (1 + \theta_0) \oG^+_{1c} N_c}{\al_+} + c_0 (1-\theta_0)
\left[
\mu_N + \frac{\left( \la_1^2 \hN_{2\la \la} - c_0 \la_1 \hN_{2 c \la} + N_c \right)\oG^+_{1c}}{2 \sqrt{\al_+}}
\right]
\end{array}
\right)
\eeq
(\ref{d:oG1+}), with $\mu_N$ as defined in (\ref{d:muN}) -- where it should be noted that we have used that
\[
\tmu_2 = \vnab \Ma(\vmuta) \cdot \tilde{\vmu}_2 = \frac{\oG^+_{1c} N_c}{2 \al_+}
\]
(\ref{d:Vtc}) -- at leading order -- and the subsequent simplifications
\[
Y_0 = \sqrt{2} \frac{\al_+^{\frac{3}{4}}}{\sqrt{N_c}}, \; \; \tY_1 = -\frac12 \sqrt{2} \frac{\al_+^{\frac{1}{4}} \oG^+_{1c}}{\sqrt{N_c}}
\]
of (\ref{e:tY0}), (\ref{e:tY1}) to obtain the expressions (\ref{e:Rclm}) for $r_{1,2}(c_0,\mu_N)$. To solve (\ref{d:MR}), we first need to evaluate the determinant $D(c_0,\la_1)$ of $\mathbb{M}(c_0,\la_1)$,
\beq
\label{d:D-1}
D(c_0,\la_1) =
2 \la_1 \left[ 2 \la_1^2 \hN_{2\la \la}^2 + c_0\la_1 \hN_{2\la \la} \hN_{2c \la} - c_0^2 \hN_{2 c \la}^2 \right] +
2 N_c \left[ 2 \la_1 (1-\theta_0) \hN_{2\la \la} - c_0 (1+\theta_0) \hN_{2 c \la} \right].
\eeq
Thus, $D(c_0) \equiv 0$ if $\la_1 = \laho^{1,2} = 0$ (since then also $\theta_0 = -1$ (\ref{e:thetaO1})). In fact, in this case (\ref{d:MR}) reduces to just one equation for $\tla_2$ and $\tth_1$
\beq
\label{d:2x2reduction}
c_0 \hN_{2 c \la} \tla_2 - N_c \, \tth_1 = c_0
\left[
\mu_N + \frac{\oG^+_{1c} N_c}{2 \sqrt{\al_+}}
\right]
\eeq
and one thus needs to go deeper in the perturbation analysis to determine $\tla_2$ -- see below.
\\ \\
{\bf (i) The next order corrections to $\laho^{2,\pm}(\cho(\tmu_2))$.}
Since for $\la_1 \neq 0$,
\beq
\label{e:N2}
\hN_{2 \la \la}^2 \la_1^2 = c^2_0 \hN_{2c \la}^2 - 2 N_c \hN_{2 \la \la} = \hN_{2 \la \la} \, \N(2; c_0^2)
\eeq
(\ref{e:quartla1}), (\ref{d:Nz}), we can eliminate the $\la_1^2$ terms from $D(c_0,\la_1)$ (and conclude that indeed $D(c_0,\la_1) = \D(c_0,\la_1)$),
\beq
\label{d:D-2}
\begin{array}{rl}
D(c_0,\la_1) & =
2 \la_1 \left[ c_0\la_1 \hN_{2\la \la} \hN_{2c \la} + c_0^2 \hN_{2 c \la}^2 - 4 N_c \hN_{2 \la \la}\right] +
2 N_c \left[ 2 \la_1 (1-\theta_0) \hN_{2\la \la} - c_0 (1+\theta_0) \hN_{2 c \la} \right]
\\[2mm]
& = 2 c_0 \frac{\hN_{2c \la}}{\hN_{2\la \la}}\left[c^2_0 \hN_{2c \la}^2 - 2 N_c \hN_{2 \la \la}\right] + 2 \la_1 c_0^2 \hN_{2 c \la}^2 -4 \la_1 (1+\theta_0) N_c \hN_{2\la \la} - 2 c_0 (1+\theta_0) N_c \hN_{2 c \la}
\\[2mm]
& = 2 \la_1 \left[c_0^2 \hN_{2 c \la}^2 - 2 (1+\theta_0) N_c \hN_{2\la \la}\right] + 2 c_0 \hN_{2 c \la} \left[ c_0^2 \frac{\hN_{2c \la}^2}{\hN_{2\la \la}} - (3 + \theta_0) N_c \right]
\\[2mm]
& = 2 \la_1 \hN_{2\la \la} \N(2 + 2\theta_0) + 2 c_0 \hN_{2 c \la} \, \N(3 + \theta_0) = \D(c_0, \la_1)
\end{array}
\eeq
with $\D(c_0, \la_1)$ and $\N(z)$ as defined in (\ref{d:ABCD}), (\ref{d:Nz}), where we note that $\D(c_0, \la_1) \to 0$ as $\la_1 \to 0$ since $\N(3 + \theta_0) \to \N(2) \to 0$ as $\la_1 \to 0$ (\ref{d:theta0cl}), (\ref{e:N2}). Assuming that det$(\mathbb{M}(c_0,\la_1)) \neq 0$, i.e. that $\D(c_0,\la_1) \neq 0$, yields (after eliminating all $\la_1^2$ terms as above),
\[
\tla_2(c_0,\mu_N) = \A(c_0,\la_1)(\mu_N - \B(c_0,\la_1)) + \C(c_0,\la_1),
\]
as in (\ref{d:la2hom-next}) -- with $\A(c_0,\la_1)$, $\B(c_0,\la_1)$ and $\C(c_0,\la_1)$ as defined in (\ref{d:ABCD}).
\\ \\
{\bf (ii) The next order correction to $\laho^{1,2}(\cho(\tmu_2))$.} Naturally, the analysis of the case $\laho^{1,2}(\cho(\tmu_2))$ automatically includes $\laho^{1,1}(\cho(\tmu_2)) \equiv 0$. Hence, as we noticed (and used) many times before: $\tla_2 = 0$ must be a solution of (\ref{d:2x2reduction}), which yields $\tth_2$ is determined uniquely and that thus
\beq
\label{d:theta12hom-next}
\theta_{\rm hom}^{1,2} (c_0, \mu_N) = \eps^2 \frac{c_0}{N_c} \left[\mu_N + \frac{\oG^+_{1c} N_c}{2 \sqrt{\al_+}}\right] |\log \eps| + \O(\eps)
\eeq
(recall that $N_c > 0$ (\ref{d:Nc})). Thus, we need to derive an equation for $\tla_2$ that has both $\tla_2 = 0$ as (trivial) solution as well as the (non-zero) leading order expression for $\laho^{1,2}(c_0)$: we need to obtain a quadratic equation involving $\tla_2$ and thus need to go to (at least) the $\O(\eps^4 |\log \eps|^2)$-level in the perturbation analysis -- note that this is similar to the $\O(\eps^2)$-analysis in the preceding Lemma where we determined 2 $\O(\eps)$ eigenvalues. Such a perturbation analysis can be done -- in principle -- although one should not underestimate the extent of the technical issues involved. For instance, it will be necessary to extend Lemma \ref{l:Lah} beyond Lemma \ref{l:Lah-sharp}, since we will need approximations as in (\ref{d:approxsvX-sharp}) for inhomogeneous terms $h(x)$ that behave at leading order like $h(X) = \oh_{0,i,j} X^i E_+^{-j}(X)$, $i \in \ZZ$, $i > 0$, for $X \gg 1$ (cf. (\ref{d:oh01j})).
\\ \\
We refrain from going deeper into the computational details of this analysis. However, we can make three basic observations. First, we know by Corollary \ref{c:hombifs} that the limit $c_0 \to 0$ corresponds to the bifurcation into traveling waves at
$\tmu_2 = \tmu_{\rm hom}^{\rm TW}$ (\ref{d:Ctmuhet}) at which the traveling pulses emerge from the stationary pulse. Since clearly $\laho^{2,\pm}(c_0) \nrightarrow 0$ as $c_0 \to 0$ in general (\ref{e:laho-TW}), it follows that $\laho^{1,2}(c_0) \rightarrow 0$ as $c_0 \to 0$. Second, it follows by symmetry (\ref{d:symm}) -- see also Theorem \ref{t:exhoms}(ii) -- that the stability of the traveling pulse solution $(\Uho(X),\Vho(X))$ of (\ref{e:RDE-S1}) with $\cho(\tmu_2) (= c_0 + \O(\eps)) > 0$ must be the same as that of its counterpropagating counterpart. Third, the term $\mu_N$ only appears with pre-factor  $c$  ($= c_0 + \eps c_1 + \O(\eps)^2$) in the approximation analysis. Therefore, we may conclude that $\laho^{1,2}(\cho(\tmu_2))$ can indeed be expressed as in (\ref{d:la12hom-next}) -- see also Corollary \ref{c:Q00} for a result on the sign of $\Q(c_0^2, c_0^2 \mu_N)$. \hfill $\Box$
\\ \\
Like $\laho^{1,2}(\cho(\tmu_2))$ (\ref{d:theta12hom-next}), also the eigenvalues $\laho^{2, \pm}(\cho(\tmu_2))$ should be symmetric under the $c_0 \to -c_0$ symmetry (\ref{d:symm}) -- since the counterpropagating pair of traveling pulses of Theorem \ref{t:exhoms}(ii) must have identical stability characteristics. This is a priori not obvious from (\ref{d:la2hom-next}). This is due to the fact that -- unlike for $\laho^{1,2}(\cho(\tmu_2))$ -- $\theta_0 \neq 1$ (\ref{e:thetaO1}), so that the $X \to -X, c_0 \to -c_0$ symmetry (\ref{d:symm}) is not built into the set-up of the analysis -- in fact $\theta_0(-c,\la) = 1/\theta_0(c,\la)$ (\ref{d:theta0cl}), see also Fig. \ref{f:EigenFctsHom}. Although leading order expression (\ref{d:la2hom}) for $\laho^{2,\pm}(\cho(\tmu_2))$ is symmetric under $c_0 \to -c_0$, the asymmetry of the set-up already shows up in leading order approximations (\ref{e:bvbvXOeps}), that clearly are not symmetric under $c_0 \to -c_0$. However, it is symmetric under $c_0 \to -c_0, \la_1 \to -\la_1$, as is (\ref{e:thetalaOeps2}): instead of the individual eigenvalues $\laho^{2, +}(c_0)$ and $\laho^{2, -}(c_0)$, the pair $\{ \laho^{2, +}(c_0), \laho^{2, -}(c_0) \}$ is invariant under $c_0 \to -c_0$ (i.e. in (\ref{d:la2hom}), $\laho^{2, +}(-\cho(\tmu_2)) = \laho^{2, +}(-c_0) = \laho^{2, -}(c_0) = \laho^{2, -}(\cho(\tmu_2))$ at leading order in $\eps$). This property persists at the next order(s). Since $\theta_0(-c_0,-\la_1) = \theta_0(c_0,\la_1)$ (\ref{d:theta0cl}), it follows that,
\[
\left(
\begin{array}{c}
m_{11}(-c_0,-\la_1)
\\
m_{21}(-c_0,-\la_1)
\end{array}
\right)
=
-\left(
\begin{array}{c}
m_{11}(c_0,\la_1)
\\
m_{21}(c_0,\la_1)
\end{array}
\right)
, \; \;
\left(
\begin{array}{c}
m_{12}(-c_0,-\la_1)
\\
m_{22}(-c_0,-\la_1)
\end{array}
\right)
=
\left(
\begin{array}{c}
m_{21}(c_0,\la_1)
\\
m_{22}(c_0,\la_1)
\end{array}
\right)
\]
(\ref{e:Mcl}), so that $\D(-c_0, -\la_1) = - \D(c_0, \la_1)$ (\ref{d:ABCD}). Moreover,
\[
\left(
\begin{array}{c}
r_1(-c_0, -\la_1, \mu_N)
\\
r_2(-c_0, -\la_1, \mu_N)
\end{array}
\right)
=
-
\left(
\begin{array}{c}
r_1(c_0, \la_1, \mu_N)
\\
r_2(c_0, \la_1, \mu_N)
\end{array}
\right)
\]
(\ref{e:Rclm}). Together, these imply by (\ref{d:MR}) that $\tla_2(-c_0, -\la_1) = \tla_2(c_0, \la_1)$, so that we indeed may conclude that,
\[
\begin{array}{rl}
\laho^{2, +}(-\cho(\tmu_2)) & = \eps \la_{1, {\rm h}}(-c_0) +
\eps^2 \tla_2(-c_0, \la_{1, {\rm h}}(-c_0))|\log \eps| + \O(\eps^2)
\\
& =
- \eps \la_{1, {\rm h}}(c_0) +
\eps^2 \tla_2(c_0, -\la_{1, {\rm h}}(c_0))|\log \eps| + \O(\eps^2)
\\
&
= \laho^{2, -}(\cho(\tmu_2))
\end{array}
\]
(\ref{d:la2hom-next}) (and we note that this is confirmed by the observations $\A(-c,-\la)=\A(c,\la)$, $\B(-c,-\la)=\B(c,\la)$ and
$\C(-c,-\la)=\C(c,\la)$ (\ref{d:ABCD})). Thus, we may also conclude that if Re$\la_{1, {\rm h}}(c_0) = 0$, i.e. if $\la_{1, {\rm h}}(c_0) = i \la_{1, {\rm h}}^i(c_0) \in i \RR$, then
\beq
\label{e:Relaho2pm}
{\rm Re}(\laho^{2,\pm}(\cho)) =
\eps^2\left[{\rm Re}(\A(c_0, i \la_{1, {\rm h}}^i)) \mu_N - {\rm Re}\left(\A(c_0, i \la_{1, {\rm h}}^i) \B(c_0, i \la_{1, {\rm h}}^i) - \C(c_0, i \la_{1, {\rm h}}^i)\right) \right]|\log \eps|
\eeq
(\ref{d:ABCD}), (\ref{d:la2hom-next}) at leading order. Note that if this is not the case -- i.e. if $\la_{1, {\rm h}}(c_0) \in \RR$, then $\laho^{2,+}(\cho(\tmu_2)) > 0$ and $\O(\eps)$: as we already observed, the homoclinic pulses of Theorem \ref{t:nearstabtravfrontspulses}(ii) are unstable in this case. We use Corollary \ref{c:hombifs} to determine precise conditions on the parameters so that stability may be possible, i.e. under which $\la_{1, {\rm h}}(c_0) = i \la_{1, {\rm h}}^i(c_0) \in i \RR$. First, we recall the existence condition $c_0^2 \tM_{2cc} + \tM_2 > 0$ (Theorem \ref{t:exhoms}(ii)) and that necessarily $\tN_{2 \la \la} > 0$ (\ref{d:la2hom}). Assuming that the parameters are in the correct range -- see below -- we define $c_0^{\rm m}$ and $\tmu_{\rm hom}^{\rm m}$ -- so that $C_{\rm hom}(\tmu_{\rm hom}^{\rm m}) = (c_0^{\rm m})^2$ (\ref{d:Ctmuhet}) -- as the critical value of $c_0$ at which the (leading order) real eigenvalues $\laho^{2,\pm}(c_0) > 0$ merge and become the purely imaginary pair of complex conjugate eigenvalues (at leading order), or vice versa,
\beq
\label{d:c0etmue}
c_0^{\rm m} = \sqrt{\frac{2 \sqrt{\al_+} \tM_2}{N_r - 2 \sqrt{\al_+} \tM_{2cc}}}, \; \;
\tmu_{\rm hom}^{\rm m} = \frac{\oG^+_{1c} \tM_2 N_r}{2 \sqrt{\al_+} \left(N_r - 2 \sqrt{\al_+} \tM_{2cc} \right)}
\eeq
(\ref{d:la2hom}), (\ref{d:Ctmuhet}) where $N_r = \hN_{2c \la}^2/\tN_{2 \la \la}$ has been defined in (\ref{d:Nz}) and we recall that (for stability) necessarily $N_r> 0$. By considering the 3 cases in Corollary \ref{c:hombifs} for which traveling pulses exist, we find for the pulses traveling with $c_0 > 0$ that $\la_{1, {\rm h}}(c_0) \in i\RR$ only in the following cases,
\beq
\label{e:laho2pm=Im}
\begin{array}{llclll}
{\bf 1.} & N_r < 2 \sqrt{\al_+} \tM_{2cc}: & {\bf (i)} & \tM_2 > 0, \; \tM_{2cc} > 0: & & c_0 > 0
\\
& & {\bf (ii)} & \tM_2 > 0, \; \tM_{2cc} < 0: & & 0 < c_0  < c_0^{\rm h-c}
\\
& & {\bf (iii)} & \tM_2 < 0, \; \tM_{2cc} > 0: & & c_0 > c_0^{\rm m} \; (> c_0^{\rm h-c})
\\
{\bf 2.} & N_r > 2 \sqrt{\al_+} \tM_{2cc}: & {\bf (i)} & \tM_2 > 0, \; \tM_{2cc} > 0: & & 0 < c_0 < c_0^{\rm m}
\\
& & {\bf (ii)} & \tM_2 > 0, \; \tM_{2cc} < 0: & & 0 < c_0  < c_0^{\rm m} \; (< c_0^{\rm h-c})
\\
& & {\bf (iii)} & \tM_2 < 0, \; \tM_{2cc} > 0: & & \emptyset
\end{array}
\eeq
(and note that we automatically have the same results for their counterpropagating counterparts by (\ref{d:symm})). Thus, the transition between real and purely imaginary eigenvalues (at leading order) only takes place in cases 1(iii), 2(i) and 2(ii); $\la_{1, {\rm h}}(c_0) = i \la_{1, {\rm h}}^i(c_0)$ for all allowed $c_0$-values in 1(i) and 1(ii), while $\la_{1, {\rm h}}(c_0) \in \RR$ always in 2(iii). Note that all 3 nontrivial eigenvalues associated to the stability of the traveling pulses are very close to $0$ at the transition from real to purely imaginary, i.e. that,
\beq
\label{e:lahomsatmerge}
\laho^{1, 2}(\cho(\tmu_{\rm hom}^{\rm m})), \;
\laho^{2, +}(\cho(\tmu_{\rm hom}^{\rm m})), \;
\laho^{2, -}(\cho(\tmu_{\rm hom}^{\rm m})) = \O(\eps^2 |\log \eps|),
\eeq
so that we need to perform a higher order analysis to determine leading order approximations of $\laho^{2, \pm}(\cho(\tmu_{\rm hom}^{\rm m}))$ (since $\D(c_0, \la_1) \to 0$ as $\la_1 \to 0$, as noted in the proof of Lemma \ref{l:asymptevhom-hot}) -- see section \ref{ss:Bifs}.
\\ \\
The upcoming theorem combines all asymptotic analysis of this (sub)section and establishes stability results on both types of traveling localized structures. Since the instability results are obvious from the expressions for the asymptotically small eigenvalues associated to the spectral problems -- as given in Lemmas \ref{l:asymptevhet}, \ref{l:asymptevhom} and \ref{l:asymptevhom-hot} -- we only formulate results on the (potential) stability of the traveling fronts -- Theorem \ref{t:nearstabtravfrontspulses}(i) -- and of the traveling pulses -- Theorem \ref{t:nearstabtravfrontspulses}(ii).
\begin{theorem}
\label{t:nearstabtravfrontspulses} Consider the traveling localized patterns $(\Uhe(X),\Vhe(X))$ and $(\Uho(X),\Vho(X))$ of (\ref{e:DS}) as constructed in, and under the conditions of, Theorems \ref{t:exhets}(ii) -- with $\vmu_1$ as in (\ref{d:vmutsi}) with $\tsi = 1$ -- and \ref{t:exhoms}(ii) -- with $\tmu_2$ and $\mu_N$ as defined in (\ref{d:tmu2}) and (\ref{d:muN}). Consider all eigenvalue curves $\la_j^\ast(\rho)$ associated to the family of Sturm-Liouville operators $\L_{\rho}^{\ast}(X)$ -- i.e. $\L_{\rho}(X)$ as defined in (\ref{d:Lrho}) with unperturbed homoclinic connection $v_0(X)$ replaced by (unperturbed) heteroclinic connection $\va(X)$. Assume that $(\la_0^\ast)'(0) = 1/\tau$, $(\la_0^\ast)'(\rho) <1/\tau$ for $\rho > 0$ and that $(\la_j^\ast)'(\rho) \leq 1/\tau$ for $\rho \geq 0$ for all $j > 0$. Moreover, assume that conditions (\ref{e:normhypsaddle}) and (\ref{e:stabtriv}) hold and that all non-degeneracy conditions formulated in Lemmas \ref{l:asymptevhet}, \ref{l:asymptevhom} and \ref{l:asymptevhom-hot} are satisfied.
\\
{\bf (i) Stable traveling fronts.} For $\tmu_1 \in \RR$ (and $\O(1)$ w.r.t. $\eps$) such that $\lahe^2(\vmu_1) < 0$ (\ref{d:la2het}), $(\Uhe(X),\Vhe(X))$ is spectrally stable as traveling front in (\ref{e:DS}).
\\
{\bf (ii) Stable traveling pulses.} For $(\tmu_2, \mu_N) \in \RR^2$ (and $\O(1)$ w.r.t. $\eps$) and $c_0 = c_0(\tmu_2)$ determined by $c_0 = \pm \sqrt{C_{\rm hom}(\tmu_2)}$ (\ref{d:Ctmuhet}) in the ranges given by (\ref{e:laho2pm=Im}) such that
${\rm Re}(\laho^{2,\pm}(c_0)) < 0$  ((\ref{e:Relaho2pm}) with $\la_{1, {\rm h}}^i(c_0) = {\rm Im}\la_{1, {\rm h}}(c_0)$ given by (\ref{d:la2hom})) -- which necessarily imposes $\tN_{2 \la \la}(\vmuta) > 0$ -- and $\laho^{1,2}(c_0) < 0$ (\ref{d:la12hom-next}), both patterns of the symmetric pair $(\Uho(X),\Vho(X))$ are (spectrally) stable as counterpropagating traveling pulses in (\ref{e:DS}) -- with speeds $\cho(\tmu_2) = \pm c_0 + \O(\eps)$.
\end{theorem}
\noindent
Naturally, the approximations of the asymptotically small eigenvalues obtained in Lemmas \ref{l:asymptevhet}, \ref{l:asymptevhom} and \ref{l:asymptevhom-hot} provide explicit open stability conditions on the main parameters $\tmu_1$ and $(\tmu_2, \mu_N)$. This is obvious for the heteroclinic fronts in (i), but also for the pulses of (ii): for fixed $\tmu_2$ there is a half-line of $\mu_N$ values given by
\beq
\label{d:munHopf}
\mu_N \lessgtr \frac{{\rm Re}\left(\A(c_0, i \la_{1, {\rm h}}^i) \B(c_0, i \la_{1, {\rm h}}^i) - \C(c_0, i \la_{1, {\rm h}}^i)\right)}{{\rm Re}(\A(c_0, i \la_{1, {\rm h}}^i))} \stackdef \mu_N^{\rm Hopf},
\eeq
under the additional assumption that ${\rm Re}(\A(c_0, i \la_{1, {\rm h}}^i)) \neq 0$. Thus (under this assumption), there is a large open region of $\mu_N$-values for which ${\rm Re}(\laho^{2,\pm}(c_0)) < 0$ -- for any given value of $\tmu_2$. The traveling pulse $(\Uho(X),\Vho(X))$ is stable for $\mu_N$ such that also $\Q(c_0^2, c_0^2 \mu_N) < 0$ (\ref{d:la12hom-next}). See also Corollary \ref{c:Q00}, in which the sign of $\Q(c_0^2, c_0^2 \mu_N)$ is directly and explicitly coupled to the bifurcation into traveling waves of Corollary \ref{c:hombifs} (for $|c_0| \ll 1$). Note that since a stable nearly heteroclinic homoclinic pulse $(\Uho(X),\Vho(X))$ necessarily must have ${\rm Im}(\laho^{2,\pm}(c_0)) \neq 0$, it follows that if such a pulse destabilizes as $\mu_N$ crosses through $\mu_N^{\rm Hopf}$ (\ref{d:munHopf}) -- i.e. if $\Q(c_0^2, c_0^2 \mu_N) < 0$ near $\mu_N = \mu_N^{\rm Hopf}$ -- then $(\Uho(X),\Vho(X))$ is destabilized by a Hopf bifurcation -- see Figs. \ref{f:HetHomBifsIntro}(b) and \ref{f:hethombifs-stab}(e,f).
\\ \\
{\bf Proof of Theorem \ref{t:nearstabtravfrontspulses}.} Naturally, the asymptotically small eigenvalues are controlled by Lemmas  \ref{l:asymptevhet}, \ref{l:asymptevhom} and \ref{l:asymptevhom-hot}, also in the sense that it follows from the analysis in these lemmas that there cannot be any other asymptotically small eigenvalues than the ones presented in these lemmas. Thus, if the conditions in (i) are (ii) are satisfied, the patterns $(\Uhe(X),\Vhe(X))$ and $(\Uho(X),\Vho(X))$ can only be destabilized by $\O(1)$ eigenvalues (as happens in (the proof of) Theorem \ref{t:nearstabstat}(ii) for the case $\Ma < 0$). However, this cannot happen by the conditions on $(\la_j^\ast)'(\rho)$ formulated in the statement of Theorem \ref{t:nearstabtravfrontspulses}. This can be concluded by arguments along exactly the same lines as those in the proofs of Theorem \ref{t:InstabRegHom}, Corollary \ref{c:1laast} and Theorem \ref{t:nearstabstat}.
\\ \\
The only difference is that the above conditions on $\la_j^\ast(\rho)$ replace the more explicit ones of Corollary \ref{c:1laast}(i) and \ref{t:nearstabstat}(i-b). Note that the latter conditions are too strong for the present case since they exclude the possibility that $(\la_0^\ast)'(0) = 1/\tau$, i.e $\Ma = 0$ (see the proofs of Theorem \ref{t:InstabRegHom} and \ref{t:nearstabstat}), which is the condition underlying the existence of traveling fronts and pulses -- see Theorems \ref{t:exhets} and \ref{t:exhoms}. Note also that we do not want to go deeper into the technical details and thus refrain from (deriving and) stating more explicit conditions like those in Corollary \ref{c:1laast}(i) and \ref{t:nearstabstat}(i-b) that do allow for $(\la_0^\ast)'(0) = 1/\tau$. As was already noticed in Remark \ref{r:larhosin} and immediately below the statement of Theorem \ref{t:nearstabstat}: the present conditions are weaker than those in Corollary \ref{c:1laast}(i) and \ref{t:nearstabstat}(i-b) -- but still stronger than strictly necessary. Both types of conditions prohibit the curves $\la_j^\ast(\rho)$ from intersecting the line $\la = \rho/\tau$ (for $\rho \geq 0$), and thus guarantee that there cannot be any $\O(1)$ unstable eigenvalues. \hfill $\Box$
\\ \\
Although we sofar only considered traveling fronts and pulses in this subsection, we can now zoom in on the special case $c_0 = 0$ of the preceding analysis and establish conditions under which the standing front $(\Uhe(X),\Vhe(X))$ of Corollary \ref{c:hetbifs}(i) (Fig. \ref{f:hetbifs}(a)) and the standing pulse $(\Uho(X),\Vho(X))$ of Corollary \ref{c:hombifs} (i), (ii) (Fig. \ref{f:hombifs}(a), (b)) may be spectrally stable for asymptotically small values of $|\Ma(\vmu)|$ -- the sofar missing case as was already announced at the beginnings of sections \ref{ss:S-NearHet} and \ref{sss:neartravstab}.
\begin{corollary}
\label{c:stabstandfrontpulseMsmall}
Consider the standing fronts $(\Uhe(X),\Vhe(X))$ of (\ref{e:DS}) of Corollary \ref{c:hetbifs}(i) with $|\vmu - \vmuta|= \O(\eps)$ and $\tmu_1$ as defined in (\ref{d:tmu1}) and the standing pulse $(\Uho(X),\Vho(X))$ of (\ref{e:DS}) of Corollary \ref{c:hombifs} (i), (ii) with $|\vmu - \vmuta|= \O(\eps^2 |\log \eps|)$ and $\tmu_2$ as defined in (\ref{d:tmu2}). Assume that the same conditions on the eigenvalue curves $\la_j^\ast(\rho)$ as formulated in Theorem \ref{t:nearstabtravfrontspulses} hold, that (\ref{e:normhypsaddle}) and (\ref{e:stabtriv}) hold and that all non-degeneracy conditions formulated in Lemmas \ref{l:asymptevhet}, \ref{l:asymptevhom} and \ref{l:asymptevhom-hot} are satisfied.
\\
{\bf (i) Standing fronts.} $(\Uhe(X),\Vhe(X))$ is spectrally stable as standing front pattern of (\ref{e:RDE-S1}) for $\tmu_1 \tN_{2 \la \la}(\vmuta) < 0$ and unstable for $\tmu_1 \tN_{2 \la \la}(\vmuta) > 0$: it is (de)stabilized by the transcritical bifurcation of Corollary \ref{c:hetbifs}(i).
\\
{\bf (ii-a) Stable standing pulses.} If $\tN_{2 \la \la}(\vmuta) > 0$ and $\oG^+_{1c} < 0$, then $(\Uho(X),\Vho(X))$ is spectrally stable as standing pulse pattern of (\ref{e:RDE-S1}) for $\tmu_2 < 0$ (and $\O(1)$ w.r.t. $\eps$) and it destabilizes at the bifurcation into traveling waves of Corollary \ref{c:hombifs}(i), (ii) as $\tmu_2$ increases through $0$.
\\
{\bf (ii-b) Unstable standing pulses.} The stationary pulse $(\Uho(X),\Vho(X))$ is unstable for any $\tmu_2 = \O(1)$ if $\tN_{2 \la \la}(\vmuta) < 0$ or $\oG^+_{1c} > 0$.
\end{corollary}
\noindent
We note that $\oG^+_{1c}$ can be bounded by the conditions (\ref{e:normhypsaddle}) and (\ref{e:stabtriv}) of Lemma \ref{l:sigmaess} that control the essential spectrum associated to the stability of $(\Uho(X),\Vho(X))$,
\[
\oG_{1c}^+ = 1- \tau \frac{f'(\oV_+)\oG^+_u}{\oF^+_u} = \frac{(\oF^+_u)^2 + \tau \oF^+_v \oG^+_u}{(\oF^+_u)^2} \; \;
\left\{
\begin{array}{ll}
< \frac{\oF^+_u + \tau \oG^+_v}{\oF^+_u} & (> 0)
\\
> \frac{-\tau(\oF^+_u \oG^+_v - \oF^+_v \oG^+_u)}{(\oF^+_u)^2} & (< 0)
\end{array}
\right.
\]
Nevertheless, this does not yield decisive insight in the sign of $\oG^+_{1c}$. More importantly, we note that the result Corollary \ref{c:stabstandfrontpulseMsmall}(ii) does not necessarily contradict the instability result of Theorem \ref{t:nearstabstat}(ii), since it is assumed in Theorem \ref{t:nearstabstat}(ii) that $\Ma(\vmu) \neq 0$ and $\O(1)$ w.r.t. $\eps$. Nevertheless, the stability result of Corollary \ref{c:stabstandfrontpulseMsmall}(ii) does imply that `something must happen' to the standing pulses as $|\Ma(\vmu)|$ decreases from $\O(1)$ to $\O(\eps^2 |\log \eps|)$ with $\tN_{2 \la \la}(\vmuta) > 0$ and $\oG^+_{1c} < 0$ (\ref{d:oG1+}) -- see section \ref{ss:Bifs}.
\\
\\
{\bf Proof of Corollary \ref{c:stabstandfrontpulseMsmall}.}
\\
{\bf (i) The standing fronts.} It follows by straightforward substitution of $c_0 =0$ in (the proof of) Lemma \ref{l:asymptevhet} that $\lahe^2(\vmu_1) = (\tmu_1/\hN_{2\la\la}) \eps + \O(\eps^2)$ (\ref{d:la2het}), (\ref{d:tmu1}) which immediately yields (i).
\\ \\
{\bf (ii) The standing pulses.} The leading order results on the 4 eigenvalues associated to the stability of standing pulse $(\Uho(X),\Vho(X))$ follow directly from (the proof of) Lemma \ref{l:asymptevhom}. In fact, since the pulse is stationary it follows by symmetry (\ref{d:symm}) that $\theta = -1$ (for eigenfunctions that are odd in $X$) or $\theta = 1$ (even) in (\ref{e:thetalaOeps2}) -- as in section \ref{sss:nearstatstab}. This yields,
\[
\la_1^2 = 0 \; ({\rm odd}), \; \; \la_1^2 \hN_{2 \la \la} + 2 \sqrt{\al_+} \tM_2 = 0 \; ({\rm even}),
\]
where we recall that $N_c(0) = \sqrt{\al_+} \tM_2(\vmuta) > 0$ (\ref{d:Nc}): the standing pulse $(\Uho(X),\Vho(X))$ is unstable for $\hN_{2 \la \la} < 0$. The leading order corrections to $\la_{\rm hom}^{2,+}(0)$ also follow directly by taking $c_0 = 0$ in (\ref{d:la2hom-next}). Since,
\[
\A(0,\la) = 0, \;
\B(0,\la) = \frac12 \tM_2 \oG^+_{1c}, \;
\C(0,\la) = \frac{\oG^+_{1c}}{\al_+ \hN_{2 \la \la}}, \;
\D(0,\la) = - 4 \la \sqrt{\al_+} \tM_2 \hN_{2 \la \la}
\]
(\ref{d:ABCD}) -- where we have used that $\N(z,0) = - z \sqrt{\al_+} \tM_2$ (\ref{d:Nz}) -- it follows that,
\[
\laho^{2,\pm}(0) =  \pm i \eps \sqrt{\frac{2 \sqrt{\al_+} \tM_2}{\hN_{2 \la \la}}} +
\eps^2 \frac{\oG^+_{1c}}{\al_+ \hN_{2 \la \la}} |\log \eps| + \O(\eps^2).
\]
Thus, the sign of Re$(\laho^{2,\pm}(0))$ is indeed determined by $\oG^+_{1c}$. Finally, we notice that by considering $|\tmu_2| \gg 1$ we have the following leading order approximation for $\bv(\Xh;0;\eps \tla)$,
\[
\bv(X_h; 0; \eps \tla) = \eps \frac{\tla}{\sqrt{2 \al_+ \tM_2}} \left[-\tla \hN_{2 \la \la} + \eps \tmu_2 |\log \eps| \right]
\]
(\ref{e:bvbvX-next}), which yields that  for  $|\tmu_2| \gg 1$ (but not asymptotically large w.r.t. $\eps$)
\beq
\label{e:laho12mularge}
\laho^{1,2}(0) = \eps^2 \frac{\tmu_2}{\hN_{2 \la \la}} |\log \eps|
\eeq
(at leading order) -- which concludes the proof of part (ii) of the corollary. \hfill $\Box$
\\
\\
Leading order approximation (\ref{e:laho12mularge}) of $\laho^{1,2}(0)$ for $|\tmu_2| \gg 1$ also (re-)confirms the bifurcation into traveling waves that takes place as $\tmu_2$ passes through 0 (Corollary \ref{c:hombifs}, Figs. \ref{f:hombifs}(a,b) and \ref{f:hethombifs-stab}(d,e)). Hence, it also yields information on the stability of the bifurcating traveling waves near the bifurcation, i.e. as $|c_0| \ll 1$ -- see again Fig. \ref{f:hethombifs-stab}(d,e). In other words, it follows that we have obtained explicit information on the sign of $\laho^{1,2}(c_{\rm hom}(\tmu_2))$ as $|\tmu_{\rm hom}^{\rm TW} - \tmu_2| \to 0$ (\ref{d:Ctmuhet}), and thus for $|c_0| \ll 1$ on (the sign of) $\Q(c_0^2, c_0^2 \mu_N)$ as defined in (\ref{d:la12hom-next}) of Lemma \ref{l:asymptevhom-hot}.
\begin{corollary}
\label{c:Q00}
Let all conditions as formulated in Lemma \ref{l:asymptevhom-hot} hold, let $\laho^{1,2}(c_{\rm hom}(\tmu_2))$, $\tmu_{\rm hom}^{\rm TW}$ be as defined in Lemma \ref{l:asymptevhom-hot}(ii), Corollary \ref{c:hombifs} and let $c_0$ be the leading order approximation of $c_{\rm hom}(\tmu_2)$. Then for $\tmu_2$ sufficiently close to $\tmu_{\rm hom}^{\rm TW}$,
\[
{\rm sign} \left[c_0^2 \Q(c_0^2,c_0^2 \mu_N)\right] = {\rm sign} \left[\Q(0,0)\right] = {\rm sign} \left[(\tmu_{\rm hom}^{\rm TW} - \tmu_2) \hN_{2\la \la}(\vmuta)\right].
\]
\end{corollary}
\noindent
Instead of formulating further (long, technical) corollaries on the precise impact of the stability results obtained in the present section -- i.e. equivalents of Corollaries \ref{c:hetbifs} and \ref{c:hombifs} that include explicit information on the stability of the (bifurcating) patters involved -- we chose to only graphically illustrate our stability results in Fig. \ref{f:hethombifs-stab}, where we complete the bifurcation diagrams of the cases previously considered in Figs. \ref{f:hetbifs} and \ref{f:hombifs} by the obtained stability insights (see also Fig. \ref{f:HetHomBifsIntro} in the Introduction).
\begin{remark}
\label{r:hN2cl2tM2cc}
\rm
The fact that $\hN_{2c\la}(\vmu_{\rm t}^\ast) = 2\tM_{2cc}(\vmu_{\rm t}^\ast)$ (\ref{e:hN2tM}) can in principle be derived directly from their definitions, (\ref{d:tM2c2}) and (\ref{d:hN2s}). A priori, one may think that such a simple relation cannot be correct since the expression for $\hN_{2c\la}(\vmu_{\rm t}^\ast)$ (\ref{d:hN2s}) involves through (\ref{d:hu1cl}), (\ref{d:hI2s}), (\ref{d:hJ2s}) integrals over the function $\hv_{1c}$ -- the unique bounded solution of $\La v =  \hJ_{1c}$ (\ref{e:OepsEigPb}), (\ref{e:bv1lc}) -- while this expression does not show up in (\ref{d:tM2c2}). However, it can be checked that $\hv_{1c}$ only appears in a very specific configuration in the (integral) expression for $\tN_{2c\la}$ and as a consequence can be removed (using (\ref{d:tv1}) and the fact that the operator $\La$ is selfadjoint),
\[
\int_{-\infty}^{\infty} \left(1-\tau \frac{\fa' \Ga_u}{\Fa_u}\right) \vax \hv_{1c} \, dX =
- \int_{-\infty}^{\infty} (\La \tv_1) \hv_{1c} \, dX =
- \int_{-\infty}^{\infty} \tv_1 (\La \hv_{1c}) \, dX =
- \int_{-\infty}^{\infty} \hJ_{1c} \tv_1 \, dX
\]
with $\hJ_{1c}$ as in (\ref{d:hJ1c}). Nevertheless, a direct transformation of $\hN_{2c\la}(\vmu_{\rm t}^\ast)$ into $2\tM_{2cc}(\vmu_{\rm t}^\ast)$ is quite a technical enterprise -- we refrain from going into the details.
\end{remark}

\subsection{Heteroclinic interfaces and homoclinic stripes}
\label{ss:IntStr}

Finally, we consider patterns in 2-dimensional $(X,Y)$-space that vary in $X$ and are trivially extended in the $Y$-direction. Naturally, we are only interested in the stability of interfaces and/or stripes that can be stable as solutions of the 1-dimensional version (\ref{e:RDE-S1}) of (\ref{e:RDE})/(\ref{e:RDE-S}), i.e. the (nearly) heteroclinic patterns $(U(X,Y,t), V(X,Y,t)) = (\Uhe(X),\Vhe(X))$ and $(U(X,Y,t), V(X,Y,t)) = (\Uho(X),\Vho(X))$ that may be stable as functions of only $X$ as established by Theorem \ref{t:nearstabstat}(i-b), Theorem \ref{t:nearstabtravfrontspulses} and Corollary \ref{c:stabstandfrontpulseMsmall}.
\\ \\
As in the proofs of Theorems \ref{t:nearstabstat} and \ref{t:nearstabtravfrontspulses}, we need to distinguish between asymptotically small and $\O(1)$ spectrum. Moreover, as before the $\O(1)$ spectrum is determined by the family of Sturm-Liouville operators $\L_\rho^\ast(X)$,
\[
\L_\rho^\ast(X) \bv_0 = \left[\La(X) - \frac{\rho f'(\va(X)) G_u(f(\va(X)),\va(X))}{\rho - F_u(f(\va(X),\va(X))}\right] \bv_0 = (\la + L^2) \bv_0
\]
(\ref{e:O1EigPb}), (\ref{d:La}), (\ref{d:Lrho}). Thus, for given $L$ fixed, the $\O(1)$ eigenvalues are determined by the intersections $\la^\ast_j(\rho) - L^2 = \rho/\tau$, where we recall that $\la^\ast_j(\rho)$ are the eigenvalues of $\L_\rho^\ast(X)$. However, in the stability results of Theorems \ref{t:nearstabstat}(i-b) and \ref{t:nearstabtravfrontspulses} it is assumed that all curves $\la^\ast_j(\rho)$ are strictly below the line $\la = \rho/\tau$ (for $\rho \geq 0$) -- except for one point at the origin (since $\la^\ast_0(0) = 0$). Hence, under the (stability) conditions of Theorems \ref{t:nearstabstat}(i-b) and \ref{t:nearstabtravfrontspulses}, the downward shifted curves $\la^\ast_j(\rho) - L^2$ cannot intersect the line $\la = \rho/\tau$ (for $\rho \geq 0$). As a consequence we may immediately conclude that there cannot be (`new') unstable $\O(1)$ spectrum: the interface/stripe patterns $(U(X,Y,t), V(X,Y,t)) = (\Uhe(X),\Vhe(X))/(\Uho(X),\Vho(X))$ of Theorems \ref{t:nearstabstat}(i-b) and \ref{t:nearstabtravfrontspulses} and Corollary \ref{c:stabstandfrontpulseMsmall} can only be destabilized by long wavelength perturbations in $Y$, i.e. by $|L| \ll 1$.
\\ \\
In the upcoming analysis, a central role will be played by the unique solution $\tv_{L}$ of
\beq
\label{d:tvL}
\La(X) \tv_{L} = \vax,
\eeq
that satisfies $\lim_{X \to - \infty} \tv_{L} = 0$ and $\tv_{L}(0) = 0$. It follows by Lemma \ref{l:Lah} that the leading order behavior of $\tv_{L}(X)$ for $X \gg 1$ is determined by $\|\vax\|_2^2 = \int_{-\infty}^\infty \vax^2 (X; \vmu) \, dX$.

\begin{theorem}
\label{t:stabstatinterfaces}
Let $(U(X,Y,t), V(X,Y,t)) = (\Uhe(X),\Vhe(X))$ be a heteroclinic interface as established by Theorem \ref{t:exhets} and let all conditions hold as formulated in Theorems \ref{t:exhets}, \ref{t:nearstabstat}(i-b), \ref{t:nearstabtravfrontspulses}(i) and Corollary \ref{c:stabstandfrontpulseMsmall}(ii).
\\
{\bf (i) $\Ma(\vmu) = \O(1)$: stationary interfaces.} Introduce $L_1$ by $L = \sqrt{\eps} L_1$. For $L_1 = \O(1)$ given, there is 1 asymptotically small eigenvalue,
\[
\la_{\rm het}(L_1; \vmu) = - \eps \frac{\|\vax\|_2^2}{\Ma(\vmu)} L_1^2 + \O(\eps^2).
\]
If the (stationary) front $(\Uhe(X),\Vhe(X))$ is stable on $\RR^1$, then the interface $(\Uhe(X),\Vhe(X))$ is stable as solution of (\ref{e:RDE})/(\ref{e:RDE-S}) on $\RR^2$.
\\
{\bf (ii) $\Ma(\vmu) = \O(\eps)$: standing and traveling interfaces.} Let $\vmu = \vmuta + \eps \vmu_1$ (cf. (\ref{d:vmutsi}) with $\tsi = 1$) and introduce $L_2$ by $L = \eps L_2$. For $L_2$ given, there are 2 asymptotically small eigenvalues $\la_{\rm het}^{1,2}(L_2; \vmu_1) = \eps \la_1 + \O(\eps^2 |\log \eps|)$ and $\la_1$ determined as solution of
\beq
\label{e:eqla1L2-interface}
\la_1^2 \hN_{2\la\la} (\vmu_{\rm t}^\ast) - \la_1 \left[\vec{\nabla} \Ma(\vmu_{\rm t}^\ast) \cdot \vmu_1 + c_0 \hN_{2c\la} (\vmu_{\rm t}^\ast) \right] -  \|\vax\|_2^2 L_2^2 = 0
\eeq
(\ref{d:la2het}), (\ref{e:la1stravfront}) -- with $c_0 = c_0(\tmu_1)$ the leading order approximation of the speed $\che(\vmu_1)$ at which the interface travels (and $c_0 =0$ for the standing interfaces). If $(\Uhe(X),\Vhe(X))$ is stable as solution of (\ref{e:RDE-S1}) for $X \in \RR$, then the interface is stable as solution of (\ref{e:RDE})/(\ref{e:RDE-S}) on $\RR^2$ if additionally $\hN_{2 \la \la} < 0$; it is unstable for $\hN_{2 \la \la} > 0$.
\end{theorem}
\noindent
Note that it follows from (\ref{e:eqla1L2-interface}) that for $|L_2|$ small (but still $\O(1)$ w.r.t. $\eps$),
\beq
\label{e:lahet2L2small}
\la_{\rm het}^1(L_2; \vmu_1) = \eps \left[- \frac{\|\vax(\vmuta)\|_2^2}{\la_{\rm het}^{2}(0; \vmu_1) \hN_{2 \la \la}(\vmuta)} L_2^2 + \O(L_2^4)\right] + \O(\eps^2|\log \eps|)
\eeq
(\ref{d:la2het}), where $\la_{\rm het}^{2}(0; \vmu_1) < \la_{\rm het}^{1}(0; \vmu_1) = 0$ since $(\Uhe(X),\Vhe(X))$ is assumed to be stable on $\RR$. Thus, the transition in (ii) from stable to unstable as $\hN_{2 \la \la}(\vmuta)$ changes sign a priori has the character of a sideband instability. However, the situation is more complex: $(\Uhe(X),\Vhe(X))$ needs to remain stable on $\RR$ as $\hN_{2 \la \la}(\vmuta)$ passes through 0. Hence, by (\ref{d:la2het}), $\vec{\nabla} \Ma(\vmu_{\rm t}^\ast) \cdot \vmu_1 + c_0(\vmu_1) \hN_{2c\la}(\vmu_{\rm t}^\ast)$ needs to change sign simultaneously with $\hN_{2 \la \la}(\vmuta)$ -- which indeed changes the sign of the pre-factor of $L_2^2$ in (\ref{e:lahet2L2small}). However, a sign change in the denominator causes a singularity: the transition from stable to unstable interfaces cannot `automatically' be seen as a sideband mechanism, one first needs to analyse the impact of (making) $|\hN_{2 \la \la}|$ asymptotically small -- see section \ref{ss:Bifs}.
\\ \\
{\bf Proof of Theorem \ref{t:stabstatinterfaces}.} We first consider case (i), i.e. we set $L = \sqrt{\eps} L_1$ and conclude that a term `$-L_1^2 \bv_0$' needs to be added to the right hand side of the second equation in the $\O(\eps)$ spectral problem (\ref{e:ODE-lin-1}). Thus, the $L_1$-dependent version of (\ref{e:OepsEigPb}) reads
\[
\La \bv_1 = \la_1 \left[1-\tau \frac{\fa' \Ga_u}{\Fa_u}\right]\vax + L_1^2 \vax
\]
where we have eliminated the term $c_0 \hJ_{1c}$ since we consider stationary fronts in (i). Thus, it follows that $\bv_1 = -\la_1 \tv_1 + L_1^2 \tv_L$ (\ref{e:bv1lc}), (\ref{d:tvL}). Since necessarily $\lim_{X \to \infty} \bv_1(X) = 0$ in this heteroclinic case, it follows by Lemma \ref{l:Lah} that indeed $\la_1 \Ma + \|\vax\|_2^2 L_1^2 = 0$. The statement of Theorem \ref{t:stabstatinterfaces}(i) follows by noticing that $\Ma > 0$ by the assumption that $(\Uhe(X),\Vhe(X))$ is stable on $\RR$ (Theorem \ref{t:nearstabstat}(i-a)).
\\ \\
The proof of part (ii) follows along exactly the same lines. We notice that including the $Y$-direction with $L = \eps L_2$ introduces an additional inhomogeneous term `$+ L_2^2 \vax$' to equation (\ref{e:Oeps2EigPb}) for $\bv_2$. As a consequence, (\ref{e:la1stravfront}) in the proof of Lemma \ref{l:asymptevhet} indeed becomes (\ref{e:eqla1L2-interface}) and the statement of the Theorem follows. \hfill $\Box$
\\ \\
The additional condition $\hN_{2 \la \la} < 0$ in Theorem \ref{t:nearstabtravfrontspulses} for the stability of traveling interfaces with $\Ma(\vmu)$ asymptotically small strongly suggests that it is unlikely that the -- traveling or standing -- stable homoclinic pulse solutions on $\RR$ of Theorem  \ref{t:nearstabtravfrontspulses}(ii) and Corollary \ref{c:stabstandfrontpulseMsmall}(ii) may extend to stable homoclinic stripes on $\RR^2$: both stability results require $\hN_{2 \la \la} > 0$. This is indeed the case.
\begin{theorem}
\label{t:unstabtravstripes}
Let all conditions of Theorem \ref{t:nearstabtravfrontspulses}(ii) and Corollary \ref{c:stabstandfrontpulseMsmall}(ii-a) hold so that the traveling or standing pulses $(\Uho(X),\Vho(X))$ are spectrally stable as 1-dimensional pattern in (\ref{e:RDE-S1}). The pattern $(U(X,Y,t), V(X,Y,t)) = (\Uho(X),\Vho(X))$ is unstable as traveling or standing homoclinic stripe solution of (\ref{e:RDE})/(\ref{e:RDE-S}) on $\RR^2$.
\end{theorem}
\noindent
{\bf Proof of Theorem \ref{t:unstabtravstripes}.} We first consider the (general) case that includes $c_0$. As in Theorem \ref{t:stabstatinterfaces} we introduce $L_2$ by setting $L = \eps L_2$ and as in the proof of Theorem \ref{t:stabstatinterfaces} this yields the appearance of new terms -- with pre-factor $L_2^2$ -- in the asymptotic analysis. More specifically, it can be checked in a straightforward fashion that for given $L_2 \neq 0$, (\ref{e:thetalaOeps2}) in the proof of Lemma \ref{l:asymptevhom} generalizes to,
\[
\begin{array}{rcl}
-\la_1^2 \hN_{2 \la \la} + c_0 \la_1 \hN_{2c \la} + \|\vax\|_2^2 L_2^2
& = & \theta \left[-\la_1^2 \hN_{2 \la \la} - c_0 \la_1 \hN_{2c \la} + \|\vax\|_2^2 L_2^2 \right]
\\[2mm]
- 2 N_c -\la_1^2 \hN_{2 \la \la} + c_0 \la_1 \hN_{2c \la} + \|\vax\|_2^2 L_2^2
& = & \theta \left[2 N_c + \la_1^2 \hN_{2 \la \la} + c_0 \la_1 \hN_{2c \la} - \|\vax\|_2^2 L_2^2 \right]
\end{array}
\]
which yields the following $L_2$-dependent version of (\ref{e:quartla1}),
\beq
\label{e:quartla1L2}
\hN_{2 \la \la}^2 \la_1^4 + \left[2 N_c \hN_{2 \la \la} - c^2_0 \hN_{2c \la}^2 - 2 \hN_{2 \la \la} \|\vax\|_2^2 L_2^2 \right] \la_1^2
- \left[ 2 N_c \|\vax\|_2^2 L_2^2 - \|\vax\|_2^4 L_2^4 \right] = 0
\eeq
where we note that since
\[
(2 N_c \hN_{2 \la \la} - c^2_0 \hN_{2c \la}^2)^2 + 4 c_0^2 \hN^2_{2 c \la} \hN_{2 \la \la} L_2^2 > 0
\]
for all $L_2$, this parabolic equation in $\Lambda_1 = \la_1^2$ only has real solutions. Moreover, since
$2 N_c \|\vax\|_2^2 L_2^2 -  \|\vax\|_2^4 L_2^4 > 0$ for $L_2^2 < 2N_c/\|\vax\|_2^2$ -- recall that $N_c > 0$ (\ref{d:Nc}) -- it follows that one of the solutions of (\ref{e:quartla1L2}), $\Lambda_1 = \Lambda_1^{+}(L_2)$, must be positive. Hence, there is a solution $\la^{+,+}_1(L_2)$ of (\ref{e:quartla1L2}) with  $\la^{+,+}_1(L_2) = + \sqrt{\Lambda_1^{+}(L_2)} > 0$ (for $L_2^2 < 2N_c/\|\vax\|_2^2$): $(\Uho(X),\Vho(X))$ is destabilized as traveling stripe solution of (\ref{e:RDE})/(\ref{e:RDE-S}) by perturbations with $|L_2|$ sufficiently small (but still $\O(1)$ w.r.t. $\eps$).
\\ \\
As in (the proof of) Corollary \ref{c:stabstandfrontpulseMsmall}(ii), the instability result for the standing stripes is embedded in the general case by setting $c_0 = 0$. In fact, we notice that parabolic equation (\ref{e:quartla1L2}) can be straightforwardly factored in the case of standing stripes
\[
\hN_{2 \la \la}^2 \left(\la_1^2 - \frac{\|\vax\|_2^2}{\hN_{2 \la \la}} L_2^2\right)\left(\la_1^2 + 2 \frac{\sqrt{\al_+} \tM_2}{\hN_{2 \la \la}} - \frac{\|\vax\|_2^2}{\hN_{2 \la \la}} L_2^2\right) = 0
\]
(\ref{d:Nc}) which directly yields instability. \hfill $\Box$
\\ \\
Finally, we note that we did not need the higher order $\O(\eps^2|\log \eps|)$ correction results of Lemma \ref{l:asymptevhom-hot} to deduce the instability of the stripes $(\Uho(X),\Vho(X))$ on $\RR^2$, while we needed this lemma to establish the stability of the traveling and standing pulses on $\RR$ (Theorem \ref{t:nearstabtravfrontspulses} and Corollary \ref{c:stabstandfrontpulseMsmall}(ii)): the $L_2$-induced instability sets in with spectrum of $\O(\eps)$ magnitude, the $\O(\eps^2|\log \eps|)$ effects cannot counteract this.

\section{Discussion}
\label{s:Disc}
In the preceding analysis, we have established a number of fundamental results on the existence and stability of stationary and (uniformly) traveling localized patterns in the general class of systems (\ref{e:RDE}) that only vary in the slow spatial coordinate $X$, i.e. that do not exhibit fast jumps or spikes. We have shown that slow pulses that correspond to homoclinic orbits in the reduced slow flow on the slow manifold do persist -- either as stationary or as traveling pattern -- but must be unstable (Theorems \ref{t:E-Pers} and \ref{t:InstabRegHom}). However, pulses and fronts that merge with a heteroclinic cycle of the reduced slow flow in the limit $\eps \to 0$ may be stable. Our results on the existence and stability of these fronts and pulses as patterns in $\RR^1$ are illustrated by the bifurcation diagrams of Fig. \ref{f:HetHomBifsIntro} and its more extended version Fig. \ref{f:hethombifs-stab} that are based on the existence results of Theorems \ref{t:exhets} and \ref{t:exhoms} and the (spectral) stability insights of Theorems \ref{t:nearstabstat}, \ref{t:nearstabtravfrontspulses} and Corollary \ref{c:stabstandfrontpulseMsmall} (and several additional results formulated in terms of lemmas and corollaries). The (in)stability of their extensions into interfaces and (homoclinic) stripes for $(x,y) \in \RR^2$ is established in Theorems \ref{t:stabstatinterfaces} and \ref{t:unstabtravstripes}. Although quite extensive, the present analysis naturally also generates further questions and/or further research topics. In this section, we mention some of these, it is divided in 3 subsections: one on observations concerning some of the choices made, one on further bifurcations and one on possible projects for future work.

\subsection{Observations}
\label{ss:Obs}
{\bf The impact of $\tau$.} On the one hand, the parameter $\tau$ as introduced in (\ref{e:RDE}) can be seen as the $(m+1)$-th component of the $(m+1)$-dimensional family of parameters $(\tau, \vec{\mu})$ of (\ref{e:RDE}). However, $\tau$ does have a direct interpretation: it measures the relative rate of evolution (in time) of the 2 components $U(x,y,t)$ and $V(x,y,t)$ of (\ref{e:RDE}). Moreover, it plays a decisive role in the (non-)appearance of traveling patterns. This follows from the simple observation that the conditions under which traveling patterns bifurcate -- $M_{\rm hom}(\vec{\mu}_{\rm t}) = 0$ ((\ref{d:M0vmuW0}) in Theorem \ref{t:E-Pers}) for potential $\W_0$ with wells of unequal depth and $\Ma(\vmuta) = 0$ (\ref{d:Ma}) in the wells of equal depth case (Theorems \ref{t:exhets} and \ref{t:exhoms}) -- clearly cannot be satisfied if $\tau$ is `too small'. In other words,
{\it if the rate of evolution (in time) of component $U(x,y,t)$ compared to that of $V(x,y,t)$ is sufficiently fast, then there are no slow traveling localized patterns} (of the type analyzed in this paper).
\\ \\
In the stability analysis (of the stationary patterns), the case `$\tau$ is small' has a similar simplifying interpretation (see especially Fig. \ref{f:Lambdajrho}): for small $\tau$ the line $\{\la = \rho/\tau\}$ approaches verticality, so that the nonlinear eigenvalue problem (\ref{e:O1EigPb}) approaches the $\rho = 0$ scalar Sturm-Liouville problem associated to (\ref{e:RedPDE}): the impact of taking $\tau$ (too) small is similar to having a(n almost) vertical slow manifold $\M_\eps$ (i.e. to having $|f'(v)|$ (too) small).
\\ \\
{\bf Traveling patterns with $|c| \gg 1$.} At the beginning of the paper, we chose the magnitude of $c$ such that it appears as leading order term in the (fast)$u$-equation (cf. (\ref{e:ODE})): the most natural choice to enable orbits to jump from one slow manifold to another. However, our analysis focuses on orbits that do not jump, moreover there are large classes of systems with only one (normally hyperbolic) slow manifold $\M_\eps$ that may or may not have a `return mechanism' (see section \ref{ss:Projects}): this point of view makes our choice of the magnitude of $c$ less natural. In fact, the existence results indicate that the branches of traveling patterns may extend beyond $\O(1)$ values of $c$ (see especially Figs. \ref{f:hetbifs} and \ref{f:hombifs}(a,c)). However, even only allowing $c$ to be logarithmically large, i.e. choosing $c_0 = \O(|\log \eps|^\si)$ for some $\si > 0$, has a nontrivial impact. Without going into the details, this can be illustrated in the existence problem for nearly heteroclinic pulses.
\\ \\
In the proof of Theorem \ref{t:exhoms}, the existence of these pulses is established by closing the gap $\Delta W_{\rm h}(\tmu_2)$ (\ref{e:gapOeps2inhomasymm}). For $|c| \gg 1$, the leading order contribution to the next order corrections of (\ref{e:gapOeps2inhomasymm}), i.e. of its $\O(\eps^2)$ term, will come from $v_3(X)$ through $\tv_{3ccc}$ (\ref{d:expv34}), i.e. it will be of the form $\tM_{3ccc}(\vmuta) c_0^3$ (where $\tM_{3ccc}(\vmuta)$ determines the leading order growth of $\tv_{3ccc}$ (cf. Lemma \ref{l:Lah})). Thus, for $c_0 = \tilde{c}_0 |\log \eps|^\si$, the condition $\Delta W_{\rm h}(\tmu_2) = 0$ becomes, at leading order,
\[
\tmu_2 \frac{|\log \eps|^{1-\si}}{\tM_\al |\tilde{c}_0|} - \frac{\oG_{1c}^+}{\al_+^{3/2}} \tM_\al |\tilde{c}_0| |\log \eps|^{1+\si} + \tM_{3ccc} \, \tilde{c}_0^2 |\log \eps|^{2} = 0
\]
with $\tM_\al = \sqrt{\al_+ \tM_{2cc}/2}$ (assuming $\tM_{2cc} > 0$, Corollary \ref{c:hombifs}, Fig. \ref{f:hombifs}(a,c)). Clearly, a transition -- and thus a potential bifurcation -- occurs at $\si = 1$. By scaling $\tmu_2 = |\log \eps|^{2} \tilde{\mu}_{2,2}$ -- so that $\|\vmu - \vmuta \| = \O(\eps^2 |\log \eps|^3)$ (\ref{d:tmu2}) -- we obtain as leading order equation for $\tilde{c}_0 = \tilde{c}_0(\tilde{\mu}_{2,2})$,
\[
\tM_\al  \frac{\oG_{1c}^+}{\al_+^{3/2}} |\tilde{c}_0|^2 - \tM_{3ccc} \, |\tilde{c}_0|^3 =  \frac{\tmu_{2,2}}{\tM_\al},
\]
which shows that a saddle node bifurcation occurs at $\tmu_{2,2}^{\rm SN} = (2\tM_{2cc}^2 (\oG_{1c}^+)^3)/(27 \al_+^{5/2} \tM_{3ccc}^2)$ if sign$(\oG_{1c}^+) =$ sign$(\tM_{3ccc})$. For instance, for the cases considered in Figs. \ref{f:hetbifs} and \ref{f:hethombifs-stab} there are 2 solutions $|\tilde{c}_0| > 0$, and thus 4 speeds $\tilde{c}_0 \in \RR$, for $\tmu_{2,2}^{\rm SN} < \tmu_{2,2} < 0$ and none for $\tmu_{2,2} < \tmu_{2,2}^{\rm SN}$: the branches of traveling pulses in Figs. \ref{f:hombifs}(a,c)/\ref{f:hethombifs-stab}(d,f) bend back for $\vmu$ at an $\O(\eps^2 |\log \eps|^3)$ distance from $\vmuta$ (if $\tM_{3ccc} < 0$). Thus, the unstable branch of  Fig. \ref{f:hethombifs-stab}(d) could be stabilized by the saddle node bifurcation.
\\ \\
Naturally, similar bifurcations may occur for the traveling pulses and fronts of Theorems \ref{t:E-Pers} and \ref{t:exhoms}, and at other ranges/magnitudes of $|c_0| \gg 1$. We refrain from going further into the details here (see also section \ref{ss:Projects}).

\subsection{Bifurcations}
\label{ss:Bifs}
{\bf The transition from real to complex eigenvalues for nearly heteroclinic pulses.} At $\tmu_2 = \tmu_{\rm hom}^{\rm m}$ (\ref{d:c0etmue}), the pair of eigenvalues $\la_{\rm hom}^{2, \pm}$ associated to the stability of a nearly heteroclinic pulse merges (at leading order, (\ref{d:la2hom})) and passes from being real to complex, or vice versa, so that the spectral stability problem has 3 nontrivial eigenvalues of $\O(\eps^2 |\log \eps|)$ (\ref{e:lahomsatmerge}). The Hopf bifurcation at $\mu_N = \mu_N^{\rm Hopf}$ (\ref{d:munHopf}) that (de)stabilizes the pulse can be tuned (by $\mu_N$) independently from $\tmu_{\rm hom}^{\rm m}$ -- since $\tmu_2$ and $\mu_N$ can be varied independent from each other (see definition (\ref{d:muN}) in Lemma \ref{l:asymptevhom-hot}). In Figs. \ref{f:HetHomBifsIntro}(b) and \ref{f:hethombifs-stab}(e,f), it is (implicitly) assumed that at $\tmu_2 = \tmu_{\rm hom}^{\rm m}$, $\mu_N$ is such that the $\O(\eps^2 |\log \eps|)$ term (\ref{d:la2hom-next}) of the pair $\la_{\rm hom}^{2, \pm}$ has positive real part (\ref{d:munHopf}), so that the (de)stabilization of the pulse is driven by the Hopf bifurcation (naturally, under the additional assumption that $\la_{\rm hom}^{1,2} < 0$ (\ref{d:la12hom-next})). However, this is only the case if indeed the real, $\O(\eps^2 |\log \eps|)$, part of the pair $\la_{\rm hom}^{2, \pm}$ is positive at $\tmu_{\rm hom}^{\rm m}$. If this is not the case, i.e. if $\mu_N$ is on the other side of $\mu_N^{\rm Hopf}$ compared to the cases sketched in Figs. \ref{f:HetHomBifsIntro}(b) and \ref{f:hethombifs-stab}(e,f), then the (de)stabilization of the nearly heteroclinic pulse sets in at -- or better: $\O(\eps |\log \eps|)$ close to -- $\tmu_{\rm hom}^{\rm m}$. In other words, in this case the nature of the (de)stabilization of the pulse is determined by the `dynamics' of the eigenvalues $\la_{\rm hom}^{2, \pm}$ as they merge $\O(\eps^2 |\log \eps|)$ close to $\la = 0$.
\\ \\
Although this (de)stabilization may seem to have the character of a Bogdanov-Takens bifurcation, it is not of co-dimension 2 (due to symmetry (\ref{d:symm})) . The exact nature of the bifurcation -- and thus of the emerging new localized patterns in (\ref{e:RDE}) generated by it -- can be studied by a higher order analysis (which we refrain from going into here).
\\ \\
{\bf Zooming in on $\hN_{2 \la \la}(\vmuta)$ passing through $0$.} Next to $\Ma(\vmu)$ (\ref{d:Ma}), $\hN_{2 \la \la}(\vmu)$ (\ref{d:hN2ll}) perhaps is the most important Melnikov-type expression/function in the preceding analysis: for $\Ma(\vmuta)$ asymptotically small, the sign of $\hN_{2 \la \la}(\vmuta)$ plays a decisive role in determining whether fronts or pulses may be stable -- Lemmas \ref{l:asymptevhet}, \ref{l:asymptevhom}, \ref{l:asymptevhom-hot}, Theorem \ref{t:nearstabtravfrontspulses} and Corollary \ref{c:stabstandfrontpulseMsmall} -- and it determines whether long wavelength perturbations in the $Y$ direction may destabilize an interface or (homoclinic) stripe associated to a 1-dimensional stable front or pulse -- Theorems \ref{t:stabstatinterfaces} and \ref{t:unstabtravstripes}.
\\ \\
Unfolding the passage of $\hN_{2 \la \la}(\vmuta)$ through $0$ is especially relevant for understanding the apparent contradiction between Theorem \ref{t:nearstabstat}(ii) -- that states that stationary pulses are unstable for $\Ma(\vmu) \neq 0$ and $\O(1)$ -- and Corollary \ref{c:stabstandfrontpulseMsmall}(iia) -- by which stationary pulses may be stable for $\Ma(\vmu) = \O(\eps^2 |\log \eps|)$ and $\tN_{2 \la \la}(\vmuta) > 0$, $\oG^+_{1c} < 0$. Clearly, these statements do not `overlap' -- and thus not contradict each other -- but the combination of Theorem \ref{t:nearstabstat}(ii) and Corollary \ref{c:stabstandfrontpulseMsmall}(iia) does imply that if $\tN_{2 \la \la}(\vmuta) > 0$ (an) additional bifurcation(s) must take place as $\|\vmu - \vmuta\|$ increases from $\O(\eps^2 |\log \eps|)$ to $\O(1)$, possibly similar -- and related -- to the new bifurcations induced by taking $|c| \gg 1$ (cf. section \ref{ss:Obs}). However, these bifurcations are not expected to occur for $\tN_{2 \la \la}(\vmuta) < 0$. Preliminary analysis indicates that considering the case of $\hN_{2 \la \la}(\vmuta)$ asymptotically small indeed will shed light on this.
\\ \\
A similar observation can be made about the transition from stable interfaces -- that must have $\tN_{2 \la \la}(\vmuta) < 0$ -- to unstable interfaces -- with $\tN_{2 \la \la}(\vmuta) > 0$ -- in Theorem \ref{t:stabstatinterfaces}. Since $\tN_{2 \la \la}(\vmuta)$ appears in the denominator of the critical eigenvalue $\la_{\rm het}^1(L_2)$ associated to long wavelength perturbations (\ref{e:lahet2L2small}), this transition cannot 'automatically' be of sideband type. To understand the exact nature of the transition from stable to unstable stripes -- and perhaps to even recover a small region of stable homoclinic stripes in parameter space (Theorem \ref{t:unstabtravstripes}) -- it is thus necessary to zoom in on the passage of $\hN_{2 \la \la}(\vmuta)$ through $0$.

\subsection{Projects}
\label{ss:Projects}

{\bf Systems with only one normally hyperbolic slow manifold.} Naturally, our focus on slow traveling waves to (\ref{e:RDE}) is especially relevant for systems that only have one normally hyperbolic slow manifold and no `return mechanism'. The latter may be because the fast reduced flow -- given by $u_{\xi \xi} + c \tau u_{\xi} + F(u, v_0)$, $v_0 \in \RR$ (cf. (\ref{e:ODE})) -- only has one critical point (associated to the slow manifold $\M_\eps(c)$), as for the family of models introduced in Remark \ref{r:larhosin}. However, it is also possible that there may be fast reduced homoclinic orbits `attached' to $\M_\eps(c)$ -- that thus may `return' orbits that jump away from $\M_\eps(c)$ (see below) -- but that these jump are irrelevant or impossible from the modelling point of view (as for the savanna grass/woodland ecosystem model (\ref{e:RDE-exvLangevelde})). In these cases, the only simple traveling patterns to (\ref{e:RDE}) -- i.e. solutions of the type $(U,x,y,t), V(x,y,t)) = (u(x-ct),v(x-ct))$ -- are those that remain on $\M_\eps(c)$ for all $\xi$/$X=\eps \xi$. However, as was already noticed in section \ref{ss:Obs}, for these models our choice $c = \O(1)$ looses its foundation.
\\ \\
Originally, the stability analysis of localized pulses -- or spikes -- in 2-component singularly-perturbed reaction-diffusion equations on $\RR$ centered around slow-fast-slow orbits in Gray-Scott and Gierer-Meinhardt models that -- in the terminology of the present paper -- followed the unstable manifold of a (linear) saddle point on a vertical slow manifold $\M_\eps \equiv \M_0$ (thus with $f'(u) \equiv 0$), jumped through the fast field following a fast reduced homoclinic  orbit and `touched down' again on $\M_\eps$ at the stable manifold of the same saddle \cite{DGK98,DGK01,IWW01,WW03}. The fast jump -- the spike -- has a positive `fast reduced eigenvalue' associated to it, which a priori suggests that the pattern cannot be stable. However, this potential instability can be removed by the interaction between the fast and slow components of the full pulse (this is called `the resolution of the NLEP paradox' in \cite{DGK01}). To build a bridge between the present analysis and the existing literature on slow-fast-slow homoclinic pulses in singularly perturbed reaction-diffusion equations -- see \cite{CW09,DV15,KWW09,SD17,Ward18} and the references therein -- it would be very natural to study slow-fast-slow orbits in (\ref{e:RDE}) that take off from, and subsequently touch down on, $\M_\eps$ meanwhile following a homoclinic orbit of the fast reduced flow. Therefore, it is a natural future project to study localized patterns in (\ref{e:RDE}) that are built from the slow patterns studied here interspersed with a fast homoclinic spike. Such a study could be set up by introducing $\tG_{\rm f}(U,V)$ and $\mu_{m+1} \in \RR$ by $G(U,V) - G(f(V),V) = \mu_{m+1}\tG_{\rm f}(U,V)$ so that (\ref{e:RDE}) transforms into,
\beq
\label{e:RDE-return}
\left\{	
\begin{array}{rcrcrcl}
\tau U_t &=& \Delta U & + & F(U,V;\vmu) & &\\
V_t &=& \frac{1}{\varepsilon^2} \Delta V & + & G(f(V),V;\vmu) & + &\mu_{m+1} \tG_{\rm f}(U,V;\vmu)
\end{array}
\right.
\eeq
with $\tG_{\rm f}(f(V),V) = 0$: $\mu_{m+1}\tG_{\rm f}(U,V)$ represents the impact of the dynamics away from the slow manifold $\M_\eps$. Written in this form, the (generalized) Gierer-Meinhardt system -- that has $\M_\eps \equiv \M_0 = \{u = 0, p = 0\}$ -- reads,
\beq
\label{e:RDE-genGM}
\left\{	
\begin{array}{rcrcrcl}
U_t &=& \Delta U & - & U & + & U^{\mu_2} V^{\mu_3}\\
V_t &=& \frac{1}{\varepsilon^2} \Delta V & - & \mu_1 V & + & \frac{1}{\eps} U^{\mu_4} V^{\mu_5}
\end{array}
\right.
\eeq 
\cite{DGK01,WW03}. Thus, $\mu_{m+1} = \O(1/\eps)$ in the Gierer-Meinhardt system: in the setting of Gierer-Meinhardt and Gray-Scott systems, the fast component $U$ needs to have a sufficiently strong impact in the slow $V$-equation in order to have stable pulses, or equivalently, the flow on $\M_\eps$ needs to be `super slow' \cite{DGK01}.  In \cite{DV15}, a much more general `slowly nonlinear' system of equations is studied, that however does have a vertical slow manifold $\M_\eps \equiv \M_0 = \{u = 0, p = 0\}$ so that there cannot be stable slow pulses. It is shown that also in this setting (the equivalent of) $\mu_{m+1}$ needs to be $\O(1/\eps)$  to have stable slow-fast-slow pulses. Studying the relations and transitions between slow and slow-fast-slow pulses in (\ref{e:RDE-return}) and especially the impact of $\mu_{m+1}$ on the stability characteristics of these pulses will deepen our understanding of the (de)stabilizing effects of fast jumps in systems like (\ref{e:RDE}).
\\ \\
{\bf Periodic patterns.} From the point of view of the physical (eco)system modeled by (\ref{e:RDE}), (uniformly traveling, or stationary) spatially periodic patterns, or wave trains, are at least as relevant as the localized patterns considered here. The existence problem can be studied very much along the lines of the present paper (in fact, Theorem \ref{t:E-Pers} already settles the problem for the case of a potential $\W_0(v)$ (\ref{d:HW0}) with unequal wells). In the stability analysis, the role of Sturm-Liouville theory can be taken over by the (classical) theory of Hill's equations \cite{MW66} in combination with the concept of $\gamma$-eigenvalues \cite{Gar93} -- that has been worked out as explicit instrument to study the stability of periodic patterns in systems like (\ref{e:RDE}) in \cite{dRDR16,vdPD05}. By associating the spectral stability problem to a family of Hill's equations parameterized by $\rho \in \RR$ -- as in section \ref{sss:famSL} for the homoclinic pulses -- an instability result like Theorem \ref{t:InstabRegHom} can again be established for the case of $\W_0(v)$ with unequal wells.
\\ \\
Like for localized patterns, the case of $\W_0(v)$ with wells of equal depth promises to be much more interesting. Traveling wave trains may be expected to be stable under circumstances similar to those formulated in Theorem \ref{t:nearstabtravfrontspulses}. However, it is expected to be quite challenging to unravel the bifurcational structure associated to the (de)stabilizations of these wave trains. For stationary -- and thus reversible -- patterns, the spectral curves $\{\la_j(\gamma): \gamma \in \S^1\}_{j=1}^{\infty}$ associated to their stability collapse into branches with endpoints associated to $\gamma = \pm 1$ (since a $\gamma$-eigenvalue $\la$ must also be a $\bar{\gamma}$-eigenvalue in the reversible case \cite{DRvdS12,vdPD05}). Therefore, it may be expected for stationary patterns that the bifurcation that (de)stabilizes the homoclinic pulses will be replaced by a similar bifurcation associated to a $\gamma = \pm 1$ endpoint. For a Hopf bifurcation, this would happen by the 2 kinds of Hopf bifurcations that together drive the `Hopf dance mechanism at the boundary of the Busse balloon' that generically appears near the homoclinic limit of a family of periodic patterns in singularly perturbed reaction-diffusion equations \cite{DRdRV18,DRvdS12}. However, this is no longer the case for traveling patterns: the curves $\{\la_j(\gamma): \gamma \in \S^1\}_{j=1}^{\infty}$ will open up into smooth images of $\S^1$. We are not aware of analytical studies of the nature of the bifurcations at the boundary of the Busse balloon -- i.e. the region of stable periodic patterns -- beyond the (reversible) Hopf dance mechanism. Such a study -- that should also include the alternative (de)stabilization scenario associated to he transition from real to complex eigenvalues (section \ref{ss:Bifs}) -- is possible along the lines set out in the present work.
\\ \\
{\bf Spatial ecology.} As explained in section \ref{ss:Eco}, the underlying motivation for the present analysis came from explicit ecosystem models in the literature. To stimulate further cross-fertilization between mathematics and ecology, it is natural -- and interesting and relevant -- to apply the here developed insights in the setting of an explicit model. Especially since it has a unique normally hyperbolic slow manifold without a return mechanism, it is appealing to consider the grass/woodland savanna ecosystem model of \cite{dKetal08,GvdVvL17,vLetal03} in the form (\ref{e:RDE-exvLangevelde}) or one of its more extended versions. Although the case in which the potential associated to its reduced slow flow (\ref{e:RedSF-exvLangevelde}) has the required double well with wells of equal depth structure will be too specific from the ecological point of view, it is a natural starting point, since the associated grass/woodland interface is commonly observed in savanna ecosystems. Therefore, we may proceed from this starting point by looking for traveling localized fronts/interfaces in the neighborhood (in parameter space) of this special case. It is expected that this neighborhood can be extended to an ecologically relevant dimensions by considering speeds $|c| \gg 1$ -- see section \ref{ss:Obs}. The simulations of \cite{GvdVvL17} indicate that the interfaces may be unstable with respect to perturbations along the interface and that the associated bifurcations generate families of spatial patterns (which can also be observed in various savannas). The recent combined ecological and mathematical insights of \cite{BDER20,Basetal18} indicate that this `multistability' may increase the resilience of the ecosystem, it is therefore of both of ecological and mathematical interest to further investigate this model by a combination of analytical and numerical methods.

\appendix

\section{A refinement of Lemma \ref{l:Lah}}
\label{a:refined}
The upcoming lemma considers exactly the same setting as that of Lemma \ref{l:Lah} on the inhomogeneous problem (\ref{e:Lah}). The accuracy of the approximations of solution $v(X)$ of (\ref{e:Lah}) necessary for the analysis varies with the exponential growth/decay rate of the inhomogeneous term $h(X)$ as $X \to \infty$ -- as is made more precise by assumptions (\ref{d:oh01j}) and the introduction of $j \in \ZZ$. Thus, the level of detail of the approximations in (\ref{d:approxsvX-sharp}) varies with $j$ in the upcoming lemma. For any $j$ -- i.e. for all relevant growth/decay rates of $h(X)$ in (\ref{e:Lah}) -- the accuracy can be in principle be improved in a straightforward fashion by the (standard) methods of the upcoming proof.
\begin{lemma}
\label{l:Lah-sharp}
Consider (\ref{e:RedSF}) with $\W_0(v)$ (\ref{d:HW0}) a double well potential with equal (non-degenerate) wells at $v = \oV_\pm$ and a local maximum at $\oV_c \in (\oV_-,\oV_+)$ and let $\va(X)$ be the increasing heteroclinic connection between the saddle points $(\oV_\pm,0)$ with $\va(0) = \oV_c$. Let $\vb(X) = \vax(X)$ and $\vu(X)$ be the 2 independent solutions of the associated homogeneous problem $\La(X) v = 0$ (\ref{d:La}),(\ref{d:vbu}). Then, for $X \gg 1$,
\beq
\label{e:vaXgg1-sharp}
\begin{array}{ccc}
\va(X) & = & \oV_{+} - \be_{+} e^{-\sqrt{\al_+} X} + \frac{\be_+^2 \ga_+}{6 \al_+} e^{-2 \sqrt{\al_+} X} + \O(E_+^3(X))
\\
\vb(X) & = & \be_{+} \sqrt{\al_+} e^{-\sqrt{\al_+} X} - \frac{\be_+^2 \ga_+ }{3 \sqrt{\al_+}} e^{-2 \sqrt{\al_+} X} + \O(E_+^3(X))
\\
\vu(X) & = & \frac{1}{2 \al_+ \be_+} e^{\sqrt{\al_+} X} + \frac{\ga_+}{2 \al_+^2} + \O(E_+(X))
\end{array}
\eeq
(cf. (\ref{d:muEnpm}), (\ref{e:vaXgg1}), (\ref{e:propsvbu})) with
\beq
\label{d:gam+}
\gamma_+ = \oG_{vv}^+ + 2f'(\oV_+)\oG_{uv}^+ + (f'(\oV_+))^2 \oG_{uu}^+ + f''(\oV_+)\oG_{u}^+
\eeq
(\ref{d:oFpmetc}). Assume for the inhomogeneous term $h(X)$ of (\ref{e:Lah}) that there are $\oh_{0,j} \neq 0 $, $\oh_{1,j}$ such that,
\beq
\label{d:oh01j}
\lim_{X \to \infty} h(X) e^{- j \sqrt{\al_+} X} = \oh_{0,j}, \; \;
\lim_{X \to \infty} \left(h(X) e^{- j \sqrt{\al_+} X} - \oh_j\right) e^{- (j-1) \sqrt{\al_+} X} = \oh_{1,j}
\eeq
for some $j \in \ZZ$ (i.e. $h(X) = \oh_{0,j} E_+^{-j}(X) + \oh_{1,j} E_+^{-(j-1)}(X) + \O(E_+^{-(j-2)}(X))$ for $X \gg 1$) and let $v(X)$ be the solution of (\ref{e:Lah}) such that $\lim_{X \to -\infty} v(X) = 0$ and $v(0) = 0$. Then, the leading order approximations of $v(X)$ for $X \gg 1$ are given by,
\beq
\label{d:approxsvX-sharp}
\begin{array}{rl}
j \leq -2: & v(X) = \frac{M_h}{2 \al+ \be_+} E_+^{-1} + \frac{\ga_+ M_h}{2 \al_+^2} + \O(E_+)
\\
j=-1: & v(X) = \frac{M_h}{2 \al+ \be_+} E_+^{-1} + \frac{\ga_+ M_h}{2 \al_+^2} - \frac{\oh_{0,-1}}{2 \sqrt{\al_+}} X E_+ - \left[\be^+ \sqrt{\al_+} M^+_{{\rm u},-1} + \frac{\oh_{0,-1}}{4 \al_+} \right] E_+ + \O(M_h E_+, X E_+^2)
\\
j=0: & v(X) = \frac{M_h}{2 \al_+ \be_+} E_+^{-1} + \left[\frac{\ga_+ M_h}{2 \al_+^2} - \frac{\oh_{0,0}}{\al_+} \right] + \left[\frac{\be_+ \ga_+ \oh_{0,0}}{2 \al_+ \sqrt{\al_+}} + \frac{\oh_{1,0}}{2 \sqrt{\al_+}} \right] X E_+ + \O(E_+)
\\
j=1: & v(X) = \frac{\oh_{0,1}}{2 \sqrt{\al_+}} X E_+^{-1} + \left[\frac{M^-_h + M^+_{{\rm b},1}}{2 \al_+ \be_+} - \frac{\oh_{0,1}}{4 \al_+} \right] E_+^{-1} + \frac{\be_+ \ga_+ \oh_{0,1}}{2 \al_+ \sqrt{\al_+}} X + \O(1)
\\
j \geq 2: & v(X) = \frac{\oh_{0,j}}{(j^2 - 1) \al_+} E_+^{-j} + \O(E_+^{-(j-1)})
\end{array}
\eeq
with $M_h$ as in (\ref{e:vgg1}), $M_h = M_h^- + M_h^+$ with $M_h^- = \int_{-\infty}^0 h(X) \vb(X) \, dX$ and,
\beq
\label{d:Mpubpm1}
\begin{array}{rlr}
j = -1: &
M^+_{{\rm u},-1} = \int_0^\infty \left[ h \vu - \frac{\oh_{0,-1}}{2 \al_+ \be_+} \right] \, dX, & |M^+_{{\rm u},-1}| < \infty
\\[2mm]
j = 1: &
M^+_{{\rm b},1} = \int_0^\infty \left[ h \vb - \oh_{0,1} \be_+ \sqrt{\al_+} \right] \, dX, & |M^+_{{\rm b},1}| < \infty
\end{array}
\eeq
\end{lemma}
\noindent
Note that we have introduced the notation $\O(M_h E_+)$ in (\ref{d:approxsvX-sharp}) for $j = -1$, since in the paper typically $M_h = \Ma(\vmu)$ (\ref{d:Ma}) and $\Ma(\vmu) = \O(\eps^{\tsi})$ (\ref{d:vmutsi}) in the nearly heteroclinic analysis.
\\ \\
{\bf Proof.} The next order approximations of (\ref{e:vaXgg1-sharp}) -- compared to (\ref{e:vaXgg1}), (\ref{e:propsvbu}) -- follow by a direct local analysis of (\ref{e:RedSF}) near the saddle $(\oV_+,0)$. Since it plays a central role in the proofs of Theorem \ref{t:exhoms} and Lemma \ref{l:asymptevhom-hot}, we derive the approximation (\ref{d:approxsvX-sharp}) for $j = -1$ in full detail. All other cases proceed along exactly the same lines: we refrain from going into the details of establishing (\ref{d:approxsvX-sharp}) for $j \neq -1$.
\\ \\
Solution $v(X)$ of (\ref{e:Lah}) is given by (\ref{e:homsolinh}) and we thus first consider the two terms in (\ref{e:homsolinh}) separately. It follows by (\ref{d:approxsvX-sharp}) and assumptions (\ref{d:oh01j}) that
\[
\lim_{X \to \infty} h(X) \vu(X) = \frac{\oh_{0,-1}}{2 \al_+ \be_+},
\]
and that $h(X) \vu(X) - \frac{\oh_{0,-1}}{2 \al_+ \be_+}$ decays as $E_+(X)$ as $X \to \infty$. Hence, $M^+_{{\rm u},-1}$ as defined in (\ref{d:Mpubpm1}) indeed converges and
\[
\begin{array}{rcl}
\int_0^X h(\tX) \vu(\tX) d\tX
& = & \frac{\oh_{0,-1}}{2 \al_+ \be_+} X + \int_0^X \left[h(\tX) \vu(\tX) - \frac{\oh_{0,-1}}{2 \al_+ \be_+} \right] d\tX
\\
& = & \frac{\oh_{0,-1}}{2 \al_+ \be_+} X + M^+_{{\rm u},-1} - \int_X^\infty \left[h(\tX) \vu(\tX) - \frac{\oh_{0,-1}}{2 \al_+ \be_+} \right] d\tX
\\
& = & \frac{\oh_{0,-1}}{2 \al_+ \be_+} X + M^+_{{\rm u},-1} + \O(E_+(X))
\end{array}
\]
for $X \gg 1$. Hence, for $X \gg 1$,
\beq
\label{e:AppA-t1}
\left[\int_0^X h \vu d\tX \right] \vb(X) =  \frac{\oh_{0,-1}}{2 \sqrt{\al_+}} X E_+(X) + \be_+ \sqrt{\al_+} M^+_{{\rm u},-1} E_+(X) + \O(X E_+^2(X)).
\eeq
Similarly,
\[
h(X) \vb(X) = h(X) \vax(X) =  \be_+ \sqrt{\al_+} \oh_{0,-1} E^2_+(X) + \O(E^3_+(X))
\]
for $X \gg 1$, so that,
\[
\begin{array}{rcl}
\int_{-\infty}^X h(\tX) \vb(\tX) d\tX
& = & \int_{-\infty}^\infty h(\tX) \vb(\tX) d\tX - \int_X^{\infty} h(\tX) \vb(\tX) d\tX
\\
& = & M_h - \be_+ \sqrt{\al_+} \oh_{0,-1} \int_X^{\infty} \left[e^{-2\sqrt{\al_+} \tX} + \O(E^3_+(\tX)) \right] d \tX
\\
& = & M_h - \frac12 \be_+ \oh_{0,-1} E^2_+(X) + \O(E^3_+(X)),
\end{array}
\]
and
\beq
\label{e:AppA-t2}
\left[\int_{-\infty}^X h \vb d\tX \right] \vu(X) = \frac{M_h}{2 \al_+ \be_+} E_+^{-1}(X) + \frac{\ga_+ M_h}{2 \al_+^2} - \frac{\oh_{0,-1}}{4 \al_+} E_+(X) + \O(M_h E_+(X), E_+^2(X))
\eeq
for $X \gg 1$ (\ref{e:vgg1}), (\ref{d:oh01j}), (\ref{d:approxsvX-sharp}). The combination of (\ref{e:AppA-t1}) and (\ref{e:AppA-t2}) yields (\ref{d:approxsvX-sharp}) for $j = -1$. \hfill $\Box$

\section{Explicit expressions for the higher order spectral problems}
\label{a:explicit}
Substituting all relevant expansions up to all $\O(\eps^2)$-terms into (\ref{e:ODE-lin}), we find as first line the equation that determines $\bu_2$ as function of $\bv_2$,
\beq
\label{e:Oeps2Specline1}
\begin{array}{rl}
\tau \la_1 \bu_1 + \tau \la_2 \bu_0 = \bu_{0,XX} + \tau c_0 \bu_{1,X} + \tau c_1 \bu_{0,X} + &
\Fa_u \bu_2 + c_0 \tF^\ast_{uu} \bu_1 + (u_2 \Fa_{uu} + v_2 \Fa_{uv})\bu_0 + c_0^2 \tF^\ast_{uuu} \bu_0 \\
+ &
\Fa_v \bv_2 + c_0 \tF^\ast_{vv} \bv_1 + (u_2 \Fa_{uv} + v_2 \Fa_{vv})\bv_0 + c_0^2 \tF^\ast_{vvv} \bv_0
\end{array}
\eeq
with $\tF^\ast_{uu}$ and $\tF^\ast_{vv}$ defined in (\ref{d:tFauuvv}), and
\beq
\label{d:tFauuuvvv}
\tF^\ast_{uuu} = \frac12 \tu_1^2 \Fa_{uuu} + \tu_1\tv_1 \Fa_{uuv} + \frac12 \tv_1^2 \Fa_{uvv}, \; \;
\tF^\ast_{vvv} = \frac12 \tu_1^2 \Fa_{uuv} + \tu_1\tv_1 \Fa_{uvv} + \frac12 \tv_1^2 \Fa_{vvv}
\eeq
with $\tv_1$ and $\tu_1$ as defined in (\ref{e:tv1}) and (\ref{d:tu1}). To arrive at a (relative) transparent expression like (\ref{e:Oeps2EigPb}) in which all $c_0$-, $c_1$-, $\la_1$-, $\la_2$-dependence is factored out, we first need to determine such an expression for $u_2$ and $v_2$ -- the $\O(\eps^2)$ terms of the existence problem. It follows from (\ref{e:Lav2}) that $v_2$ can be written as
\beq
\label{d:tv2s}
v_2 = c_1 \tv_1 + c_0^2 \tv_{2cc} + \tv_2,
\eeq
with again $\tv_1$ as in (\ref{e:tv1}), and $\tv_{2cc}$ and $\tv_2$ uniquely determined by the fact that $\vs \subset W^u((\oV_-,0))$ (\ref{d:usvsexp}) and initial conditions $\tv_{2cc}(0) = \tv_2(0) = 0$ (cf. Lemma \ref{l:Lah}). Hence, by (\ref{e:u12}) we have that also
\beq
\label{d:tu2s}
u_2 = c_1 \tu_1 + c_0^2 \tu_{2cc} + \tu_2,
\eeq
with $\tu_1$ as in (\ref{d:tu1}) and
\beq
\label{e:tu2c}
\tu_{2cc} = \fa' \tv_{2cc} + \frac12 \fa'' \tv_1^2- \tau \left[\qa (\tf_1^\ast)' \tv_1 + \tf_1^\ast \tv_{1,X} \right] - \tau^2 \tF_{2cc}^\ast , \; \; \tu_2 = \fa' \tv_2 - \tF_2^\ast,
\eeq
where $\tf_1^\ast$, $\tF_{2cc}^\ast$ and $\tF_2^\ast$ are defined in the usual way, cf. (\ref{d:fFGast}). It thus follows from (\ref{e:Oeps2Specline1}) that
\beq
\label{e:bu2inbv2}
\bu_2 = \fa'\bv_2 + \la_2 \frac{\tau\fa'\vax}{\Fa_u} - c_1 \frac{\tau (\fa'\vax)_X + (\tF^\ast_{uu} \fa' + \tF^\ast_{vv})\vax}{\Fa_u} + \la_1^2 \frac{\tau \hu_{1\la}}{\Fa_u} - c_0 \la_1 \hI_{2c\la} - c_0^2 \hI_{2cc} - \hI_2
\eeq
(cf. (\ref{e:bu1inbv1})) with,
\beq
\label{d:hI2s}
\begin{array}{ccl}
\hI_{2c\la} & = & \frac{1}{\Fa_{u}}\left[\tau(\hu_{1\la, X}-\hu_{1c}) + \tF^\ast_{uu}\hu_{1\la} - \tF^\ast_{vv}\tv_1 \right] \\[2mm]
\hI_{2cc} & = & \frac{1}{\Fa_{u}}\left[\tau\hu_{1c, X}+\tF^\ast_{uu}\hu_{1c} + \tF^\ast_{vv}\hv_{1c} + \left((\tu_{2cc}\Fa_{uu} + \tv_{2cc}\Fa_{uv} + \tF^\ast_{uuu}) \fa' + (\tu_{2cc}\Fa_{uv} + \tv_{2cc}\Fa_{vv} + \tF^\ast_{vvv}) \right) \vax \right] \\[2mm]
\hI_2 & = & \frac{1}{\Fa_{u}}\left[(\tu_2 \Fa_{uu} + \tv_2 \Fa_{uv}) \fa' + (\tu_2 \Fa_{uv} + \tv_2 \Fa_{vv})\right] \vax
\end{array}
\eeq
Substitution of (\ref{e:bu2inbv2}) and all other relevant expressions into the second line of the equation obtained from (\ref{e:ODE-lin}) at the $\O(\eps^2)$ level, yields (\ref{e:Oeps2EigPb}) with,
\beq
\label{d:hJ2s}
\begin{array}{ccl}
\hJ_{2c\la} & = & \Ga_u \hI_{2c\la} + \hv_{1c}+ \tv_{1,X} - \tG^\ast_{uu}\hu_{1\la} + \tF^\ast_{vv}\tv_1, \\
\hJ_{2cc} & = & \Ga_u \hI_{2cc} - \hv_{1c, X} - \tG^\ast_{uu}\hu_{1c} - \tG^\ast_{vv}\hv_{1c} - \left((\tu_{2cc}\Ga_{uu} + \tv_{2cc}\Ga_{uv} + \tG^\ast_{uuu}) \fa' + (\tu_{2cc}\Ga_{uv} + \tv_{2cc}\Ga_{vv} + \tG^\ast_{vvv}) \right) \vax, \\
\hJ_2 & = & \Ga_u \hI_{2} - \left((\tu_2 \Ga_{uu} + \tv_2 \Ga_{uv}) \fa' + (\tu_2 \Ga_{uv} + \tv_2 \Ga_{vv})\right) \vax
\end{array}
\eeq

\end{document}